\documentclass[11pt]{book}

\usepackage{amsmath, amsthm,amsfonts,amssymb,mathtools,mathdots}

\usepackage{enumitem, imakeidx, framed}
\PassOptionsToPackage{svgnames}{xcolor}
\usepackage{tikz, fancyhdr, stmaryrd, cmll}
\usepackage[a4paper]{geometry}
\usepackage[font=footnotesize]{caption}
\usepackage[british]{babel}

\newtheorem{thm}{Theorem}[chapter]
\newtheorem{cor}[thm]{Corollary}
\newtheorem{prop}[thm]{Proposition}
\newtheorem{lem}[thm]{Lemma}
\newtheorem{requir}{Desideratum}
\newtheorem{conj}[thm]{Conjecture}
\theoremstyle{definition}
\newtheorem{dfn}[thm]{Definition}
\theoremstyle{remark}
\newtheorem{remark}[thm]{Remark}
\newtheorem{exm}[thm]{Example}
\newtheorem{cons}[thm]{Construction}

\usetikzlibrary{arrows, decorations.markings, shapes.geometric, decorations.pathmorphing, decorations.pathreplacing, intersections, patterns,calc, backgrounds}

\pgfdeclarelayer{bg}
\pgfdeclarelayer{mid}
\pgfdeclarelayer{fg}
\pgfsetlayers{bg,mid,fg,main}

\tikzset{
	0c/.style={circle, draw, fill, inner sep=1.5pt},
	1c/.style={->, thick, shorten <=2pt, shorten >=2pt},
	1cboth/.style={<->, thick, shorten <=2pt, shorten >=2pt},
	1clong/.style={->, thick},
	1cthin/.style={->, shorten <=4pt, shorten >=4pt},
	1cdot/.style={->, dashed, thick, shorten <=2pt, shorten >=2pt},
	1cinc/.style={right hook->, thick, shorten <=2pt, shorten >=2pt},
	1cincl/.style={left hook->, thick, shorten <=2pt, shorten >=2pt},
	follow/.style={->, >=stealth, very thick, shorten <=3pt, shorten >=3pt, color=magenta},
	2c/.style={double, thick, shorten <=6pt, shorten >=8pt, decoration={markings,mark=at position -6pt with {\arrow[scale=1.75]{>}}}, preaction={decorate}},
	2cdot/.style={double, dashed, thick, shorten <=10pt, shorten >=10pt, decoration={markings,mark=at position -8pt with {\arrow[scale=1.75]{>}}}, preaction={decorate}},
	3c1/.style={thick, double, double distance=3pt, shorten <=9pt, shorten >=11pt},
    	3c2/.style={thick, shorten <=9pt, shorten >=10pt},
	3c3/.style={shorten <=9pt, shorten >=10pt, decoration={markings,mark=at position -8pt with {\arrow[scale=3]{>}}},preaction={decorate}},
	4c1/.style={thick, double, double distance=4pt, shorten <=1pt, shorten >=2.75pt},
	4c2/.style={thick, double, double distance=1pt, shorten <=1pt, shorten >=1.25pt, decoration={markings,mark=at position -.05pt with {\arrow[scale=3,ultra thin]{>}}},preaction={decorate}},
	edge/.style={line width=.8pt, color=black},
	edgedot/.style={densely dotted, line width=.8pt, color=black},
	edgethdot/.style={densely dotted, line width=.4pt, color=gray},
	edgeth/.style={line width=.4pt, color=gray!60},
	edgethin/.style={line width=.8pt, color=gray!60},
	edgedotdark/.style={densely dotted, line width=.8pt, color=gray!80},
	dot/.style={circle, draw=black, line width=.8pt, fill=white, inner sep=1.7pt},
	dotth/.style={circle, draw=gray!60, fill=gray!60, inner sep=1.5pt},
	dotwh/.style={circle, draw=gray!60, line width=.4pt, fill=white, inner sep=1.7pt},
	dotwhite/.style={circle, draw=black, line width=.8pt, fill=white, inner sep=1.8pt},
	dotdark/.style={circle, draw, fill=black, inner sep=1.5pt},
	dotgrey/.style={circle, draw=black, line width=.8pt, fill=gray!60, inner sep=1.8pt},
	trian/.style={regular polygon,regular polygon sides=3,shape border rotate=0,fill=white, line width=.8pt, draw=black, inner sep=1.8pt},
	trianh/.style={regular polygon,regular polygon sides=3,shape border rotate=0,fill=white, draw=gray!60, line width=.4pt, inner sep=1.8pt},
	trib/.style={regular polygon,regular polygon sides=3,shape border rotate=0,fill=black, draw, inner sep=1.5pt},
	tribh/.style={regular polygon,regular polygon sides=3,shape border rotate=0,fill=gray!60, draw=gray!60, inner sep=1.5pt},
	tribco/.style={regular polygon,regular polygon sides=3,shape border rotate=180,fill=black, draw, inner sep=1.5pt},
	cover/.style={circle, draw=gray!10, fill=gray!10, inner sep=3.5pt},
	coverc/.style={circle, draw=gray!80, line width=.4pt, inner sep=3pt},
	coverch/.style={circle, draw=gray!30, line width=.4pt, inner sep=3pt},
	coverb/.style={circle, draw=gray!80, line width=.4pt, fill=gray!10, inner sep=3.5pt},
	every node/.style={scale=.8},
}

\newcommand\cp[1]{*_{#1}}
\newcommand\comp[1]{\triangleright_{#1}}
\newcommand\bord[2]{\partial_{#1}^{#2}}
\newcommand\idd[1]{\varepsilon_{#1}}

\newcommand\opp[1]{{#1}^\mathrm{op}}
\newcommand\coo[1]{{#1}^\mathrm{co}}
\newcommand\oppn[2]{{#1}^{\mathrm{op}(#2)}}
\newcommand\cone[2]{C^{#1}(#2)}

\newcommand\rimp[2]{#1\!\multimap\!#2}
\newcommand\limp[2]{#2\!\multimapinv\!#1}
\newcommand\rcimp[2]{#1 \diagdown #2}
\newcommand\lcimp[2]{#2 \diagup #1}

\newcommand\cat[1]{\mathbf{#1}}
\newcommand\omegacat{\omega\cat{Cat}}
\newcommand\adc{\mathbf{ADC}}
\newcommand\loopfree{\adc_\uparrow}

\newcommand\globpos{\cat{GlobPos}_\subset}
\newcommand\rpol{\cat{RPol}}

\newcommand\sbord[2]{\Delta_{#1}^{#2}}
\newcommand\skel[2]{#2^{(#1)}}
\newcommand\dmn[1]{\mathrm{dim}(#1)}
\newcommand\augm[1]{(#1,d)}
\newcommand\lfgen[1]{\langle #1 \rangle}
\newcommand\transp[1]{\lfloor #1 \rfloor}
\newcommand\invrs[1]{#1^{-1}}
\newcommand\cutt[1]{\mathrm{cut}_{#1}}
\newcommand\dagg[1]{#1^\dagger}

\newcommand\braket[2]{\langle\, #1 \,|\, #2 \,\rangle}
\newcommand\bra[1]{\langle\, #1\,|}
\newcommand\ket[1]{|\, #1 \, \rangle}
\newcommand\ketbra[2]{\ket{#1} \! \bra{#2}}
\newcommand\upto{\underset{\,\diamond}{=}}
\newcommand\ghz{\ket{\textit{GHZ}}}
\newcommand\wstate{\ket{W}}
\newcommand\qint[1]{[#1]_q}

\newcommand\gbullet{{\!\text{
	\begin{tikzpicture}
		\node[dotgrey] at (0,0) {};
	\end{tikzpicture}}}
}
\newcommand\wbullet{{\!\text{
	\begin{tikzpicture}
		\node[dotwhite] at (0,0) {};
	\end{tikzpicture}}}
}
\newcommand\bbullet{{\!\text{
	\begin{tikzpicture}
		\node[dotdark] at (0,0) {};
	\end{tikzpicture}}}
}

\newcommand\restr[2]{{
  \left.\kern-\nulldelimiterspace 
  #1 
  \vphantom{\big|} 
  \right|_{#2}
  }}

\makeindex[title=Index]

\begin{document}

\pagestyle{empty}
\begin{titlepage}
\begin{center}
\mbox{}\\[6pt]
\Huge \textbf{The algebra of entanglement \\[6pt] and the geometry of composition} \\[100pt] 

\includegraphics[width=0.2\columnwidth]{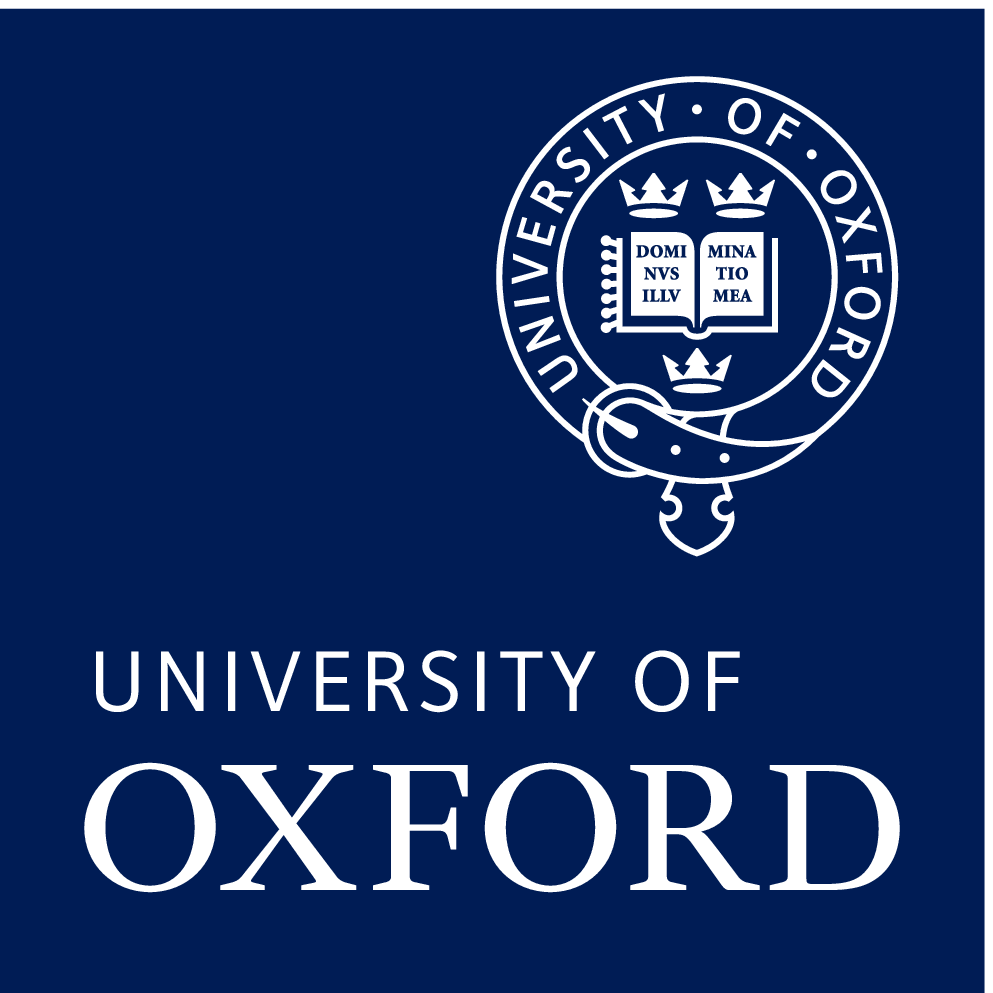} \\[90pt]

\LARGE Amar Hadzihasanovic \\[6pt]
\Large Wolfson College \\[3pt] 
University of Oxford \\[75pt]

A thesis submitted for the degree of \\[3pt]
\emph{Doctor of Philosophy} \\[12pt]
Trinity 2017
\end{center}
\end{titlepage}

\thispagestyle{empty} 

 \cleardoublepage 
\selectlanguage{british}
\begin{center}\textbf{Abstract}\end{center}

\vspace{10pt}

\begin{quotation}
\noindent \small 
String diagrams turn algebraic equations into topological moves. These moves have often-recurring shapes, involving the sliding of one diagram past another. In the past, this fact has mostly been considered in terms of its computational convenience; this thesis investigates its deeper reasons.

In the first part of the thesis, we individuate, at its root, the dual nature of polygraphs --- freely generated higher categories --- as presentations of higher algebraic theories, and as combinatorial descriptions of ``directed spaces'': CW complexes whose cells have their boundary subdivided into an input and an output region. Operations of polygraphs modelled on operations of topological spaces, including an asymmetric tensor product refining the topological product of spaces, can then be used as the foundation of a compositional universal algebra, where sliding moves of string diagrams for an algebraic theory arise from tensor products of sub-theories. Such compositions come automatically with higher-dimensional coherence cells. We provide several examples of compositional reconstructions of higher algebraic theories, including theories of braids, homomorphisms of algebras, and distributive laws of monads, from operations on polygraphs.

In this regard, the standard formalism of polygraphs, based on strict $\omega$-categories, suffers from some technical problems and inadequacies, including the difficulty of computing tensor products. As a solution, we propose a notion of regular polygraph, barring cell boundaries that are degenerate, that is, not homeomorphic to a disk of the appropriate dimension. We develop the theory of globular posets, based on ideas of poset topology, in order to specify a category of shapes for cells of regular polygraphs. We prove that these shapes satisfy our non-degeneracy requirement, and show how to calculate their tensor products. Then, we introduce a notion of weak unit for regular polygraphs, allowing us to recover weakly degenerate boundaries in low dimensions. We prove that the existence of weak units is equivalent to the existence of cells satisfying certain divisibility properties --- an elementary notion of equivalence cell --- which prompts new questions on the relation between units and equivalences in higher categories.

In the second part of the thesis, we turn to specific applications of diagrammatic algebra to quantum theory. First, we re-evaluate certain aspects of quantum theory from the point of view of categorical universal algebra, which leads us to define a 2-dimensional refinement of the category of Hilbert spaces. Then, we focus on the problem of axiomatising fragments of quantum theory, and present the ZW calculus, the first complete diagrammatic axiomatisation of the theory of qubits.

The ZW calculus has specific advantages over its predecessors, the ZX calculi, including the existence of a computationally meaningful normal form, and of a fragment whose diagrams can be interpreted physically as setups of fermionic oscillators. Moreover, the choice of its generators reflects an operational classification of entangled quantum states, which is not well-understood for states of more than 3 qubits: our result opens the way for a compositional understanding of entanglement for generic multipartite states. We conclude with preliminary results on generalisations of the ZW calculus to quantum systems of arbitrary finite dimension, including the definition of a universal set of generators in each dimension.
\end{quotation}
\cleardoublepage

\frontmatter \pagestyle{plain} 
\chapter{Preface}
\thispagestyle{plain}

This thesis is a report of some results, and the state of my thoughts on a variety of subjects, after four years spent as a doctoral student in the Oxford Quantum Group. The overarching theme is the relationship between algebra, or computation, and geometry, or topology; another recurring theme is compositionality. 

That algebra is useful for geometry is a notion by now ingrained in mathematical thinking at any level. The converse is somewhat less mainstream, driven recently to the spotlight by the contributions of higher/homotopical algebra and higher category theory to a number of fields, including physics via topological quantum field theory. Even so, in many cases geometry inspires at most a change of environment for algebraic structures, where the algebraic concepts are interpreted in a less rigid context --- ``weaker'', in higher-categorical jargon --- yet by themselves remain atomic and unvaried. That is, we may have ``homomorphisms up to homotopy'', but the concept of homomorphism is still fundamentally algebraic. Is it?

Being already interested, as a student, in category theory and mathematical physics, at some point in 2012 I saw a variant of
\begin{equation} \label{eq:stringadjunct}
\input{img/c0_adjunction.tex}
\end{equation} 
in a paper of Bob Coecke's on categorical quantum mechanics, and I was introduced to string diagrams. The fact that adjunctions, the omnipresent concept of categorical algebra, is encapsulated in graphical equations such as (\ref{eq:stringadjunct}) struck me immediately; the following year I applied to become Bob's student, largely because of the appeal that the diagrammatic approach, central to categorical quantum mechanics, had on me.

Later on, thanks to Dan Marsden, I was exposed to 
\begin{equation} \label{eq:stringhomo}
\input{img/c0_homomorphism.tex}
\end{equation} 
and related equations, capturing the notion of algebra homomorphism in great generality; similar pictures describe modules and their intertwiners. They endow the basic ingredients of algebra with specific geometries, comprising planar isotopy (\ref{eq:stringadjunct}) and ``sliding rules'' such as (\ref{eq:stringhomo}), and hinting at a connection at the most fundamental level. Sliding rules are particularly common in diagrammatic algebra, and after a while I was compelled to ask: why these specific pictures? Are there general combinatorial-geometric processes by which equations of a certain shape arise? 

With time, this question, and others spawned by its partial answers, have become the main focus of my research; they are the subject of the first three chapters, to which the \emph{geometry of composition} in the title refers. The key to understanding came from an entirely separate line of applications of topology to algebra, namely, homotopical methods in rewriting theory, and in particular their recent reworking into the theory of polygraphs, due to Albert Burroni and others. A polygraph is a presentation of a higher category or algebraic theory, but it is also a combinatorial model of a CW complex, made up of abstract ``directed cells'': it is my claim, in this thesis, that the specific efficacy of string diagrams and their higher-dimensional generalisations rests on this analogy.

On the other hand, the \emph{algebra of entanglement} refers to certain algebraic theories, relevant to quantum information theory, on which I worked before and besides focussing on the questions above. Entanglement is a property of quantum systems that translates, diagrammatically, into connectedness, hence a topological property; yet Bob and Aleks Kissinger linked it to the properties of certain algebras in the category of Hilbert spaces (Section \ref{sec:entanglement}). It seemed an ideal testing ground for the dual algebraic-topological nature of string diagrams, and it led me to develop the ZW calculus, the first complete diagrammatic axiomatisation of the theory of qubits, which is the main subject of the last two chapters.

In fact, the bulk of my work on the ZW calculus came before and was a catalyst for the broader scope of the first part. Even from this short summary, then, it should be clear that the structure of the thesis does not reflect chronology in any meaningful way. The path to the current shape of this work has been (unsurprisingly) tortuous, full of wrong intuitions and all-out changes of mind; more will come in the future, no doubt. So, before leaving you to the post-processed flow of arguments, I will attempt a quick recollection of my actual thought process through its main turns, and thus provide an alternative summary of the thesis. I hope that the knowledge of some abandoned ideas may complement the surviving ones.

Early in 2014, I followed Bob's suggestion to look into the ``GHZ/W calculus'', as an early sketch of the ZW calculus was called. I had recently learned about its precursor, the ZX calculus (Section \ref{sec:frobenius}), whose completeness for the stabiliser fragment of quantum mechanics was proved by Miriam Backens around that time, and what I especially liked about it was the undirectedness: the generators were symmetrical with respect to the swapping and transposition of wires, so the string diagrams could effectively be manipulated as labelled graphs. The GHZ/W calculus did not initially share this property, and I decided that I had to make it undirected as well.

My intuition, at the time, was that these graph-like string diagrams reveal a deeper, undirected geometry of processes, which gets lost in the ``bureaucracy'' of categorical algebra, with its need for direction and composition functors and structural morphisms; in this I was influenced by Bob's comments in this direction, and those of Jean-Yves Girard, and others, on proof nets being a kind of de-bureaucratised proof theory. In other words, the geometry is fundamental, while the categorical algebra is imposed on it to express compositions combinatorially.

In the early summer of 2014 I discovered a normal form for ZW diagrams, and in August 2014 I had a complete axiomatisation, together with a proof of completeness for what is called the vanilla ZW calculus in Section \ref{sec:fragments}. The proof for the general $R$-labelled case, Theorem \ref{thm:completenesszw}, came only a couple of years later, not because it is substantially different or more complicated, but because I had not tried. The case $R = \mathbb{C}$ of Theorem \ref{thm:completenesszw} is completeness for the theory of qubits, what had escaped the ZX calculi for nearly a decade: it means that any calculation on pure quantum states of a number of qubits that can be performed with vector calculus, or matrix calculus, can also be performed diagrammatically, using the axioms of the ZW calculus.

Anyway, I felt that this result vindicated my ``undirected'' perspective. The next step, I thought, would be to fully disclose the hidden underlying topology of the ZW calculus, by somehow making all the axioms ``trivial'' for a certain topological interpretation of the diagrams. A fragment of the theory was about commutative Frobenius algebras, and those already had a topological interpretation as cobordisms of 2-dimensional manifolds, so instead I focussed on the bialgebra axioms, such as
\begin{equation} \label{eq:stringbialg}
\input{img/c0_bialgebra.tex}
\end{equation} 
In September, I came up with the following picture:
\begin{equation} \label{eq:stringbialg3d}
\input{img/c0_bialgebra3d.tex}
\end{equation} 
where the two sides of (\ref{eq:stringbialg}) appear as the intersections of two perpendicular branching surfaces, as they slide past each other.

I showed it to several people, including Jamie Vicary, who immediately observed that the same representation could be applied to different kinds of naturality equations --- the sliding rules that I mentioned earlier. I decided that I needed to make sense of (\ref{eq:stringbialg3d}), formally, and I soon realised two things:
\begin{enumerate}
	\item that the two perpendicular planes of (\ref{eq:stringbialg3d}) correspond to copies of the planar theories of monoids and comonoids, yet their interaction results in a three-dimensional theory, and
	\item that whatever the interaction is, its objects have to be fundamentally ``directed'', for applied to ``undirected'' branching surfaces, as the ones that could arise from the theory of Frobenius algebras, it would produce ``Frobenius bialgebras'' which are degenerate.
\end{enumerate}
The second was a blow to my ambition of purely undirected foundations, but I would not give up for a while longer. The first indicated that objects of some fixed dimension would not be closed under the interaction, which to me suggested that I should learn about higher categories in arbitrary dimension.

I spent the end of 2014 and the first half of 2015 reading about any notion of higher category that I could find, searching for one that would suit my purpose. At some point I learned about the Crans-Gray tensor product of strict $\omega$-categories, which would turn out to be, essentially, what I was looking for; I did not realise it at first, both because (\ref{eq:stringbialg3d}) actually comes from a \emph{quotiented} tensor product (a kind of smash product), and because everybody around me seemed to think that strict higher categories are wrong and useless.

Eventually, I decided that it would be simpler to implement the interaction at the level of a ``shape category'', and then extend it to its presheaves, so I studied various combinatorics of shapes. Somehow, I still believed that undirected shapes would turn out to be fundamental, which led me to study combinatorial topology and develop a core of what has become Section \ref{sec:globularposet}, but also to abandon that line of research for a long time, when a mixture of directed and undirected, as I had vaguely in mind, did not seem to work out.

In July 2015, I took a break from this when I presented the ZW calculus at LICS in Kyoto and QPL in Oxford. Because the ZX calculus is symmetrical with respect to the exchange of its two component algebras (the Z and the X in the name), and the original paper by Bob and Aleks on GHZ/W also focussed on the similarities, with one point of divergence, between two algebras, in its early versions I tried to make the ZW calculus as symmetrical as possible, developing the axiomatisation of the two algebras parallelly; this is still evident in the axioms published in \cite{hadzihasanovic2015diagrammatic}. However, following a talk by Alessandro Tosini at QPL, I realised that the ``pure fragment'' of the ZW calculus was suitable for describing fermionic quantum computation, and I shifted to a more asymmetrical view, where the Z algebra is a late addition to the fermionic core of the calculus. This is how it is presented in Chapter \ref{chap:zwcalculus}.

One of the reasons I insisted on ``undirected category theory'' was the lack of a satisfactory placement of dagger categories, one of the fundamental concepts of categorical quantum mechanics, within general categories with dualities. In August, I tried a different approach, attempting to derive the dagger of Hilbert spaces from some other categorical construction. I worked out a generalisation of Peter Selinger's CPM construction, related to the Grothendieck construction, which made daggers derivable from inverses of isomorphisms in some relevant cases, and answered a question of Chris Heunen's about the relation between the categories of Hilbert spaces and of partial injections. It also lacked good functorial properties, and left me unconvinced overall; but it contained the idea of decomposing contractive maps of Hilbert spaces as pairs of isometries, which is the basis of a more persuasive reconstruction in Section \ref{sec:abelian}. This result, Proposition \ref{prop:hilbtrunc}, presents the category of Hilbert spaces and short linear maps as the truncation --- a kind of lower-dimensional projection --- of a bicategory of cospans of isometries, whose dagger is purely combinatorial.

In the autumn, as I taught some classes for Samson Abramsky's course on \emph{Categories, Proofs and Processes}, I returned to my favourite subject, categorical logic and proof theory, whose ideas inform the entire thesis, even though they are not addressed explicitly. In the ``topology-first, shapes-only'' approach I had by then adopted, I found it troubling that equation (\ref{eq:stringadjunct}) could be interpreted in a $*$-autonomous category --- a categorical model of multiplicative linear logic --- but the inputs and outputs of each cell would be composed according to different monoidal structures.

I found a solution which involved representing the monoidal structures by universal cells, internal to a wider structure with a purely combinatorial-geometric composition; then, when in January 2016, invited by Claudia Faggian,  I visited the Laboratoire PPS in Paris, I learnt from Paul-Andr\'e Melli\`{e}s that I had just re-discovered Claudio Hermida's \emph{coherence via universality}.

At the same time, I learnt that not only strict, but \emph{freely generated} strict $\omega$-categories were still popular in France, under the name of polygraphs. I finally felt justified in studying them, leaving the weak coherent structure to be handled subsequently through internal representability conditions; and in a couple of months, thanks in particular to the beautiful works of Richard Steiner on the subject, I came up with the basic framework of Chapter \ref{chap:interacting} --- minus some important technical details that I initially glossed over, but are now fully worked out in Chapter \ref{chap:polygraphs}, such as a proof of the fact that the tensor product of two polygraphs is still a polygraph (Theorem \ref{thm:tensorpoly}).

This is a kind of compositional universal algebra, where higher algebraic theories, in the guise of polygraphs, are composed together and manipulated as if they were topological spaces, in order to produce different algebraic theories. Basic algebraic notions and constructions are recovered as directed refinements of topological ones: homomorphisms of algebras from cylinders (Definition \ref{dfn:cylinder}), actions from cones (Definition \ref{dfn:cone}); topological moves for braids come from cubes (Example \ref{exm:braidings}), and products of copies of the theory of monoids underlie both distributive laws of monads (Example \ref{exm:tensorpro}) and the equational theory of bialgebras (Section \ref{sec:dimension}). 

The fact that this form of ``directed topology'' seemed to have such a richer space of possibilities definitively dissuaded me from my earlier intuition that the undirected should come first. I also finally had an explanation of equation (\ref{eq:stringbialg3d}), see Section \ref{sec:dimension}: a conceptual one, at least. Technically, I was soon facing the limitations of strict $\omega$-categories, both the reasons why they are often dismissed by category theorists, and others that were specific to my purposes. I formulated a personal list of \emph{desiderata} as a guideline for a version of polygraphs more adequate for compositional universal algebra; that list has now become Section \ref{sec:strictness}.

At first, I thought an answer could be achieved by cubical methods, if anything because of the role that cubical sets play in Marco Grandis's directed algebraic topology, that I wanted to connect to the homotopical methods in polygraph-based rewriting theory. This resulted in several months of work that have left no trace in this thesis. 

Meanwhile, I had become convinced that the main obstacle laid in the modelling of weak units, as suggested by a conjecture of Carlos Simpson's. Through the papers of Joachim Kock and Andr\'e Joyal, I found out about a notion of weak unit that could be formulated in a regularised version of polygraphs, where cells with degenerate boundaries are barred. Interestingly, in low dimensions, I was able to prove that the existence of such weak units is equivalent to the existence of cells satisfying certain divisibility properties: an elementary notion of equivalence cell, independent of any prior notion of unit (Theorem \ref{thm:0representable}). Thus, at least in this framework, equivalences can be seen as fundamental, and units as derived. The relation between identities and equivalences has been subject to much scrutiny in the last few years, due to homotopy type theory and the Univalent Foundations programme, which makes this result particularly thought-provoking.

Around August 2016, my focus became to define what ``cells with non-degenerate boundaries'' should be in arbitrary dimension. At first, I considered a category of shapes ``modelled on Batanin cells at different resolutions'', but that went nowhere. Then, I realised that the reasons why, one year before, I had dropped the combinatorial topology approach did not apply anymore; the result is the basic theory of globular posets, as presented in Section \ref{sec:globularposet}. I defined a category of cell shapes that are non-degenerate in a sense made precise by Theorem \ref{thm:homeodisk}: the $k$-boundaries of $n$-cells have geometric realisations homeomorphic to $k$-disks, for all $k \leq n$. Moreover, the shapes are closed under tensor products (Theorem \ref{thm:globpostensor}), which are also easy to calculate.

I called presheaves on this category of shapes \emph{regular polygraphs}. The combined theory of regular polygraphs and representability conditions for weak coherent structure is fully worked out, in a low-dimensional case, in Section \ref{sec:weakness}. Some of it was intended to be a part of a survey on the polycategorical semantics of linear logic, that I started in December jointly with Alex Kavvos, but we froze in 2017, as we intensified work on our respective theses.

Finally, I developed the basic idea for the qudit generalisations of the ZW calculus, Section \ref{sec:qudits}, during QPL 2016 in Glasgow, based on an earlier suggestion by Prakash Panangaden that the calculus may have a connection to anyonic algebra. Very recently, there has been renewed excitement about the ZW calculus, after it was used to derive missing axioms of the ZX calculus, so more developments may be coming soon.

This concludes my quick recollection; I leave you to the actual content, with a few indications. I chose a conversational style, with ample discussions of ideas and results, rather than a ``wall of theorems''. This also means that technical notions are introduced only as soon as they are needed, rather than collected into a single, homogeneous section. The downside is that remarks made at some point may be referenced at any other point, and that some notational conventions change as the subjects progress, and certain distinctions stop being relevant; I tried to avoid slipping from consistency into pedantry. For these reasons, the thesis is better read from cover to cover.

This is also the only part in which I use the first-person singular; coupled with the conversational style, I felt it would distract from the mathematical content. Instead, I will disappear behind an inclusive first-person plural. 

The results are my own, except where otherwise stated. Acknowledgements are at the end.

\vspace{15pt}
\hfill Oxford, June and September 2017
 \thispagestyle{empty}
\renewcommand\thepage{}
\tableofcontents
\newpage \thispagestyle{empty} \renewcommand\thepage{\arabic{page}}

\mainmatter \pagestyle{fancy}
\chapter{Polygraphs, categories, diagrams} \label{chap:polygraphs}
\thispagestyle{plain}

\noindent\emph{In this chapter:}
\begin{itemize}
	\item[$\triangleright$] Abstracting from the definition of CW complex, through the combinatorics of globular sets and $\omega$-categories, we arrive at the definition of polygraph. We introduce some basic operations of $\omega$-categories: skeleta, truncations, disjoint unions. --- \emph{Section \ref{sec:topological}}
	\item[$\triangleright$] The tensor product of $\omega$-categories is a non-symmetric refinement of the product of topological spaces. We show how to calculate it directly for a class of suitably loop-free polygraphs, then extend it to arbitrary $\omega$-categories. Finally, we prove that the tensor product of two polygraphs is a polygraph. --- \emph{Section \ref{sec:loopfree}}
	\item[$\triangleright$] We discuss the relation between pasting diagrams and string diagrams, two common devices for picturing cells in higher categories. We complete the toolkit of polygraph-building operations with quotients and opposites. --- \emph{Section \ref{sec:pasting}}
\end{itemize}

\section{From topological to directed cells} \label{sec:topological}
This chapter is meant to be an introduction to higher categories from a combinatorial-topological point of view. Rather than treating them as generalised categories --- which comes with an expectation of large, organisational, container objects: think of the ``infinity-category of infinity-categories'' --- we treat them as \emph{generalised spaces}, ``small'' objects to be studied on their own, without any meta-theoretical presumption. 

With that in mind, we make a number of choices. We keep our treatment elementary, favouring explicit, inductive definitions over ones that are more structural and open to generalisation, yet rely on more abstract algebraic machinery. Moreover, we pick one route --- one that is particularly suited to the later applications --- and follow it, only commenting on alternative paths when it seems to substantially enrich comprehension. Finally, we insist on the ``generalised space'' analogy, and always favour definitions and constructions that have a topological flavour, over ones which may be more familiar to category theorists. Categories in arbitrary dimension are fundamental, low-dimensional ones are only particular examples, and if we show a bias towards the latter, it will only be for want of pictorial representations.

We assume a basic knowledge of category theory and of algebraic topology; \cite{riehl2016category} and \cite{may1999concise} are recommended sources. When here, and in the rest of the thesis, some peripheral notions are mentioned without definition, a good starting point is usually the $n$Lab wiki, \texttt{https://ncatlab.org}.

The focus is on the objects called \emph{computads} by Street \cite{street1976limits}, independently rediscovered by Burroni under the name of \emph{polygraphs} \cite{burroni1993higher}, by which they are best known in the higher-dimensional rewriting community \cite{metayer2003resolutions,lafont2007algebra,guiraud2016polygraphs}. To explain what a polygraph is, we can start from the notion of CW complex; see \cite[Chapter 10]{may1999concise}, or another textbook of basic algebraic topology.

A CW complex is a topological space obtained by successively adjoining \emph{$n$-cells}, that is, copies of $D^n$ --- the topological $n$-disk --- for increasing $n$, attaching them along their boundaries $\partial D^n$, homeomorphic to the $(n-1)$-sphere $S^{n-1}$. Such a space $X$ comes with a stratification in the $n$-dimensional spaces obtained at each step of the construction, the \emph{$n$-skeleta} $\skel{n}{X}$:\index{CW complex}\index{skeleton!of a CW complex}
\begin{itemize}
	\item $\skel{0}{X}$ is just a set of points, that is, a discrete topological space;
	\item $\skel{1}{X}$ is obtained by attaching the boundary points of a number of copies of $D^1 \equiv I$, the topological interval, to the points in $\skel{0}{X}$: it is a \emph{topological graph};
\end{itemize}
and so on. How an $n$-cell $x: D^n \to X$ is attached is determined by a \emph{gluing map} $\partial x: \partial D^n \to \skel{n-1}{X}$. We say that the $n$-cells used in constructing the space $X$ are the $n$-dimensional \emph{generators} of $X$, and write $G_n(X)$ for their set.

A reason for the interest in CW complexes in algebraic topology is that they are at least partially combinatorial, hence more amenable to calculations; yet any topological space is weakly equivalent to a CW complex, that is, has homotopy groups isomorphic to those of a CW complex in all dimensions, so topology ``up to weak equivalence'' can be reduced to the study of CW complexes. Moreover, CW complexes are precisely the cofibrant objects of the classical Quillen-Serre model structure on topological spaces \cite[Section II.3]{quillen1967homotopical}.

Polygraphs are one way of distilling the combinatorial structure of CW complexes, and generalising it in the process. Instead of topological $n$-cells with their point-set gluing maps, we want a set of ``abstract $n$-dimensional generators'', and a combinatorial specification of gluing. For 0-cells, nothing needs to be done: they are just a set. The gluing of a 1-cell $f: D^1 \to X$ is determined by the image of the two boundary points of $D^1$, which we call $+$ and $-$, in the 0-skeleton; hence, by a pair of 0-cells $f(+)$, $f(-)$. We can thus describe our 1-dimensional CW complex $\skel{1}{X}$ with a diagram
\begin{equation} \label{eq:graph}
\begin{tikzpicture}[baseline={([yshift=-.5ex]current bounding box.center)}]
	\node[scale=1.25] (0) at (-2,0) {$X_0$};
	\node[scale=1.25] (1) at (0.3,0) {$G_1(X)$};
	\draw[1c] (-.35,.15) to node[auto,swap] {$\bord{}{+}$} (-1.7,.15);
	\draw[1c] (-.35,-.15) to node[auto] {$\bord{}{-}$} (-1.7,-.15);
\end{tikzpicture}
\end{equation}
of functions of sets, where, for all $f \in G_1(X)$ and $\alpha \in \{+,-\}$, we define $\bord{}{\alpha}(f) := f(\alpha)$.

\begin{itemize}
	\item[] \textbf{Notation.} The variables $\alpha, \beta$ will always implicitly range over $\{+,-\}$; moreover, let $\alpha\beta := +$ when $\alpha = \beta$, and $-$ when $\alpha \neq \beta$.
\end{itemize}
	
On the other hand, this description presupposes an arbitrary choice, where for each cell we pick one of the points in the boundary of $D^1$ to be $-$, and the other to be $+$. This is the same as a choice of direction on the 1-cell, say from $f(-)$ to $f(+)$; in fact, what we gave is the definition of a \emph{directed graph}. If $a = \bord{}{-}f$, $b = \bord{}{+}f$, we picture $f$ as
\begin{equation*}
\begin{tikzpicture}
	\node[0c] (0) at (0,0) [label=left:$a$] {};
	\node[0c] (1) at (1.5,0) [label=right:$b$] {};
	\draw[1c] (0) to node[auto] {$f$} (1);
	\node[scale=1.25] at (2.1,-.1) {,};
\end{tikzpicture}
\end{equation*}
and also write $f: a \to b$.

On to 2-cells: in many cases, we can require that the gluing map $\partial x: \partial D^2 \to \skel{1}{X}$ cover a finite loop in the graph $\skel{1}{X}$, with constant time per edge with respect to the standard metric on the unit circle; such a map is combinatorially described by the finite sequence of edges that it traverses. Since our edges are given with a direction, we have to specify whether $f$ is traversed with the chosen orientation (write $f$) or the opposite one (write $f^*$, and $\bord{}{\alpha}f^* := \bord{}{-\alpha}f$). Moreover, $f$ followed by $f^*$, or vice versa, can be eliminated from the sequence, without affecting the homotopy type of the resulting space.

This can be rephrased as follows. Call $f, g$ \emph{composable} if $\bord{}{+}f = \bord{}{-}g$. Construct the \emph{free groupoid}
\begin{equation} \label{eq:free_gpd}
\begin{tikzpicture}[baseline={([yshift=-.5ex]current bounding box.center)}]
	\node[scale=1.25] (0) at (-2,0) {$X_0$};
	\node[scale=1.25] (1) at (0,0) {$X_1$};
	\draw[1c] (-.35,.15) to node[auto,swap] {$\bord{}{+}$} (-1.7,.15);
	\draw[1c] (-.35,-.15) to node[auto] {$\bord{}{-}$} (-1.7,-.15);
\end{tikzpicture}
\end{equation}
on the graph (\ref{eq:graph}): that is, consider the set of sequences
\begin{equation*} 
\begin{tikzpicture}[baseline={([yshift=-.5ex]current bounding box.center)}]
	\node[0c] (0) at (0,0) [label=left:$a_1$] {};
	\node[0c] (1) at (1.5,0) [label=above:$a_2$] {};
	\node[0c] (2) at (3,0) [label=above:$a_n$] {};
	\node[0c] (3) at (4.5,0) [label=right:$a_{n+1}$] {};
	\draw[1c] (0) to node[auto] {$f_1$} (1);
	\draw[1cdot] (1) to (2);
	\draw[1c] (2) to node[auto] {$f_n$} (3);
\end{tikzpicture}
\end{equation*}
of composable edges --- here $n \in \mathbb{N}$, and the dotted arrow stands for a sequence of arbitrary length --- where each $f_i$ is either an edge $f$, or an edge with opposite orientation $f^*$, for some $f \in G_1(X)$. Then, $X_1$ is the quotient of this set by the equivalence relation generated by
\begin{equation*}
\begin{tikzpicture}
\begin{scope}
	\node[0c] (0) at (0,0) [label=left:$a$] {};
	\node[0c] (1) at (1.5,0) [label=above:$b$] {};
	\node[0c] (2) at (3,0) [label=right:$a$] {};
	\draw[1c] (0) to node[auto] {$f$} (1);
	\draw[1c] (1) to node[auto] {$f^*$} (2);
	\node[scale=1.25] at (4,0) {$\sim$};
	\node[0c] (0) at (4.8,0) [label=above:$a$] {};
	\node[scale=1.25] at (5.2,-.1) {,};
\end{scope}
\begin{scope}[shift={(6,0)}]
	\node[0c] (0) at (0,0) [label=left:$b$] {};
	\node[0c] (1) at (1.5,0) [label=above:$a$] {};
	\node[0c] (2) at (3,0) [label=right:$b$] {};
	\draw[1c] (0) to node[auto] {$f^*$} (1);
	\draw[1c] (1) to node[auto] {$f$} (2);
	\node[scale=1.25] at (4,0) {$\sim$};
	\node[0c] (0) at (4.8,0) [label=above:$b$] {};
	\node[scale=1.25] at (5.2,-.1) {,};
\end{scope}
\end{tikzpicture}
\end{equation*}
for all $f \in G_1(X)$, freely extended to longer sequences containing either side. The boundary functions are extended to finite sequences by $\bord{}{-}(f_1,\ldots,f_n) := \bord{}{-}f_1$, and $\bord{}{+}(f_1,\ldots,f_n) := \bord{}{+}f_n$, and pass to the quotient. Then, the gluing of a 2-cell is specified by an element $f$ of $X_1$ such that $\bord{}{+}f = \bord{}{-}f$.

This requires an arbitrary choice of a ``starting point'': any cyclic permutation of the sequence of edges characterises the same gluing map. Indeed, we might as well opt fosr a more convenient setup: subdivide the boundary $S^1$ of a 2-cell $x$ into two copies of $D^1$, $\bord{}{-}x$ and $\bord{}{+}x$, and specify the gluing with a \emph{pair} of sequences of edges that have the same initial and final point. Overall, the extension to the 2-skeleton is described by a diagram
\begin{equation} \label{eq:2_extension}
\begin{tikzpicture}[baseline={([yshift=-.5ex]current bounding box.center)}]
	\node[scale=1.25] (0) at (-2,0) {$X_0$};
	\node[scale=1.25] (1) at (0,0) {$X_1$};
	\node[scale=1.25] (2) at (2.3,0) {$G_2(X),$};
	\draw[1c] (-.35,.15) to node[auto,swap] {$\bord{}{+}$} (-1.7,.15);
	\draw[1c] (-.35,-.15) to node[auto] {$\bord{}{-}$} (-1.7,-.15);
	\draw[1c] (1.65,.15) to node[auto,swap] {$\bord{}{+}$} (.3,.15);
	\draw[1c] (1.65,-.15) to node[auto] {$\bord{}{-}$} (.3,-.15);
\end{tikzpicture}
\end{equation}
satisfying the \emph{globularity} condition $\bord{}{\alpha}\bord{}{+} = \bord{}{\alpha}\bord{}{-}$. We picture $x$  as \index{globularity}
\begin{equation*}
\input{img/c1_1_2cell.tex}
\end{equation*} 
where $\bord{}{-}x$ is the sequence $f_1, \ldots, f_n$, and $\bord{}{+}x$ is the sequence $g_1,\ldots,g_m$. 

Formally, this step does not rely on (\ref{eq:free_gpd}) being a free groupoid: a \emph{free category} would work as well, or even an entirely different ``free algebra'', describing a notion of cell that has no direct point-set analogue. In other words, what is needed is any algebraic notion of composition with respect to which, at each step of the construction, we can generate all possible composites of the newly added and of the pre-existing cells, so that they can function as boundaries of higher-dimensional ones.

From (\ref{eq:2_extension}), iterating the construction, we would first generate
\begin{equation} \label{eq:2cat}
\begin{tikzpicture}[baseline={([yshift=-.5ex]current bounding box.center)}]
	\node[scale=1.25] (0) at (-2,0) {$X_0$};
	\node[scale=1.25] (1) at (0,0) {$X_1$};
	\node[scale=1.25] (2) at (2,0) {$X_2,$};
	\draw[1c] (-.35,.15) to node[auto,swap] {$\bord{}{+}$} (-1.7,.15);
	\draw[1c] (-.35,-.15) to node[auto] {$\bord{}{-}$} (-1.7,-.15);
	\draw[1c] (1.65,.15) to node[auto,swap] {$\bord{}{+}$} (.3,.15);
	\draw[1c] (1.65,-.15) to node[auto] {$\bord{}{-}$} (.3,-.15);
\end{tikzpicture}
\end{equation}
where $X_2$ contains all possible composite 2-cells, then extend it to
\begin{equation} \label{eq:2globular}
\begin{tikzpicture}[baseline={([yshift=-.5ex]current bounding box.center)}]
	\node[scale=1.25] (0) at (-2,0) {$X_0$};
	\node[scale=1.25] (1) at (0,0) {$X_1$};
	\node[scale=1.25] (2) at (2,0) {$X_2$};
	\node[scale=1.25] (3) at (4.3,0) {$G_3(X),$};
	\draw[1c] (-.35,.15) to node[auto,swap] {$\bord{}{+}$} (-1.7,.15);
	\draw[1c] (-.35,-.15) to node[auto] {$\bord{}{-}$} (-1.7,-.15);
	\draw[1c] (1.65,.15) to node[auto,swap] {$\bord{}{+}$} (.3,.15);
	\draw[1c] (1.65,-.15) to node[auto] {$\bord{}{-}$} (.3,-.15);
	\draw[1c] (3.65,.15) to node[auto,swap] {$\bord{}{+}$} (2.3,.15);
	\draw[1c] (3.65,-.15) to node[auto] {$\bord{}{-}$} (2.3,-.15);
\end{tikzpicture}
\end{equation}
where the boundary of the generating 3-cells has in turn been subdivided into two copies of $D^2$, and again $\bord{}{\alpha}\bord{}{+} = \bord{}{\alpha}\bord{}{-}$. 

Observe that, while the sets in (\ref{eq:2cat}) contain all composites of their elements --- so the diagram can be endowed with an internal composition operation --- this is not the case for (\ref{eq:2globular}): in the procedure, we are moving back-and-forth between two different categories, one of which is independent of the chosen notion of composition. This is the category of globular sets \cite[Section 1.4]{leinster2004higher}, defined as follows.

\begin{dfn}
Let $\cat{G}$ be the category generated by the diagram \index{G@$\cat{G}$}
\begin{equation*}
\begin{tikzpicture}[baseline={([yshift=-.5ex]current bounding box.center)}]
	\node[scale=1.25] (0) at (-2,0) {$0$};
	\node[scale=1.25] (1) at (0,0) {$1$};
	\node[scale=1.25] (2) at (2.2,-.1) {$\ldots$};
	\node[scale=1.25] (3) at (4.4,0) {$n$};
	\node[scale=1.25] (4) at (6.6,-.1) {$\ldots,$};
	\draw[1c] (-1.7,.15) to node[auto] {$\imath_+$} (-.35,.15);
	\draw[1c] (-1.7,-.15) to node[auto,swap] {$\imath_-$} (-.35,-.15);
	\draw[1c] (.3,.15) to node[auto] {$\imath_+$} (1.65,.15);
	\draw[1c] (.3,-.15) to node[auto,swap] {$\imath_-$} (1.65,-.15);
	\draw[1c] (2.7,.15) to node[auto] {$\imath_+$} (4.05,.15);
	\draw[1c] (2.7,-.15) to node[auto,swap] {$\imath_-$} (4.05,-.15);
	\draw[1c] (4.7,.15) to node[auto] {$\imath_+$} (6.05,.15);
	\draw[1c] (4.7,-.15) to node[auto,swap] {$\imath_-$} (6.05,-.15);
	\node[scale=1.25] (5) at (8,0) {$n \in \mathbb{N},$};
\end{tikzpicture}
\end{equation*}
with the relation $\imath_+\imath_\alpha = \imath_-\imath_\alpha$ for all consecutive pairs of morphisms. A \emph{globular set} is a presheaf $X: \cat{G}^\mathrm{op} \to \cat{Set}$; we denote it by the $\cat{Set}$-diagram \index{globular set}
\begin{equation*}
\begin{tikzpicture}[baseline={([yshift=-.5ex]current bounding box.center)}]
	\node[scale=1.25] (0) at (-2,0) {$X_0$};
	\node[scale=1.25] (1) at (0,0) {$X_1$};
	\node[scale=1.25] (2) at (2.2,-.1) {$\ldots$};
	\node[scale=1.25] (3) at (4.4,0) {$X_n$};
	\node[scale=1.25] (4) at (6.6,-.1) {$\ldots,$};
	\draw[1c] (-.35,.15) to node[auto,swap] {$\bord{}{+}$} (-1.7,.15);
	\draw[1c] (-.35,-.15) to node[auto] {$\bord{}{-}$} (-1.7,-.15);
	\draw[1c] (1.65,.15) to node[auto,swap] {$\bord{}{+}$} (.3,.15);
	\draw[1c] (1.65,-.15) to node[auto] {$\bord{}{-}$} (.3,-.15);
	\draw[1c] (4.05,.15) to node[auto,swap] {$\bord{}{+}$} (2.7,.15);
	\draw[1c] (4.05,-.15) to node[auto] {$\bord{}{-}$} (2.7,-.15);
	\draw[1c] (6.05,.15) to node[auto,swap] {$\bord{}{+}$} (4.7,.15);
	\draw[1c] (6.05,-.15) to node[auto] {$\bord{}{-}$} (4.7,-.15);
	\node[scale=1.25] (5) at (8,0) {$n \in \mathbb{N}.$};
\end{tikzpicture}
\end{equation*}
The elements of $X_n$ are called the \emph{$n$-cells} of $X$; for all cells $x$, we call $\bord{}{-}x$ the \emph{input} boundary, and $\bord{}{+}x$ the \emph{output} boundary of $x$. We also write $x: \bord{}{-}x \to \bord{}{+} x$. \index{boundary!in $\omega$-categories} \index{input|see {boundary}} \index{output|see {boundary}} \index{$\bord{}{+}, \bord{}{-}$}

Moreover, for all $n$-cells $x$, all $k < n$, and all $\alpha \in \{+,-\}$, let
\begin{equation*}
	\bord{k}{\alpha}x \; := \underbrace{(\bord{}{\alpha}\ldots\bord{}{\alpha})}_{n-k} x\,;
\end{equation*} 
we call $\bord{k}{-}x$ the input $k$-boundary, and $\bord{k}{+}x$ the output $k$-boundary of $x$. By definition, $\bord{n-1}{\alpha}x = \bord{}{\alpha}x$. 

Given two globular sets $X$, $Y$, a morphism $f: X \to Y$ is a morphism of presheaves on $\mathbf{G}$; it is given by a sequence of functions $\{f_n: X_n \to Y_n\}_{n\in \mathbb{N}}$ that commute with the boundary operations. We write $\mathbf{Glob}$ for the category of globular sets and morphisms. \index{Glob@$\cat{Glob}$}
\end{dfn}

A finite diagram of the form
\begin{equation*}
\begin{tikzpicture}[baseline={([yshift=-.5ex]current bounding box.center)}]
	\node[scale=1.25] (0) at (-2,0) {$X_0$};
	\node[scale=1.25] (1) at (0,0) {$X_1$};
	\node[scale=1.25] (2) at (2.2,-.1) {$\ldots$};
	\node[scale=1.25] (3) at (4.4,0) {$X_n$};
	\draw[1c] (-.35,.15) to node[auto,swap] {$\bord{}{+}$} (-1.7,.15);
	\draw[1c] (-.35,-.15) to node[auto] {$\bord{}{-}$} (-1.7,-.15);
	\draw[1c] (1.65,.15) to node[auto,swap] {$\bord{}{+}$} (.3,.15);
	\draw[1c] (1.65,-.15) to node[auto] {$\bord{}{-}$} (.3,-.15);
	\draw[1c] (4.05,.15) to node[auto,swap] {$\bord{}{+}$} (2.7,.15);
	\draw[1c] (4.05,-.15) to node[auto] {$\bord{}{-}$} (2.7,-.15);
\end{tikzpicture}
\end{equation*}
satisfying $\bord{}{\alpha}\bord{}{+} = \bord{}{\alpha}\bord{}{-}$ can be extended to a globular set by setting $X_m := \emptyset$ for all $m > n$.

At each step of our combinatorial construction of a CW complex, we are extending such a finite diagram with a set of higher-dimensional generators; then, we freely generate and include all their composites. To make this formally precise, we need to introduce an algebraic notion of composition for abstract $n$-cells. In \cite{batanin1998computads}, a variety of options is collected under the umbrella of \emph{finitary monads} on the category of globular sets. The standard notion of polygraph corresponds to the simplest (arguably) combinatorics, given by strict, globular $\omega$-categories.

\begin{dfn} \index{composable cells}
Let $X$ be a globular set. Given two $n$-cells $x$, $y$ of $X$, and $1 \leq k \leq n$, we say that $x$ and $y$ are \emph{$k$-composable}, and write $x \comp{k} y$, if $\bord{k-1}{+}x = \bord{k-1}{-}y$.

We write $X_n \comp{k} X_n \subseteq X_n \times X_n$ for the set of pairs of $k$-composable $n$-cells of $X$.
\end{dfn}

\begin{dfn} \index{omega-category@$\omega$-category} \index{omegacat@$\omegacat$}
An \emph{$\omega$-category} is a globular set $X$, together with \emph{unit} and \emph{$k$-composition} operations
\begin{equation*}
	\idd{}: X_n \to X_{n+1}, \qquad \cp{k}: X_n \comp{k} X_n \to X_n
\end{equation*}
for all $n \in \mathbb{N},$ and $1 \leq k \leq n,$ satisfying the following conditions: for all $k$-cells $x$, $n > k$, let
\begin{equation*}
	\idd{n}x := \underbrace{(\idd{}\ldots\idd{})}_{n-k} x\,;
\end{equation*}
\begin{enumerate}
	\item for all $n$-cells $x$, $y$, and $k < n$, whenever the left-hand side is defined,
	\begin{align*}
		& \bord{}{\alpha}(\varepsilon x) = x, \\
		& \bord{}{-}(x \cp{n} y) = \bord{}{-} x, \quad \bord{}{+}(x \cp{n} y) = \bord{}{+} y, \\
		& \bord{}{\alpha}(x \cp{k} y) = \bord{}{\alpha} x \cp{k} \bord{}{\alpha} y\,;
	\end{align*}
	\item whenever the left-hand side is defined, for all $n$, and $k < n$,
	\begin{align*}
		& \idd{}(x \cp{n} y) = \idd{}x \cp{n} \idd{}y, \\
		& (x \cp{n} y) \cp{n} z = x \cp{n} (y \cp{n} z), & \text{(associativity)}, \\
		& x \cp{k} \idd{n}(\bord{k}{+}x) = x, \quad \idd{n}(\bord{k}{-}x) \cp{k} x = x, & \text{(unitality)}, \\
		& (x \cp{n} x') \cp{k} (y \cp{n} y') = (x \cp{k} y) \cp{n} (x' \cp{k} y'), & \text{(interchange)}.
	\end{align*}
\end{enumerate} 

Given two $\omega$-categories $X$, $Y$, a \emph{map} $f: X \to Y$ is a morphism of the underlying globular sets that commutes with units and compositions. A map $f$ is an \emph{inclusion}, written $f: X \hookrightarrow Y$, if each component $f_n$ is an injection of sets. There is a category $\omegacat$ whose objects are $\omega$-categories, and morphisms are maps. \index{inclusion!of $\omega$-categories}
\end{dfn}

\begin{dfn}
Let $X$ be an $\omega$-category, $n \in \mathbb{N}$, and $x$ be an $n$-cell of $X$. We define a number $\dmn{x}$, the \emph{dimension} of $x$, by induction on $n$: \index{dimension!in $\omega$-categories}
\begin{itemize}
	\item if $x \in X_0$, $\dmn{x} := 0$;
	\item if $x \in X_n$, if $x = \idd{}y$ for some $y$, then $\dmn{x} := \dmn{y}$, else $\dmn{x} := n$.
\end{itemize}
\end{dfn}

\begin{itemize}
	\item[] \textbf{Terminology.} Let $X$ be an $\omega$-category. When we say that $x$ is an \emph{$n$-cell} of $X$, we mean $x \in X_n$. When we say that $x$ is \emph{$n$-dimensional}, we mean $\dmn{x} = n$. 
\end{itemize}

The unit operations are only a combinatorial device to allow the composition of $n$-dimensional and $k$-dimensional cells when $n, k$ do not coincide: for all $n$, the symbols $x$ and $\idd{n} x$ should be seen as names for the \emph{same} $k$-dimensional cell, where $k = \dmn{x}$. Graphically, $k$-compositions are depicted as pastings of cells along their $(k-1)$-boundary, and units as the lower-dimensional cells they originate from, by \emph{pasting diagrams}: \index{pasting diagram}
\begin{equation} \label{eq:compositionpast}
\input{img/c1_1_composition_axioms.tex}
\end{equation} 
\begin{equation}
\input{img/c1_1_unitality.tex}
\end{equation} 
\begin{equation} \label{eq:interpast}
\input{img/c1_1_interchange.tex}
\end{equation} 
The fact that we can univocally interpret pasting diagrams of more than two cells, such as the middle terms of equations (\ref{eq:compositionpast}) and (\ref{eq:interpast}), is a consequence of the associativity, unitality and interchange axioms; we can see these axioms as ``geometric'' equations, ensuring that $\omega$-categories are a reasonable algebra of composition. On the other hand, a generic $\omega$-category will also have ``algebraic'', or ``non-geometric'' equations between pairs of its cells, that is, ones that are not derivable from the axioms of composition.

There is a functor $U: \omegacat \to \mathbf{Glob}$, taking an $\omega$-category to its underlying globular set. As shown in \cite{batanin1998monoidal}, this functor is finitary monadic, which implies that $\omegacat$ has all limits and colimits, because $\mathbf{Glob}$, a topos of presheaves, does. We refer to \cite{batanin1998monoidal, berger2002cellular} for some combinatorial descriptions of the monad on $\mathbf{Glob}$ whose algebras are $\omega$-categories. 

\begin{dfn} \label{dfn:skeleton} \index{skeleton!of an $\omega$-category} \index{truncation}
Let $n \in \mathbb{N}$. An \emph{$n$-category} is an $\omega$-category $X$ such that $X_m$ is the image of $X_{m-1}$ through $\idd{}$ for all $m > n$. Since all the information about $X$ is contained in the $X_k$ with $k \leq n$, we represent it as a diagram ending at $X_n$.

For each $n$, $n$-categories and their maps form a full subcategory $\imath: n\cat{Cat} \hookrightarrow \omegacat$; the inclusion functors have both a left and a right adjoint, 
\begin{equation*}
	\tau_{\leq n} \dashv \imath \dashv \sigma_{\leq n},
\end{equation*}
called the \emph{$n$-truncation} and the \emph{$n$-skeleton}, respectively. Explicitly, for $X$ an $\omega$-category, and $n \in \mathbb{N}$, the $n$-category $\sigma_{\leq n} X$ has the underlying globular set
\begin{equation*}
\begin{tikzpicture}[baseline={([yshift=-.5ex]current bounding box.center)}]
	\node[scale=1.25] (0) at (-2,0) {$X_0$};
	\node[scale=1.25] (1) at (0,0) {$X_1$};
	\node[scale=1.25] (2) at (2.2,-.1) {$\ldots$};
	\node[scale=1.25] (3) at (4.4,0) {$X_n$};
	\draw[1c] (-.35,.15) to node[auto,swap] {$\bord{}{+}$} (-1.7,.15);
	\draw[1c] (-.35,-.15) to node[auto] {$\bord{}{-}$} (-1.7,-.15);
	\draw[1c] (1.65,.15) to node[auto,swap] {$\bord{}{+}$} (.3,.15);
	\draw[1c] (1.65,-.15) to node[auto] {$\bord{}{-}$} (.3,-.15);
	\draw[1c] (4.05,.15) to node[auto,swap] {$\bord{}{+}$} (2.7,.15);
	\draw[1c] (4.05,-.15) to node[auto] {$\bord{}{-}$} (2.7,-.15);
\end{tikzpicture}
\end{equation*}
simply obtained by interrupting the sequence of the $X_m$ at $n$, and the same compositions and units as $X$; whereas $\tau_{\leq n} X$ has the underlying $n$-globular set 
\begin{equation*}
\begin{tikzpicture}[baseline={([yshift=-.5ex]current bounding box.center)}]
	\node[scale=1.25] (0) at (-2,0) {$X_0$};
	\node[scale=1.25] (1) at (0,0) {$X_1$};
	\node[scale=1.25] (2) at (2.2,-.1) {$\ldots$};
	\node[scale=1.25] (3) at (4.5,0) {$X_n^\tau,$};
	\draw[1c] (-.35,.15) to node[auto,swap] {$\bord{}{+}$} (-1.7,.15);
	\draw[1c] (-.35,-.15) to node[auto] {$\bord{}{-}$} (-1.7,-.15);
	\draw[1c] (1.65,.15) to node[auto,swap] {$\bord{}{+}$} (.3,.15);
	\draw[1c] (1.65,-.15) to node[auto] {$\bord{}{-}$} (.3,-.15);
	\draw[1c] (4.05,.15) to node[auto,swap] {$\bord{}{+}$} (2.7,.15);
	\draw[1c] (4.05,-.15) to node[auto] {$\bord{}{-}$} (2.7,-.15);
\end{tikzpicture}
\end{equation*}
where $X_n^\tau$ is the coequaliser of 
\begin{equation*}
\begin{tikzpicture}[baseline={([yshift=-.5ex]current bounding box.center)}]
	\node[scale=1.25] (0) at (-2,0) {$X_n$};
	\node[scale=1.25] (1) at (0.2,0) {$X_{n+1}\,.$};
	\draw[1c] (-.35,.15) to node[auto,swap] {$\bord{}{+}$} (-1.7,.15);
	\draw[1c] (-.35,-.15) to node[auto] {$\bord{}{-}$} (-1.7,-.15);
\end{tikzpicture}
\end{equation*}
The coequaliser induces well-defined composition operations on $X_n^\tau$, since for all $n$-cells $x$, $x'$, $y$, $y'$, if $x \sim x'$ is witnessed by $z: x \to x'$,
\begin{equation*}
	y \comp{k} x \comp{k} y' \quad \text{ if and only if } \quad y \comp{k} x' \comp{k} y',
\end{equation*}
and in that case the $(n+1)$-cell $\idd{} y \cp{k} z \cp{k} \idd{} y'$ is a witness of $y \cp{k} x \cp{k} y' \sim y \cp{k} x' \cp{k} y'$.

For all $\omega$-categories $X$, we will also use $n$-truncation of $X$ and $n$-skeleton of $X$ to denote the unit $X \to \imath \tau_{\leq n}X$ and the counit $\imath \sigma_{\leq n}X \hookrightarrow X$, respectively, in $\omegacat$. We will write $\skel{n}{X} := \imath \sigma_{\leq n} X$.
\end{dfn}

\begin{remark}
Our notation and use of the term ``truncation'' is consistent with \cite{kapranov1991infty, simpson2009homotopy}; on the other hand, \cite{batanin1998computads} uses it to refer to our $\sigma_{\leq n}$. 
\end{remark}

By definition, if $X$ is an $n$-category, any map $f: X \to Y$ factors uniquely through a map $X \to \skel{n}{Y}$, which we still call $f$.

For all $n \in \mathbb{N}$, the $n$-skeleton $\skel{n}{X} \hookrightarrow X$ factors through the $(n+1)$-skeleton, and $X$ is a colimit of the sequence of the $\skel{n}{X} \hookrightarrow \skel{n+1}{X}$. Hence, any $\omega$-category can be described by the sequence of its $n$-skeleta. 

\begin{cons}\label{cons:freeomega} \index{omega-category@$\omega$-category!free}
The left adjoint $F: \mathbf{Glob} \to \omegacat$ of $U$, taking a globular set $X$ to the \emph{free} $\omega$-category on $X$, can be inductively described as follows. Given a globular set $X$, let $\skel{0}{(FX)} := X_0$. Suppose $\skel{n-1}{(FX)}$ has been defined, together with boundary-preserving inclusions $X_k \subseteq FX_k$, for all $k < n$. Define $FX^\star_{n}$, together with boundary functions $\bord{}{\alpha}: FX^\star_n \to FX_{n-1}$, as the set freely generated by the following clauses: 
\begin{itemize}
	\item for all $x \in X_n$, $x: \bord{}{-}x \to \bord{}{+}x \in FX^\star_n$;
	\item for all $y \in FX_{n-1}$, $\epsilon (y): y \to y \in FX^\star_n$;
	\item for all $(x, y) \in FX^\star_n \comp{n} FX^\star_n$, $(x \star_n y): \bord{}{-}x \to \bord{}{+} y \in FX^\star_n$;
	\item for all $k < n$, and $(x, y) \in FX^\star_n \comp{k} FX^\star_n$, $(x \star_k y) : \bord{}{-}x \cp{k} \bord{}{-}y \to \bord{}{+}x \cp{k} \bord{}{+}y \in FX^\star_n$.
\end{itemize}
Define $FX_n$ to be the quotient $FX_n^\star/\sim$, where $\sim$ is generated by
\begin{align*}
	& \epsilon(x \cp{j} y) \sim (\epsilon(x) \star_j \epsilon(y)), \\
	& ((x \star_n y) \star_n z) \sim (x \star_n (y \star_n z)),  \\
	& (x \star_k \underbrace{\epsilon(\idd{}\ldots\idd{}}_{n-k}y)) \sim x \sim (\underbrace{\epsilon(\idd{}\ldots\idd{}}_{n-k}z) \star_k x), \\
	& ((x \star_j x') \star_k (y \star_j y')) \sim ((x \star_k y) \star_j (x' \star_k y')),
\end{align*}
for all compatible $x, x', y, y', z$, and all $j < k \leq n$. Compositions and units are the obvious ones.

Given a morphism $f: X \to Y$ of globular sets, $Ff: FX \to FY$ is the map whose components $Ff_n$ are extended from the $f_n$ by induction on the syntax of elements of $FX_n$, $FY_n$, in the obvious way.
\end{cons}

Many simple, yet important examples of $\omega$-categories are free on a globular set, including the simplest version of a combinatorial, oriented $n$-disk.
\begin{dfn} \index{globe!standard}
For all $n \in \mathbb{N}$, the \emph{standard $n$-globe} is defined as the $n$-category $G^n := F(\mathrm{Hom}_\mathbf{G}(-,n))$, free on the Yoneda embedding of the object $n$ of $\cat{G}$ into $\mathbf{Glob}$.
\end{dfn}
More explicitly, $G^n$ is the $n$-category with a single $n$-dimensional cell $\top$, and, for all $k < n$, exactly two $k$-dimensional cells $k^+$, $k^-$, with $\bord{k}{\alpha}(\top) = k^\alpha$. We write $\partial G^n$ for the $(n-1)$-skeleton of $G^n$, which results from removing its $n$-dimensional cell. 

From the Yoneda lemma and the adjunction $F \dashv U$, it follows that maps $x: G^n \to X$ are in bijection with cells $x \in X_n$. Maps $\partial G^n \to X$, on the other hand, are in bijection with pairs $x^-, x^+ \in X_{n-1}$ such that $\bord{}{\alpha}x^- = \bord{}{\alpha}x^+$, that is, pairs of $(n-1)$-cells that can form the boundary of an $n$-cell. For all $x: G^n \to X$, there is a map $\partial x : \partial G^n \to X$, uniquely determined by the fact that it maps $k^\alpha$ to $\bord{k}{\alpha} x$, for all $k < n$.

\begin{dfn} \index{union!disjoint}
Given a family $\{X_i\}_{i \in I}$ of $\omega$-categories, its \emph{disjoint union} $\coprod_{i \in I} X_i$ is its coproduct in $\omegacat$.
\end{dfn}

There is nothing unexpected about disjoint unions: sets of $n$-cells and operations are all induced ``pointwise'' by coproducts in $\mathbf{Set}$. We write $X + Y$ for the disjoint union of two $\omega$-categories. 

We now have all the ingredients for a formal definition of polygraphs, modelled on the definition of a CW complex.

\begin{dfn} \index{polygraph} \index{generator}
A \emph{polygraph} is an $\omega$-category $X$ together with, for all $n > 0$, a subset $G_n(X)$ of $X_n$, such that
\begin{equation*}
\begin{tikzpicture}[baseline={([yshift=-.5ex]current bounding box.center)}]
	\node[scale=1.25] (0) at (0,2) {$\displaystyle \coprod_{G_n(X)} \partial G^n$};
	\node[scale=1.25] (1) at (3.5,2) {$\displaystyle\coprod_{G_n(X)} G^n$};
	\node[scale=1.25] (2) at (0,0) {$\skel{n-1}{X}$};
	\node[scale=1.25] (3) at (3.5,0) {$\skel{n}{X}$};
	\draw[1c] (0) to node[auto,swap] {$\displaystyle\coprod_{G_n(X)} \partial x$} (2);
	\draw[1c] (1) to node[auto] {$\displaystyle\coprod_{G_n(X)} x$} (3);
	\draw[1cinc] (0) to (1);
	\draw[1cinc] (2) to (3);
	\draw[edge] (2.5,0.2) to (2.5,0.8) to (3.3,0.8);
\end{tikzpicture}
\end{equation*}
is a pushout diagram in $\omegacat$. The cells of $G_n(X)$ are called the \emph{$n$-dimensional generators} of $X$. A polygraph $X$ is an \emph{$n$-polygraph} if it is an $n$-category.
\end{dfn}
By the essential uniqueness of pushouts, given $\skel{n-1}{X}$ and the maps $\{\partial x: \partial G^n \to \skel{n-1}{X}\}_{x \in G_n(X)}$ --- or equivalently, the graph
\begin{equation*}
\begin{tikzpicture}[baseline={([yshift=-.5ex]current bounding box.center)}]
	\node[scale=1.25] (0) at (-2.2,0) {$X_{n-1}$};
	\node[scale=1.25] (1) at (0.3,0) {$G_n(X),$};
	\draw[1c] (-.35,.15) to node[auto,swap] {$\bord{}{+}$} (-1.7,.15);
	\draw[1c] (-.35,-.15) to node[auto] {$\bord{}{-}$} (-1.7,-.15);
\end{tikzpicture}
\end{equation*}
with $\bord{}{\alpha}x := \partial x((n-1)^\alpha$) --- we can reconstruct $\skel{n}{X}$ up to isomorphism. This allows us to give an inductive reconstruction of $X$, by free syntactic extensions in the style of Construction \ref{cons:freeomega}.
\begin{cons} \label{cons:polygraph}
Let $\skel{0}{X} := G_0(X)$. Suppose $\skel{n-1}{X}$ has been defined, and we have maps $\{\partial x: \partial G^n \to \skel{n-1}{X}\}_{x \in G_n(X)}$. Define $X_n^\star$, together with boundary functions $\bord{}{\alpha}: X^\star_n \to X_{n-1}$, as the set freely generated by the following clauses:
\begin{itemize}
	\item for all $x \in G_n(X)$, $x: \bord{}{-}x \to \bord{}{+}x \in X^\star_n$;
	\item for all $y \in X_{n-1}$, $\epsilon (y): y \to y \in X^\star_n$;
	\item for all $(x, y) \in X^\star_n \comp{n} X^\star_n$, $(x \star_n y): \bord{}{-}x \to \bord{}{+} y \in X^\star_n$;
	\item for all $k < n$, and $(x, y) \in X^\star_n \comp{k} X^\star_n$, $(x \star_k y) : \bord{}{-}x \cp{k} \bord{}{-}y \to \bord{}{+}x \cp{k} \bord{}{+}y \in X^\star_n$.
\end{itemize}
Then $X_n$ is the quotient $X_n^\star/\sim$, where $\sim$ is the same as in Construction \ref{cons:freeomega}.
\end{cons}
We can immediately see, from this description, that a free $\omega$-category on a globular set is a special type of polygraph, where boundaries of generators are generators.

As we mentioned, unlike polygraphs, general $\omega$-categories have ``non-geometric'' equations between cells. Nevertheless, we can replace, in a suitable sense, any $\omega$-category with a polygraph: the idea is to treat each non-geometric equation as arising from the truncation of some higher-dimensional cell, in the sense of Definition \ref{dfn:skeleton}, and try to ``invert'' the truncation by adding enough cells. 

This leads to the concept of \emph{polygraphic resolution}, introduced by M\'etayer in \cite{metayer2003resolutions}.

\begin{dfn} \index{resolution}
A map $p: E \to X$ of $\omega$-categories is a \emph{trivial fibration} if it satisfies the \emph{right-lifting property} with respect to the set $\{\partial G^n \hookrightarrow G^n \}_{n \in \mathbb{N}}$: that is, for all commutative diagrams
\begin{equation*}
\begin{tikzpicture}[baseline={([yshift=-.5ex]current bounding box.center)}]
	\node[scale=1.25] (0) at (0,2) {$\partial G^n$};
	\node[scale=1.25] (1) at (2,2) {$E$};
	\node[scale=1.25] (2) at (0,0) {$G^n$};
	\node[scale=1.25] (3) at (2,0) {$X$,};
	\draw[1c] (0) to node[auto] {$f$} (1);
	\draw[1c] (2) to node[auto,swap] {$x$} (3);
	\draw[1cinc] (0) to (2);
	\draw[1c] (1) to node[auto] {$p$} (3);
\end{tikzpicture}
\end{equation*}
there exists a map $\tilde{x}: G^n \to E$ making 
\begin{equation*}
\begin{tikzpicture}[baseline={([yshift=-.5ex]current bounding box.center)}]
	\node[scale=1.25] (0) at (0,2) {$\partial G^n$};
	\node[scale=1.25] (1) at (2,2) {$E$};
	\node[scale=1.25] (2) at (0,0) {$G^n$};
	\node[scale=1.25] (3) at (2,0) {$X$};
	\draw[1c] (0) to node[auto] {$f$} (1);
	\draw[1c] (2) to node[auto,swap] {$x$} (3);
	\draw[1cinc] (0) to (2);
	\draw[1c] (1) to node[auto] {$p$} (3);
	\draw[1c] (2) to node[auto] {$\tilde{x}$} (1);
\end{tikzpicture}
\end{equation*}
commute. A trivial fibration is called a \emph{resolution} if $E$ is a polygraph.
\end{dfn}
The notion of resolution captures the intended meaning of ``replacing equations with higher cells''. First of all, $p$ is necessarily surjective on cells of $X$. Moreover, if $x = y$ is a non-geometric equation of $n$-cells of $X$, then $x = p(\tilde{x})$, $y = p(\tilde{y})$ for some pair of necessarily distinct $n$-cells of $E$; the unit $(n+1)$-cell $\idd{}x$, considered as a cell $x \to y$, will then lift to an $(n+1)$-cell $z: \tilde{x} \to \tilde{y}$ of $E$.

\begin{cons} \label{cons:resolution} Given an $\omega$-category $X$, construct a polygraph $RX$, together with a map $\varphi_X: RX \to X$ as follows. Let $RX_0 := X_0$, and $\varphi_{X,0}: RX_0 \to X_0$ be the identity. 

Suppose $\skel{n}{RX}$ and $\skel{n}{\varphi_X}: \skel{n}{RX} \to \skel{n}{X}$ have been defined. Let the set $G_{n+1}(RX)$ of $(n+1)$-generators of $RX$ be 
\begin{equation*}
	\{(x, x^-, x^+) \,|\, x: \varphi_{X,n}(x^-) \to \varphi_{X,n}(x^+)\} \subseteq X_{n+1} \times RX_{n} \times RX_{n},
\end{equation*}
with boundary functions $\bord{}{\alpha}(x,x^-,x^+) := x^\alpha$. Then, define $\varphi_{X,n+1}$ to be the free extension of the first projection $G_{n+1}(RX) \to X_{n+1}, (x,x^-, x^+) \mapsto x$ to $RX_{n+1}$.

By construction, $\varphi: RX \to X$ is a resolution of $X$.
\end{cons}
\begin{dfn} \index{resolution!standard}
For all $\omega$-categories $X$, the map $\varphi_X: RX \to X$ is called the \emph{standard resolution} of $X$.
\end{dfn}

As shown in \cite[Section 4]{metayer2003resolutions}, $R$ can be extended to a functor $\omegacat \to \mathbf{Pol}$, where $\mathbf{Pol}$ is the category whose objects are polygraphs, and maps $f: X \to Y$ are maps of the underlying $\omega$-categories that map generators to generators, that is, are induced by sequences of functions $\{f_n: G_n(X) \to G_n(Y)\}_{n \in \mathbb{N}}$. The functor $R$ is right adjoint to the forgetful functor $\mathbf{Pol} \to \omegacat$, and the standard resolution is in fact the counit of this adjunction, shown to be monadic in \cite{metayer2016strict}.

In \cite{lafont2010folk}, Lafont, M\'etayer, and Worytkiewicz incorporated the definition of trivial fibrations into the wider context of a model structure on $\omegacat$. In the model structure, a trivial fibration is a particular weak equivalence, in the following sense.

\begin{dfn} \label{dfn:equivalence} \index{equivalence!in $\omega$-categories} \index{weak inverse}
Let $X$ be an $\omega$-category. The class of \emph{equivalences} in $X$ is defined coinductively as follows: an $n$-cell $f: x \to x'$ is an equivalence if there exists an $n$-cell $\bar{f}: x' \to x$, and equivalence $(n+1)$-cells $e: f \cp{n} \bar{f} \to \idd{} x$, $h: \idd{} x' \to \bar{f} \cp{n} f$. If $f$ is an equivalence, then any $n$-cell $\bar{f}$ satisfying the definition is called a \emph{weak inverse} of $f$.

A map $f: X \to Y$ of $\omega$-categories is a \emph{weak equivalence} if, for all $n$-cells $x, x'$ of $X$ such that either $n=0$, or $n>0$ and $\bord{}{\alpha} x = \bord{}{\alpha} x'$, and for all $(n+1)$-cells $y': f(x) \to f(x')$ of $Y$, there exists an $(n+1)$-cell $y: x \to x'$ of $X$ and an equivalence $e_y: f(y) \to y'$. Two $\omega$-categories $X, Y$ are \emph{weakly equivalent} if there exists a weak equivalence $f: X \to Y$.
\end{dfn} \index{equivalence!weak}

In this context, one obtains an elegant characterisation of polygraphs, parallel to the characterisation of CW complexes in the Quillen-Serre model structure on topological spaces, as shown by M\'etayer in \cite{metayer2008cofibrant}.
\begin{dfn}
An $\omega$-category $X$ is \emph{cofibrant} if, for all trivial fibrations $p: E \to Y$, and all maps $f: X \to Y$, there exists a map $\tilde{f}: X \to E$ making the following diagram commute:
\begin{equation*}
\begin{tikzpicture}[baseline={([yshift=-.5ex]current bounding box.center)}]
	\node[scale=1.25] (1) at (2,2) {$E$};
	\node[scale=1.25] (2) at (0,0) {$X$};
	\node[scale=1.25] (3) at (2,0) {$Y$.};
	\draw[1c] (2) to node[auto,swap] {$f$} (3);
	\draw[1c] (1) to node[auto] {$p$} (3);
	\draw[1c] (2) to node[auto] {$\tilde{f}$} (1);
\end{tikzpicture}
\end{equation*}
\end{dfn}
\begin{thm}
An $\omega$-category $X$ is cofibrant if and only if it admits the structure of a polygraph.
\end{thm}
\begin{proof}
This is \cite[Theorem 7.4]{metayer2008cofibrant}.
\end{proof}
Being a trivial fibration, a resolution is a weak equivalence for the model structure, in particular, for all $X$, the polygraph $RX$ can be taken as a \emph{cofibrant replacement} of $X$. We will not delve much further into these aspects; we refer to \cite{riehl2009categorical} for an introduction to the general theory of model categories, and to the original source \cite{lafont2010folk} for more details on this example.

We now have a notion of generalised space made of directed cells. However, we still lack one fundamental feature of topological spaces: the ability to reason \emph{compositionally} --- to assemble smaller spaces into a larger one, whose features depend on those of its components, and the operations by which they have been put together. So far, we have only used disjoint unions, analogous to their topological counterpart, and not particularly interesting. 

Far more interesting is the product of topological spaces, on which several important constructions rely. This has, indeed, a counterpart in $\omegacat$, which is \emph{not} the categorical product --- the latter we denote by $X \with Y$, to avoid confusion: this is the $\omega$-category whose $n$-cells are pairs $(x, y)$ of an $n$-cell of $X$, and an $n$-cell of $Y$, and all operations are defined pointwise. It is, instead, a different, non-symmetric monoidal structure, which we define in the next sections. \index{product!of $\omega$-categories}

\section{Loop-free polygraphs and the tensor product} \label{sec:loopfree}
The tensor product of $\omega$-categories, variably called \emph{Crans-Gray} or \emph{lax Gray} tensor product, has been approached in different ways in the literature. 

Arguably, the simplest definition is obtained when $\omega$-categories are defined as \emph{cubical sets}, rather than globular sets with structure; see \cite{grandis2003cubical} for a survey. In \cite{al2002multiple}, Al-Agl, Brown, and Steiner proved that so-called cubical $\omega$-categories with connections are equivalent to globular $\omega$-categories, yet the equivalence is non-trivial, and by working directly in the cubical picture, one pays the price of simplicity in defining the tensor product with complexity in handling redundant ``thin cells'', see \cite{steiner2006thin}. A direct, yet highly complex globular definition was given by Crans in \cite{crans1995pasting}, and another, simplicial in nature, by Verity in \cite{verity2008complicial}. 

Instead, we will follow Steiner's \cite{steiner2004omega}: in brief, we will define a subcategory $\loopfree$ of $\omegacat$ on which the tensor product can be readily calculated, and such that any $\omega$-category is the colimit of a diagram in $\loopfree$; then extend the tensor product along such colimits to the whole of $\omegacat$. 

The main point is that the tensor product can be calculated directly on those polygraphs that can be uniquely reconstructed from their (oriented) \emph{incidence poset}, that is, from the assignment, to each generator, of the set of generators lying in its input and output boundary. Presentations of $\omega$-categories in this style have been studied in \cite{johnson1989combinatorics, power1991pasting, street1991parity, steiner1993algebra, crans1997presentations}, among others.

In \cite{kapranov1991combinatorial}, Kapranov and Voevodsky observed that the same information could conveniently be encapsulated in a \emph{chain complex}: replace the set of $n$-cells lying in the boundary of an $(n+1)$-cell with their sum in the free abelian group $\mathbb{Z}X_{n}$ on the set of all $n$-cells, the sign of each summand specifying whether the cell is in the input or in the output boundary. This is the formalism adopted by Steiner, and the one we will use. \index{abelian group!free}

The tensor product of $\omega$-categories has never been explicitly considered in relation to polygraphs and higher-dimensional rewriting (although M\'etayer makes implicit use of a special instance --- a cylinder construction --- in \cite{metayer2003resolutions}). In this section, we connect these strands of literature and fill in some gaps, working towards the main new result, Theorem \ref{thm:tensorpoly}: the tensor product of two polygraphs $X, Y$ is a polygraph, whose $n$-generators are pairs of a $k$-generator of $X$ and an $(n-k)$-generator of $Y$, for $0 \leq k \leq n$.

\begin{dfn} \index{augmented directed complex} \index{adc@$\adc$}
Let $X = \coprod_{n \in \mathbb{N}} X_n$ be a finite graded set, $\dmn{x} := n$ if $x \in X_n$. An \emph{augmented directed complex} $\augm{X}$ with basis $X$ is a chain complex, bounded below, of free abelian groups
\begin{equation*}
\begin{tikzpicture}[baseline={([yshift=-.5ex]current bounding box.center)}]
	\node[scale=1.25] (0) at (-2,0) {$\mathbb{Z}$};
	\node[scale=1.25] (1) at (0.2,0) {$\mathbb{Z}X_0$};
	\node[scale=1.25] (2) at (2.4,-.1) {$\ldots$};
	\node[scale=1.25] (3) at (4.75,0) {$\mathbb{Z}X_n$};
	\node[scale=1.25] (4) at (7,-.1) {$\ldots,$};
	\draw[1c] (-.35,0) to node[auto,swap] {$e$} (-1.7,0);
	\draw[1c] (2.05,0) to node[auto,swap] {$d$} (.7,0);
	\draw[1c] (4.2,0) to node[auto,swap] {$d$} (2.85,0);
	\draw[1c] (6.6,0) to node[auto,swap] {$d$} (5.25,0);
	\node[scale=1.25] (5) at (8.4,0) {$n \in \mathbb{N}.$};
\end{tikzpicture}
\end{equation*}
A morphism $f: \augm{X} \to \augm{Y}$ of augmented directed complexes is a morphism of chain complexes such that $f(\mathbb{N}X_n) \subseteq \mathbb{N}Y_n$ for all $n \in \mathbb{N}$. Augmented directed complexes and their morphisms form a category $\adc$.
\end{dfn}

\begin{cons} There is a functor $P: \adc \to \omegacat$, defined as follows. Given an augmented directed complex $\augm{X}$, $P\augm{X}$ is the $\omega$-category whose set of $n$-cells is
\begin{equation*}
	P\augm{X}_n := \{(x_0^-, x_0^+, \ldots, x_{n-1}^-, x_{n-1}^+, x_n) \,|\, x_k^\alpha \in \mathbb{N}X_k, \, e x_0^\alpha = 1, \text{ and } d x_{k+1}^\alpha = x_k^+ - x_k^- \},
\end{equation*}
where $x_n^\alpha := x_n$, boundary and unit operations are defined by
\begin{align*}
	\bord{}{\alpha}(x_0^-, x_0^+, \ldots, x_{n-1}^-, x_{n-1}^+, x_n) & := (x_0^-, x_0^+, \ldots, x_{n-2}^-, x_{n-2}^+, x_{n-1}^\alpha), \\
	\idd{}(x_0^-, x_0^+, \ldots, x_{n-1}^-, x_{n-1}^+, x_n) & := (x_0^-, x_0^+, \ldots, x_n, x_n, 0),
\end{align*}
and $k$-composition is defined, for $k$-composable $n$-cells $x = (x_0^-, x_0^+, \ldots, x_{n-1}^-, x_{n-1}^+, x_n)$ and $y = (y_0^-, y_0^+, \ldots, y_{n-1}^-, y_{n-1}^+, y_n)$, by
\begin{equation*}
	x \cp{k} y := (x_0^-, y_0^+, \ldots, x_{k-1}^-, y_{k-1}^+, x_k^- + y_k^-, x_k^+ + y_k^+, \ldots, x_n + y_n).
\end{equation*}
It is straightforward to verify that $P\augm{X}$ satisfies the axioms of an $\omega$-category. Moreover, any morphism $f: \augm{X} \to \augm{Y}$ of augmented directed complexes induces a map $Pf: P\augm{X} \to P\augm{Y}$ of $\omega$-categories, simply by
\begin{equation*}
	Pf(x_0^-, x_0^+, \ldots, x_{n-1}^-, x_{n-1}^+, x_n) := (f(x_0^-), f(x_0^+), \ldots, f(x_{n-1}^-), f(x_{n-1}^+), f(x_n)).
\end{equation*}
\end{cons}

The idea, of course, is that $x_k^\alpha$ should be the sum of all $k$-dimensional generators in the input or output $k$-boundary of an $n$-cell. Not all augmented directed complexes, however, admit such a combinatorial interpretation: for instance, given a cell of $P\augm{X}$, $x_0^\alpha$ could be a sum of more than one element of $X_0$, while any cell in a polygraph should have a single generator in its input and output $0$-boundary. 

Nevertheless, by restricting to a full subcategory of $\adc$, we can make sure that such situations do not arise, as will be described in what follows and defined in Definition \ref{dfn:loopfree}. 

\begin{dfn} \index{covering relation}
Let $X$ be a finite poset with order relation $\leq$. For all elements $x, y \in X$, we say that $y$ \emph{covers} $x$ if $x < y$ and, for all $y' \in X$, if $x < y' \leq y$, then $y' = y$. 

The directed graph $HX$ with $X$ as set of vertices, and an edge $c_{y,x}: y \to x$ for all pairs $y, x$ such that $y$ covers $x$, is called the \emph{Hasse diagram} of $X$. We can reconstruct $X$ from $HX$ by taking its transitive closure. \index{Hasse diagram}
\end{dfn}

Given an augmented directed complex $\augm{X}$, we can make $X$ into a poset --- the \emph{incidence poset} of $\augm{X}$ --- as follows.
\begin{dfn} \label{dfn:incidence} \index{incidence poset}
Let $(X,d)$ be an augmented directed complex, and $y \in \mathbb{Z}X_n$; there exist unique integers $\{m_y(x)\}_{x \in X_n}$ such that $y = \sum_{x \in X_n} m_y(x)\,x$. Given an element $x \in X_n$, we write $x \in y$ if $m_y(x) \neq 0$.

Let $HX$ be the following directed graph: the vertices of $HX$ are the elements of $X$; for all $x, y \in X$, there is an edge $c_{y,x}: y \to x$ whenever $x \in dy$. The transitive closure of $HX$ defines an order relation on $X$, whose Hasse diagram is $HX$.
\end{dfn}

For all elements $x \in \mathbb{Z}X_n$, $n > 0$, there is a unique decomposition of the form
\begin{equation*}
	d x = d^+ x - d^- x,
\end{equation*}
such that $d^\alpha x \in \mathbb{N}X_{n-1}$, $\alpha \in \{+,-\}$, and for all generators $y$, if $y \in d^\alpha x$ then $y \notin d^{-\alpha} x$. For all $k < n$, we also define
\begin{equation*}
	d_k^\alpha x \; := \underbrace{(d^\alpha\ldots d^\alpha)}_{n-k} x.
\end{equation*}
Since $dd = 0$, we have $d^+(d^+x) - d^-(d^+x) = d(d^+ x) = d(d^- x) = d^+(d^-x) - d^-(d^-x)$; by uniqueness of the decomposition, it follows that $d^\alpha(d^+ x) = d^\alpha(d^- x)$. By induction, it is always the case that $d(d_k^+ x) = d(d_k^- x) = d_{k-1}^+ x - d_{k-1}^- x$.

We can add information about the decomposition of boundaries into input and output to the Hasse diagram of $X$.

\begin{dfn} \label{dfn:orientation} \index{orientation}
Let $X$ be a finite poset. An \emph{orientation} on $X$ is a labelling of edges of $HX$ with elements of $\{+,-\}$, that is, a function $o: HX_1 \to \{+,-\}$, where $HX_1$ is the set of edges of $HX$. 

Given an orientation $o$ on $X$, we define a directed graph $HX^o$ as follows: the vertices of $HX^o$ are the elements of $X$, and for all $c_{y,x}: y \to x$ in $HX$,
\begin{itemize}
	\item if $o(c_{y,x}) = +$, there is an edge $c^+_{y,x}: y \to x$ in $HX^o$;
	\item if $o(c_{y,x}) = -$, there is an edge $c^-_{y,x}: x \to y$ in $HX^o$.
\end{itemize}
\end{dfn}

\begin{dfn} \label{dfn:loopfree} \index{augmented directed complex!unital} \index{augmented directed complex!strongly loop-free} \index{adcl@$\loopfree$}
Let $\augm{X}$ be an augmented directed complex. If $X$ is the incidence poset of $\augm{X}$, it admits an orientation defined by $o(c_{y,x}) = +$ if $x \in d^+ y$, and $o(c_{y,x}) = -$ if $x \in d^- y$. We say that $\augm{X}$ is \emph{strongly loop-free} if $HX^o$ is acyclic.

We say that $\augm{X}$ is \emph{unital} if, for all $x \in X_n$,
\begin{equation*}
	 \lfgen{x} := (d_0^- x, d_0^+ x, \ldots, d_{n-1}^- x, d_{n-1}^+ x, x)
\end{equation*}
is an $n$-cell in $P\augm{X}$, or, equivalently, if, for all $\alpha$, it holds that $e(d_0^\alpha x) = 1$.

We write $\loopfree$ for the full subcategory of $\adc$ on unital, strongly loop-free augmented directed complexes.
\end{dfn}

This characterisation of loop-freeness is inspired by acyclic matchings in discrete Morse theory, see \cite[Chapter 11]{kozlov2008combinatorial}. The idea is that the orientation of boundaries induces a ``flow'', described by the graph $HX^o$, on the generators of $P\augm{X}$: from the input boundary of a cell, into the cell, into the output boundary. If $X$ is strongly loop-free, the flow never returns to a cell after leaving it.

\begin{exm}
All standard $n$-globes are strongly loop free. Let $X$ be the incidence poset of $G^2$. The following are representations of the labelled Hasse diagram of $X$, of the graph $HX^o$, and of the flow induced on cells of $G^2$:
\begin{equation*}
\input{img/c1_2_2globe_bis.tex}
\end{equation*} 
\end{exm}

\begin{thm} \label{thm:loopsubcat}
Let $P: \loopfree \to \omegacat$ be the restriction of $P$ to unital, strongly loop-free augmented directed complexes. Then $P$ is full and faithful. Moreover, for all objects $\augm{X}$ of $\loopfree$, $P\augm{X}$ admits the structure of a polygraph whose set of $n$-dimensional generators is $\{\lfgen{x} \,|\, x \in X_n\}$.
\end{thm}
\begin{proof}
Follows from the results in \cite[Section 5]{steiner2004omega}.
\end{proof}

Thus, we can treat $\loopfree$ as a full subcategory of $\omegacat$, and say that a polygraph $X$ is \emph{loop-free} if it belongs to it. We will freely employ notions defined for augmented directed complexes, such as the incidence poset, when discussing loop-free polygraphs. \index{polygraph!loop-free}

As shown in \cite{steiner2007simple}, loop-free polygraphs contain an important subcategory $\Theta$, whose objects are called simple $\omega$-categories in \cite{makkai2001duality}, and Batanin cells in \cite{joyal1997disks}. These can be described as follows. Let $X$ be a finite set, together with a function $\mathrm{dim}: X \to \mathbb{N}$ and a total ordering $(x_0, \ldots, x_m)$ of its elements, satisfying the following conditions: \index{omega-category@$\omega$-category!simple} \index{Batanin cell} \index{theta@$\Theta$}
\begin{enumerate}
	\item $\dmn{x_0} = \dmn{x_m} = 0$;
	\item $\dmn{x_{k+1}} = \dmn{x_k} \pm 1$, for all $k < m$.
\end{enumerate}
Then, $X$ becomes a graded set with $X_n := \{x \in X \,|\, \dmn{x} = n \}$, and an augmented directed complex $\augm{X}$ in the following way. Let $ex = 1$ for all $x \in X_0$, and, if $n > 0$, $x \in X_n$, and $x = x_k$ in the ordering, let
\begin{equation*}
	dx := x_{k+j} - x_{k-i},
\end{equation*}
where $i, j > 0$ are chosen so that $\dmn{x_{k+j}} = \dmn{x_{k-i}} = n-1$, and for all $i' < i, j' < j$, $\dmn{x_{k+j'}}, \dmn{x_{k-i'}} \geq n$. It is straightforward to check that $\augm{X}$ is well-defined, unital and strongly loop-free; by construction, the total ordering of elements of $X$ is the transitive closure of the graph $HX^o$.

Clearly, $\augm{X}$ is determined up to isomorphism by the sequence of natural numbers $(\dmn{x_k})_{k\leq m}$. In fact, by the properties of the sequence, it suffices to remember its local minima and maxima: for instance, the sequence
\begin{equation*}
	(0, 1, 2, 3, 2, 3, 2, 1, 0, 1, 2, 1, 0)
\end{equation*}
can be encoded as the sequence of alternating maxima and minima $(3, 2, 3, 0, 2)$. 

The elements of $\Theta$ can therefore be characterised by such sequences; they also admit a number of alternative characterisations, in terms of planar trees, see \cite{batanin1998monoidal}, or iterated wreath products, see \cite{berger2007iterated}. Their importance lies in the fact that they classify compositions of an arbitrary number of cells in an $\omega$-category, in the sense that maps $X \to Y$, where $X$ is an object of $\Theta$, are in bijection with composable diagrams of cells of $Y$. In particular, a sequence $(\dmn{x_k})_{k\leq m}$ with $p$ maxima, $x_{i_1}, \ldots, x_{i_p}$, classifies a composition of $p$ cells, of dimensions $\dmn{x_{i_j}}$, along the boundaries specified by the local minima between each consecutive pair of local maxima. 

\begin{exm} Consider the sequence $(0,1,2,1,2,1,2,1,0,1,0,1,2,1,0)$, alternatively encoded by its sequence of maxima and minima $(2,1,2,1,2,0,1,0,2)$. We can depict the corresponding loop-free polygraph as the pasting diagram
\begin{equation*}
\input{img/c1_2_thetacat.tex}
\end{equation*} 
classifying the composition $(x \cp{2} x' \cp{2} x'') \cp{1} y \cp{1} z$, for composable $2$-cells $x, x', x'', z$, and $1$-cells $y$.
\end{exm}

\begin{remark}
For all $n > 0$, the standard $n$-globe is a simple $\omega$-category, the one encoded by the sequence $(0, 1, \ldots, n, n-1, \ldots 0)$ with a global maximum at $n$.
\end{remark}

\begin{thm} \label{thm:bergercell}
The inclusion of subcategories $\imath: \Theta \hookrightarrow \omegacat$ is dense, that is, each $\omega$-category $X$ is the colimit of the projection functor $\imath/X \to \omegacat$ from the comma category $\imath/X$.
\end{thm}
\begin{proof}
This is proved in \cite[Theorem 1.12]{berger2002cellular}.
\end{proof}

\begin{remark}
In fact, Berger proved, more specifically, that $\omegacat$ is equivalent to the category of presheaves $\opp{\Theta} \to \mathbf{Set}$ satisfying a Segal-like condition; we refer to the source for more details.
\end{remark}

Since $\Theta$ is a full subcategory of $\loopfree$, it follows, \emph{a fortiori}, that each $\omega$-category $X$ is the colimit of a diagram of loop-free polygraphs. The reason for working with $\loopfree$ rather than $\Theta$ is the existence of the following monoidal structure, based on the tensor product of chain complexes \cite[Section 12.3] {may1999concise}.
\begin{dfn} \index{tensor product!of augmented directed complexes}
Let $\augm{X}, \augm{Y}$ be augmented directed complexes with bases $X, Y$. The \emph{tensor product} $\augm{X} \otimes \augm{Y}$ of $X$ and $Y$ is the augmented directed complex $\augm{X \times Y}$ with basis $X \times Y$, graded as
\begin{equation*}
	(X \times Y)_n := \coprod_{k=0}^n X_k \times Y_{n-k},
\end{equation*}
so that $\mathbb{Z}(X \times Y)_n = \bigoplus_{k=0}^n \mathbb{Z}X_k \otimes \mathbb{Z}X_{n-k}$, whose boundary morphisms are defined, for $x \in \mathbb{Z}X_k$, $y \in \mathbb{Z}Y_{n-k}$, by
\begin{align*}
	e(x \otimes y) & := ex \, ey, \qquad\qquad\qquad\qquad\qquad\;\;\; n = k = 0, \\
	d(x \otimes y) & := \begin{cases}
		dx \otimes y, & n = k > 0, \\
		x \otimes dy, & n > k = 0, \\
		dx \otimes y + (-1)^{\dmn{x}} x \otimes dy, & n > k > 0.
	\end{cases}
\end{align*}
The tensor product extends to morphisms in the obvious way, to define a functor $\adc \times \adc \to \adc$, colimit-preserving separately in each variable. Together with the unit $\augm{1}$ --- the augmented directed complex with a single 0-dimensional basis element $\bullet$, and $e\bullet = 1$ --- it defines a monoidal structure on $\adc$.
\end{dfn}

Clearly, $\augm{X} \otimes \augm{Y}$ is unital when both $\augm{X}$ and $\augm{Y}$ are. The incidence poset of $\augm{X} \otimes \augm{Y}$ is just the product of incidence posets $X \times Y$, and $H(X \times Y)$ is the cartesian product of the graphs $HX$ and $HY$ (this is not to be confused with the categorical product; see \cite[Section 2.6]{hahn2013graph}). The graph $H(X \times Y)^o$ is ``almost'' the cartesian product of $HX^o$ and $HY^o$, with the difference that copies of $HY^o$ in the fibres $\{x\} \times HY^o$ have reversed edges when $x$ is odd-dimensional; this is acyclic when $HX^o$ and $HY^o$ are. 

Since $(1,d)$ is unital and strongly loop-free, the tensor product restricts to a monoidal structure on $\loopfree$. We are now ready to define the tensor product of $\omega$-categories, as in \cite[Section 7]{steiner2004omega}.

\begin{dfn} \index{tensor product!of $\omega$-categories}
Let $X$, $Y$ be $\omega$-categories. The \emph{tensor product} $X \otimes Y$ of $X$ and $Y$ is defined, up to natural isomorphism, as the colimit in $\omegacat$ of the functor
\begin{equation*}
\begin{tikzpicture}[baseline={([yshift=-.5ex]current bounding box.center)}]
	\node[scale=1.25] (0) at (-2,0) {$P/X \times P/Y$};
	\node[scale=1.25] (1) at (2.1,0) {$\loopfree \times \loopfree$};
	\node[scale=1.25] (2) at (5.7,0) {$\loopfree$};
	\node[scale=1.25] (3) at (8.4,0) {$\omegacat$,};
	\draw[1c] (0) to node[auto] {$\pi_X \times \pi_Y$} (1);
	\draw[1c] (1) to node[auto] {$\otimes$} (2);
	\draw[1cinc] (2) to node[auto] {$P$} (3);
\end{tikzpicture}
\end{equation*}
where $\pi_X$, $\pi_Y$ are projection functors from the comma categories $P/X$, $P/Y$. 
\end{dfn}
The tensor product naturally extends to a functor $\omegacat \times \omegacat \to \omegacat$, colimit-preserving separately in each variable. By Theorem \ref{thm:loopsubcat} and Theorem \ref{thm:bergercell}, it determines a monoidal structure on $\omegacat$, extending the monoidal structure on $\loopfree$, whose unit is $1$, the terminal $\omega$-category with a single $0$-dimensional cell.

Since each $\omega$-category can be seen as a presheaf $\opp{\Theta} \to \mathbf{Set}$, and the tensor product preserves colimits, the monoidal structure is part of a biclosed monoidal structure on $\omegacat$, in a standard way: that is, for all $\omega$-categories $X$ and $Y$, $\rimp{X}{Y}$ is defined as the $\omega$-category, unique up to isomorphism, such that
\begin{equation*}
	\mathrm{Hom}(-,\rimp{X}{Y}) \simeq \mathrm{Hom}(X \otimes -, Y),
\end{equation*}
and similarly $\limp{X}{Y}$ is defined as the $\omega$-category, unique up to isomorphism, such that
\begin{equation*}
	\mathrm{Hom}(-,\limp{X}{Y}) \simeq \mathrm{Hom}(- \otimes X, Y).
\end{equation*}

It remains to be shown that the tensor product of two polygraphs is still a polygraph. Conceptually, this is completely analogous to the fact that the product of two CW complexes $X$ and $Y$ has a canonical CW structure, whose $n$-generators are products of a $k$-generator of $X$, and an $(n-k)$-generator of $Y$, for $k \leq n$.

For technical purposes, it is convenient to expand the combinatorial notion of directed cell used in the definition of polygraph.

\begin{dfn} \label{dfn:loopfree-globe} \index{globe!loop-free}
A loop-free polygraph $X$ is an \emph{$n$-globe} if its incidence poset has a greatest element $\top$ with $\dmn{\top} = n$.
\end{dfn}
Let $\partial X$ be the $(n-1)$-skeleton of an $n$-globe $X$. Then $X$ is obtained from $\partial X$ by gluing the $n$-cell $\top$ to $\partial X$, that is,
\begin{equation} \label{eq:generalglobe}
\begin{tikzpicture}[baseline={([yshift=-.5ex]current bounding box.center)}]
	\node[scale=1.25] (0) at (-1,2) {$\partial G^n$};
	\node[scale=1.25] (1) at (2,2) {$G^n$};
	\node[scale=1.25] (2) at (-1,0) {$\partial X$};
	\node[scale=1.25] (3) at (2,0) {$X$};
	\draw[1cinc] (0) to (1);
	\draw[1cinc] (2) to (3);
	\draw[1c] (0) to node[auto,swap] {$\partial \top$} (2);
	\draw[1c] (1) to node[auto] {$\top$} (3);
	\draw[edge] (1,0.2) to (1,0.8) to (1.8,0.8);
\end{tikzpicture}
\end{equation}
is a pushout diagram in $\omegacat$. Hence, in the definition of a polygraph $X$, we can as well ask that, for all $x \in G_n(Y)$, there exist maps $\hat{x}: G_x \to Y$, such that $G_x$ is \emph{any} $n$-globe, $\hat{x}(\top) = x$, and
\begin{equation*}
\begin{tikzpicture}[baseline={([yshift=-.5ex]current bounding box.center)}]
	\node[scale=1.25] (0) at (0,2) {$\displaystyle \coprod_{G_n(X)} \partial G_x$};
	\node[scale=1.25] (1) at (3.5,2) {$\displaystyle\coprod_{G_n(X)} G_x$};
	\node[scale=1.25] (2) at (0,0) {$\skel{n-1}{X}$};
	\node[scale=1.25] (3) at (3.5,0) {$\skel{n}{X}$};
	\draw[1c] (0) to node[auto,swap] {$\displaystyle\coprod_{G_n(X)} \partial \hat{x}$} (2);
	\draw[1c] (1) to node[auto] {$\displaystyle\coprod_{G_n(X)} \hat{x}$} (3);
	\draw[1cinc] (0) to (1);
	\draw[1cinc] (2) to (3);
	\draw[edge] (2.5,0.2) to (2.5,0.8) to (3.3,0.8);
\end{tikzpicture}
\end{equation*}
is a pushout diagram; by composing with the diagrams (\ref{eq:generalglobe}), we reobtain the original definition. 

\begin{lem} \label{lem:tensorskeleton}
For all $\omega$-categories $X, Y$, and $n \in \mathbb{N}$,
\begin{equation*}
	\skel{n}{(X \otimes Y)} \simeq \mathrm{colim}_{k\leq m \leq n}\,\skel{k}X \otimes \skel{m-k} Y,
\end{equation*}
the colimit being over all inclusions $\skel{i}{X} \otimes \skel{j}{Y} \hookrightarrow \skel{i'}{X} \otimes \skel{j'}{Y}$.
\end{lem}
\begin{proof}
This is true by definition for loop-free polygraphs. For general $\omega$-categories $X$, the $n$-skeleton of $X$ is the colimit of the simple $n$-categories over $X$, hence also of the loop-free $n$-polygraphs over $X$. In the diagram defining $X \otimes Y$, the loop-free $n$-polygraphs are precisely those of the form $U \otimes V$ for some loop-free $k$-polygraph $U$ and $(n-k)$-polygraph $V$. Since colimits commute, and are preserved by $\otimes$, we can factorise as required.
\end{proof}

\begin{remark}
For all $x \in X$, $y \in Y$, we have $\idd{}x \otimes y = \idd{}(x \otimes y) = x \otimes \idd{}y$: this can be verified directly for loop-free polygraphs, and carries over to general $\omega$-categories since they are colimits of diagrams of unit-preserving maps. 

Therefore, from Lemma \ref{lem:tensorskeleton}, we can deduce that if $X$ is an $n$-category and $Y$ an $m$-category, then $X \otimes Y$ is an $(n+m)$-category.
\end{remark}

\begin{thm} \label{thm:tensorpoly}
Let $X, Y$ be polygraphs. Then $X \otimes Y$ is a polygraph, whose set of $n$-generators, for all $n \in \mathbb{N}$, is
\begin{equation*}
	G_n(X \otimes Y) = \coprod_{k=0}^n G_k(X) \times G_{n-k}(Y).
\end{equation*}
\end{thm}
\begin{proof}
The claim follows from the definition for loop-free polygraphs, and is obvious for $0$-polygraphs. It suffices, then, as an inductive step, to consider the extensions $\skel{n-1}{(X \otimes Y)} \hookrightarrow \skel{n}{(X \otimes Y)}$, for $n > 0$.

Using the fact that tensor products preserve colimits, for all $k \leq n$, we have a diagram
\begin{equation} \label{eq:tensorpoly}
\begin{tikzpicture}[baseline={([yshift=-.5ex]current bounding box.center)}]
	\node[scale=1.25] (0) at (0,2) {$\displaystyle \coprod \partial G^k \otimes G^{n-k}$};
	\node[scale=1.25] (1) at (4,2) {$\displaystyle \coprod G^k \otimes G^{n-k}$};
	\node[scale=1.25] (0b) at (8,2) {$\displaystyle \coprod G^k \otimes \partial G^{n-k}$};
	\node[scale=1.25] (2) at (0,0) {$\skel{k-1}{X} \otimes \skel{n-k}{Y}$};
	\node[scale=1.25] (2b) at (8,0) {$\skel{k}{X} \otimes \skel{n-k-1}{Y}$};
	\node[scale=1.25] (3) at (4,0) {$\skel{k}{X} \otimes \skel{n-k}{Y}$};
	\draw[1c] (0) to node[auto,swap] {$\displaystyle\coprod \partial x \otimes y$} (2);
	\draw[1c] (1) to node[auto, pos=0.2] {$\displaystyle\coprod x \otimes y$} (3);
	\draw[1c] (0b) to node[auto] {$\displaystyle\coprod x \otimes \partial y$} (2b);
	\draw[1cinc] (0) to (1);
	\draw[1cincl] (0b) to (1);
	\draw[1cincl] (2b) to (3);
	\draw[1cinc] (2) to (3);
	\draw[edge] (3,0.3) to (3,0.8) to (3.8,0.8);
	\draw[edge] (5,0.3) to (5,0.8) to (4.2,0.8);
	\node[scale=1.25] at (9.7,-.2) {,};
\end{tikzpicture}
\end{equation}
where the coproducts are over $x \in G_k(X)$, $y \in G_{n-k}(Y)$, and both squares are pushouts. Now, the tensor products in the top row only involve loop-free polygraphs, so they can be calculated directly.

In particular, $G^k \otimes G^{n-k}$ is an $n$-globe, and it is easy to verify that $\partial(G^k \otimes G^{n-k})$ is the join of $\partial G^k \otimes G^{n-k}$ and $G^k \otimes \partial G^{n-k}$ as subobjects of $G^k \otimes G^{n-k}$. Hence, the maps $\partial x \otimes y$ and $x \otimes \partial y$ induce a unique map from $\partial(G^k \otimes G^{n-k})$ to the join of $\skel{k-1}{X} \otimes \skel{n-k}{Y}$ and $\skel{k}{X} \otimes \skel{n-k-1}{Y}$ as subobjects of $\skel{n}{(X \otimes Y)}$ --- which is in in fact a subobject of $\skel{n-1}{(X \otimes Y)}$.

By considering the diagram composed of all diagrams (\ref{eq:tensorpoly}) for $k = 0,\ldots,n$, by Lemma \ref{lem:tensorskeleton} and commutation of colimits, we obtain a pushout square
\begin{equation*}
\begin{tikzpicture}[baseline={([yshift=-.5ex]current bounding box.center)}]
	\node[scale=1.25] (0) at (0,2) {$\displaystyle \coprod_{k=0}^n \coprod \partial (G^k \otimes G^{n-k})$};
	\node[scale=1.25] (1) at (5,2) {$\displaystyle \coprod_{k=0}^n \coprod G^k \otimes G^{n-k}$};
	\node[scale=1.25] (2) at (0,0) {$\skel{n-1}{(X \otimes Y)}$};
	\node[scale=1.25] (3) at (5,0) {$\skel{n}{(X \otimes Y)}$};
	\draw[1c] (0,1.6) to node[auto,swap] {$\displaystyle\coprod \coprod \partial (x \otimes y)$} (2);
	\draw[1c] (5,1.6) to node[auto] {$\displaystyle\coprod \coprod x \otimes y$} (3);
	\draw[1cinc] (0) to (1);
	\draw[1cinc] (2) to (3);
	\draw[edge] (4,0.2) to (4,0.8) to (4.8,0.8);
	\node[scale=1.25] at (7,-.2) {.};
\end{tikzpicture}
\end{equation*}
Since the $G^k \otimes G^{n-k}$ are $n$-globes, this proves the claim.
\end{proof}

\section{From pasting diagrams to string diagrams} \label{sec:pasting}

Besides disjoint unions and tensor products, another way of obtaining polygraphs from polygraphs is by quotients. Given a polygraph $X$ with set of generators $G(X) = \coprod_{n \in \mathbb{N}} G_n(X)$, and an equivalence relation $E$ on $G(X)$, we can try to extend it to cells of $X$, by retracing the syntactic construction of $X$: for all $(x, x') \in E$, with $\dmn{x}, \dmn{x'} \leq n$, let
\begin{align*}
	& (x,x') \in E^\star & \text{if } \dmn{x} = \dmn{x'} = n, \\
	& (x, \epsilon(\idd{n-1}x')) \in E^\star & \text{if } \dmn{x} = n, \dmn{x'} < n, \\
	& (\epsilon(\idd{n-1}x),\epsilon(\idd{n-1}x)) \in E^\star & \text{if } \dmn{x}, \dmn{x'} < n,
\end{align*}
then extend to $k$-composites involving $x$ or $x'$, and so on. 

This will define an actual equivalence relation on $X$ in $\omegacat$ if and only if, at each step, all relations $(x,x') \in E^\star$ imposed on $n$-cells are such that $(\bord{}{\alpha} x, \bord{}{\alpha} x') \in E^\star$ in $X_{n-1}$; in that case, we say that $E$ is an equivalence relation on the polygraph $X$. Pairs of equivalent cells then define a sub-$\omega$-category $E^\star \hookrightarrow X \with X$.

\begin{dfn} \index{quotient!of polygraphs}
Let $E$ be an equivalence relation on a polygraph $X$. The \emph{quotient} $X/E$ of $X$ by $E$ is the coequaliser in $\omegacat$ of 
\begin{equation*}
\begin{tikzpicture}[baseline={([yshift=-.5ex]current bounding box.center)}]
	\node[scale=1.25] (0) at (0,0) {$E^\star$};
	\node[scale=1.25] (1) at (2,0) {$X,$};
	\draw[1c] (0.35,.15) to node[auto] {$\pi_1$} (1.7,.15);
	\draw[1c] (0.35,-.15) to node[auto,swap] {$\pi_2$} (1.7,-.15);
\end{tikzpicture}
\end{equation*}
where $\pi_1, \pi_2$ are composites of $E^\star \hookrightarrow X \with X$ with the two projections. 
\end{dfn}

The $\omega$-category $X/E$ is clearly still a polygraph, with set of generators $G(X)/E$, and $\dmn{[x]} = \min \{\dmn{x} \,|\, x \in [x] \}$, for all equivalence classes $[x] \in G(X)/E$.

A particular case is that in which $A \hookrightarrow X$ is a sub-polygraph, that is, an inclusion in $\mathbf{Pol}$, that is, a map of polygraphs induced by an inclusion $G(A) \subseteq G(X)$ of their generators; and $E$ is the equivalence relation $G(A) \times G(A)$ (hence, $E^\star = A \with A$). In the quotient of $X$ by $E$, the sub-polygraph $A$ is reduced to a single 0-cell. We introduce the following special notation.

\begin{dfn} 
Let $X$ be a polygraph, $A \hookrightarrow X$ a sub-polygraph. The \emph{quotient} $X/A$ of $X$ by $A$ is the quotient of $X$ by $G(A) \times G(A)$.
\end{dfn}

By Theorem \ref{thm:tensorpoly}, if we want to calculate the tensor product of two polygraphs, it suffices to calculate the tensor products of their generators. In practice, a good strategy is the following. Given generators $x \in G(X)$ and $y \in G(Y)$, let $U_x \hookrightarrow X, U_y \hookrightarrow Y$ be the smallest sub-polygraphs containing $x$ and $y$, respectively. Then we can:
\begin{enumerate}
	\item find loop-free $n$-globes $G_x$ and $G_y$ such that $\hat{x}: G_x \to U_x$, $\hat{y}: G_y \to U_y$ are quotient maps of polygraphs by equivalence relations $E_x$, $E_y$, respectively;
	\item explicitly calculate $G_x \otimes G_y$ in $\loopfree$;
	\item quotient by $E_x \times E_y$, to obtain a description of $U_x \otimes U_y = U_{x \otimes y}$.
\end{enumerate}
In fact, the informal substitution of cells with loop-free globes is at the very essence of the use of pasting diagrams in reasoning about higher categories.

\begin{exm} \label{exm:cylinder}
Let $X$ be the 2-polygraph with the following presentation: $X$ has a single $0$-generator $x$, a single $1$-generator $a: x \to x$, and a single $2$-generator $m: a \cp{1} a \to a$. Picturing $m$ as a pasting diagram, we would draw
\begin{equation} \label{diag:multiply}
\input{img/c1_3_multiply.tex}
\end{equation} 
where the three occurrences of $x$ and the three occurrences of $a$ in the boundary of $m$ are drawn as spatially distinct objects. In fact, if we forget the labelling, we can relabel the pasting diagram as a loop-free 2-globe
\begin{equation} \label{diag:multiply_loopfree}
\input{img/c1_3_multiply_loopfree.tex}
\end{equation} 
the labels in diagram (\ref{diag:multiply}) indicate what identifications of cells are needed to obtain $m$ as a quotient.

Let us calculate the tensor product $\vec{I} \otimes X$, where $\vec{I} := G^1$ is the ``directed interval'', that is, the standard $1$-globe; we could directly calculate an expression for its 3-generator, whose boundary contains all the others, but for this first example we will proceed step by step. The 0-generators are just $0^- \otimes x$ and $0^+ \otimes x$. There are three 1-generators: $\top \otimes x: 0^- \otimes x \to 0^+ \otimes x$, $0^- \otimes a$, and $0^+ \otimes a$. \index{i@$\vec{I}$}

To calculate an expression for the latter two, we can replace $U_a$ with a standard 1-globe $G_a$, with the unique 1-cell $a_1: x_1 \to x_2$. Then, in $\{0^\alpha\} \otimes G_a$, we have a 1-cell $0^\alpha \otimes a_1: 0^\alpha \otimes x_1 \to 0^\alpha \otimes x_2$; quotienting to $\vec{I} \otimes X$, we find $0^\alpha \otimes a: 0^\alpha \otimes x \to 0^\alpha \otimes x$.

The 2-generators are either products of a 0-generator and a 2-generator, or of two 1-generators. The former are of the form $0^\alpha \otimes m$, for which we can replace $U_m$ with the 2-globe $G_m$ of diagram (\ref{diag:multiply_loopfree}); proceeding as before, we find $0^\alpha \otimes m: (0^\alpha \otimes a) \cp{1} (0^\alpha \otimes a) \to (0^\alpha \otimes a)$. 

The latter only include $\top \otimes a$, for which we can calculate $\vec{I} \otimes G_a$. In terms of augmented directed complexes, we find
\begin{align*}
	d(\top \otimes a_1) & = d\top \otimes a_1 - \top \otimes da_1 = \\
		& = 0^+ \otimes a_1 - 0^- \otimes a_1 - \top \otimes x_2 + \top \otimes x_1;
\end{align*}
in the induced polygraph, $\bord{}{+}\lfgen{\top \otimes a_1}$ is
\begin{align*}
	& \;\; (0^- \otimes x_1, 0^+ \otimes x_2,  \top \otimes x_1 + 0^+ \otimes a_1) = \\
	= & \;\; \lfgen{\top \otimes x_1} \cp{1} \lfgen{0^+ \otimes a_1},
\end{align*}
and $\bord{}{-}\lfgen{\top \otimes a_1}$ is
\begin{align*}
	& \;\; (0^- \otimes x_1, 0^+ \otimes x_2,  0^- \otimes a_1 + \top \otimes x_2) = \\
	= & \;\; \lfgen{0^- \otimes a_1} \cp{1} \lfgen{\top \otimes x_2}.
\end{align*}
As a pasting diagram, this is
\begin{equation*}
\input{img/c1_3_square.tex}
\end{equation*} 
where we wrote $x y$ for $x \otimes y$. Passing to the quotient, we find $\top \otimes a: (0^- \otimes a) \cp{1} (\top \otimes x) \to (\top \otimes x) \cp{1} (0^+ \otimes a)$. 

Finally, $\vec{I} \otimes X$ has a single 3-generator $\top \otimes m$, for which we can calculate $\vec{I} \otimes G_m$. In the corresponding augmented directed complex,
\begin{align} \begin{split} \label{eq:topboundary}
	d(\top \otimes m_1) = & \;\; d\top \otimes m_1 - \top \otimes d m_1 = \\
	= & \;\; 0^+ \otimes m_1 - 0^- \otimes m_1 - \top \otimes a_3 + \top \otimes a_1 + \top \otimes a_2.
\end{split} \end{align}
Then, in the induced polygraph, $\bord{}{+}\lfgen{\top \otimes m_1}$ is
\begin{align*}
	& \;\; (0^- \otimes x_1, 0^+ \otimes x_3, 0^- \otimes a_1 + 0^- \otimes a_2 + \top \otimes x_3, \top \otimes x_1 + 0^+ \otimes a_3, \\ & \qquad \qquad \top \otimes a_1 + \top \otimes a_2 + 0^+ \otimes m_1) = \\
	= & \;\; (\idd{}\lfgen{0^- \otimes a_1} \cp{1} \lfgen{\top \otimes a_2}) \cp{2} (\lfgen{\top \otimes a_1} \cp{1} \idd{}\lfgen{0^+ \otimes a_2}) \cp{2} (\idd{}\lfgen{\top \otimes x_1} \cp{1} \lfgen{0^+ \otimes m_1}),
\end{align*}
and $\bord{}{-}\lfgen{\top \otimes m_1}$ can similarly be calculated to be $(\lfgen{0^- \otimes m_1} \cp{1} \idd{}\lfgen{\top \otimes x_3}) \cp{2} \lfgen{\top \otimes a_3}$. 

While from \cite[Section 5]{steiner2004omega} it is possible to extract an algorithmic procedure to calculate these expressions, it can be helpful, in low-dimensional cases, to consider the pasting diagrams for the top-dimensional generators in the boundary of a cell, and reason diagrammatically to understand how they fit together. In this case, knowing the shape of the 2-generators in the output boundary of $\top \otimes m_1$ --- which by equation (\ref{eq:topboundary}) are $\top \otimes a_1, \top \otimes a_2$, and $0^+ \otimes m_1$ --- we can reconstruct the pasting diagram
\begin{equation*}
\input{img/c1_3_multiply_cylinder.tex}
\end{equation*} 
for $\bord{}{+} (\top \otimes m_1)$, where we omitted the labels of 0-cells to avoid clutter.

Passing to the quotient, we obtain the following picture for $\top \otimes m$ in pasting diagrams:
\begin{equation} \label{eq:multiply_cylinder}
\input{img/c1_3_multiply_cylinder_2.tex}
\end{equation} 
\end{exm}
In the remainder, we will omit such explicit calculations; in Chapter \ref{chap:interacting}, we will see how, in some common cases, it is possible to ``guess'', by diagrammatical reasoning, the shape of the tensor product of two cells.

A common way of eliminating loops, when picturing cells by pasting diagrams, is by introducing ``filler'' cells, which are meant to be quotiented down to lower-dimensional cells. For instance, a 2-cell $y: \idd{} x \to a$, whose input is a 0-dimensional cell, can be represented as a quotient of any of the following loop-free 2-globes:
\begin{equation*}
\input{img/c1_3_loopfree_unit.tex}
\end{equation*} 
The first diagram is just a standard 2-globe. The second is a cubical cell; the possibility of transforming globes into cubes by introducing enough degenerate cells, and, conversely, of using connections to transform cubes into ``combinatorial globes'', where non-degenerate $n$-cells only appear in the first $n$ dimensions, is the basis of the equivalence theorem of \cite{al2002multiple}. 

In the third diagram, we started from the second, cubical picture, and shrank it --- it appears as the smaller, central square in the diagram --- by filling the surrounding space with unit cells. Because of the combinatorial nature of pasting diagrams, the size of cells does not have any intrinsic meaning, and because unit cells can be freely introduced without affecting the quotient, it is clear that we can ``fill all available space'' with them. In particular, in a 2-dimensional pasting diagram,
\begin{itemize}
	\item any 2-dimensional cell can be shrunk to a point;
	\item any 1-dimensional cell can be expanded vertically, then shrunk horizontally to a line;
	\item any 0-dimensional cell can be expanded to fill a region of the plane.
\end{itemize}
In this limit, we can picture the 2-cell $y$ as a diagram
\begin{equation} \label{diag:unit_string}
\input{img/c1_3_unit_string.tex}
\end{equation} 
with the convention that points represent 2-dimensional cells, lines 1-dimensional cells, and bounded regions of space 0-dimensional cells. 

As another example, the 2-cell (\ref{diag:multiply}) becomes (labels omitted on the first diagram)
\begin{equation} \label{diag:multiply_str}
\input{img/c1_3_multiply_string.tex}
\end{equation} 
If we are not interested in the specific 0-cells appearing in the diagram, but only in what pairs are identified, we can encode the same information by using different colours for different regions of the plane. For example, the loop-free 2-globe (\ref{diag:multiply_loopfree}) can be represented as the diagram
\begin{equation} \label{diag:multiply_string}
\input{img/c1_3_multiply_string_colour.tex}
\end{equation} \index{string diagram}

Diagrams such as (\ref{diag:unit_string}) and (\ref{diag:multiply_string}) are called \emph{string diagrams}. Although they are related to earlier, informal notations such as Petri nets \cite{sobocinski2010representations} and Feynman diagrams \cite{blute2010proof}, their popularisation as a general abstract notation is usually credited to Penrose \cite{penrose1971applications}, and their introduction in category theory to Kelly \cite{kelly1972many}. Having been used informally by category theorists for many years, they were first studied as a formal language in its own right by Joyal and Street in \cite{joyal1991geometry}.

Whereas pasting diagrams are closer to the intuition of the analogy between polygraphs and CW complexes, string diagrams can be better suited to a process-theoretic interpretation of 2-categories, where a 2-cell $y: a_1 \cp{1} \ldots \cp{1} a_n \to b_1 \cp{1} \ldots \cp{1} b_m$ is seen as a process taking $n$ objects of types $a_1, \ldots, a_n$ as input, and outputting $m$ objects of types $b_1, \ldots, b_m$; see for example \cite[Chapter 3]{coecke2017picturing}. In fact, a good way of reading a string diagram, in analogy with Feynman diagrams, is as a timeline of the evolution of a number of objects --- say, particles --- where lines are the particles' worldlines. For example, diagram (\ref{diag:unit_string}) represents the creation, at point $y$, of a particle of type $a$; and diagram (\ref{diag:multiply_string}) represents the fusion, at point $m_1$, of two particles of types $a_1$ and $a_2$ into a particle of type $a_3$.

Moreover, as geometric objects, string diagrams can be seen as certain topological graphs, with an orientation on their edges, and an embedding in the plane compatible with the orientation \cite{joyal1991geometry}. Now, graphs have their own set of ``natural'' topological moves, such as the bending of wires; and there are other topological moves associated with the embedding of graphs into spaces other than the plane: for example, the \emph{braiding} of wires when the embedding is into a 3-dimensional manifold. 

It turns out that many interesting classes of higher categories are characterised by the fact that their set of cells is closed under certain topological moves on the corresponding string diagrams, and various kinds of topological equivalences on the diagrams induce equalities or equivalences of cells. A survey of such results is \cite{selinger2011survey}. 

Added to the fact that graphs may be computationally better behaved for purposes of automated rewriting --- see for example \cite{bonchi2016rewriting} --- this makes string diagrams a powerful reasoning tool. We will employ them extensively in the rest of the thesis, introducing various sets of rules as we need them.

\begin{remark} \label{remark:strings}
There exist several formalisations of string diagrams as mathematical objects, more or less restricted in purpose; the references already listed are a good starting point. 

A common statement is that string diagrams are ``Poincar\'e duals'' of pasting diagrams. However, we find this slightly misleading, and prefer the intuition beneath our informal ``definition'', exemplified by diagram (\ref{diag:multiply_str}):
\begin{center}
	\emph{String diagrams are pasting diagrams filled up with unit cells}.
\end{center}
The reason is that, while so far we have been dealing with ``strict'' unit cells, that are to be quotiented down to a lower-dimensional cell, one often wants to work with ``weak unit cells'', that \emph{represent} lower-dimensional cells, yet whose elimination is non-trivial. For instance, the non-trivial braidings of strings one obtains when reasoning with 3-dimensional string diagrams are trivialised in the quotient when the ``surrounding space'' is strictly degenerate, but not when it consists of weak units.  

If string diagrams are just special pasting diagrams, it suffices to slightly adjust the definition: instead of filling up with strict units, fill up with weak units. These weak unit cells do not have any special, structural status: in particular, weak unit $n$-cells are genuine $n$-dimensional cells, so the Poincar\'e dual definition fails to generalise. 

Rather, weak units are characterised by their properties, and these properties are reflected in the intuitive use of the ``surrounding space'' of string diagrams. For example, the possibility of arbitrarily decomposing the surrounding space into cells corresponds to an idempotency property of units; and the fact that we can create and eliminate space on the edges corresponds to a cancellability property. In fact, until one of these properties is needed in a computation, we may as well treat the surrounding space as an \emph{arbitrary} cell. We will discuss this more in detail in Chapter \ref{chap:directed}.
\end{remark}

Some operations on polygraphs have a simple effect on the string-diagrammatic representations of their cells. For example, given the string diagram for a cell $x$ of $X$, and a sub-polygraph $A \hookrightarrow X$, in the string diagram for the quotient $[x]$ of $x$ in $X/A$ all generators in $G(A)$ are merged together into a single-coloured region of the plane. 

\begin{exm}
Let $X$ be the loop-free 2-globe (\ref{diag:multiply_loopfree}), and let $A$ be the sub-polygraph $U_{a_2}$. The operation of quotienting $X$ by $A$ can be pictured as follows in string diagrams:
\begin{equation*}
\input{img/c1_3_multiply_quotient.tex}
\end{equation*} 
\end{exm}

Unlike the operations we have introduced so far, which have an undirected, topological counterpart, the following only makes sense for directed cells. Let $\mathbb{N}^+ := \mathbb{N} - \{0\}$.

\begin{dfn} \index{opposite!$\omega$-category}
Let $X$ be an $\omega$-category, $S \subseteq \mathbb{N}^+$. The \emph{$S$-opposite $\omega$-category} of $X$ is the $\omega$-category $\oppn{X}{S}$ so defined: for each $n \in \mathbb{N}$, let $\oppn{X}{S}_n := \{\oppn{x}{S} \,|\, x \in X_n\}$; boundary and unit operations are defined, at each $n$-cell $x$ of $X$, by
\begin{align*}
	\bord{}{\alpha}(\oppn{x}{S}) & := \begin{cases}
		\oppn{(\bord{}{-\alpha} x)}{S}, & n \in S, \\
		\oppn{(\bord{}{\alpha} x)}{S}, & n \not\in S,
	\end{cases} \\
	\idd{} \oppn{x}{S} & := \oppn{(\idd{}x)}{S}.
\end{align*}
If $x, y$ are $n$-cells of $X$, and $k < n$, then $\oppn{x}{S} \comp{k} \oppn{y}{S}$ in $\oppn{X}{S}$ if and only if either $y \comp{k} x$ and $k \in S$, or $x \comp{k} y$ and $k \not\in S$. In either case,
\begin{equation*}
	\oppn{x}{S} \cp{k} \oppn{y}{S} := \begin{cases}
		\oppn{(y \cp{k} x)}{S}, & k \in S, \\
		\oppn{(x \cp{k} y)}{S}, & k \not\in S.
	\end{cases}
\end{equation*}
For each $S \subseteq \mathbb{N}^+$, the assignment $\oppn{(-)}{S}$ extends to a functor $\omegacat \to \omegacat$, defined by $\oppn{f}{S}: \oppn{x}{S} \mapsto \oppn{f(x)}{S}$ for each map $f: X \to Y$ and cell $x$ of $X$.

If $S$ is a finite set $\{x_1, \ldots, x_n\}$, let $\oppn{X}{x_1,\ldots,x_n} := \oppn{X}{S}$. We write $\opp{X} := \oppn{X}{1}$, $\coo{X} := \oppn{X}{2}$, and $X^- := \oppn{X}{\mathbb{N}^+}$.
\end{dfn}

Informally, $\oppn{X}{S}$ is $X$ with the direction of $n$-cells reversed, for all $n \in S$.

\begin{prop}
Let $X$ be an $\omega$-category, $S_1, S_2 \subseteq \mathbb{N}^+$. Then $\oppn{(\oppn{X}{S_1})}{S_2} = \oppn{(\oppn{X}{S_2})}{S_1} = \oppn{X}{S}$ with $S := (S_1 \cup S_2) - (S_1 \cap S_2)$. In particular, $\oppn{(\oppn{X}{S})}{S} = \oppn{X}{\emptyset} = X$.
\end{prop}
\begin{proof}
Obvious.
\end{proof}

\begin{prop}
Let $X$ be a polygraph, $S \subseteq \mathbb{N}^+$. Then $\oppn{X}{S}$ is a polygraph whose set of $n$-generators is $\oppn{G_n(X)}{S}$.
\end{prop}
\begin{proof}
Simple induction on the syntactic definition of $X$, as in Construction \ref{cons:freeomega}.
\end{proof}

While for arbitrary $S$, $\oppn{(-)}{S}$ does not interact with the tensor product in any obvious manner, it does when $S = \mathbb{N}^+$.

\begin{prop}
Let $X$, $Y$ be $\omega$-categories. Then $(X \otimes Y)^- \simeq X^- \otimes Y^-$.
\end{prop}
\begin{proof}
Let $\augm{X}$ be a loop-free polygraph, seen as an augmented directed complex. Then $\augm{X}^- = \augm{X^-}$, where $X^- := \{x^- \,|\, x \in X\}$, and, for all $x \in X_n$,
\begin{align*}
	e(x^-) & := ex, & n = 0, \\
	d(x^-) & := - (dx)^-, & n > 0.
\end{align*}
Let $\augm{Y}$ be another loop-free polygraph; we want to show that the assignment $(x \otimes y)^- \mapsto x^- \otimes y^-$, for all $x \in X$, $y \in Y$, defines an isomorphism $(\augm{X} \otimes \augm{Y})^- \simeq \augm{X^-} \otimes \augm{Y^-}$.

When both $\dmn{x}$ and $\dmn{y}$ are 0, this is obviously well-defined. If $\dmn{x} > \dmn{y} = 0$, then
\begin{align*}
	d((x \otimes y)^-) & \; = \; - (dx \otimes y)^- \mapsto \\
	& \mapsto \, - (dx)^- \otimes y^- = d(x^-) \otimes y^- = d(x^- \otimes y^-),
\end{align*}
and similarly if $\dmn{y} > \dmn{x} = 0$. Finally, if $\dmn{x}, \dmn{y} > 0$,
\begin{align*}
	d((x \otimes y)^-) & \; = -(dx \otimes y + (-1)^{\dmn{x}} x \otimes dy)^- \mapsto \\
	& \mapsto \, - (dx)^- \otimes y^- - (-1)^{\dmn{x}} x^- \otimes (dy)^- = \\
	& \; = \; d(x^-) \otimes y^- + (-1)^{\dmn{x}} x \otimes	d(y^-) = d(x^- \otimes y^-).
\end{align*}

So $(X \otimes Y)^- \simeq X^- \otimes Y^-$ for all pairs of loop-free polygraphs. Since $(-)^-$ is an invertible endofunctor, in particular it preserves colimits, so the isomorphism extends to tensor products of arbitrary $\omega$-categories.
\end{proof}

In string diagrams, the operations $X \mapsto \opp{X}, \coo{X}$ become horizontal and vertical reflection, respectively.

\begin{exm} Let $X$ be the loop-free 2-globe (\ref{diag:multiply_loopfree}). The following are representations of $X$, $\opp{X}$, $\coo{X}$, and $X^\mathrm{op\,co} = X^-$ in string diagrams:
\begin{equation*}
\input{img/c1_3_multiply_opposite.tex}
\end{equation*} 
\end{exm}

Together with disjoint unions, tensor products and quotients, opposites complete our elementary toolkit of polygraph-building operations. In the next chapter, we will put it to use in one of the traditional areas of application of category theory, universal algebra; and learn to recognise what tensor products look like in string diagrams.

  \thispagestyle{empty} 
\chapter{Interacting algebraic theories} \label{chap:interacting}
\thispagestyle{plain}

\noindent\emph{In this chapter:}
\begin{itemize}
	\item[$\triangleright$] We give a brief overview of categorical universal algebra from a higher-dimensional point of view, centred on the notion of PRO and on its variants. We present some basic algebraic theories, and touch on connections with higher algebra and rewriting theory. --- \emph{Section \ref{sec:univ_algebra}}
	\item[$\triangleright$] Starting from the observation that ``directed cylinders'' capture in great generality the notion of homomorphism of algebras, we show how to construct algebraic theories compositionally, using the operations described in Chapter \ref{chap:polygraphs}. Theories assembled with tensor products come automatically with higher-dimensional coherence cells. --- \emph{Section \ref{sec:sliding}}
	\item[$\triangleright$] We define a smash product operation for pointed $\omega$-categories, and show that it is related to certain interactions of algebraic theories, leading for instance to the theory of bialgebras, yet fails to capture them faithfully, due to a degeneracy in the combinatorics of $\omega$-categories. --- \emph{Section \ref{sec:dimension}}
\end{itemize}

\section{Higher categories and universal algebra} \label{sec:univ_algebra}
To a first approximation, an algebraic theory is a mathematical object described by \index{algebraic theory}
\begin{itemize}
	\item a collection of basic finitary operations $f$ that take $n$ objects as input, and output $m$ objects;
	\item axioms imposing the identity of certain composites of the operations.
\end{itemize}
In logical terms, this corresponds to a first-order theory over a signature containing only finitary function symbols, whose axioms are all identities, universally quantified over their variables. The class of structures axiomatised by such a theory is called a \emph{variety}, and is the primary object of study of universal algebra.

Recalling the process-theoretic interpretation of 2-categories that we mentioned in the previous chapter, in relation to string diagrams, we see that the same information can be given in the following form.
\begin{dfn}
The \emph{oriented 1-sphere} is the 1-polygraph $\vec{S}^1$ with a single 0-generator $\bullet$, and a single 1-generator $a: \bullet \to \bullet$.

Let $X$ be a 2-category. Then $X$ is a \emph{PRO} if $\skel{1}{X} \simeq \vec{S}^1$. A \emph{presentation} of $X$ is a polygraph $X'$ such that $X \simeq \tau_{\leq 2} X'$. \index{PRO} \index{presentation}

If $X$, $Y$ are PROs, with 1-generators $a_X$, $a_Y$, a \emph{map} $f: X \to Y$ is a map of 2-categories such that $f(a_X) = a_Y$. PROs and their maps form a category $\cat{PRO}$.
\end{dfn}
The translation between algebraic theories and PROs is the following. If $X$ is a PRO, then each 1-cell of $X$ is uniquely of the form
\begin{equation*}
	[0] := \idd{}\bullet \qquad \text{ or } \qquad [n] := \underbrace{a \cp{1} \ldots \cp{1} a}_{n}, \quad n>0.
\end{equation*}
Then, given a presentation $X'$ of $X$, we can interpret each 2-generator $f: [n] \to [m]$ of $X'$ as a basic operation with $n$-ary input, and $m$-ary output, and each 3-generator $e: f_1 \to f_2$ of $X'$ as an axiom $f_1 = f_2$ between composites of basic operations. Because $X = \tau_{\leq 2} X'$, any such relation becomes an equality in $X$, and any non-geometric equality in $X$ is induced by such an axiom.

We will employ the two viewpoints interchangeably; when we say that a PRO $X$ is presented by certain 2-generators and axioms $\{f_i = g_i\}_{i\in I}$, we mean that the 3-polygraph with 3-generators $\{e_i: f_i \to g_i\}_{i\in I}$, or equivalently $\{\bar{e}_i: g_i \to f_i\}_{i\in I}$, is a presentation of $X$.

\begin{prop}
Every PRO has a presentation.
\end{prop}
\begin{proof}
If $X$ is a PRO, we have by definition a polygraphic structure on $\skel{1}{X}$. We can extend this to a polygraphic resolution of $X$, as in Construction \ref{cons:resolution}; we obtain a polygraph $X'$ such that $X$ is weakly equivalent to $\tau_{\leq 2}X'$. But a weak equivalence between 2-categories is a bijection on 2-cells, and by construction we also have an isomorphism of 1-skeleta; thus, $X$ is actually isomorphic to $\tau_{\leq 2} X'$, and $X'$ is a presentation of $X$.
\end{proof}

While the traditional logical foundations of universal algebra presuppose a strict dichotomy between \emph{syntax}, the theory, and \emph{semantics}, the models, the category-theoretic formulation internalises it, transforming it into a simpler dichotomy between domain and codomain.
\begin{dfn} \index{model}
Let $X$ be a PRO. A \emph{model} of $X$ in the $\omega$-category $Y$, also called an \emph{$X$-algebra}, is a map $X \to Y$.
\end{dfn}
Foundations permitting, we reobtain the set-theoretic notion of model when $Y$ is the 2-category $\textit{Set}_\times$ with a single 0-cell, whose 1-cells are finite sequences $(A_1, \ldots, A_n)$ of sets, picked from a suitably small ``set of sets'', with concatenation of sequences as 1-composition, and whose 2-cells $f: (A_1, \ldots, A_n) \to (B_1, \ldots, B_m)$ are functions $f: A_1 \times \ldots \times A_n \to B_1 \times \ldots \times B_m$.

When one is only interested in set-theoretic models, where operations are interpreted as functions, an operation $f$ with $m$-ary output is, of course, the same as $m$ operations $f_1, \ldots, f_m$ with unary output; the same holds for models in all 2-category $Y$ where, like in $\textit{Set}_\times$, 1-composition is the categorical product in a 1-category. We rephrase this more precisely.

\begin{cons} \label{remark:lawvere} Given an $\omega$-category $Y$, and 0-cells $x, y$ of $Y$, there is an $\omega$-category $Y(x,y)$ whose 0-cells $\transp{a}$ are 1-cells $a: x \to y$ of $Y$, 1-cells $\transp{f}: \transp{a} \to \transp{b}$ are 2-cells $f: a \to b$ of $Y$, and so on; composition is defined by $\transp{w} \cp{n} \transp{z} = \transp{w \cp{n+1} z}$, when $w \comp{n+1} z$ in $Y$. When $Y$ is an $n$-category, $Y(x,y)$ is an $(n-1)$-category. 
\end{cons} \index{Lawvere theory} \index{Law@$\cat{Law}$}

Then, if $Y$ is a 2-category with a single 1-cell $\bullet$, when we say that 1-composition is a categorical product in $Y$, we mean that
\begin{enumerate}
	\item $\idd{} \bullet$ is terminal in the 1-category $Y(\bullet,\bullet)$, and
	\item for all 1-cells $a, b$ of $Y$, if $d_a: a \to \idd{} \bullet$, $d_b: b \to \idd{} \bullet$ are the unique 2-cells from $a, b$ to $\idd{} \bullet$ in $Y$, then $\idd{} a \cp{1} d_b: a \cp{1} b \to a$ and $d_a \cp{1} \idd{} b: a \cp{1} b \to b$ form a product diagram in $Y(\bullet,\bullet)$.
\end{enumerate}
We can demand that this hold already in the PRO $X$. This was, in fact, a requirement in the first formulation of categorical universal algebra, given in Lawvere's PhD thesis \cite{lawvere1963functorial}; the following is equivalent to Lawvere's definition.

\begin{dfn}
Let $X$ be a PRO. $X$ is a \emph{Lawvere theory} if its 1-composition is a categorical product. We write $\cat{Law}$ for the subcategory of $\cat{PRO}$ whose objects are Lawvere theories, and maps $f: X \to Y$ induce product-preserving maps $X(\bullet, \bullet) \to Y(\bullet, \bullet)$.
\end{dfn}
By a well-known criterion for monoidal structures to be cartesian \cite[Theorem 4.31]{heunen2017lectures}, Lawvere theories can be characterised algebraically among PROs, as the ones that contain 2-cells \index{swap} \index{copy operation} \index{discard operation}
\begin{equation*}
\input{img/c2_1_cartesian.tex}
\end{equation*} 
satisfying the axioms (all labels left implicit)
\begin{equation*}
\input{img/c2_1_cartesian_ax1.tex}
\end{equation*} 
\begin{equation*}
\input{img/c2_1_cartesian_ax2.tex}
\end{equation*} 
and, for all 2-cells $f: [n] \to [m]$ of $X$, including $s$, $c$, and $d$ themselves,
\begin{equation} \label{eq:cartesian_nat}
\input{img/c2_1_cartesian_nat.tex}
\end{equation} 
\begin{equation} \label{eq:cartesian_nat2}
\input{img/c2_1_cartesian_nat2.tex}
\end{equation} 
where the thinner, lighter lines and points indicate the repetition of a pattern for a suitable number of times, in this case $n$ times in the input, and $m$ times in the output of $f$ (see \cite{kissinger2014pattern} for a more formal version of string diagrams with repeated patterns).

To obtain an intuition about these, it is useful to think about their intended interpretations in $\textit{Set}_\times$, where, for some set $A$,
\begin{align*}
	s: A \times A \to A \times A, \;\qquad (x,y) & \; \mapsto (y,x), \\
	c: A \to A \times A, \qquad\qquad  x & \; \mapsto (x,x), \\
	d: A \to \{*\}, \qquad \qquad x & \; \mapsto *.
\end{align*}
We have just defined an important example of an algebraic theory.
\begin{dfn}
The theory of \emph{cocommutative comonoids} $\coo{\textit{CMon}}$ is the Lawvere theory presented by the 2-generators $s$, $c$, $d$, and the axioms $sym$, $\coo{un}_L$, $\coo{un}_R$ and $\{nat^f_g \,|\, f, g = s,c,d\}$. \index{comonoid}
\end{dfn}
If $X$, $Y$ are Lawvere theories, and $f: X \to Y$ is a map of PROs, an equivalent condition for $f$ to be a map of Lawvere theories is that $f$ send the swap, copy and discard 2-cells of $X$ to those of $Y$. It follows that $\coo{\textit{CMon}}$ is the initial object of $\cat{Law}$.

We can prove a couple of derived equations of this theory diagrammatically, as an example:
\begin{equation*}
\input{img/c2_1_coassociativity.tex}
\end{equation*} 
\begin{equation*}
\input{img/c2_1_cocommutativity.tex}
\end{equation*} 
The latter is a dual form of the Eckmann-Hilton argument, that a pair of monoids satisfying an interchange equation collapses to a single commutative monoid \cite{eckmann1962group}. \index{Eckmann-Hilton argument}

An intermediate notion between a PRO and a Lawvere theory is obtained when only the swap generator, and related axioms, are kept.
\begin{dfn} \index{PROP}
A \emph{PROP} is a PRO $X$ together with a distinguished swap 2-cell $s: [2] \to [2]$, satisfying the axioms $sym$ and $\{nat_s^f \,|\, f \in X_2\}$.

We write $\cat{PROP}$ for the category whose objects are PROPs, and maps $f: X \to Y$ send the swap cell of $X$ to the swap cell of $Y$.
\end{dfn}
Any Lawvere theory is a PROP with its swap cell. The initial object of $\cat{PROP}$ is $\textit{Sym}$, the PROP presented by the 2-generator $s$ and the axioms $sym$, $nat^s_s$. This is equivalent to the PROP of permutations, whose 2-cells are of the form $\sigma: [n] \to [n]$ for some $n \in \mathbb{N}$, and permutation $\sigma \in S_n$, the group of permutations of $n$ elements. Because every PROP $X$ contains $Sym$, for all $n, m \in \mathbb{N}$, the group $S_n$ acts from the left, and $S_m$ acts from the right by 2-composition on $X([n], [m])$. \index{Sym@$\textit{Sym}$}

The string diagrams for a PROP $X$ can be seen as topological graphs with an \emph{immersion}, rather than an embedding, into the plane. The interpretation of a cell depends both on the graph structure, and the choice of immersion: hence, while $\sigma \cp{2} f \cp{2} \tau$ and $f$ have the same underlying graph for all $f: [n] \to [m]$, all $\sigma \in S_n$, and all $\tau \in S_m$, in a model they can be mapped to different cells. 

A still weaker notion is obtained when permutation groups are replaced by braid groups.
\begin{dfn} \index{PROB} \index{braiding}
A \emph{PROB} is a PRO $X$ together with distinguished braiding 2-cells
\begin{equation*}
\input{img/c2_1_braiding.tex}
\end{equation*} 
satisfying the axioms, for all 2-cells $f: [n] \to [m]$ (labels left implicit),
\begin{equation*}
\input{img/c2_1_braiding_ax.tex}
\end{equation*} 
\begin{equation} \label{eq:braiding_nat}
\input{img/c2_1_braiding_nat.tex}
\end{equation} 
We write $\cat{PROB}$ for the category whose objects are PROBs, and maps $f: X \to Y$ send the braiding cells of $X$ to the braiding cells of $Y$.
\end{dfn}
The swap map of a PROP satisfies the axioms of a braiding with $b = \bar{b}$; therefore, there is an inclusion of categories $\cat{PROP} \hookrightarrow \cat{PROB}$, and both come with forgetful functors to $\cat{PRO}$. 

Let $rei_3 := nat_{b,L}^b = nat_{b,R}^b$, that is,
\begin{equation*}
\input{img/c2_1_reidemeister3.tex}
\end{equation*}
The initial PROB is $\textit{Brd}$, the PROB presented by the 2-generators $b$, $\bar{b}$ and the axioms $rei_2$, $rei'_2$, $rei_3$. These axioms are known in knot theory as the second and third Reidemeister moves; we refer to \cite{kassel2008braid}, or any other monograph on the subject, for more details. The third Reidemeister move $rei_3$ is also known in physics as the Yang-Baxter equation \cite{perk2006yang}. The string diagrams for a PROB can be seen as topological graphs embedded into 3-dimensional space, then projected onto a plane. \index{Reidemeister moves} \index{Brd@$\textit{Brd}$}

The following are some basic duality results.
\begin{prop}
Let $X$ be a PRO. Then:
\begin{enumerate}[label=(\alph*)]
	\item $\opp{X}$ and $\coo{X}$ are both PROs, and they are PROBs or PROPs if and only if $X$ is;
	\item $\opp{X}$ is a Lawvere theory if and only if $X$ is.
\end{enumerate}
\end{prop}
\begin{proof}
The first point is obvious for PROs, and follows for PROBs and PROPs from the fact that $\textit{Brd}$ and $\textit{Sym}$ are self-dual in both senses. The second point follows from $\coo{\textit{CMon}}$ being isomorphic to $\textit{CMon}^\mathrm{co \,op}$, but not to $\textit{CMon}$. 
\end{proof}

So, every time we define an algebraic theory $X$, we obtain three more by duality.
\begin{dfn} \index{monoid} \index{Mon@$\textit{Mon}$, $\textit{CMon}$}
The theory of \emph{commutative monoids} is the PROP $\textit{CMon}$, the 2-opposite of $\coo{\textit{CMon}}$.

The theory of \emph{comonoids} is the sub-PRO $\coo{\textit{Mon}}$ of $\coo{\textit{CMon}}$, whose 2-generators are $c$ and $d$, satisfying the axioms $\coo{un}_L$, $\coo{un}_R$, and
\begin{equation*}
\input{img/c2_1_coassociativity_bis.tex}
\end{equation*} 
The theory of \emph{monoids} is its 2-opposite $\textit{Mon}$. Because of its importance, we name its 2-generators: \index{multiplication} \index{unit!of a monoid}
\begin{equation*}
\input{img/c2_1_monoid.tex}
\end{equation*} 
If $X$, $Y$ are two PROs, we can form their \emph{union} $X \uplus Y$, as the pushout \index{union!of PROs}
\begin{equation*}
\begin{tikzpicture}[baseline={([yshift=-.5ex]current bounding box.center)}]
	\node[scale=1.25] (0) at (-1,2) {$\vec{S}^1$};
	\node[scale=1.25] (1) at (2,2) {$X$};
	\node[scale=1.25] (2) at (-1,0) {$Y$};
	\node[scale=1.25] (3) at (2,0) {$X \uplus Y$};
	\draw[1cinc] (0) to (1);
	\draw[1cinc] (0) to (2);
	\draw[1cinc] (2) to (3);
	\draw[1cinc] (1) to (3);
	\draw[edge] (1,0.2) to (1,0.8) to (1.8,0.8);
\end{tikzpicture}
\end{equation*}
in $\omegacat$; by construction, $X \uplus Y$ is still a PRO. Given presentations $X'$ of $X$, and $Y'$ of $Y$, $X \uplus Y$ is presented by the disjoint union of the generators and axioms of $X'$ and $Y'$. The theory of \emph{Frobenius algebras} is the PRO $\textit{Frob}$ obtained as a quotient of $\textit{Mon} \uplus \coo{\textit{Mon}}$ by the additional axioms
\begin{equation} \label{eq:frobenius}
\input{img/c2_1_frobenius.tex}
\end{equation} 
The theory $SFrob$ of \emph{special Frobenius algebras} has the further axiom \index{Frobenius algebra} \index{Frob@$\textit{Frob}, \textit{CFrob}, \textit{SCFrob}$}  
\begin{equation*}
\input{img/c2_1_frobenius_special.tex}
\end{equation*} 
These have PROP analogues $\textit{CFrob}$ and $\textit{SCFrob}$, based on $\textit{CMon} \uplus \coo{\textit{CMon}}$.
\end{dfn}

\begin{remark} \label{remark:spider} \index{spider presentation}
The theories of monoids, comonoids and Frobenius algebras have a useful alternative presentation, --- the ``spider'' presentation, in the terminology of \cite{coecke2008interacting}. For Frobenius algebras, this has infinite generators of the form $s_n^m: [n] \to [m]$, for $m, n \in \mathbb{N}$, and axioms $s^1_1 = \idd{} a$ and, for all $s_n^m, s_p^q$,
\begin{equation*}
\input{img/c2_1_spider.tex}
\end{equation*} 
\begin{equation*}
\input{img/c2_1_spider2.tex}
\end{equation*} 
where $n' = n+p-1$, and $m' = m+q-1$. By restricting to 2-generators of the form $s_n^1$ or $s_1^m$,  one obtains presentations of $\textit{Mon}$ and $\coo{\textit{Mon}}$, respectively.

For special Frobenius algebras, we allow more than one shared edge on the left hand side of a $cut$ axiom, for example
\begin{equation} \label{eq:spider_special}
\input{img/c2_1_spider_special.tex}
\end{equation} 
while for commutative monoids, comonoids, Frobenius algebras, we add axioms
\begin{equation*}
\input{img/c2_1_spider_commutative.tex}
\end{equation*} 
and in this case, any of the $cut$ axioms implies all the others.

The PROP of commutative Frobenius algebras is equivalent to the PROP of \emph{2-cobordisms}, in terms of which 2-dimensional topological quantum field theories are defined; see \cite{kock2003frobenius} for a review of the subject.
\end{remark}

A rich strand of categorical universal algebra, that we have not addressed so far, is based on monads, instead of PROs, as embodiments of algebraic theories. 

There is no clear hierarchy of generality between the two: on the one hand, the theory of monads on suitably small categories can be internalised in the framework of PROs, since monads are models of $\textit{Mon}$ in the 2-category $\textit{Cat}$ of categories, functors and natural transformations. On the other hand, if monads are taken as an ``external'', fundamental notion, the same information contained in a PRO can be encoded as a monad in a certain bicategory, as shown in \cite{lack2004composing}. We will say a little more about this in Section \ref{sec:sliding}. 

Another strand is based on \emph{operads}, see for example \cite{markl2007operads}. However, it is not our goal to give an exhaustive survey of categorical approaches to universal algebra; for that, we refer to better-suited sources such as \cite{hyland2007category, markl2008operads, berger2012monads}.

One advantage of the PRO framework over others is that the condition that a PRO be a 2-category is a forced one to begin with; to obtain a higher-dimensional generalisation, it suffices to remove it.
\begin{dfn}
An \emph{$\omega$PRO} is an $\omega$-category $X$ such that $\skel{1}{X} \simeq \vec{S}^1$. If $X$, $Y$ are $\omega$PROs, with 1-generators $a_X$, $a_Y$, a map $f: X \to Y$ is a map of $\omega$-categories such that $f(a_X) = a_Y$. \index{omegaPRO@$\omega$PRO} 
\end{dfn}
In $\omega$PROs, identities between 2-cells can be replaced with 3-cells, which in turn can have 4-cells between them, and so on; presentations of PROs are a particular example of $\omega$PROs. There are at least two areas of mathematics where this is of interest.
\begin{enumerate}
\item In \emph{higher algebra}, the focus is on algebras arising in homotopy theory, or more generally in categories with weak equivalences, so all axioms are required to hold not exactly, but ``up to higher homotopies''; \cite{stasheff1970hspaces, boardman1973homotopy} are seminal works, and \cite{lurie2017higher} is a recent, extensive development of the subject. In this case, the higher cells that replace identities are usually equivalences; it is then common to speak of \emph{pseudo}-algebras. A key word in this field is \emph{coherence}: roughly, the possibility of ``strictifying'' some, or all of the equivalences of a pseudo-algebra to \emph{identities} of another pseudo-algebra which is equivalent to the original one in a suitable sense. See \cite[Section 5]{maclane1976topology} for an early historical overview of coherence problems. \index{coherence}

\item In \emph{higher-dimensional rewriting theory}, the focus is on computational properties of presentations of algebraic theories, such as the possibility of effectively determining whether two operations, expressed as formal composites of basic operations, are equal --- a higher-dimensional version of the classical word problem for monoids. This is the application for which polygraphs were (re)invented in \cite{burroni1993higher}; see \cite{mimram2014towards} for a more recent overview. The higher cells that replace identities are \emph{not} equivalences, and are interpreted as rewrites, then rewrites of rewrites, and so on; one speaks, sometimes, of \emph{lax} algebras. The key words are \emph{confluence}, \emph{termination}, \emph{convergence}, \emph{normal form}; see \cite[Section 1.2]{guiraud2009identities} for their definition in the framework of polygraphs.
\end{enumerate}
The two are not unrelated: for instance, coherence problems for pseudo-algebras can in some cases be solved by reducing them to word problems, see for example \cite{mellies2003maclane}.

The definition of higher algebraic theories is often a product of ``categorification'' \cite{baez1998}, a kind of reverse-engineering: one takes an algebraic theory, and makes an informed guess of what higher-dimensional cells its identities may come from. While there are known ways of finding obstructions to convergence of a given presentation, such as the analysis of its critical pairs, which can suggest what higher cells need to be introduced, with coherence the issue can be more subtle; not having canonical higher-dimensional cells to start with, one is reduced to an unsystematic, case-by-case analysis. Furthermore, direct calculations with higher-dimensional cells can be unwieldy, although in recent times computational aids, such as the \emph{Globular} proof assistant \cite{vicary2016globular}, have been developed.

In the next sections, we will show a potential pathway to improving this situation: by building up higher algebraic theories from simpler ones, using the operations introduced in Chapter \ref{chap:polygraphs}, we automatically obtain a selection of higher-dimensional generators.

\section{Homomorphisms, sliding rules and tensor products} \label{sec:sliding}

\begin{itemize}
\item[]	\textbf{Note.} This section and the following are partially based on material previously published in \cite{hadzihasanovic2017topological}.
\end{itemize}
Algebraic theories come with a notion of homomorphism, a ``structure-preserving map'', such that algebras and their homomorphisms form a category. In the PRO framework, an $X$-algebra in $Y$ is simply a map $f: X \to Y$; thus, by the biclosed structure of $\omegacat$, there are two obvious candidates for an $\omega$-category of $X$-algebras in $Y$: 
\begin{equation*}
\rimp{X}{Y} \quad \text{ and } \quad \limp{X}{Y}.
\end{equation*}
They both come with ``forgetful'' maps to $Y$, in the first case
\begin{equation*}
\begin{tikzpicture}[baseline={([yshift=-.5ex]current bounding box.center)}]
	\node[scale=1.25] (0) at (0,0) {$\rimp{X}{Y} \; \simeq \; 1 \otimes (\rimp{X}{Y})$};
	\node[scale=1.25] (1) at (4.9,0) {$X \otimes (\rimp{X}{Y})$};
	\node[scale=1.25] (2) at (8,0) {$Y$,};
	\draw[1c] (0) to node[auto] {$\bullet \otimes \mathrm{id}$} (1);
	\draw[1c] (1) to node[auto] {$e^R_{X,Y}$} (2);
\end{tikzpicture}
\end{equation*}
where $e^R_{X,Y}$ is the canonical evaluation map; the second case is analogous.

Let us consider $\limp{X}{Y}$. A 1-cell in $\limp{X}{Y}$ --- one candidate for an $X$-algebra homomorphism --- corresponds to a map $\vec{I} \to \limp{X}{Y}$, which in turn corresponds to a map $\vec{I} \otimes X \to Y$. In the other case, homomorphisms correspond to maps $X \otimes \vec{I} \to Y$. 

\begin{remark} 
Continuing the parallel between polygraphs and topological spaces, a homomorphism of $X$-algebras, in any of the two versions, is the directed analogue of a homotopy of maps $X \to Y$.
\end{remark}

In other words, to understand these notions of homomorphism, it suffices to understand the following.
\begin{dfn} \label{dfn:cylinder} \index{cylinder}
Let $X$ be an $\omega$-category. The \emph{left cylinder} and \emph{right cylinder} of $X$ are the $\omega$-categories $\vec{I} \otimes X$ and $X \otimes \vec{I}$, respectively.
\end{dfn}

Their structure can be intuitively understood as follows: $\vec{I} \otimes X$ contains two copies of $X$, $\{0^-\} \otimes X$ and $\{0^+\} \otimes X$; for all $n$-cells $x$ of $X$, it contains an $(n+1)$-cell $\top \otimes x$, which has $0^- \otimes x$ in its input boundary, and $0^+ \otimes x$ in its output boundary --- it ``mediates'' between the two copies at the $n$-dimensional level. The right cylinder is only slightly more complicated, because the direction of the $(n+1)$-cells $x \otimes \top$ with respect to the $x \otimes 0^\alpha$ alternates with the dimension of $x$.

In the literature on PROs and Lawvere theories \cite{lack2004composing, hyland2007category}, a more restrictive notion of homomorphism is usually considered, where models in $Y$ are ``located'' on the same 0-cell of $Y$, and the mediating 1-cell is a unit. These correspond to the following constructions.

\begin{dfn} \index{cylinder!reduced}
Let $X$ be an $\omega$PRO. The \emph{reduced left cylinder} and \emph{right cylinder} of $X$ are the quotient $\omega$-categories $\vec{I} \otimes X/\vec{I} \otimes \{\bullet\}$ and $X \otimes \vec{I}/\{\bullet\} \otimes \vec{I}$, respectively.
\end{dfn}

In fact, we have already explicitly calculated a left cylinder in Example \ref{exm:cylinder}. Generalising the calculation to the presentation of a PRO $X$ with 2-generators $f: [n] \to [m]$, we can picture the 3-cells $\top \otimes f$ as
\begin{equation} \label{eq:left_cylinder}
\input{img/c2_2_left_cylinder.tex}
\end{equation} 
where the grey region corresponds to $0^- \otimes \bullet$ and the yellow region to $0^+ \otimes \bullet$. In this case, because $a$ is the only odd-dimensional generator of $X$, in the right cylinder only $a \otimes \top$ swaps directions, so $f \otimes \top$ can be pictured as
\begin{equation} \label{eq:right_cylinder}
\input{img/c2_2_right_cylinder.tex}
\end{equation} 
Pictures for reduced cylinders are then easily obtained my merging regions: the quotient of $\top \otimes f$ in the reduced left cylinder, for example, is
\begin{equation*}
\input{img/c2_2_reduced_left.tex}
\end{equation*} 
The cells (\ref{eq:left_cylinder}) and (\ref{eq:right_cylinder}) are examples of \emph{sliding rules}: cells that, pictured in string diagrams, can be interpreted as ``sliding'' a diagram past the boundary between two regions of the plane. 

Another way of understanding these pictures is through the following 3-dimensional representation, shown in the case $f: [2] \to [1]$, which is meant to convey the intuition of the product of a 2-dimensional string diagram and the ``1-dimensional string diagram'' 
\begin{tikzpicture} 
	\draw[edge, line width=4pt, color=gray!10] (0,0) to (1,0);
	\draw[edge, line width=4pt, color=yellow] (1,0) to (2,0);
	\node[dotdark] at (1,0) {};
\end{tikzpicture} corresponding to $\top: 0^- \to 0^+$:
\begin{equation*}
\input{img/c2_2_cube.tex}
\end{equation*} 

In \cite{hinze2016equational}, Hinze and Marsden gave a survey of string-diagrammatic reasoning in 2-category theory, where sliding rules are omnipresent, associated to various sorts of ``naturality'' conditions. We can reinterpret this as an abundance of cylinders, or tensor products in general; that is, we can see the tensor product of $\omega$-categories as the meeting point of the algebraic notion of naturality, or structure-preserving morphism, and its topological correlate, the sliding rules.

\begin{exm}
The two notions of homomorphism of monads, that is, left and right cylinders $\vec{I} \otimes \textit{Mon}$, $\textit{Mon} \otimes \vec{I}$ in the 2-category $\textit{Cat}$, correspond to the notions, respectively, of Kleisli law and Eilenberg-Moore law between two monads. As shown in \cite[Section 4]{hinze2016equational}, it is possible to recover the notion of Eilenberg-Moore algebra for a monad $T$ on the category $\cat{C}$ as an Eilenberg-Moore law between $T$ and the identity monad on the terminal category.
\end{exm}

Some of the axioms of PROs that we have already encountered can be reinterpreted as sliding rules: most obviously, $nat^f_{b,L}$ in diagram (\ref{eq:braiding_nat}) is obtained from the cell $\top \otimes f$, diagram (\ref{eq:left_cylinder}), by the quotient defined by
\begin{align*}
	0^\alpha \otimes x & \; \mapsto \; x, \quad \text{for all cells $x$, and $\alpha \in \{+,-\}$}, \\
	\top \otimes \bullet & \; \mapsto \; a, \\
	\top \otimes a & \; \mapsto \; b; 
\end{align*}
the same goes for $nat^f_s$, replacing $b$ with $s$. Similarly, $nat^f_{b,R}$ can be obtained from $f \otimes \top$, diagram (\ref{eq:right_cylinder}).

This suggests that we can use tensor products together with quotients to construct axioms for algebraic theories. 

\begin{exm} \label{exm:braidings}
For $0 \leq k \leq n$, the polygraph 
\begin{equation*}
	\bigotimes^n \vec{S}^1 := \underbrace{\vec{S}^1 \otimes \ldots \otimes \vec{S}^1}_n
\end{equation*} 
has $n \choose k$ $k$-dimensional generators, corresponding to all possible tensor products of $k$ copies of $a$ and $(n-k)$ copies of $\bullet$. Because $\vec{S}^1 \simeq \opp{(\vec{S}^1)}$, we have in this particular case well-defined ``permutation maps'' $\sigma: \bigotimes^n \vec{S}^1 \to \bigotimes^n \vec{S}^1$, sending a generator $x_1 \otimes \ldots \otimes x_n$ to $x_{\sigma(1)} \otimes \ldots \otimes x_{\sigma(n)}$, for some permutation $\sigma$.

Let $(\bigotimes^n \vec{S}^1)_P$ be the coequaliser of all permutations; this is a polygraph with a single $k$-dimensional generator $a_k$ for all $k \leq n$. Clearly, $(\bigotimes^n \vec{S}^1)_P$ is an $\omega$PRO. Moreover, for all $k \leq n$, the $k$-skeleton of $(\bigotimes^n \vec{S}^1)_P$ is isomorphic to $(\bigotimes^k \vec{S}^1)_P$, so we can define an $\omega$PRO $(\bigotimes^\infty \vec{S}^1)_P$ as the colimit of the sequence of inclusions of $n$-skeleta $(\bigotimes^n \vec{S}^1)_P \hookrightarrow (\bigotimes^{n+1} \vec{S}^1)_P$.

Let us examine the low-dimensional generators of $(\bigotimes^\infty \vec{S}^1)_P$. The 2-dimensional generator is a cell $a_2: a \cp{1} a \to a \cp{1} a$, which we picture as a braiding:
\begin{equation*}
\input{img/c2_2_braiding2d.tex}
\end{equation*}
The 3-dimensional generator is a cell
\begin{equation*}
\input{img/c2_2_braiding3d.tex}
\end{equation*}
which can be seen as a presentation of the axiom $nat_b^b$ of $\textit{Brd}$, that is, the third Reidemeister move. The 4-dimensional generator can be pictured as follows:
\begin{equation*}
\input{img/c2_2_braiding4d.tex}
\end{equation*}
This is a presentation of the so-called Zamolodchikov tetrahedron equation; see for example \cite{kapranov1994zamolodchikov, baez2004higher}. From the point of view of rewriting theory, this can be seen as a coherence condition between different sequences of $rei_3$ moves; however, through the interface with the topology of braids given by string diagrams, it can also be seen as a combinatorial presentation of an isotopy of surfaces. The higher-dimensional generators give further coherences, or further isotopies.
\end{exm}

\begin{exm}
In general, given a number of $\omega$PROs $X_1, \ldots, X_n$, their tensor product comes with an inclusion $\bigotimes^n \vec{S}^1 \hookrightarrow \bigotimes^n_{i=1} X_i$; we can then define $(\bigotimes^n_{i=1} X_i)_P$ to be the pushout 
\begin{equation*}
\begin{tikzpicture}[baseline={([yshift=-.5ex]current bounding box.center)}]
	\node[scale=1.25] (0) at (-1.5,2) {$\bigotimes^n \vec{S}^1$};
	\node[scale=1.25] (1) at (2,2) {$\bigotimes^n_{i=1} X_i$};
	\node[scale=1.25] (2) at (-1.5,0) {$(\bigotimes^n \vec{S}^1)_P$};
	\node[scale=1.25] (3) at (2,0) {$(\bigotimes^n_{i=1} X_i)_P$};
	\draw[1cinc] (0) to (1);
	\draw[1c] (0) to (2);
	\draw[1cinc] (2) to (3);
	\draw[1c] (1) to (3);
	\draw[edge] (1,0.2) to (1,0.8) to (1.8,0.8);
\end{tikzpicture}
\end{equation*}
in $\omegacat$. Because the 1-skeleton of $(\bigotimes^n \vec{S}^1)_P$ is isomorphic to $\vec{S}^1$, it follows that $(\bigotimes^n_{i=1} X_i)_P$ is an $\omega$PRO. 

When $X$ is a PROB, it actually comes with a map $(\vec{S}^1 \otimes \vec{S}^1)_P \to X$, sending $a_2$ to the braiding cell of $X$, which can be pre-composed to obtain a map $\vec{S}^1 \otimes \vec{S}^1 \to X$. Then, $X \otimes Y$ comes with a map $\bigotimes^4 \vec{S}^1 \to X \otimes Y$, and we can take $(X \otimes Y)_{P,2}$ to be its pushout along the quotient $\bigotimes^4 \vec{S}^1 \to (\bigotimes^4 \vec{S}^1)_P$. This has the effect of identifying the braiding cells of $X$ and $Y$; it is then easy to verify that the 2-truncation of $(X \otimes Y)_{P,2}$ is itself a PROB. When either $X$ or $Y$ is a PROP, $(X \otimes Y)_{P,2}$ also is.
\end{exm}

To obtain a presentation of $\textit{Brd}$ from $(\bigotimes^n \vec{S}^1)_P$, for $n \geq 3$, it suffices to add an inverse for the 2-dimensional generator. We can obtain this as a particular instance of a general procedure for turning a cell into an equivalence. 

In the following, when we speak of ``attaching an $n$-cell $x: x^- \to x^+$ to an $\omega$-category $X$'', we mean taking the pushout
\begin{equation*}
\begin{tikzpicture}[baseline={([yshift=-.5ex]current bounding box.center)}]
	\node[scale=1.25] (0) at (-1,2) {$\partial G^n$};
	\node[scale=1.25] (1) at (2,2) {$G^n$};
	\node[scale=1.25] (2) at (-1,0) {$X$};
	\node[scale=1.25] (3) at (2,0) {$X'$};
	\draw[1cinc] (0) to (1);
	\draw[1c] (0) to node[auto,swap] {$\partial x$} (2);
	\draw[1cinc] (2) to (3);
	\draw[1c] (1) to node[auto] {$x$} (3);
	\draw[edge] (1,0.2) to (1,0.8) to (1.8,0.8);
\end{tikzpicture}
\end{equation*}
in $\omegacat$, where $\partial x((n-1)^\alpha) := x^\alpha$.

\begin{cons}
Let $X$ be an $\omega$-category, and $C$ a set of non-degenerate cells of $X$. We define $X\{\invrs{C}\}$ to be the $\omega$-category obtained from $X$ by attaching consecutively, for all $x: x^- \to x^+$ in $C$, with $\dmn{x} = n$,
\begin{enumerate}
	\item an $n$-cell $\invrs{x}: x^+ \to x^-$;
	\item a pair of $(n+1)$-cells, $e(x) : x \cp{n} \invrs{x} \to \idd{} x^-$, $h(x) : \idd{} x^+ \to \invrs{x} \cp{n} x$.
\end{enumerate}
Then, we define a sequence $\{X_n \hookrightarrow X_{n+1}\}_{n \in \mathbb{N}}$ of inclusions of $\omega$-categories, as follows. Let $X_0 := X$, $C_0 := C$; then, for all $n > 0$,
\begin{align*}
	X_n & := X_{n-1}\{\invrs{C_{n-1}}\}, \\
	C_n & := \{e(x), h(x) \,|\, x \in C_{n-1}\}.
\end{align*}
We define $X[\invrs{C}]$ to be the colimit of this sequence. By construction, all cells in $C$ are equivalences in $X[\invrs{C}]$. 
\end{cons}

When $C = \{x\}$ for some cell $x$, we write $X[\invrs{x}]$ for $X[\invrs{C}]$.
\begin{remark}
If $X$ is a polygraph and $C$ a set of generators, then $X[\invrs{C}]$ is still a polygraph, with an infinite set of additional generators
\begin{equation*}
	\invrs{C_0} + \coprod_{n>0} (C_n + \invrs{C_n}),
\end{equation*}
where $\invrs{C_n} := \{\invrs{x} \,|\, x \in C_n\}$.
\end{remark}

\begin{remark} \label{remark:nequivalence} \index{equivalence!standard}
We can define the ``standard $n$-equivalence'' $Eq^n$ to be the polygraph $G^n[\invrs{\top}]$, obtained from the standard $n$-globe by weakly inverting its $n$-dimensional generator; equivalence $n$-cells in an $\omega$-category $X$ are in bijection with maps $Eq^n \to X$. Then, we can alternatively obtain $X[\invrs{C}]$ in a single step as the pushout of $\omega$-categories
\begin{equation*}
\begin{tikzpicture}[baseline={([yshift=-.5ex]current bounding box.center)}]
	\node[scale=1.25] (0) at (0,2) {$\displaystyle \coprod_{x \in C} G^{\dmn{x}}$};
	\node[scale=1.25] (1) at (3.5,2) {$\displaystyle\coprod_{x \in C} Eq^{\dmn{x}}$};
	\node[scale=1.25] (2) at (0,0) {$X$};
	\node[scale=1.25] (3) at (3.5,0) {$X[\invrs{C}].$};
	\draw[1c] (0) to node[auto,swap] {$\displaystyle\coprod_{x \in C} x$} (2);
	\draw[1c] (1) to (3);
	\draw[1cinc] (0) to (1);
	\draw[1cinc] (2) to (3);
	\draw[edge] (2.5,0.2) to (2.5,0.8) to (3.3,0.8);
\end{tikzpicture}
\end{equation*}
\end{remark} 

\begin{prop}
For all $n \geq 3$, $(\bigotimes^n \vec{S}^1)_P[\invrs{a_2}]$ is a presentation of $\textit{Brd}$.
\end{prop}
Thus, we have a first example of a compositional reconstruction of an algebraic theory --- a simple one, but the point of interest is the inductive, natively higher-dimensional character: there is no reason why we should truncate at the 2-dimensional level. By weakly inverting any cell that needs to be an equivalence after tensoring with each copy of $\vec{S}^1$, and before quotienting, we can obtain pseudo-versions of $\textit{Brd}$ that also include coherence cell between the second and third Reidemeister move, and their higher-dimensional analogues.

The case of cylinders of homomorphisms invites us to reconsider the boundaries of an algebraic theory: while $\vec{I} \otimes X$, $X \otimes \vec{I}$ are not $\omega$PROs, their models still have a clear algebraic interpretation. It is worth, then, to also consider constructions that yield ``multi-sorted'' theories (more than one 1-cell), or ``multi-coloured'' theories (more than one 0-cell). 

\begin{exm} \label{exm:tensorpro}
When $X$ and $Y$ are PROs, $X \otimes Y$ contains a copy of $X$ and a copy of $Y$ (let $x := x \otimes \bullet$, $y := \bullet \otimes y$) together with an additional braiding-like generator $b := a_X \otimes a_Y : a_Y \cp{1} a_X \to a_X \cp{1} a_Y$, and ``half-braiding'' cells of the same shape as $nat^f_{b,R}$, diagram (\ref{eq:braiding_nat}), for all 2-cells $f$ of $X$, and $nat^g_{b,L}$ for all 2-cells $g$ of $Y$.

In the particular case $X = Y = \textit{Mon}$, models of $\textit{Mon} \otimes \textit{Mon}$ in $\textit{Cat}$ are pairs of monads together with a \emph{distributive law} \cite{beck1969distributive}. The main point is that the axioms of $\textit{Mon} \otimes \textit{Mon}$ are the minimal conditions ensuring that \index{distributive law}
\begin{equation*}
\input{img/c2_2_distributive.tex}
\end{equation*}
satisfy the equations for the multiplication and unit operations of a monoid; thus, a distributive law in $\textit{Cat}$ corresponds to a way of composing two monads and obtaining another monad.

Lack's article \cite{lack2004composing} introduced a framework for composing PROs and PROPs based on distributive laws, through the identification of PROs with monads in a bicategory of spans of monoids. This has since sparked a fruitful line of research; see for example \cite{bonchi2014interacting, bonchi2014hopf, duncan2016interacting}.

While the concept of ``composing algebraic theories'' in this thesis was inspired by Lack's approach, the two seem only vaguely related from a technical standpoint. The applications also appear to be different: the monad composition seems to be particularly useful to derive axioms when concrete models of the component theories are available, suggesting a ``concrete'' distributive law between them; it is also strictly 2-dimensional, that is, it does not produce any higher-dimensional coherence equations. On the other hand, tensor products are easier to calculate for polygraphs, hence for presentations of algebraic theories, so the ``syntactic'' perspective is favoured over the ``semantic''. 
\end{exm}

Next, we want to show how associativity-like axioms can be interpreted as quotiented sliding rules. 
\begin{dfn} \label{dfn:cone} \index{cone} \index{future cone|see {cone}} \index{past cone|see {cone}} \index{conemath@$\cone{+}{-}, \cone{-}{-}$}
Let $X$ be an $\omega$-category. The \emph{future} and \emph{past left cone} of $X$ are the quotient $\omega$-categories $\cone{+}{X} := \vec{I} \otimes X/\{0^+\} \otimes X$ and $\cone{-}{X} := \vec{I} \otimes X/\{0^-\} \otimes X$, respectively.
\end{dfn}
There are similar notions of future and past right cone, obtained as quotients of right cylinders. If $X$ is a PRO, cones can be understood diagrammatically starting from cylinders: for all 2-cells $f: [n] \to [m]$, $\cone{+}{X}$ contains a 3-cell
\begin{equation*}
\input{img/c2_2_future_cone.tex}
\end{equation*} 
and $\cone{-}{X}$ a 3-cell
\begin{equation*}
\input{img/c2_2_past_cone.tex}
\end{equation*} 
In the case $X = \textit{Mon}$, models of $\cone{+}{X}$ in an 2-category $Y$ are 1-cells of $Y$ together with a left action, and a right co-action, respectively, of a monoid in $Y$. Right cones similarly yield right actions and left co-actions. Thus, cones on $X$ can be seen as presentations of generalised theories of ``actions and co-actions of $X$-algebras''. \index{action} \index{co-action}

\begin{exm} \index{oriental} 
Iterated future cones on the terminal $\omega$-category yield Street's oriented simplices, or orientals \cite{street1987algebra}; this was proven in detail in \cite{buckley2016orientals}.

Because tensor products preserve colimits, it is possible to obtain the $n$-oriental directly as a quotient of the ``directed $n$-cube'' $\bigotimes^n \vec{I}$ by the sub-$\omega$-category of all cells of the form $x \otimes 0^+ \otimes x'$. This quotient can be followed by the identification of all $k$-dimensional generators, for all $k \leq n$, similarly to Example \ref{exm:braidings}. We write $\textit{As}^n$ for the $\omega$PRO so obtained, and $as_k$ for its $k$-dimensional generator, $k > 1$; again, we have inclusions $\textit{As}^n \hookrightarrow \textit{As}^{n+1}$ of $n$-skeleta for all $n \in \mathbb{N}$, and we can define an $\omega$PRO $\textit{As}^\infty$ as the colimit of this sequence. \index{as@$\textit{As}^n$}

The generators of $\textit{As}^\infty$ can be understood graphically by taking string diagrams for $(\bigotimes^\infty \vec{S}^1)_P$, and interrupting strands of wire whenever they reach an underpassing. The 2-dimensional generator can be pictured as
\begin{equation*}
\input{img/c2_2_asso2d.tex}
\end{equation*} 
The 3-dimensional generator is
\begin{equation*}
\input{img/c2_2_asso3d.tex}
\end{equation*} 
which is a presentation of the associativity axiom of the theory of monoids. The 4-dimensional generator is
\begin{equation*}
\input{img/c2_2_asso4d.tex}
\end{equation*} 
This is a presentation of Mac Lane's pentagon equation \cite{maclane1963natural}; the way we drew it is meant to enhance the intuition of associativity axioms as ``quotiented sliding rules''. 

The $\textit{As}^n$, for $n \geq 3$, are presentations of the theory of semigroups (monoids without the unit), and higher-dimensional lax generalisations. They can be promoted to presentations of the theory of monoids as follows: take the 2-polygraph $U$ with a single 2-generator $u: [0] \to [1]$, so that $\cone{+}{U}$ contains the 3-generator
\begin{equation*}
\input{img/c2_2_leftunit.tex}
\end{equation*} 
which can be quotiented to an $\omega$PRO $X$ presenting the left unitality axiom $un_L$. Then, $\opp{X}$ presents the right unitality axiom $un_R$, and $\cone{+}{\opp{X}}$ contains a 4-generator of shape 
\begin{equation*}
\input{img/c2_2_triangle.tex}
\end{equation*} 
which can be quotiented to a presentation of Mac Lane's triangle axiom, also containing a left unitality cell and an associativity cell. We can then iterate cones and quotients of generators with the same shape as in the construction of $\textit{As}^n$ in order to generate higher coherence cells.

If we weakly invert all 3-generators, we obtain presentations of the higher algebraic theory of pseudomonoids \cite{street2004frobenius}. If $\textit{Cat}_\times$ is the 3-category whose 1-cells are sequences of suitably small categories, 2-cells $f: (C_1, \ldots, C_n) \to (D_1, \ldots D_n)$ are functors $f: C_1 \times \ldots \times C_n \to D_1 \times \ldots \times D_n$, and 3-cells are natural transformations, then pseudomonoids in $\textit{Cat}_\times$ are small monoidal categories.
\end{exm}

We can observe, in this and other examples, that in the compositional approach to universal algebra it is the ``multi-coloured'' theories that come first, and have simpler combinatorics: so, for instance, actions of 1-cells on 1-cells are fundamental, and monoids are derived as the special case where all the 1-cells involved, and all the operations of the same arity are equal.
\begin{exm} 
We can derive the theory of Frobenius algebras as a special case of the theory of 1-cells with a two-sided action of a monoid, compatible with the two-sided co-action of a comonoid, as follows. Quotient the cylinder $\vec{I} \otimes \textit{Mon} \otimes \coo{\textit{Mon}}$ by the sub-$\omega$-categories $\{0^+\} \otimes \textit{Mon}$ and $\{0^-\} \otimes \coo{\textit{Mon}}$, where $\textit{Mon}$, $\coo{\textit{Mon}}$ are seen as subtheories of $\textit{Mon} \otimes \coo{\textit{Mon}}$ as in Example \ref{exm:tensorpro}. 

The quotient then contains the subtheories $\cone{+}{\textit{Mon}}$ of left actions of monoids and $\cone{-}{\coo{\textit{Mon}}}$ of right co-actions of comonoids, together with a 3-cell 
\begin{equation*}
\input{img/c2_2_frob_left.tex}
\end{equation*} 
obtained as the quotient of $\top \otimes a \otimes a$; the dotted wires mark the 1-cells that have been quotiented out. Identifying as usual all 2-cells of the same shape --- in particular, $0^+ \otimes a \otimes a$ and $0^- \otimes a \otimes a$ become units --- we obtain a presentation of the axiom $frob_L$, diagram (\ref{eq:frobenius}).

A presentation of $frob_R$ is similarly constructed by using right cylinders, or directly as the 1-opposite of the theory just defined; the union of the two is a presentation of $\textit{Frob}$.
\end{exm}

\section{Dimension-raising interactions} \label{sec:dimension}
We have reinterpreted, and reconstructed many of the axioms seen in Section \ref{sec:univ_algebra} as sliding rules, obtained by tensor products and quotients; but the axioms $nat^f_d$ of diagram (\ref{eq:cartesian_nat}), and $nat^f_c$ of diagram (\ref{eq:cartesian_nat2}) still escape this characterisation. They have something in common with the sliding rules we have seen so far --- they also posit the possibility of pulling certain generators through other generators --- but they do not simply consist in ``sliding diagrams past a boundary''.

Let us focus on a specific theory, whose characteristic axioms are of this kind. We also define an important extension for future use. 

\begin{dfn} \label{dfn:bialg} \index{bialgebra} \index{Bialg@$\textit{Bialg}$} \index{Hopf algebra} \index{Hopf@$\textit{Hopf}$}
The theory of \emph{bialgebras} is the PROB $\textit{Bialg}$ containing $\textit{Mon} \uplus \coo{\textit{Mon}}$, and satisfying the additional axioms
\begin{equation*}
\input{img/c2_3_bialgebra.tex}
\end{equation*} 
\begin{equation*}
\input{img/c2_3_bialgebra2.tex}
\end{equation*} 
together with the usual axioms of PROBs. The theory of \emph{Hopf algebras} is the PROB $\textit{Hopf}$ containing $\textit{Bialg}$, with an additional generator $i: [1] \to [1]$, called the \emph{antipode}, satisfying the additional axioms
\begin{equation*}
\input{img/c2_3_hopf.tex}
\end{equation*} 
together with the usual axioms of PROBs.
\end{dfn}

Each of the axioms of bialgebras involves the interaction of a generator of $\textit{Mon}$ and of a generator of $\coo{\textit{Mon}}$, so it is reasonable that they should come from a composition of these two theories. However, $nat^m_c$ also involves a braiding cell, standing in the middle between cells of $\textit{Mon}$ and cells of $\coo{\textit{Mon}}$, not obviously belonging to either side. This ambiguity is also reflected in the standard algebraic explanation of this axiom, which says that the multiplication operation is a comonoid homomorphism, \emph{or equivalently}, the copy operation is a monoid homomorphism: should we privilege one point of view over the other? Moreover, both $\textit{Mon}$ and $\coo{\textit{Mon}}$ are simple PROs, that is, ``planar'' theories, with no braidings: why should their composition be a PROB?

A way of resolving the ambiguity is picturing the axioms with string diagrams in 3-dimensional space:
\begin{equation} \label{eq:bialgebra3d}
\input{img/c2_3_bialgebra3d.tex}
\end{equation} 
Here, the wires of the $\textit{Mon}$ operations and of the $\coo{\textit{Mon}}$ operations open up in orthogonal planes; and the braiding, rather than being an algebraic generator, arises from the geometry of the surrounding space.

So, what we are after is a composition of planar theories $X, Y$ --- PROs with a 1-dimensional ``elementary object'' $a$ --- that yields a three-dimensional theory, whose ``elementary object'' is 2-dimensional. Because of the orthogonality suggested by diagram (\ref{eq:bialgebra3d}), an obvious candidate is the 2-cell $a \otimes a$ in $X \otimes Y$.

We will now describe a simple such composition, and show some compelling evidence that it is conceptually correct; unfortunately, due to the limitations of the framework that we have adopted so far, it will also turn out to be technically wrong.

\begin{dfn} \index{omega-category@$\omega$-category!pointed}
A \emph{pointed $\omega$-category} $(X, \bullet)$ is an $\omega$-category $X$ together with a distinguished 0-cell $\bullet$, called its \emph{basepoint}. A map $f: (X, \bullet_X) \to (Y, \bullet_Y)$ of pointed $\omega$-categories is a map $f: X \to Y$ such that $f(\bullet_X) = \bullet_Y$. Pointed $\omega$-categories and their maps form a category $\omegacat_\bullet$.

Given two pointed $\omega$-categories $(X, \bullet_X), (Y, \bullet_Y)$, their \emph{wedge sum} is the pointed $\omega$-category $(X \lor Y, \bullet)$, where $X \lor Y := X + Y/\{\bullet_X\} + \{\bullet_Y\}$, and $\bullet$ is the 0-cell resulting from the identification of $\bullet_X$ and $\bullet_Y$. \index{wedge sum}

There is an inclusion of $\omega$-categories $X \lor Y \hookrightarrow X \otimes Y$, with $X \mapsto X \otimes \{\bullet_Y\}$ and $Y \mapsto \{\bullet_X\} \otimes Y$. The \emph{smash product} of $(X, \bullet_X)$ and $(Y, \bullet_Y)$ is the pointed $\omega$-category $(X \owedge Y, \bullet)$, where $X \owedge Y := X \otimes Y/X \lor Y$, and $\bullet$ is the image of $X \lor Y$ in the quotient. \index{smash product}
\end{dfn}

Of course, all these operations have topological analogues of the same name, and defined in the same way, as long as all the earlier notions are translated into their undirected counterparts.

The wedge sum and smash product define monoidal structures on $\omegacat_\bullet$, with $(1, \bullet)$ and $(1 + 1, \bullet)$, respectively, as units. Moreover, if the underlying $\omega$-categories of the components are polygraphs, their wedge sum and smash product are also polygraphs: the generators of $X \owedge Y$ are the basepoint, and the images $x \owedge y$ of the generators $x \otimes y$ of $X \otimes Y$ such that $x \neq \bullet_X$ and $y \neq \bullet_Y$.

In the following, we will always consider $\omega$PROs to be pointed with their unique 0-cell as basepoint.

\begin{exm} \index{cylinder!reduced}
Let $\vec{I}_\bullet$ be $\vec{I} + \{\bullet\}$, with basepoint $\bullet$, and let $X$ be an $\omega$PRO. Then the reduced left and right cylinders of $X$ are equal to $\vec{I}_\bullet \owedge X$ and $X \owedge \vec{I}_\bullet$, respectively. In general, we can define reduced cylinders for all pointed $\omega$-categories.
\end{exm}

\begin{exm} \label{exm:suspension}
Given a pointed $\omega$-category $(X, \bullet)$, let its \emph{reduced suspension} be $\Sigma X := X \owedge \vec{S}^1$. \index{reduced suspension}

For all $n$-cells $x$ of an $\omega$-category $X$, $n > 0$, $\Sigma X$ contains an $(n+1)$-cell $x \owedge a$, and $\bord{}{\alpha}(x \owedge a) = \bord{}{\alpha} x \owedge a$. If in particular $X$ is a polygraph with a single 0-cell, such as the presentation of a PRO, the net effect of the suspension is raising the dimension of all other generators by 1. 

In that case, moreover, $\vec{S}^1 \owedge X \simeq X^- \owedge \vec{S}^1$, because $\bord{}{\alpha}(a \owedge x) = a \owedge \bord{}{-\alpha} x$.
\end{exm}

\begin{exm}
Let us calculate $\textit{Mon} \owedge \textit{Mon}$. The simplest way to understand this smash product is to start from string-diagrammatic representations of the generators of $\textit{Mon} \otimes \textit{Mon}$, and then delete all cells of the form $x \otimes \bullet$ or $\bullet \otimes x$.

The only 2-dimensional generator is $a \owedge a$. The 3-dimensional generators are $a \owedge m$, $a \owedge u$, $m \owedge a$, $u \owedge a$. In the terminology of \ref{exm:suspension}, these are either suspensions of the generators of $\textit{Mon}$, which will therefore have the same arity as the original operations, or suspensions of their 2-opposites, which will have the opposite arity. 

String-diagrammatically, this derives from the sliding of a generator past a boundary changing the number of their intersections:
\begin{equation*}
\input{img/c2_3_suspension1.tex}
\end{equation*} 
\begin{equation*}
\input{img/c2_3_suspension2.tex}
\end{equation*} 
\begin{equation*}
\input{img/c2_3_suspension3.tex}
\end{equation*} 
\begin{equation*}
\input{img/c2_3_suspension4.tex}
\end{equation*} 
Another way of picturing these is by drawing surface diagrams, as in \cite{dunn2017surface}, tracing a graph of the sliding of string diagrams of one copy of $\textit{Mon}$ through string diagrams of the other copy, and then only remembering the intersections.

The 4-generators of $\textit{Mon} \owedge \textit{Mon}$ are $m \owedge m$, $m \owedge u$, $u \owedge m$, and $u \owedge u$, whose string-diagrammatic presentation should not come as a surprise: $m \owedge m$ is
\begin{equation*}
\input{img/c2_3_bialgebra_mm.tex}
\end{equation*} 
which in three-dimensional string diagrams becomes
\begin{equation*}
\input{img/c2_3_bialgebra3d_mm.tex}
\end{equation*} 
Similarly, the other three generators yield versions of the other three axioms of bialgebras. 
\end{exm}

\begin{exm} 
Generalising this example, for all generators $f: [n] \to [m]$ of a PRO $X$, the cell $c \owedge \coo{f}$ in $\coo{\textit{Mon}} \owedge \coo{X}$ has a string-diagrammatic presentation of the same shape as the axiom $nat^f_c$ of diagram (\ref{eq:cartesian_nat2}); for instance, $\coo{\textit{Mon}} \owedge \textit{Mon}$ has 4-cells corresponding to the axioms of commutative comonoids. 
\end{exm}

Again, nothing prevents us from applying the same reasoning to higher-dimensional $\omega$PROs: for example, given a presentation $X$ of $\textit{Mon}$ with higher-dimensional cells presenting associativity and unitality axioms, $X \owedge X$ will have higher coherence cells between these and the bialgebra axioms.

But while this indicates that, say, $\textit{Mon} \otimes \textit{Mon}$ does contain all the information needed to construct the axioms of bialgebras, its quotient $\textit{Mon} \owedge \textit{Mon}$ is in fact too degenerate to recover it --- something that we failed to notice when writing \cite{hadzihasanovic2017topological}. This is due to the following, well-known degeneracy of $\omega$-categories.

\begin{prop} \label{prop:eckmann}
Let $X$ be an $\omega$-category with a single 0-cell $\bullet$ and no 1-dimensional cells. Then, for all $n$-cells $x, y$ of $X$, $n > 0$,  $x \cp{1} y = y \cp{1} x$.
\end{prop}
\begin{proof}
This is just a version of the Eckmann-Hilton argument: \index{Eckmann-Hilton argument}
\begin{align*}
	x \cp{1} y & = (x \cp{2} \idd{n} \bullet) \cp{1} ({\idd{n}\bullet} \cp{2} y) = \qquad \qquad \qquad \text{(interchange)} \\
		& = (x \cp{1} \idd{n} \bullet) \cp{2} ({\idd{n} \bullet} \cp{1} y) = x \cp{2} y = \\
		& = ({\idd{n}\bullet} \cp{1} x) \cp{2} (y \cp{1} \idd{n} \bullet) = \qquad \qquad \qquad \text{(interchange)} \\
		& = ({\idd{n} \bullet} \cp{2} y) \cp{1} (x \cp{2} \idd{n} \bullet) = y \cp{1} x. \qedhere
\end{align*} 
\end{proof}

In $\textit{Mon} \owedge \textit{Mon}$, let 
\begin{equation*}
	[n] := \underbrace{(a\owedge a) \cp{1} \ldots \cp{1} (a\owedge a)}_{n}, \quad n>0.
\end{equation*}
Then, the 3-generators $c' := (a \owedge m): [1] \to [2]$, and $m' := (m \owedge a): [2] \to [1]$ should correspond, respectively, to the copy and multiplication operations of the theory of bialgebras. However, 
\begin{align*}
	(m' \cp{3} c') \cp{1} \idd{}[1] & = (m' \cp{3} c') \cp{1} ({\idd{}[1]} \cp{3} {\idd{}[1]}) = & \text{(interchange)} \\
		& = (m' \cp{1} {\idd{}[1]}) \cp{3} (c' \cp{1} {\idd{}[1]}) = & \text{(Proposition \ref{prop:eckmann})} \\
		& = ({\idd{}[1]} \cp{1} m') \cp{3} (c' \cp{1} {\idd{}[1]}) = & \text{(interchange)} \\
		& = ({\idd{}[1]} \cp{3} c') \cp{1} (m' \cp{3} {\idd{}[1]}) = c' \cp{1} m',
\end{align*}
which, in string diagrams projected onto the plane, is
\begin{equation*}
\input{img/c2_3_degeneracy.tex}
\end{equation*} 
clearly, not a part of the theory of bialgebras.

So, the naive quotient of $X \lor Y$ down to a 0-cell is not the correct answer to our question. In particular, $a \otimes \bullet$ and $\bullet \otimes a$ should not become strictly degenerate; rather, they should become \emph{weak units}, of the kind we informally discussed in Section \ref{sec:pasting}. 

A minimal set of conditions for a 1-dimensional cell $i$ to be a weak unit, in an otherwise strict 3-category $X$ with a single 0-cell, that ensures that $X$ will contain braiding 3-cells was given in \cite{joyal2007weak, joyal2013coherence}; we rephrase it in our notation.
\begin{dfn} \index{unit!Joyal-Kock}
Let $X$ be a 3-category with a single 0-cell $\bullet$. For all 1-cells $x$ of $X$, there are maps $- \cp{1} x, x \cp{1} - : X(\bullet,\bullet) \to X(\bullet,\bullet)$ sending an $n$-cell $y$ of $X$ to $\idd{n} x \cp{1} y$ and $y \cp{1} \idd{n} x$, respectively. 

Let $X(\bullet, \bullet)_\mathrm{reg}$ be the sub-2-category of $X(\bullet,\bullet)$ whose 0-cells are 1-dimensional (non-degenerate) cells of $X$. A 1-dimensional cell $i$ of $X$ is a \emph{Joyal-Kock weak unit} if
\begin{enumerate}
	\item $i$ is pseudo-idempotent, that is, there is an equivalence $e: i \cp{1} i \to i$, and
	\item the maps $- \cp{1} i, i \cp{1} -: X(\bullet,\bullet)_\mathrm{reg} \to X(\bullet,\bullet)_\mathrm{reg}$ are weak equivalences of 2-categories. 
\end{enumerate}
For all 1-cells $x$ of $X$, we write $\textit{End}(x)$ for the sub-2-category of $X(\bullet,\bullet)$ whose unique 0-cell is $x$. \index{end@$\textit{End}(-)$}
\end{dfn}

We summarise the main points of the construction. In \cite[Theorem A]{joyal2013coherence}, it is shown that, when $i$ is a weak unit, $e$ is automatically pseudo-associative; that is, there is a map $f: \textit{As}^4 \to X$ where $f(as_1) = i$, $f(as_2) = e$, $f(as_3)$ is an equivalence, and $f(as_4)$ is an identity. If a weak inverse $\bar{e}$ of $e$ is chosen so that is also adjoint to $e$ (this is always possible, by a standard result on 2-categories), $\bar{e}$ will also be coherently co-associative; moreover, $e$ and $\bar{e}$ automatically also satisfy pseudo-Frobenius axioms. 

\begin{remark} \label{remark:spideridem} 
In particular, the part of the 2-skeleton of $X$ relative to the pseudo-idempotency of the weak unit can be given, up to weak equivalence, a ``spider presentation'' as in Remark \ref{remark:spider}, with generators $s_n^m$ for $n, m >0$.
\end{remark}

Moreover, for all 2-cells $x: i \to i$ of $X$, there are equivalence 3-cells
\begin{equation*}
\input{img/c2_3_unitor.tex}
\end{equation*} 
\begin{equation*}
\input{img/c2_3_unitor_rev.tex}
\end{equation*} 
satisfying certain coherence equations, see \cite[Lemma 1.8]{joyal2007weak}. The braiding of two 2-cells $x, y: i \to i$ is the composite of the equivalence 3-cells
\begin{equation*}
\input{img/c2_3_braiding_jk.tex}
\end{equation*} 
In the three-dimensional analogue of a PRO, besides the 2-cells relative to the combinatorics of the weak unit, there should be a unique generator $a: i \to i$; in fact, we can accept a generator $a: i \cp{1} i \to i \cp{1} i$ instead, because $\bar{e} \cp{2} - \cp{2} e: \textit{End}(i \cp{1} i) \to \textit{End}(i)$ is a weak equivalence of 2-categories. 

\begin{prop} \label{prop:prob_from_jk}
Let $X$ be a 3-category such that
\begin{enumerate}
	\item $\skel{1}{X} \simeq \vec{S}^1$, with 1-generator $i$,
	\item there is a 2-cell $a: i \cp{1} i \to i \cp{1} i$, such that the sub-2-category of $X$ generated by $a$ is a polygraph, and
	\item $i$ is a Joyal-Kock weak unit.
\end{enumerate}
Then the 2-category obtained by restricting the 1-cells of $\textit{End}(i)$ to composites of $a$ is a PROB.
\end{prop}
\begin{proof}
Follows from \cite[Proposition 1.3]{joyal2007weak}.
\end{proof}

We can think of the following strategy for obtaining an improved smash product of algebraic theories. Suppose $X$, $Y$ are presentations of PROs, whose 2-generators have arities $[n] \to [m]$ with $n, m > 0$. Then, in the tensor product $X \otimes Y$, we can
\begin{enumerate}
	\item identify $a \otimes \bullet$ and $\bullet \otimes a$, and call the resulting 1-generator $i$;
	\item identify all generators $f \otimes \bullet$, $\bullet \otimes g$ such that $f$, $g$ have the same arity $[n] \to [m]$, and call the resulting generator $s_n^m$;
	\item add enough cells to make $i$ a weak unit for the 3-truncation of the polygraph so constructed, with $s_n^m$ a part of the spider presentation of Remark \ref{remark:spideridem}.
\end{enumerate}
The resulting 3-category, together with the 2-cell $a \otimes a$, would satisfy the conditions of Proposition \ref{prop:prob_from_jk}.

While it could be an interesting problem to formalise this construction --- especially the ``weak unit completion'' --- and see if it can be made functorial in some sense, it also seems quite limited and \emph{ad hoc}. Because operations with 0-ary inputs and outputs are barred, which even excludes presentations of $\textit{Mon}$ (but not the $\textit{As}^n$), we would need to simulate them, by introducing some kind of weak unit already at the 2-categorical level. Higher-dimensional generalisations are also non-obvious, which deprives the compositional approach of its main attraction.

Most of all, it is somewhat disappointing that the simplicity of the smash product should be abandoned due to what appears to be a failure of the $\omega$-categorical combinatorics of composition, rather than the construction itself. 

In truth, this is not the only shortcoming of the $\omega$-category framework, whose critical points --- many of them centred on the degeneracy of Proposition \ref{prop:eckmann} --- we have carefully avoided before this section. In the next chapter, we will examine these issues, and begin delineating a more general strategy for overcoming them, while preserving the main features of our approach: the combinatorial-topological perspective, and the compositionality.
  \thispagestyle{empty} 
\chapter{Towards directed spaces} \label{chap:directed}
\thispagestyle{plain}

\noindent\emph{In this chapter:}
\begin{itemize}
	\item[$\triangleright$] We review the critical points in the combinatorics of $\omega$-categories that lead to the degeneracy of smash products, and motivate the quest for weak higher categories. We formulate four \emph{desiderata} for a new notion of directed space, that could serve as an improved foundation for compositional universal algebra. --- \emph{Section \ref{sec:strictness}}
	\item[$\triangleright$] In partial response to the \emph{desiderata}, we introduce globular posets, directed analogues of incidence posets of regular CW complexes, and develop their basic theory. We prove that geometric realisations of globular posets are regular CW complexes, and that the category of globular posets and inclusions admits a monoidal structure and a monoidal functor to the category of Steiner's augmented directed complexes. We then define a regular polygraph as a space locally modelled on globes, the directed analogues of regular CW disks. --- \emph{Section \ref{sec:globularposet}}
	\item[$\triangleright$] As a proof of concept, we show how to develop the theory of weak 2-categories within regular polygraphs satisfying representability properties. We introduce adequate notions of unit cells, and prove that the existence of units is equivalent to the existence of cells satisfying certain divisibility properties; the results raise some thought-provoking questions on the relation between units and equivalences. Finally, we outline a strategy for higher dimensions. --- \emph{Section \ref{sec:weakness}}
\end{itemize}

\section{The problem with strictness} \label{sec:strictness}

At the beginning of this thesis, we adopted the perspective that polygraphs are a notion of directed space, and can be effectively handled as such, exploiting the compositionality inherent to geometric objects, even when applied to apparently un-topological fields like universal algebra and rewriting. Now, we need to partially revise this viewpoint. \index{directed space}

This is not to say that polygraphs are not, generally speaking, a notion of space: given their status in the folk model structure on $\omegacat$ of \cite{lafont2010folk}, where they play the role of CW complexes in the classical model structure on topological spaces, if one considers a space to be simply an ``object of abstract homotopy theory'', polygraphs clearly fulfill the requirement. Rather, we need to reconsider their adequateness as the foundation of ``compositional higher-dimensional algebra''.

In a way, already the complexity involved in the definition of the tensor product was troubling: such a fundamental operation for the applications of Chapter \ref{chap:interacting}, yet it forced us to take a detour through the entirely separate combinatorics of chain complexes. It is with smash products, though, that we reached an \emph{impasse}. To understand why, we need to ask: what do we expect from a good notion of directed space?

We borrowed the term from directed algebraic topology, a field of research that has found applications in concurrency theory and related areas; see the monographs \cite{grandis2009directed, fajstrup2016directed}. The main notions of directed space employed there are point-set topological spaces with some additional structure, such as local partial orderings of their points, which do not suit our focus on presentations of algebraic theories, and the combinatorial \emph{penchant} it entails. On the other hand, these directed spaces rely on an uncontroversial underlying notion of undirected space, so perhaps we need to ask, first: what is a good combinatorial notion of space? 

One very general, category-theoretic answer, originally due to Grothendieck and based on ideas that are discussed in depth by Lawvere in \cite{lawvere1992categories}, states that there is a notion of space associated to any category $\cat{S}$ of ``probes'', or ``elements of space'': an $\cat{S}$-space is a mathematical object that is ``tested'' by the probes in $\cat{S}$, that is, is determined by how the elements of space fit into it. For instance,
\begin{itemize}
	\item sets are spaces tested by points,
	\item graphs are spaces tested by vertices and edges,
	\item smooth spaces are spaces tested by copies of $\mathbb{R}^n$,
\end{itemize}
and so on. Formally, this information is given simply by a presheaf on $\cat{S}$; this can be supplemented by the information of how elements of space can be glued together to obtain another element, which corresponds to making $\cat{S}$ a site: a space is then a sheaf on this site. This leads us to a first \emph{desideratum} for our hypothetical definition.

\begin{requir}
	The category of directed spaces is a Grothendieck topos.
\end{requir}

Or, at least, it should fit into one, after relaxing some constraints on what objects qualify as spaces; consider, for example, the relation between topological spaces and simplicial sets. Given that polygraphs generalise graphs, one would expect them to fulfill this condition: a polygraph should be a space tested by all possible shapes of cells, whose borders are well-formed pasting diagrams of lower-dimensional shapes. 

However, that is not the case.
\begin{prop}
The category $\cat{Pol}$ is not cartesian closed; in particular, it is not a Grothendieck topos.
\end{prop}
\begin{proof}
This is \cite[Theorem 3.2]{makkai2008category}.
\end{proof}
In fact, the restriction of $\cat{Pol}$ to 2-polygraphs is a presheaf category, and the property is lost from dimension 3 onwards. 

A direct proof by Cheng \cite{cheng2012direct} points specifically to the simplest polygraph to which Proposition \ref{prop:eckmann} applies non-trivially --- the one with a single 0-dimensional generator, and two 2-dimensional generators $x, y$ --- as the crux of the problem: roughly, because $x \cp{1} y$ is equal to $y \cp{1} x$, given a 3-dimensional cell $z$ that has $x \cp{1} y$ as its input or output boundary, it is impossible to recover $x$ or $y$ as sub-cells of $z$; they ``cannot be told apart'' in the composite. In other words, polygraphs lack a notion of subspace that is compatible, in the way we would expect it to be, with the $\omega$-categorical boundaries and compositions.

Another reasonable expectation, given our setup, is that directed spaces should generalise CW complexes: so any combinatorial description of a CW complex should yield a particular directed space, and it should be possible to study the homotopy theory of CW complexes as that of ``undirected objects'' in the category of directed spaces, with the equivalence witnessed by a geometric realisation functor. 

\begin{requir}
	Combinatorial descriptions of CW complexes are descriptions of directed spaces. There exists a realisation functor from the category of directed spaces to the category of topological spaces, inverting this construction at least up to weak equivalence.
\end{requir}

The natural notion of undirected cell in the context of $\omega$-categories is an equivalence cell, as in Definition \ref{dfn:equivalence}, leading to the following definition of undirected object.
\begin{dfn} \index{omega-groupoid@$\omega$-groupoid} \index{omegagpd@$\omega\cat{Gpd}$}
An \emph{$\omega$-groupoid} is an $\omega$-category whose cells are all equivalences. We write $\omega\cat{Gpd}$ for the full subcategory of $\omegacat$ on $\omega$-groupoids. An $\omega$-groupoid is an $n$-groupoid if it is an $n$-category; $n$-groupoids form a full subcategory $n\cat{Gpd}$ of $\omega\cat{Gpd}$.
\end{dfn}
Given an $\omega$-category $X$, for all 0-cells $x$ of $X$, we define a sequence of monoids
\begin{equation*}
	\pi_n(X,x) := \textit{End}_{\tau_{\leq n}X}(\idd{n-1} x), \qquad n > 0,
\end{equation*}
When $X$ is an $\omega$-groupoid, these are in fact groups, because all $n$-cells in $\tau_{\leq n}X$ are invertible. We also define a set $\pi_0(X) := \tau_{\leq 0}X$. The following is \cite[Definition 4.1.1]{simpson2009homotopy}.

\begin{dfn} \index{realisation} \index{geometric realisation|see {realisation}}
A \emph{realisation functor} for $n$-groupoids is a functor $\mathfrak{R}: n\cat{Gpd} \to \cat{Top}$, such that there exist a natural transformation $r_X: \skel{0}{X} \to \mathfrak{R}X$, a natural isomorphism $\zeta^0_X: \pi_0(X) \to \pi_0(\mathfrak{R}X)$ and, for $0 < k \leq n$, natural isomorphisms $\zeta^k_{X,x} : \pi_k(X,x) \to \pi_k(\mathfrak{R}X, r_X(x))$, such that, for all 0-cells $x$ of $X$, the function $\zeta^0_X$ maps $x$ to the connected component of $r_X(x)$.
\end{dfn}
Requiring natural isomorphisms for all $n > 0$, we also obtain a notion of realisation functor of $\omega$-groupoids.

\begin{prop}
For each $n$, there exists a realisation functor for $n$-groupoids.
\end{prop}
\begin{proof}
Such a functor is constructed in \cite{kapranov1991infty} for a different definition of $n$-groupoid, shown in \cite[Theorem 2.2.1]{simpson2009homotopy} to be equivalent to the one given here.
\end{proof}
Given the incidence poset of a regular CW complex, we can build a polygraph $X$ in $\omega\cat{Gpd}$ by progressively attaching standard $n$-equivalences rather than $n$-globes, as in Remark \ref{remark:nequivalence}; from this perspective, given a reasonable realisation functor, we would expect $\mathfrak{R}X$ to be weakly equivalent to the original CW complex. Unfortunately, as shown by Simpson, this is doomed to fail.

\begin{prop} \label{prop:no2sphere}
Given any realisation functor $\mathfrak{R}$ for $\omega$-groupoids, there exists no $\omega$-groupoid $X$ such that $\mathfrak{R}X$ is weakly equivalent to the 2-sphere $S^2$.
\end{prop}
\begin{proof}
Any realisation functor for $\omega$-groupoids, precomposed with $\tau_{\leq 3}$, induces a realisation functor for 3-groupoids. This is then the content of \cite[Theorem 4.4.2]{simpson2009homotopy}.
\end{proof}
In fact, the impossibility of realising a space as $\mathfrak{R}X$ for some $\omega$-groupoid $X$ extends to all spaces whose homotopy type has non-trivial Whitehead products \cite{whitehead1941adding}, and the source of failure can be traced back to Proposition \ref{prop:eckmann}: two 2-cells $x, y: \idd{} \bullet \to \idd{} \bullet$ commute strictly in an $\omega$-category, but the commutation of two 2-cells is generally non-trivial in the homotopy sequence of a topological space, as witnessed by their Whitehead product; see the discussion at the beginning of \cite[Section 4]{simpson2009homotopy}.

This is another \emph{desideratum} failed by polygraphs and $\omega$-categories, and among the main reasons why the focus, in homotopical algebra, is not on $\omega$-categories, but \emph{weak} higher categories --- in which all, or some of the axioms of associativity, unitality, interchange hold not strictly, but ``up to homotopy'', in some sense.

Notions of weak higher category abound, and it is not our intention to give a detailed overview; \cite{leinster2004higher}, \cite{cheng2004higher} and the already cited \cite{simpson2009homotopy} are excellent surveys, with a focus on higher algebra, low-dimensional intuition, and homotopy theory, respectively. In brief, the approaches can usually be sorted into one of two families, called the algebraic and non-algebraic in \cite{cheng2004higher}; we could also frame them as the \emph{external} and \emph{internal} approach.
\begin{itemize}
\item In the algebraic, or external approach, higher algebraic theories, including some concept of homotopy and axiom up to homotopy, are taken as fundamental, and pre-existing to the notion of weak higher category. Compositions and units are operations of such a higher algebra.
\item In the non-algebraic, or internal approach, only cells of different shapes exist as structure: ``being a composite of other cells'' and ``being a unit'' are properties of cells, witnessed by other, higher-dimensional cells. The existence of enough composite and unit cells is a property of certain objects within a larger category of spaces.
\end{itemize}
These tend to be formalised, respectively, as ``presheaves with structure'' and ``presheaves with properties'' on some category $\cat{S}$ of shapes of cells. \index{shape category}

In the first approach, the combinatorics of shapes can be kept simple --- usually, the underlying presheaves are globular sets, as in the $\omega$-categorical case --- and all complications are left to the higher algebraic definitions of compositions and units. In the second approach, the combinatorics of shapes must be expressive enough, in combination with the internal notion of homotopy, to distinguish the ``structural'' cells, witnessing compositions, units and their axioms, from the rest.

Mirroring these two approaches, there are two different attitudes as to the role of polygraphs in higher category theory, and how to amend their definition in the transition to weak higher categories.
\begin{itemize}
\item Polygraphs are presentations of $\omega$-categories: thus, corresponding to a notion of weak higher category, one needs a notion of weak polygraph to present it. This is the route taken by \cite{batanin1998computads}, where categories of polygraphs --- some of them presheaf categories --- are defined relative to different monads on $\cat{Glob}$, embodying different algebras of composition.
\item Polygraphs are like the presheaves on a shape category in the non-algebraic approaches: combinatorial complexes of cells, whose shapes happen to be specified algebraically. The issue is then that the algebraic combinatorics of shapes are too unrestricted, leading to the degeneracy which makes $\cat{Pol}$ not a presheaf category. The theory of weak higher categories should be formulated internally to some restriction of $\cat{Pol}$.
\end{itemize}
It is due remarking that, while from the first point of view it would seem that the theory of polygraphs is ``not expressive enough'', from the second point of view it is even too expressive: consider the works of Steiner, who has rephrased in the language of augmented directed complexes basically all the combinatorics of shapes commonly in use \cite{steiner2004omega, steiner2007orientals, steiner2007simple, steiner2012opetopes}. 

In the same way, one should be cautious when interpreting Proposition \ref{prop:no2sphere}: it is not necessarily the case that one cannot realise topological spaces of all homotopy types as geometric realisations of polygraphs. Rather, it is the naive translation of CW complexes and their homotopy theory into the theory of polygraphs that fails.

While both approaches come with specific advantages, our interest in polygraphs as a foundational tool in higher algebra makes the choice particularly easy: it would be somewhat perverse to rely on a different foundation of higher algebra, just to be able to define our basic objects. Thus, we will pursue a non-algebraic approach, restricting the shapes of cells of polygraphs. Having explained why the restriction is needed, that is, where the unrestricted polygraphs fail, we now have to consider where they succeed, and what we do not want to lose in the process.

\begin{requir} The tensor product of directed spaces can be efficiently calculated.
\end{requir}
This, in fact, would be an improvement over polygraphs, where the efficiency is confined to the loop-free case; but given the importance of the construction, it is not an unreasonable requirement. 

If directed spaces were to form a presheaf category, there may be one canonical way of endowing them with a tensor product. If the shape category $\cat{S}$ is monoidal, the category of presheaves on $\cat{S}$ also comes with a monoidal product, the Day convolution product \cite{day1970closed}: for all $X, Y: \opp{\cat{S}} \to \cat{Set}$, this is defined by the coend \index{convolution!Day}
\begin{equation*}
	X \otimes Y := \int^{s_1, s_2} X(s_1) \times Y(s_2) \times \mathrm{Hom}_\cat{S}(-,s_1 \otimes s_2),
\end{equation*}
where $s_1, s_2$ range over the objects of $\cat{S}$.

Of the shape categories commonly in use --- standard globes, simplices, cubes, Batanin cells, opetopes --- only cubes have a monoidal structure that, indeed, corresponds to the tensor product of polygraphs, after a translation in the style of Steiner. Mainly because of this, and its usefulness in developing a directed homology theory, cubical sets are the one combinatorial notion of directed space that is used in \cite{grandis2009directed}.

On the other hand, while some work has been done on cubical sets as a setting for weak higher category theory from the algebraic side, see \cite{grandis2007higher,kachour2017aspects}, the non-algebraic side remains, as far as we know, underdeveloped. We hypothesise that this is due to a certain mismatch between the cubical shapes, where everything is perfectly symmetrical and coming in powers of two, and the need to express the composition of cells, which is an intrinsically asymmetrical, many-to-one operation, through the combinatorics of higher-dimensional cells.

In any case, all the simplest shape categories are ruled out by our final, and most restrictive \emph{desideratum}.
\begin{requir} \label{requir:string} String diagrams can be directly interpreted in a directed space.
\end{requir}
Chapter \ref{chap:interacting} only gives a sample of the invaluable aid to intuition given by string diagrams in reasoning with 2-dimensional and 3-dimensional algebraic theories, which, thanks to tools like \emph{Globular}, could be brought to higher-dimensional ones as well; that is not something that we want to give up. 

This rules out all category of shapes with an upper bound on the number of inputs or outputs of a 2-cell, that is, all of the ones commonly used except Batanin cells, which, however, failed the previous requirement. What it does not rule out is disallowing 0-ary, or more generally degenerate boundaries: under the interpretation of string diagrams as pasting diagrams (and not their duals), as in Remark \ref{remark:strings}, the intuitive fact that string diagrams, including the ``surrounding space'', are square regions of a plane means that they have non-degenerate boundaries. On the other hand, in order to model, for instance, 0-ary operations of an algebraic theory, it will be necessary to have weak unit cells of the appropriate dimensions.

Given that 2-cells with 0-dimensional boundary seem to be the source of most technical problems of $\omega$-categories, it is to be expected that some restriction to degenerate boundaries should be involved. The opetopic approach \cite{baez1998higher,cheng2004weak} addresses the issue by permitting unrestricted input boundaries (including degenerate ones), but constraining the output boundaries to consist of a single, non-degenerate cell. Our \emph{desiderata}, instead, lead us to considering symmetric constraints, which means barring degeneracy on both sides. This comes with its own set of difficulties: in the opetopic approach, weak units are obtained from cells with degenerate input, so we know that we will need an alternative construction. \index{opetope}

There is a potential candidate for a class of shapes, satisfying many of our requirements, that we have already seen in Chapter \ref{chap:polygraphs}: it is the loop-free globes of Definition \ref{dfn:loopfree-globe}. They are, after all, closed under tensor products, they have no lower bound on the number of input and output cells, and any cell in the boundary of an $n$-globe is automatically an $(n-1)$-globe. 

However, loop-freeness is a strong, global condition, and while it is harmless in 2 dimensions, in 3 dimensions it bars shapes that should be valid, as in the example, shown in \cite{power1991pasting},
\begin{equation} \label{eq:loopglobe}
\input{img/c3_1_loopglobe.tex}
\end{equation} 
where both the input and the output boundary are non-degenerate, loop-free, and regular as CW decompositions of 2-disks, yet $x$ and $\bar{x}$ form a loop in the overall 3-cell. Such a cell may appear, for example, in a presentation of the theory of Frobenius algebras: the left hand sides of the axioms $frob_L$ and $frob_R$, diagram (\ref{eq:frobenius}), which are equal in $\textit{Frob}$, are shaped as the input and output of diagram (\ref{eq:loopglobe}).

In the next section, we will try to expand this class to one that may capture our intended shapes, or at least bring us closer to the goal. The specific kind of non-degeneracy that we want from an $n$-dimensional cell is having a geometric realisation \emph{homeomorphic} to an $n$-disk, and input and output boundaries homeomorphic to $(n-1)$-disks; to achieve that, we go directly to the source of our initial analogy, adopting ideas from combinatorial algebraic topology, related, in particular, to the classification of incidence posets of regular CW complexes.

We call the directed spaces so obtained \emph{regular polygraphs}. The development of a relative non-algebraic notion of weak higher category is still in its early stages. Nevertheless, in the final section, we will show how bicategories, that is, weak 2-categories can be reconstructed in this framework, and propose a pathway to generalisation.

\section{Globular posets and regular polygraphs} \label{sec:globularposet}

In this section, we mention some notions of combinatorial topology without actively using them; we refer to \cite{wachs2006poset, kozlov2008combinatorial} for more details. We recall, however, the following fundamental definitions, which can all be found in \cite{wachs2006poset}. 

\begin{dfn} \index{dimension!in graded posets} \index{poset!graded} 
Let $X$ be a finite poset with order relation $<$, and let $X_\bot$ be $X$ extended with a least element $\bot$. We say that $X$ is \emph{graded} if, for all $x \in X$, all paths from $x$ to $\bot$ in the Hasse diagram $HX_\bot$ have the same length, that is, traverse the same number of edges. If $X$ is graded, for all $x \in X$, let $n$ be the length of paths from $x$ to $\bot$. Then, we define $\dmn{x} := n-1$, the \emph{dimension} of $x$, and write $X_n := \{x \in X \,|\, \dmn{x} = n\}$.

We say that a subset $U$ of a poset $X$ is \emph{closed} when, for all $x, y \in X$, if $y \in U$ and $x \leq y$, then $x \in U$. Given any subset $U$ of $X$, its \emph{closure} is the closed subset $\mathrm{cl}(U) := \{x \in X \,|\, \exists y \in U \, x \leq y \}$. For all $x \in X$, let $U_x := \mathrm{cl}\{x\}$. \index{poset!subsets!closed} 

We write $\dmn{U} := \mathrm{max}\{\dmn{x} \,|\, x \in U\}$; in particular, $\dmn{U_x} = \dmn{x}$. We say that $U$ is \emph{pure} if all maximal elements of $U$ have dimension $n = \dmn{U}$, or, equivalently, if $U = \mathrm{cl}(U \cap X_n)$. \index{poset!subsets!pure}

For all $x, y \in X$ such that $x \leq y$, the \emph{interval} from $x$ to $y$ is the subset $[x,y] := \{z \in X \,|\, x \leq z \leq y\}$. If $X$ is graded, all paths from $y$ to $x$ in $HX$ have length $\dmn{y} - \dmn{x}$; we call this the \emph{length} of $[x,y]$. A graded poset $X$ is \emph{thin} if all intervals $[x,y]$ of length 2 in $X_\bot$ contain precisely 4 elements, that is, they are of the form \index{poset!thin} \index{interval} 
\begin{equation*}
\begin{tikzpicture}[baseline={([yshift=-.5ex]current bounding box.center)}]
	\node[scale=1.25] (0) at (0,2) {$y$};
	\node[scale=1.25] (1) at (-1,1) {$z_1$};
	\node[scale=1.25] (1b) at (1,1) {$z_2$};
	\node[scale=1.25] (2) at (0,0) {$x$};
	\draw[1c] (0) to (1);
	\draw[1c] (0) to (1b);
	\draw[1c] (1) to (2);
	\draw[1c] (1b) to (2);
\end{tikzpicture}
\end{equation*}
in the Hasse diagram $HX_\bot$.
\end{dfn}

To proceed from the undirected to the directed case, we want to endow posets with an orientation, as in Definition \ref{dfn:orientation}.
\begin{dfn} \index{poset!oriented} \index{oriented poset|see {poset}}
Let $X$ be a finite poset with an orientation $o$. We say that $X$ is an \emph{oriented thin poset} if $X$ is thin, and the orientation is compatible in the following sense. Extend $o$ to $X_\bot$ by setting $o(c_{x,\bot}) := +$ for all $x$ of dimension 0. Then, for all intervals $[x,y]$ of length 2 in $X_\bot$, the labelling
\begin{equation} \label{eq:signed}
\begin{tikzpicture}[baseline={([yshift=-.5ex]current bounding box.center)}]
	\node[scale=1.25] (0) at (0,2) {$y$};
	\node[scale=1.25] (1) at (-1,1) {$z_1$};
	\node[scale=1.25] (1b) at (1,1) {$z_2$};
	\node[scale=1.25] (2) at (0,0) {$x$};
	\draw[1c] (0) to node[auto,swap] {$\alpha_1$} (1);
	\draw[1c] (0) to node[auto] {$\alpha_2$} (1b);
	\draw[1c] (1) to node[auto,swap] {$\beta_1$} (2);
	\draw[1c] (1b) to node[auto] {$\beta_2$} (2);
\end{tikzpicture}
\end{equation}
must satisfy $\alpha_1\beta_1 = -\alpha_2\beta_2$. 
\end{dfn}
Eventually, we want to interpret $X$ as the oriented incidence poset of a polygraph: $n$-dimensional elements of $X$ are $n$-dimensional generators, and the labelling of $HX$ specifies whether an $(n-1)$-generator is in the input or in the output boundary of an $n$-generator. 

Thinness is part of a combinatorial criterion for a poset to be a CW poset, that is, the incidence poset of a regular CW complex, see \cite[Proposition 2.2]{bjorner1984posets}: it is a local condition, which essentially imposes that the cells be manifold-like by ruling out irregular situations, such as three 1-cells in the boundary of a 2-cell meeting at a point. Compatibility of the orientation will be necessary for globularity (case $\alpha_1 \neq \alpha_2$), and composability of cells in the same half of the boundary (case $\alpha_1 = \alpha_2$). \index{CW poset}

\begin{dfn} \index{boundary!in oriented posets} \index{$\bord{}{+}, \bord{}{-}$} \index{$\sbord{}{+}, \sbord{}{-}$}
Let $X$ be an oriented thin poset, $U \subseteq X$ closed, $\dmn{U} = n$. For $\alpha \in \{+,-\}$, let
\begin{align*}
	\sbord{}{\alpha} U & := \{x \in U \,|\, \dmn{x} = n-1 \text{ and, for all $y \in U$, if $y$ covers $x$, then $o(c_{y,x}) = \alpha$} \}, \\
	\bord{}{\alpha} U & := \mathrm{cl}(\sbord{}{\alpha} U) \cup \{ x \in U \,|\, \text{for all $y \in U$, if $x \leq y$, then $\dmn{y} < n$} \}.
\end{align*}
In particular, when $U$ is pure, $\bord{}{\alpha} U = \mathrm{cl}(\sbord{}{\alpha} U)$. We will also write $\bord{}{}U := \bord{}{+}U \cup \bord{}{-} U$.
\end{dfn}
The remaining condition for a poset to be a CW poset is a version of shellability for all subsets $U_x$ \cite{bjorner1993lexico}; this is a global condition, preventing cells from having globally non-spherical (for example, toroidal) boundaries. We introduce a directed variant: a kind of sequential, pairwise composability of cells in the boundary of another cell.

\begin{dfn} \index{globe!in oriented posets}
Let $X$ be an oriented thin poset. The class of \emph{globes} in $X$ is defined inductively on dimension and number of maximal elements, as follows.
\begin{itemize}
	\item For all $x \in X$, such that $\dmn{x} = 0$, the subset $\{x\}$ is a 0-dimensional globe.

	\item For all $x \in X$, such that $\dmn{x} = n > 0$, the subset $U_x$ is an \emph{atomic} $n$-globe if $\bord{}{\alpha} U_x$ is an $(n-1)$-dimensional globe, $\alpha \in \{+,-\}$. \index{globe!atomic}

	\item Given two closed, pure, $n$-dimensional $U, U' \subseteq X$, we say that $U$ and $U'$ are \emph{mergeable} if \index{poset!subsets!mergeable}
\begin{enumerate}
	\item $U \cap U' = \bord{}{\alpha}U \cap \bord{}{-\alpha} U'$ for some $\alpha \in \{+,-\}$;
	\item $U \cap U'$ is an $(n-1)$-dimensional globe;
	\item $\bord{}{\beta} (U \cup U')$ is an $(n-1)$-dimensional globe, for $\beta \in \{+,-\}$.
\end{enumerate}
Then, a pure, closed $n$-dimensional $U \subseteq X$ is an $n$-globe if it is atomic, or if there exists a non-trivial bi-partition $\{x_{1,1}, \ldots, x_{1,p}\}$, $\{x_{2,1}, \ldots, x_{2,q}\}$ of its $n$-dimensional elements such that 
\begin{equation*}
	U_1 := \mathrm{cl}\{x_{1,1},\ldots, x_{1,p}\} \quad \text{ and } \quad U_2 := \mathrm{cl}\{x_{2,1}, \ldots, x_{2,q}\}
\end{equation*}
are mergeable $n$-globes.
\end{itemize}
\end{dfn}

\begin{dfn}
An oriented thin poset $X$ is a \emph{globular poset} if, for all $x \in X$, the subset $U_x$ is a globe in $X$. \index{poset!globular} \index{globular poset|see {poset}}
\end{dfn}

By unravelling the induction in the definition of non-atomic $n$-globes, we find that it amounts to requiring that for all $n$-dimensional elements $x \in U$, the subset $U_x$ is an atomic $n$-globe, and there exists a binary tree of bi-partitions of $U \cap X_n$, whose leaves are labelled by singletons, and such that at each branching, the closures of the sets of leaves at each branch are mergeable.

\begin{exm} \label{exm:mergeable}
The following pasting diagram represents a 2-dimensional globe:
\begin{equation*}
\input{img/c1_2_mergeable.tex}
\end{equation*} 
the following are both valid trees of bi-partitions of its 2-dimensional elements:
\begin{equation*}
\begin{tikzpicture}[baseline={([yshift=-.5ex]current bounding box.center)}, yscale=1.2]
\begin{scope}
	\node[scale=1.25] (0) at (0,2) {$\{x_1, x_2, x_3\}$};
	\node[scale=1.25] (1) at (-1,1) {$\{x_1, x_2\}$};
	\node[scale=1.25] (1b) at (1,1) {$\{x_3\}$};
	\node[scale=1.25] (2) at (-2,0) {$\{x_1\}$};
	\node[scale=1.25] (2b) at (0,0) {$\{x_2\}$};
	\draw[1c] (0) to (1);
	\draw[1c] (0) to (1b);
	\draw[1c] (1) to (2);
	\draw[1c] (1) to (2b);
	\node[scale=1.25] at (1.5,-.1) {,};
\end{scope}
\begin{scope}[shift={(5,0)}, xscale=-1]
	\node[scale=1.25] (0) at (0,2) {$\{x_1, x_2, x_3\}$};
	\node[scale=1.25] (1) at (-1,1) {$\{x_2, x_3\}$};
	\node[scale=1.25] (1b) at (1,1) {$\{x_1\}$};
	\node[scale=1.25] (2) at (-2,0) {$\{x_3\}$};
	\node[scale=1.25] (2b) at (0,0) {$\{x_2\}$};
	\draw[1c] (0) to (1);
	\draw[1c] (0) to (1b);
	\draw[1c] (1) to (2);
	\draw[1c] (1) to (2b);
	\node[scale=1.25] at (-2.5,-.1) {.};
\end{scope}
\end{tikzpicture}
\end{equation*}
\end{exm}

\begin{remark} \label{remark:stronger}
There are some stronger conditions that we could impose on $n$-globes: for instance, we could require that there exist an ordering $x_1, \ldots, x_m$ of the $n$-dimensional elements, such that, for all $j \leq m$, the subsets $\mathrm{cl}\{x_1, \ldots x_{j-1}\}$ and $U_{x_j}$ are mergeable; or that given any sub-globe $V \subseteq U$, any tree of bi-partitions of $V$ can be completed to a tree of bi-partitions of $U$. The latter is a directed variant of extendable shellability, see for example \cite{moriyama2003incremental}. Up to dimension 2 these are definitely equivalent, and it would be desirable that this extended to higher dimensions, but for now we leave it as an open problem.
\end{remark}

\begin{prop} \label{prop:basicglobe}
Let $U$ be an $n$-globe in $X$, with $n > 0$. Then:
\begin{enumerate}[label=(\alph*)]
	\item $\sbord{}{+}U$ and $\sbord{}{-}U$ are both inhabited;
	\item for all $n$-dimensional $x$, $x'$ in $U$, let $x \prec x'$ if there exists $y$ in $\sbord{}{+} U_x \cap \sbord{}{-} U_{x'}$; then, the transitive closure of $\prec$ is a connected partial order on $U \cap X_n$;
	\item any $(n-1)$-dimensional element of $U$ is covered at most by two $n$-dimensional elements, and is covered by a single $n$-dimensional element if and only if it belongs to $\sbord{}{\alpha}U$, for some $\alpha \in \{+,-\}$.
\end{enumerate}
\end{prop}
\begin{proof}
For the first point, observe that $\bord{}{\alpha} U$ is pure and $(n-1)$-dimensional, hence it is equal to the closure of $\sbord{}{\alpha} U$, necessarily inhabited.

For the second point, we proceed by induction on the number of top-dimensional elements of $U$. If $U$ is atomic, there is nothing to prove. Otherwise, let $U_1, U_2$ be the mergeable sub-$n$-globes of $U$ given by the definition; by the inductive hypothesis, their top-dimensional elements form a connected partial order with $\prec$. We have $U_1 \cap U_2 = \bord{}{\alpha} U_1 \cap \bord{}{-\alpha} U_2$ for some $\alpha \in \{+,-\}$, and because it is an $(n-1)$-dimensional globe, it is the closure of $\sbord{}{\alpha} U_1 \cap \sbord{}{-\alpha} U_2$. If $\alpha = +$, for some $x_1 \in U_1$ and $x_2 \in U_2$, it holds that $x_1 \prec x_2$, else $x_2 \prec x_1$; in both cases, the two partial orders connect.

The third point is obvious for atomic $n$-globes. If $U$ is non-atomic, given $U_1$, $U_2$ as in the definition, $U_1 \cap U_2 = \mathrm{cl}(\sbord{}{\beta}U_1 \cap \sbord{}{-\beta}U_2)$ for some $\beta$, and the elements of $\sbord{}{\beta}U_1$, $\sbord{}{-\beta}U_2$ are covered by a single element of $U_1$, $U_2$, respectively, by the inductive hypothesis; therefore, the elements of $\sbord{}{\beta}U_1 \cap \sbord{}{-\beta}U_2$ are covered by exactly two elements in $U$. A simple case distinction on the other $(n-1)$-dimensional elements completes the proof.
\end{proof}

\begin{remark}
In particular, if $X$ is a globular poset, and $x$ of an element of dimension $n > 0$,
\begin{equation*}
	\sbord{}{\alpha} U_x = \{y \in X \,|\, \text{$x$ covers $y$ and $o(c_{x,y}) = \alpha$} \},
\end{equation*}
by Proposition \ref{prop:basicglobe}.$(a)$. It follows that both the input boundary and the output boundary of $x$ contain at least one $(n-1)$-dimensional element.
\end{remark}

\begin{cons} \label{cons:merging}
Let $X$ be an oriented thin poset, and suppose that there are $n$-dimensional elements $x_1$ and $x_2$ with the following property: $U_{x_1} \cap U_{x_2} = U_y$ for some $(n-1)$-dimensional $y$, only covered by $x_1$ and $x_2$, and
\begin{equation*}
\begin{tikzpicture}[baseline={([yshift=-.5ex]current bounding box.center)}]
	\node[scale=1.25] (1) at (-1,1) {$x_1$};
	\node[scale=1.25] (1b) at (1,1) {$x_2$};
	\node[scale=1.25] (2) at (0,0) {$y$};
	\draw[1c] (1) to node[auto,swap] {$+$} (2);
	\draw[1c] (1b) to node[auto] {$-$} (2);
\end{tikzpicture}
\end{equation*}
in the labelled Hasse diagram of $X$. Let $X'$ be the poset obtained from $X$ by identifying the elements $x_1$, $x_2$, and $y$. Then $X'$ is graded: if $\tilde{x}$ is the element of $X'$ obtained from the identification, the only problem could arise if length $m$ paths from $z$ to $y$ in $HX$ got mapped to length $m$ paths from $z$ to $\tilde{x}$ in $HX'$; but because $y$ is only covered by $x_1$ and $x_2$, any such path in $HX$ passes through either $x_1$ or $x_2$, hence becomes a length $(m-1)$ path in $HX'$.

Now, define a labelling of edges of $HX'$ as follows: $o'(c_{w,z}) := o(c_{w,z})$ if $w, z \neq \tilde{x}$; $o'(c_{w,\tilde{x}}) := o(c_{w,x_i})$ if $w$ covers $x_i$ in $X$, and similarly $o'(c_{\tilde{x},z}) := o(c_{x_i,z})$ if $x_i$ covers $z$ in $X$, $i= 1,2$. We will show that $X'$ is thin, and that $o'$ is a well-defined orientation. 

By hypothesis, only $y$ is covered both by $x_1$ and by $x_2$ in $X$, so $o'(c_{\tilde{x},z})$ is well-defined for all $z \in X'$. Suppose that $w$ covers both $x_1$ and $x_2$ in $X$; then, in the labelled Hasse diagram of $X$,
\begin{equation*}
\begin{tikzpicture}[baseline={([yshift=-.5ex]current bounding box.center)}]
	\node[scale=1.25] (0) at (0,2) {$w$};
	\node[scale=1.25] (1) at (-1,1) {$x_1$};
	\node[scale=1.25] (1b) at (1,1) {$x_2$};
	\node[scale=1.25] (2) at (0,0) {$y$};
	\draw[1c] (0) to node[auto,swap] {$\alpha_1$} (1);
	\draw[1c] (0) to node[auto] {$\alpha_2$} (1b);
	\draw[1c] (1) to node[auto,swap] {$+$} (2);
	\draw[1c] (1b) to node[auto] {$-$} (2);
	\node[scale=1.25] at (1.25,0) {,};
\end{tikzpicture}
\end{equation*}
hence $o(c_{w,x_1}) = \alpha_1 = \alpha_2 = o(c_{w,x_2})$, which proves that $o'(c_{w,\tilde{x}})$ is also well-defined.

For thinness of $X'$, it suffices to consider paths $\tilde{x} \to z \to w$ in $HX'_\bot$. Any such path must come from a path $x_i \to z \to w$ in $HX_\bot$, say $i = 1$. This can be completed to
\begin{equation*}
\begin{tikzpicture}[baseline={([yshift=-.5ex]current bounding box.center)}]
	\node[scale=1.25] (0) at (0,2) {$x_1$};
	\node[scale=1.25] (1) at (-1,1) {$z$};
	\node[scale=1.25] (1b) at (1,1) {$z'$};
	\node[scale=1.25] (2) at (0,0) {$w$};
	\draw[1c] (0) to node[auto,swap] {$\alpha_1$} (1);
	\draw[1c] (0) to node[auto] {$\alpha_2$} (1b);
	\draw[1c] (1) to node[auto,swap] {$\beta_1$} (2);
	\draw[1c] (1b) to node[auto] {$\beta_2$} (2);
\end{tikzpicture}
\end{equation*}
for a unique $z'$ in $X$, with $\alpha_1\beta_1 = -\alpha_2\beta_2$. If $z' \neq y$, the diagram is faithfully mapped onto $X'$. If $z' = y$, by thinness applied to the interval $[x_2, w]$, there exists a unique $z'' \neq z, y$ such that 
\begin{equation*}
\begin{tikzpicture}[baseline={([yshift=-.5ex]current bounding box.center)}]
	\node[scale=1.25] (0) at (0,2) {$x_1$};
	\node[scale=1.25] (0b) at (2,2) {$x_2$};
	\node[scale=1.25] (1) at (-1,1) {$z$};
	\node[scale=1.25] (1b) at (1,1) {$y$};
	\node[scale=1.25] (1c) at (3,1) {$z''$};
	\node[scale=1.25] (2) at (1,0) {$w$};
	\draw[1c] (0) to node[auto,swap] {$\alpha_1$} (1);
	\draw[1c] (0) to node[auto] {$\!+$} (1b);
	\draw[1c] (0b) to node[auto,swap] {$-\!$} (1b);
	\draw[1c] (0b) to node[auto] {$\alpha_2$} (1c);
	\draw[1c] (1) to node[auto,swap] {$\beta_1$} (2);
	\draw[1c] (1b) to node[auto] {$\beta$} (2);
	\draw[1c] (1c) to node[auto] {$\beta_2$} (2);
\end{tikzpicture}
\end{equation*}
in $HX_\bot$, and $\alpha_1\beta_1 = -\beta = -\alpha_2 \beta_2$. This becomes
\begin{equation*}
\begin{tikzpicture}[baseline={([yshift=-.5ex]current bounding box.center)}]
	\node[scale=1.25] (0) at (0,2) {$\tilde{x}$};
	\node[scale=1.25] (1) at (-1,1) {$z$};
	\node[scale=1.25] (1b) at (1,1) {$z''$};
	\node[scale=1.25] (2) at (0,0) {$w$};
	\draw[1c] (0) to node[auto,swap] {$\alpha_1$} (1);
	\draw[1c] (0) to node[auto] {$\alpha_2$} (1b);
	\draw[1c] (1) to node[auto,swap] {$\beta_1$} (2);
	\draw[1c] (1b) to node[auto] {$\beta_2$} (2);
\end{tikzpicture}
\end{equation*}
in $HX'_\bot$, which completes the proof that $X'$ is an oriented thin poset.
\end{cons}

\begin{dfn} \index{simple merger}
Let $X$, $X'$ be oriented thin posets as in Construction \ref{cons:merging}. We say that $X'$ is obtained from $X$ by a \emph{simple merger}, and write $X \leadsto X'$.
\end{dfn}

The following depicts a sequence of simple mergers on a 2-globe, the coloured arrow pointing from $x_1$ to $x_2$:
\begin{equation*}
\input{img/c1_2_simplemerger.tex}
\end{equation*} 
In general, the result of a simple merger on a globular poset is not a globular poset: consider, for example, the 2-globe
\begin{equation*}
\input{img/c3_2_notsimplemerger.tex}
\end{equation*} 
where the two indicated 2-cells fulfil the conditions of a simple merger, but the result of their identification is no longer a 2-globe. However, we will see that it is possible to choose what cells to merge in a way that preserves globular posets, and induction on such sequences of simple mergers will be one of our main proof techniques.

In the following, we use equality somewhat improperly, to mean that a subset of $X$ is preserved by the mergers.
\begin{lem} \label{lem:simplemerge}
Let $X$ be an $n$-globe. Then, there exists a sequence of simple mergers $X \leadsto \ldots \leadsto \tilde{X}$ such that $\tilde{X}$ is an atomic $n$-globe, and $\bord{}{\alpha} \tilde{X} = \bord{}{\alpha} X$, $\alpha \in \{+,-\}$.
\end{lem}
\begin{proof}
We proceed by induction on $n$, and on the number $p$ of top-dimensional elements of $X$. If $n = 0$, or if $p = 1$, there is nothing to prove. 

Suppose $n = 1$, $p > 1$; then, given a tree of bi-partitions of the 1-dimensional elements of $X$, there exist two leaves $\{x_1\}$, $\{x_2\}$ that branch out of the two-element set $\{x_1, x_2\}$, so in particular $U_{x_1}, U_{x_2}$ are mergeable 1-globes. Then, $U_{x_1} \cap U_{x_2} = \sbord{}{\alpha} U_{x_1} \cap \sbord{}{-\alpha} U_{x_2} = \{y\}$ for some 0-dimensional $y$, which by Proposition \ref{prop:basicglobe}.$(c)$ is only covered by $x_1$ and $x_2$. Hence, $x_1$ and $x_2$ fulfill the conditions for a simple merger $X \leadsto X'$. 

Clearly, $\bord{}{\beta} U_{\tilde{x}} = \bord{}{\beta} (U_{x_1} \cup U_{x_2})$, so $U_{\tilde{x}}$ is an atomic 1-dimensional globe. Moreover, the identification of $x_1, x_2$, and the truncation of the corresponding final branching in the tree of bi-partitions of $X$ yield a valid tree of bi-partitions of $X'$, so $X'$ is a 1-globe with one less top-dimensional element, and the inductive hypothesis applies.

Now, suppose $n > 1$, and let $x_1, x_2$ as in the 1-dimensional case. If $V := U_{x_1} \cap U_{x_2} = \bord{}{\alpha} U_{x_1} \cap \bord{}{-\alpha} U_{x_2}$ is an atomic $(n-1)$-globe, we can proceed exactly as in the 1-dimensional case. Otherwise, we know that $V$ is an $(n-1)$-dimensional globe, and since it is a subset of both $\bord{}{\alpha} U_{x_1}$ and $\bord{}{-\alpha} U_{x_2}$, which are also globes, from Proposition \ref{prop:basicglobe}.$(c)$ we can derive that $V \cap \bord{}{\beta}(U_{x_1} \cup U_{x_2}) \subset \bord{}{}V$. By the inductive hypothesis, we can perform a sequence of simple mergers transforming $V$ into an atomic $(n-1)$-globe, without modifying $\bord{}{} V$ and thus $\bord{}{\beta} (U_{x_1} \cup U_{x_2})$. Then, we can proceed by merging $x_1$ and $x_2$ as in the 1-dimensional case, and apply the inductive hypothesis.
\end{proof}

\begin{thm} \label{thm:globularity}
Let $X$ be an $n$-globe, $n > 1$. Then, for $\alpha = \{+,-\}$,
\begin{equation*}
	\bord{}{\alpha}(\bord{}{+} X) = \bord{}{\alpha}(\bord{}{-} X).
\end{equation*}
\end{thm}
\begin{proof}
First, suppose that $X$ is atomic, that is $X = U_x$ for some $n$-dimensional $x$. Suppose that $z \in \sbord{}{\alpha}(\bord{}{+} X)$, and let $y \in \sbord{}{+} X$ be covering $z$. By thinness, the interval $[x,z]$ is of the form
\begin{equation*}
\begin{tikzpicture}[baseline={([yshift=-.5ex]current bounding box.center)}]
	\node[scale=1.25] (0) at (0,2) {$x$};
	\node[scale=1.25] (1) at (-1,1) {$y$};
	\node[scale=1.25] (1b) at (1,1) {$y'$};
	\node[scale=1.25] (2) at (0,0) {$z$};
	\draw[1c] (0) to node[auto,swap] {$+$} (1);
	\draw[1c] (0) to node[auto] {$\beta_1$} (1b);
	\draw[1c] (1) to node[auto,swap] {$\alpha$} (2);
	\draw[1c] (1b) to node[auto] {$\beta_2$} (2);
\end{tikzpicture}
\end{equation*}
for some $y'$. Suppose that $\beta_1 = +$; then $y' \in \sbord{}{+} X$, yet necessarily $\beta_2 = -\alpha$, contradicting $z \in \sbord{}{\alpha}(\bord{}{+} X)$. It follows that $\beta_1 = -$, and $\beta_2 = \alpha$; since $z$ is only covered by $y$ and $y'$, $z \in \sbord{}{\alpha}(\bord{}{-} X)$. The converse is symmetrical.

Now, suppose $X$ is not atomic. By Lemma \ref{lem:simplemerge}, we can construct an atomic $n$-globe $\tilde{X}$ such that $\bord{}{\alpha} X = \bord{}{\alpha} \tilde{X}$. By the first part, $\bord{}{\alpha}(\bord{}{+} X) = \bord{}{\alpha}(\bord{}{+} \tilde{X}) = \bord{}{\alpha}(\bord{}{-} \tilde{X}) = \bord{}{\alpha}(\bord{}{-} X)$.
\end{proof}

In fact, the proof shows that $\bord{}{+} X \cap \bord{}{-} X = \bord{}{}(\bord{}{+} X) = \bord{}{}(\bord{}{-} X)$. This result, in addition to being technically useful, enables us to verify that our non-degeneracy criterion is satisfied.

In the following, we identify the standard $n$-globe with its oriented incidence poset.
\begin{cor}
Let $X$ be an $n$-globe. Then, there exist a sequence of simple mergers $X \leadsto \ldots \leadsto G^n$ from $X$ to the standard $n$-globe, whose restriction to $\bord{}{}X$ is a sequence of simple mergers from $\bord{}{}X$ to $\bord{}{}G^n$.
\end{cor}
\begin{proof}
By Lemma \ref{lem:simplemerge}, with a sequence of simple mergers we can obtain $X'$ having a single $n$-dimensional element $\top$, and $\bord{}{\alpha} X = \bord{}{\alpha} X'$.

Then, we can do the same with the $(n-1)$-globes $\bord{}{+} X'$, $\bord{}{-} X'$, obtaining an $n$-globe $X''$ with a single $(n-1)$-dimensional element both in the input and in the output boundary of $\top$, all the while preserving the $(n-2)$-globes $\bord{}{\alpha}(\bord{}{+} X') = \bord{}{\alpha}(\bord{}{-} X')$. Iterating on successive boundaries, we end up with an $n$-globe that has only two $k$-dimensional elements $k^+, k^-$ for all $k < n$, and $\sbord{}{\alpha} U_{k^+} = \sbord{}{\alpha} U_{k^-} = \{(k-1)^\alpha\}$, that is, a standard $n$-globe.

Notice that, after the first application of Lemma \ref{lem:simplemerge}, all simple mergers only involve elements of $\bord{}{} X' = \bord{}{} X$, and the single remaining $n$-dimensional element is unaffected, so by restriction we obtain a sequence of simple mergers from $\bord{}{}X$ to $\bord{}{}G^n$.
\end{proof}

The standard $n$-globe admits no simple merger; thus, it can be characterised as the $n$-globe that is terminal with respect to this operation. 

\begin{dfn}
Let $X$ be a poset. The \emph{nerve} of $X$ is the simplicial set $NX$ whose $n$-simplices are chains $(x_1 \leq \ldots \leq x_n)$ of length $n$ in $X$, faces are defined by $(x_1 \leq \ldots \leq x_n) \mapsto (x_1 \leq \ldots \leq x_{k-1} \leq x_{k+1} \leq \ldots \leq x_n)$, and degeneracies by $(x_1 \leq \ldots \leq x_n) \mapsto (x_1 \leq \ldots \leq x_{k} \leq x_{k} \leq \ldots \leq x_n)$, for $1 \leq k \leq n$. \index{nerve}
\end{dfn}
The nerve can be composed with the usual geometric realisation functor to obtain a point-set topological space $|X|$. It is customary, in combinatorial topology, to attribute to a poset topological properties of $|X|$. By \cite[Theorem 1.7]{lundell1969topology}, if $X$ is the incidence poset of a regular CW complex $C$, then $|X|$ is homeomorphic to $C$.

\begin{thm} \label{thm:homeodisk}
Let $X$ be an $n$-globe. Then $|X|$ is homeomorphic to an $n$-disk, and $|\bord{}{}X|$ is homeomorphic to an $(n-1)$-sphere.
\end{thm}
\begin{proof}
The underlying posets of $G^n$, $\bord{}{} G^n$ are the incidence posets of the minimal regular CW decomposition of the $n$-disk and of the $(n-1)$-sphere, respectively.  Moreover, every simple merger has an underlying map of posets, so it suffices to show this map induces a homeomorphism of topological spaces. For $n = 0$, there is nothing to prove. For $n = 1$, if $X$ is atomic, it is $G^1$; otherwise, each merger induces the homeomorphic mapping of two 1-disks connected at one extremity onto a single 1-disk with the same boundary. 

For $n > 1$, if $X$ is atomic, all mergers involve $k$-dimensional elements with $k < n$, so by the inductive hypothesis $|X|$ is an $n$-disk. If $X$ is not atomic, let $x_1$, $x_2$ be the first $n$-dimensional elements to be merged. We already know that $|U_{x_1}|$ and $|U_{x_2}|$ are $n$-disks, and their intersection is an $(n-1)$-disk; again, their union is an $n$-disk, and the induced map is a homeomorphism.
\end{proof}

\begin{cor}
Let $X$ be a globular poset. Then the underlying poset of $X$ is a CW poset.
\end{cor}
\begin{proof}
By \cite[Definition 2.1]{bjorner1984posets}, it suffices that, for all $n$-dimensional $x$, $|\bord{}{}U_x|$ be homeomorphic to an $(n-1)$-sphere, which follows from the previous Theorem.
\end{proof}

\begin{remark}
Dualising the notion of simple merger, one obtains a notion similar to \emph{cellular collapses} as considered, for instance, in discrete Morse theory: an $n$-cell with only one $(n-1)$-cell in the input, and one in the output boundary is collapsed down to an $(n-1)$-cell. In fact, a consequence of Lemma \ref{lem:simplemerge} is that, when $X$ is a globe, $(X - \bord{}{} X)^\mathrm{op}$ is collapsible. 

While mergers induce homeomorphisms, collapses induce (simple) homotopy equivalences. This seems also related to the fact that shellability, which in general only decides the homotopy type of a space, as a property of the \emph{dual} of an incidence poset can decide the homeomorphism class, see \cite[Proposition 4.5]{bjorner1984posets}. There seems to be a duality at play between topology up-to-homeomorphism and topology up-to-homotopy, whose exact nature we do not know.
\end{remark}

Next, we want to assemble globular posets into a category. While there may be more interesting notions of morphism between globular posets, for our purposes we are only really interested in inclusions of globes into the boundary of other globes.
\begin{dfn} \index{inclusion!of globular posets}
Let $X, Y$ be two globular posets. An \emph{inclusion} $\imath: X \hookrightarrow Y$ is a closed embedding of posets that is compatible with the orientations, that is, $o^Y(c_{\imath(y),\imath(x)}) = o^X(c_{y,x})$ for all $y, x \in X$ such that $y$ covers $x$. An inclusion is an \emph{isomorphism} if it is also surjective. We write $\globpos$ for the category of globular posets and inclusions. \index{globpos@$\globpos$} 
\end{dfn}

\begin{remark}
Being a closed embedding of graded posets, an inclusion preserves dimensions and preserves and reflects the covering relation: if $y' = \imath(y)$ and $y'$ covers $x'$, then there is a unique $x$ such that $y$ covers $x$ and $x' = \imath(x)$. 

In particular, if $X$ is an $n$-globe with a greatest element $x$ and $\imath: X \to Y$ an embedding, then $\imath(X) = U_{\imath(x)}$, an $n$-globe in $Y$. 
\end{remark}

\begin{cons}
Let $X$ be a globular poset, and, for all $n$, define $\textit{Gl}_n(X)$ to be the set of $n$-globes in $X$. By Theorem \ref{thm:globularity}, 
\begin{equation*}
\begin{tikzpicture}[baseline={([yshift=-.5ex]current bounding box.center)}]
	\node[scale=1.25] (0) at (-2.4,0) {$\textit{Gl}_0(X)$};
	\node[scale=1.25] (1) at (0.4,0) {$\textit{Gl}_1(X)$};
	\node[scale=1.25] (2) at (3,-.1) {$\ldots$};
	\node[scale=1.25] (3) at (5.6,0) {$\textit{Gl}_n(X)$};
	\node[scale=1.25] (4) at (8.2,-.1) {$\ldots,$};
	\draw[1c] (-.35,.15) to node[auto,swap] {$\bord{}{+}$} (-1.7,.15);
	\draw[1c] (-.35,-.15) to node[auto] {$\bord{}{-}$} (-1.7,-.15);
	\draw[1c] (2.45,.15) to node[auto,swap] {$\bord{}{+}$} (1.1,.15);
	\draw[1c] (2.45,-.15) to node[auto] {$\bord{}{-}$} (1.1,-.15);
	\draw[1c] (4.85,.15) to node[auto,swap] {$\bord{}{+}$} (3.5,.15);
	\draw[1c] (4.85,-.15) to node[auto] {$\bord{}{-}$} (3.5,-.15);
	\draw[1c] (7.65,.15) to node[auto,swap] {$\bord{}{+}$} (6.3,.15);
	\draw[1c] (7.65,-.15) to node[auto] {$\bord{}{-}$} (6.3,-.15);
	\node[scale=1.25] (5) at (9.6,0) {$n \in \mathbb{N},$};
\end{tikzpicture}
\end{equation*}
is a globular set, and because inclusions of globular posets preserve globes, the construction extends to a functor $\textit{Gl}: \globpos \to \cat{Glob}$.
\end{cons}

\begin{cons} \label{cons:induced_adc} 
Let $X$ be a globular poset, with its natural grading. Then, define $e: \mathbb{Z}X_0 \to \mathbb{Z}$ by $ex := 1$ for all $x \in X_0$, and, for all $n > 0$, $d: \mathbb{Z}X_n \to \mathbb{Z}X_{n-1}$ by
\begin{equation*}
	dx := \sum_{y \in \sbord{}{+} U_x} y \; - \sum_{y' \in \sbord{}{-} U_x} y',
\end{equation*}
for all $x \in X_n$. We claim that $(X,d)$ so defined is an augmented directed complex. Writing 
\begin{equation*}
	d^+x := \sum_{y \in \sbord{}{+} U_x} y, \quad\quad d^-x := \sum_{y' \in \sbord{}{-} U_x} y',
\end{equation*}
we have that
\begin{equation*}
	dd^\alpha x = \sum_{z \in \sbord{}{+}\bord{}{\alpha} U_x} z \; - \sum_{z' \in \sbord{}{-}\bord{}{\alpha} U_x} z'
\end{equation*}
because $\bord{}{\alpha} U_x = \mathrm{cl}(\sbord{}{\alpha}U_x)$ is an $(n-1)$-globe, so by Proposition \ref{prop:basicglobe}.$(c)$ any element covered by some $y \in \sbord{}{\alpha}U_x$ is either covered by another element of $\sbord{}{\alpha}U_x$, and appears a second time with the opposite sign in $dd^\alpha x$, or is only covered by $y$, and belongs to $\sbord{}{\beta}\bord{}{\alpha} U_x$ for some $\beta \in \{+,-\}$. Then, by Theorem \ref{thm:globularity}, $dd^+ x = dd^- x$, so $dd = 0$, and $ed = 0$ is verified similarly. The construction extends to a functor $D: \globpos \to \adc$. 

By construction, if an augmented directed complex is induced by a globular poset $X$, its oriented incidence poset, Definition \ref{dfn:incidence}, is the original globular poset. In particular, if the graph $HX^o$ is acyclic, $DX$ is strongly loop-free. Moreover, $DX$ is always unital.
\end{cons}

So, a globular poset induces an augmented directed complex; we suspect that this can be strengthened as follows. In\cite{steiner1993algebra}, a \emph{directed precomplex} $X$ is, in our terminology, an oriented graded poset. For a closed subset $U$ of $X$, let
\begin{align*}
	\sbord{k}{\alpha} U & := \{x \in U \,|\, \dmn{x} = k \text{ and, for all $y \in U$, if $y$ covers $x$, then $o(c_{y,x}) = \alpha$} \}, \\
	\tilde{\partial}_k^\alpha U & := \mathrm{cl}(\sbord{k}{\alpha} U) \cup \{ x \in U \,|\, \text{for all $y \in U$, if $x \leq y$, then $\dmn{y} \leq k$} \}.
\end{align*}
Two closed subsets $U, V \subseteq X$ are then said to be $k$-composable if $U \cap V = \tilde{\partial}_{k-1}^+ U = \tilde{\partial}_{k-1}^- V$, and we define $U \cp{k} V := U \cup V$ in that case. Such compositions are compatible with the $k$-boundary operations in the same way as compositions in $\omega$-categories.
\begin{dfn}
A directed precomplex $X$ is a \emph{directed complex} if, for all $x \in X$, $\dmn{x} = n$,
\begin{enumerate}
	\item $\bord{}{\alpha}\bord{}{\beta}U_x = \tilde{\partial}_{n-2}^\alpha U_x$, and
	\item $\bord{}{\alpha}U_x$ is a composite of the $U_y$, $y < x$.
\end{enumerate}
\end{dfn}
Now, if $X$ is a globular poset, it satisfies the first condition, by a proof similar to the atomic case of Theorem \ref{thm:globularity}. Moreover, for all $n$-globes $U$ in $X$, we define 
\begin{equation*}
	\bord{k}{\alpha} U \; := \underbrace{(\bord{}{\alpha}\ldots\bord{}{\alpha})}_{n-k} U\,,
\end{equation*} 
knowing that $\bord{}{\alpha}\bord{k}{+} U = \bord{}{\alpha}\bord{k}{-} U$.

\begin{conj}
Let $X$ be a globular poset. Then, $X$ is a directed complex, and for all globes $U$ in $X$ it holds that $\tilde{\partial}_k^\alpha U = \bord{k}{\alpha} U$.
\end{conj}
What makes us quite confident that the conjecture is true is that the condition of mergeability is, at least topologically, more restrictive than $\omega$-categorical composability: at the level of the underlying CW posets, a merger of $n$-disks is homeomorphic to an $n$-disk, while a composition is in general only homotopy-equivalent. However, the translation between the two formalisms seems to involve some technical subtleties, mainly because $k$-composability depends on the single \emph{datum} of the $(k-1)$-boundary, which is formally the same for all $k$, whereas mergeability depends on the shared boundary between two $n$-globes \emph{and} on the way it fits in the rest of their boundaries, which presents itself in different ways for increasing $n$. 

Proving the conjecture would justify, \emph{a posteriori}, our naming of ``regular polygraphs'', by showing that they are, indeed, polygraphs in the $\omega$-categorical sense; on the other hand, as a notion of directed space, they can be interesting independently of its validity. 

Next, we turn our attention to our third \emph{desideratum}, concerning tensor products.
\begin{dfn} \index{tensor product!of oriented posets}
Let $X, Y$ be oriented posets. The \emph{tensor product} $X \otimes Y$ of $X$ and $Y$ is the graded poset $X \times Y$, oriented as follows: write $x \otimes y$ for an element $(x,y)$ of $X \times Y$; then, for all $x'$ covered by $x$ in $X$, $y'$ covered by $y$ in $Y$, let
\begin{align*}
	o(c_{x \otimes y, x' \otimes y}) & := o^X(c_{x,x'}), \\
	o(c_{x \otimes y, x \otimes y'}) & := (-1)^{\dmn{x}}o^Y(c_{y,y'}).
\end{align*}
\end{dfn}
\begin{thm} \label{thm:globpostensor}
Let $X, Y$ be oriented posets. Then,
\begin{enumerate}
	\item if $X, Y$ are oriented thin posets, $X \otimes Y$ is an oriented thin poset;
	\item if $X, Y$ are globular posets, $X \otimes Y$ is a globular poset.
\end{enumerate}
\end{thm}
\begin{proof}
An interval of length 2 in $X \otimes Y$, whose least element has dimension $n > 0$, has one of the following forms:
\begin{itemize}
	\item $[x', x] \otimes \{y\}$ for some interval $[x', x]$ of length 2 in $X$, and $y \in Y$;
	\item $\{x\} \otimes [y',y]$ for some interval $[y',y]$ of length 2 in $Y$, and $x \in X$;
	\item $[x'\otimes y', x \otimes y]$, where $x$ covers $x'$ in $X$ and $y$ covers $y'$ in $Y$.
\end{itemize}
In the first two cases, oriented thinness follows from the oriented thinness of $X$ and $Y$, respectively. In the last case, let $\alpha := o^X(c_{x,x'})$, $\beta := o^Y(c_{y,y'})$. Then, the interval has the form
\begin{equation*}
\begin{tikzpicture}[baseline={([yshift=-.5ex]current bounding box.center)}, scale=1.5]
	\node[scale=1.25] (0) at (0,2) {$x\otimes y$};
	\node[scale=1.25] (1) at (-1,1) {$x' \otimes y$};
	\node[scale=1.25] (1b) at (1,1) {$x \otimes y'$};
	\node[scale=1.25] (2) at (0,0) {$x' \otimes y'$};
	\draw[1c] (0) to node[auto,swap] {$\alpha$} (1);
	\draw[1c] (0) to node[auto] {$(-1)^{\dmn{x}}\beta$} (1b);
	\draw[1c] (1) to node[auto,swap] {$(-1)^{\dmn{x'}}\beta$} (2);
	\draw[1c] (1b) to node[auto] {$\alpha$} (2);
\end{tikzpicture}
\end{equation*}
in the labelled Hasse diagram of $X \otimes Y$, and $(-1)^{\dmn{x'}} = (-1)^{\dmn{x} - 1} = - (-1)^{\dmn{x}}$. The case of intervals $[\bot, x \otimes y]$ can easily be handled explicitly. This proves the first point.

Next, suppose $X$ and $Y$ are globes; we will show that $X \otimes Y$ is a globe, by simultaneous double induction on the dimensions $n$ of $X$ and $m$ of $Y$, and on the number of their top-dimensional elements. When either $n$ or $m$ is 0, $X \otimes Y \simeq X$ or $X \otimes Y \simeq Y$, obviously a globe.

Suppose that $n, m > 0$, and let $\alpha' := (-1)^n \alpha$; then 
\begin{equation*}
	\bord{}{\alpha}(X \otimes Y) = (\bord{}{\alpha}X \otimes Y) \cup (X \otimes \bord{}{\alpha'} Y).
\end{equation*}
By the inductive hypothesis, both $U_1 := \bord{}{\alpha}X \otimes Y$ and $U_2 := X \otimes \bord{}{\alpha'} Y$ are globes, and their intersection $\bord{}{\alpha}X \otimes \bord{}{\alpha'}Y$ is a globe. Thus, $U_1$ and $U_2$ satisfy the condition for $\bord{}{\alpha}(X \otimes Y)$ to be a globe.

If $X, Y$ are atomic, this is sufficient. Otherwise, suppose $X$ is non-atomic, and let $U^X_1$ and $U^X_2$ be mergeable sub-globes of $X$ as in the definition. Then $U^X_1 \otimes Y$, $U^X_2 \otimes Y$, and their intersection $(U^X_1 \cap U^X_2) \otimes Y$ are all globes by the inductive hypothesis; therefore $X \otimes Y$ is a globe. Induction on the number of top-dimensional elements of $Y$ is similar. 

Finally, suppose that $X, Y$ are globular posets. For all $x \in X, y \in Y$, we have $U_{x \otimes y} = {U^X_x \otimes U^Y_y}$, which is a globe since $U^X_x, U^Y_y$ are. Therefore, $X \otimes Y$ is a globular poset.
\end{proof}

The tensor product defines a monoidal structure on $\globpos$, whose unit is the 0-globe $1$. Moreover, the full subcategory whose objects are atomic globes is closed under the tensor product. 

\begin{prop}
The functor $D: \globpos \to \adc$ is monoidal.
\end{prop}
\begin{proof}
Simple inspection of the definitions.
\end{proof}

This means, in particular, that the tensor product of globular posets agrees with the tensor product of $\omega$-categories on the globular posets that are mapped by $D$ to loop-free polygraphs. 

We are now ready to define our directed spaces.

\begin{dfn} \index{RG@$\cat{RG}$} \index{polygraph!regular}
Let $\cat{RG}$ be a skeleton of the full subcategory of $\globpos$ whose objects are atomic globes. A \emph{regular polygraph} $X$ is a presheaf $X: \opp{\cat{RG}} \to \cat{Set}$. A \emph{map} $f: X \to Y$ of regular polygraphs is a morphism of presheaves.

The \emph{tensor product} $X \otimes Y$ of two regular polygraphs $X, Y$ is their Day convolution with respect to the tensor product of globular posets. The tensor product defines a monoidal structure on the category $\rpol$ of regular polygraphs and maps, whose unit is the Yoneda embedding of the 0-globe $1$. \index{tensor product!of regular polygraphs}
\end{dfn}
\begin{remark}
The Day convolution product is, in fact, canonically part of a biclosed structure; see \cite[Theorem 3.3]{day1970closed}.
\end{remark}

We will casually identify a globe $G$ with its Yoneda embedding $\mathrm{Hom}_\cat{RG}(-,G)$. Because $\rpol$ is a presheaf category, it is cocomplete, so it has coproducts modelling the disjoint union of regular polygraphs, and coequalisers, which can be used to define quotients; unlike the $\omega$-categorical case, however, only generators of the same dimension can be identified in a quotient. We can also define opposites.
\begin{dfn} \index{opposite!oriented poset}
Let $X$ be an oriented poset, $S \subseteq \mathbb{N}^+$. Then $\oppn{X}{S}$ is the oriented poset with the same underlying poset as $X$, and the orientation $o'$ defined by 
\begin{equation*}
o'(c_{x,y}) := \begin{cases}
		-o(c_{x,y}), & \dmn{x} \in S, \\
		o(c_{x,y}), & \dmn{x} \not\in S,
	\end{cases}
\end{equation*}
for all elements $x, y \in X$ such that $x$ covers $y$. It is easy to verify that $\oppn{X}{S}$ is a globular poset if $X$ is; in fact, $\oppn{(-)}{S}$ defines a functor $\globpos \to \globpos$, restricting to a functor $\cat{RG} \to \cat{RG}$.

Let $X$ be a regular polygraph, $S \subseteq \mathbb{N}^+$. The \emph{$S$-opposite} of $X$ is the regular polygraph $\oppn{X}{S} := X(\oppn{(-)}{S})$. 
\end{dfn}
The applicable results of the last part of Section \ref{sec:pasting} all hold, and we will use the same notational shortcuts. 

The definition of the shape category $\cat{RG}$ as a skeleton is not very satisfactory, since recognising isomorphic globular posets may be a hard problem. Our hope is that an inductive definition can be found, akin to the zoom complexes of \cite{kock2010polynomial} for the category of opetopes, inducing a recursive enumeration of isomorphism classes. We suspect that such a definition may involve a classification of the connected partial orders of Proposition \ref{prop:basicglobe}.$(b)$; we know, for instance, that for 1-globes, one obtains a linear order of atomic globes, and for 2-globes, an order with a planar Hasse diagram.

On the other hand, the definition of regular polygraphs automatically satisfies our \emph{desiderata} 1 and 3, and seems well-equipped for satisfying number 4 as well. As for \emph{desideratum} 2, we can define the geometric realisation of a regular polygraph $X$ in the standard way, as the coend
\begin{equation*}
	|X| := \int^G |G| \times X(G),
\end{equation*}
in the category of topological spaces, where $G$ ranges over the globes in $\cat{RG}$, and $|G|$ is the geometric realisation of their underlying poset. For now, this answers our question partially: if we can label the incidence poset of a regular CW complex in such a way that it becomes a globular poset, then its geometric realisation is going to be homeomorphic to the original CW complex. 

However, to be able to turn such a globular poset into an ``undirected object'' in $\rpol$ --- preferably equivalent to the ``$\cat{RG}$-nerve'' of a topological space, obtained abstractly as the right adjoint of the geometric realisation --- we need a notion of equivalence cell in a regular polygraph; and to even be able to define a homotopy sequence for a regular polygraph, in the absence of degenerate boundaries, we will need a notion of weak unit.

In the next section, we do this in the 2-truncated case, to suggest that, while much work still needs to be done, there is at least a viable path.

\section{Weakness via representability: a case study} \label{sec:weakness}

The approach that we follow is the same that underlies the opetopic definition of weak higher categories: this is what Hermida calls \emph{coherence via universality} \cite{hermida2000representable, hermida2001coherent}, where coherent higher algebraic structure is subsumed by the existence of universally characterised elements in a wider context (representability). \index{coherence!via universality}

As is the case for opetopic sets, the category $\cat{G}$ of standard globes appears as a full subcategory of the shape category $\cat{RG}$ of regular polygraphs. This allows us to associate to any regular polygraph $X$ a globular set $GX$, defined as the restriction of $X$ to $\cat{G}$; when $X$ satisfies certain properties, $GX$ can be made the carrier of certain higher algebraic structures. \index{Globu@$G(-)$} \index{globular set!induced|see {Globu@$G(-)$}}

Because of the overlap between our globes and opetopes, we will see that the modelling of weak associativity can more or less be carried over from the opetopic approach; the original part is our treatment of equivalences and weak units. 
\begin{itemize}
	\item[] \textbf{Terminology.} We assume that a representative is fixed for all isomorphism classes of $n$-globes. An $n$-cell in a regular polygraph $X$ is an element of $X(G)$ for some atomic $n$-globe $G$; in the $\omega$-categorical version of polygraphs, this would have corresponded to an $n$-dimensional generator. 
	
	We write $\mathrm{Hom}(-,-)$ for $\mathrm{Hom}_\rpol(-,-)$. Recall that, by the Yoneda lemma, for all regular polygraphs $X$ and atomic $n$-globes $G$, $X(G) \simeq \mathrm{Hom}(G,X)$.
\end{itemize}

\subsection{Regular poly-bicategories and divisible 2-cells}

As a case study, we will focus on regular 2-polygraphs --- presheaves on the full subcategory of $\cat{RG}$ on atomic globes of dimension $\leq 2$ --- endowed with an associative composition of 2-cells. The immediate advantage is that atomic 2-globes are easily classified: they are simply the shapes 
\begin{equation*}
\input{img/c3_3_2shape.tex}
\end{equation*} 
with $n, m > 0$. Thus, a regular 2-polygraph $X$ is described by a diagram
\begin{equation*}
\begin{tikzpicture}
	\node[scale=1.25] (0) at (-2,0) {$X_0$};
	\node[scale=1.25] (1) at (0,0) {$X_1$};
	\node[scale=1.25] (2) at (2.4,0) {$X_2^{(n,m)}\;,$};
	\draw[1c] (-.35,.15) to node[auto,swap] {$\bord{}{+}$} (-1.7,.15);
	\draw[1c] (-.35,-.15) to node[auto] {$\bord{}{-}$} (-1.7,-.15);
	\draw[1c] (1.65,.15) to node[auto,swap] {$\bord{(j)}{+}$} (.3,.15);
	\draw[1c] (1.65,-.15) to node[auto] {$\bord{(i)}{-}$} (.3,-.15);
	\node[scale=1.25] at (5.25,.25) {$i = 1,\ldots,n,$};
	\node[scale=1.25] at (5.3,-.25) {$j = 1,\ldots,m,$};
\end{tikzpicture}
\end{equation*}
for all $n, m \geq 1$, satisfying
\begin{align*}
	& \bord{}{-}\bord{(1)}{-} = \bord{}{-}\bord{(1)}{+}\;, \\ 
	& \bord{}{+}\bord{(n)}{-} = \bord{}{+}\bord{(m)}{+}\;, \\
	& \bord{}{+}\bord{(i)}{-} = \bord{}{-}\bord{(i+1)}{-}\;, \quad \quad \;i=1,\ldots,n-1, \\
	& \bord{}{+}\bord{(j)}{-} = \bord{}{+}\bord{(j+1)}{+}\;, \quad \quad \;j=1,\ldots,m-1.
\end{align*}
Here $X_0, X_1$ are the values of the presheaf $X$ at the unique atomic 0-globe and 1-globe, and $X_2^{(n,m)}$ at the atomic 2-globe with $n$ input cells and $m$ output cells. The induced globular set $GX$ is
\begin{equation*}
\begin{tikzpicture}
	\node[scale=1.25] (0) at (-2,0) {$X_0$};
	\node[scale=1.25] (1) at (0,0) {$X_1$};
	\node[scale=1.25] (2) at (2.4,0) {$X_2^{(1,1)}\;.$};
	\draw[1c] (-.35,.15) to node[auto,swap] {$\bord{}{+}$} (-1.7,.15);
	\draw[1c] (-.35,-.15) to node[auto] {$\bord{}{-}$} (-1.7,-.15);
	\draw[1c] (1.65,.15) to node[auto,swap] {$\bord{(1)}{+}$} (.3,.15);
	\draw[1c] (1.65,-.15) to node[auto] {$\bord{(1)}{-}$} (.3,-.15);
\end{tikzpicture}
\end{equation*}

Composable pairs of $n$-cells in a regular polygraph are classified by mergeable pairs of atomic $n$-globes, equivalently, $n$-globes with two top-dimensional elements. 

\begin{cons} We can identify a non-atomic $n$-globe $G$ with a regular polygraph as follows. Take an atomic $(n+1)$-globe $G'$ with a single output $n$-cell $x$, and such that $\bord{}{-}G' \simeq G$; then, take the sub-presheaf $\bord{}{-}G'$ of the Yoneda embedding of $G'$, where we discard the $(n+1)$-cell corresponding to the identical embedding $G = G$ and the $n$-cell corresponding to the inclusion $U_x \subseteq G$. We will identify $G$ with the regular polygraph $\bord{}{-}G'$.
\begin{itemize}
\item[] \textbf{Terminology.} An $n$-globe in a regular polygraph $X$ is a map $G \to X$, where $G$ is the regular polygraph corresponding to a generic $n$-globe $G$.
\end{itemize}
In particular, we obtain composable pairs of $n$-cells in a regular polygraph $X$ as $n$-globes $G \to X$ where $G$ has two $n$-cells.
\end{cons} 

When $n = 2$, by a case distinction, we see that 2-globes with two 2-cells fall into one of the following geometric situations, where the shared boundaries contain at least one 1-cell, and for the rest only regularity constraints apply:
\begin{equation} \label{eq:mergeable2globes}
\input{img/c3_3_mergeable2globes.tex}
\end{equation}
We note that there is a correspondence between these, and the axioms of the ``spider presentation'' of (special) Frobenius algebras, Remark \ref{remark:spider}.

Similarly, $n$-globes with three top-dimensional elements classify sequences of two compositions. For our algebraic composition on regular 2-polygraphs, they will classify associativity and interchange axioms: whenever there are two possible trees of mergeable bi-partitions of the three cells, as in Example \ref{exm:mergeable}, composing them in one order or the other leads to the same result. 

Before considering this notion of composition, however, we will work with a more restrictive one, since it leads to some clearer, more symmetrical combinatorics.
\begin{dfn} \index{globe!simple} 
Let $X$ be an oriented thin poset, and let $U$ and $U'$ be two closed, pure, $n$-dimensional subsets of $X$. We say that $U$ and $U'$ are \emph{simply mergeable} if they are mergeable, and $U \cap U'$ is an atomic $(n-1)$-globe.

By induction on the number of top-dimensional elements of $U$: a pure, closed $n$-dimensional $U \subseteq X$ is a \emph{simple $n$-globe} if it is an atomic $n$-globe, or if there are $U_1$, $U_2$ as in the definition of $n$-globe, such that $U_1, U_2$ are simply mergeable simple $n$-globes.
\end{dfn}
Alternatively, a simple $n$-globe may be characterised as one which reduces to an atomic $n$-globe, as in Lemma \ref{lem:simplemerge}, by a sequence of simple mergers involving only $n$-dimensional cells. 

The simple 2-globes with two 2-cells are also classified by diagram (\ref{eq:mergeable2globes}), where the shared boundary is now constrained to consist of a single 1-cell. However, there are fewer possible simple 2-globes with three 2-cells. The following two pictures, for instance, describe a family both of 2-globes and of simple 2-globes, depending on constraints on the shared boundaries:
\begin{equation*}
\input{img/c3_3_simple_assoc.tex}
\end{equation*} 
the following, however, only describes a family of non-simple 2-globes, even though its 2-cells may be pairwise simply mergeable:
\begin{equation*}
\input{img/c3_3_nonsimple_assoc.tex}
\end{equation*} 

Far from being an arbitrary algebraic concoction, the kind of composition classified by simple 2-globes is the correct one for sequents in proof theory, the so-called cut rule; and it is the one underlying planar polycategories, a context for the categorical semantics of proof theory, see \cite{koslowski2003monadic, pastro2004sigmapi}. In \cite{cockett2003morphisms}, Cockett, Koslowski, and Seely gave a ``multi-coloured'' generalisation of planar polycategories, on which the following is based. 
\begin{dfn}
A \emph{regular poly-$n$-category} is a regular $n$-polygraph $X$, together with
\begin{enumerate}
	\item for all simple $n$-globes $G$ with two $n$-cells, reducing to the atomic $n$-globe $G'$ after one simple merger, a composition operation 
	\begin{equation*}
		\cutt{G}: \mathrm{Hom}(G,X) \to X(G'), 
	\end{equation*}
	restricting to the identity on $\mathrm{Hom}(\bord{}{\alpha}G, X) = \mathrm{Hom}(\bord{}{\alpha}G', X)$, and such that
	\item for all simple $n$-globes $G_{1,2,3}$ with three $n$-cells $x_1, x_2, x_3$, reducing to the atomic $n$-globe $G'$ in two steps through either the simple merger of $U_{x_1}$ and $U_{x_2}$, or the simple merger of $U_{x_2}$ and $U_{x_3}$, the following diagram commutes:
\begin{equation*}
\begin{tikzpicture}[baseline={([yshift=-.5ex]current bounding box.center)}]
	\node[scale=1.25] (0) at (-2,2) {$\mathrm{Hom}(G_{1,2,3},X)$};
	\node[scale=1.25] (1) at (2,2) {$\mathrm{Hom}(G_{12,3},X)$};
	\node[scale=1.25] (2) at (-2,0) {$\mathrm{Hom}(G_{1,23},X)$};
	\node[scale=1.25] (3) at (2,0) {$X(G').$};
	\draw[1c] (0) to node[auto] {$\cutt{G_{1,2}}$} (1);
	\draw[1c] (2) to node[auto,swap] {$\cutt{G_{1,23}}$} (3);
	\draw[1c] (0) to node[auto,swap] {$\cutt{G_{2,3}}$} (2);
	\draw[1c] (1) to node[auto] {$\cutt{G_{12,3}}$} (3);
\end{tikzpicture}
\end{equation*}
\end{enumerate}
A morphism $f: X \to Y$ of regular poly-$n$-categories is a map of the underlying regular polygraphs such that
\begin{equation*}
\begin{tikzpicture}[baseline={([yshift=-.5ex]current bounding box.center)}]
	\node[scale=1.25] (0) at (-1,1.5) {$\mathrm{Hom}(G,X)$};
	\node[scale=1.25] (1) at (2,1.5) {$X(G')$};
	\node[scale=1.25] (2) at (-1,0) {$\mathrm{Hom}(G,Y)$};
	\node[scale=1.25] (3) at (2,0) {$Y(G').$};
	\draw[1c] (0) to node[auto] {$\cutt{G}$} (1);
	\draw[1c] (2) to node[auto,swap] {$\cutt{G}$} (3);
	\draw[1c] (0) to node[auto,swap] {$\mathrm{Hom}(G,f)$} (2);
	\draw[1c] (1) to node[auto] {$f_{G'}$} (3);
\end{tikzpicture}
\end{equation*}
commutes for all suitable $G \leadsto G'$. We write $\cat{RPnCat}$ for the category of regular poly-$n$-categories and morphisms.
\end{dfn}

We have taken some notational shortcuts, which we hope are self-explanatory. The endofunctors $\oppn{(-)}{S}$ on $\globpos$ induce endofunctors on $\cat{RPnCat}$: if $X$ has compositions $\cutt{G}$, $\oppn{X}{S}$ has compositions $\mathrm{cut}'_G := \cutt{\oppn{G}{S}}$, defined through the isomorphisms $\mathrm{Hom}(G, \oppn{X}{S}) \simeq \mathrm{Hom}(\oppn{G}{S},X)$.

We speak of regular poly-bicategories in the 2-dimensional case, where the $\mathrm{cut}$ operations correspond to functions \index{poly-bicategory!regular}
\begin{equation*}
\begin{tikzpicture}
	\node[scale=1.25] (0) at (-2.2,0) {$X_2^{(n,m)} {}_{\bord{(i)}{+}}\!\!\times_{\bord{(j)}{-}} X_2^{(p,q)}$};
	\node[scale=1.25] (1) at (2.5,0) {$X_2^{(n+p-1,m+q-1)},$};
	\draw[1c] (0) to node[auto] {$\cutt{i,j}$} (1);
\end{tikzpicture}
\end{equation*}
such that $i, j$ satisfy the two conditions on any side of the following square:
\begin{equation*}
\begin{tikzpicture}
	\node[scale=1.25] (0) at (-1,.6) {$i=1$};
	\node[scale=1.25] (1) at (1,.6) {$j=p$};
	\node[scale=1.25] (2) at (-1,-.6) {$i=m$};
	\node[scale=1.25] (3) at (1,-.6) {$j=1$};
	\draw[edge] (-.5,.6) to node[auto] {$(b)$} (.5,.6);
	\draw[edge] (-.5,-.6) to node[auto,swap] {$(d)$} (.5,-.6);
	\draw[edge] (-1,.4) to node[auto,swap] {$(a)$} (-1,-.4);
	\draw[edge] (1,.4) to node[auto] {$(c)$} (1,-.4);
	\node[scale=1.25] at (1.6,-.65) {.};
\end{tikzpicture}
\end{equation*}
Note that for the conditions $(a)$ and $(c)$ to be satisfied, it is necessary that $m = 1$ and $p = 1$, respectively. There are 9 associativity schemes and 8 interchange schemes that these must satisfy; we refer to \cite[Section 1.7]{pastro2004sigmapi} for an explicit classification.

\begin{remark}
When $X$ is a regular poly-bicategory, $GX$ is equipped with the strictly associative composition of 2-cells $\cutt{1,1}$.
\end{remark}

\begin{dfn}
A \emph{regular multi-bicategory} is a regular poly-bicategory $X$ such that $X_2^{(n,m)}$ is empty whenever $m > 1$. A \emph{regular polycategory} is a regular poly-bicategory with a single 0-cell. A \emph{regular multicategory} is a regular multi-bicategory with a single 0-cell.
\end{dfn} \index{multi-bicategory|see {poly-bicategory}} \index{polycategory|see {poly-bicategory}} \index{multicategory|see {poly-bicategory}} 

We will write 
\begin{equation*}
	p: (a_1, \ldots, a_n) \to (b_1, \ldots, b_m)
\end{equation*}
for a 2-cell $p \in X_2^{(n,m)}$ in a regular poly-bicategory, with $\bord{(i)}{-}p = a_i$, $\bord{(j)}{+}p = b_j$, for $i = 1,\ldots, n$, and $j=1,\ldots m$; we will omit the parentheses when there is a single 1-cell in either boundary. We will use Greek letters $\Gamma, \Delta$ for generic composable sequences of 1-cells.

The construction of a weakly associative composition of 1-cells in a poly-bicategory can be carried out exactly as in \cite{cockett2003morphisms}, which is in turn a two-sided version of \cite{hermida2000representable}, where the multi-bicategorical case is considered. We paraphrase it here, before moving on to the original treatment of weak units. 

The main idea of Hermida's \cite{hermida2000representable} is the following: suppose that we want to ``internalise'' the operation of forming sequences of composable 1-cells in a multi-bicategory $X$ as an algebra on its set of 1-cells, say $(a,b) \mapsto a \# b$. Then, we should require that 2-cells
\begin{equation} \label{eq:2celltorep}
	p: (\Gamma_1, a, b, \Gamma_2) \to d
\end{equation}
correspond bijectively to 2-cells
\begin{equation*}
	p': (\Gamma_1, a \# b, \Gamma_2) \to d.
\end{equation*}
In particular, the identity $a \# b \to a \# b$ must correspond to a 2-cell $t: (a, b) \to a \# b$. The bijective correspondence can be obtained from the existence of a unique factorisation of each 2-cell $p$ as in (\ref{eq:2celltorep}) through $t$, that is, from a universal property of the 2-cell $t$. 

Then, any 1-cell $c$ such that there exists a universal 2-cell $(a, b) \to c$ with the same property is an equally valid candidate for a composite of $a, b$. Such a 1-cell is an ``internal representative'' of the composable pair $(a,b)$ in the set of 1-cells; not working with our regularity constraint, Hermida also postulates the existence of 1-cells representing ``empty sequences'' (0-cells). The existence of representatives for arbitrary sequences is called \emph{representability}. 

Hermida then shows that any choice of representatives defines a bicategory structure \cite{benabou1967introduction} on the 2-globular set $GX$; and that performing such a choice for each representable multi-bicategory defines an equivalence between the category of bicategories and pseudofunctors, and the category of representable multi-bicategories and maps that preserve universal 2-cells. In this precise sense, the structure of a bicategory is equivalent to a property of multi-bicategories.

Beyond the one considered by Hermida, there are other universal properties that 2-cells in a multi-bicategory may satisfy. Given a 2-cell $p: (a,b) \to c$, there are three ways in which it can be composed with another 2-cell:
\begin{equation} \label{eq:2divisibility}
\input{img/c3_3_2divisibility.tex}
\end{equation} 
Note that we will usually identify pasting diagrams of 2-globes with the unique result of their composition in any order. \index{divisible!2-cell!$(c^*), (l_*), (r_*)$}

\begin{dfn} 
Let $a$, $b$ be two 1-cells in a regular poly-bicategory $X$, with compatible boundaries as required separately by each definition. 

A \emph{tensor} of $a$ and $b$ is a 1-cell $a \otimes b$, together with a 2-cell $t_{a,b}: (a,b) \to a\otimes b$ that is divisible for compositions of type $(c^*)$: that is, for all \index{tensor!in a poly-bicategory} 
\begin{equation*} 
p: (\Gamma_1, a, b, \Gamma_2) \to (\Delta),
\end{equation*}
there exists a unique 
\begin{equation*} 
\tilde{p}: (\Gamma_1, a \otimes b, \Gamma_2) \to (\Delta),
\end{equation*}
such that
\begin{equation*}
\input{img/c3_3_tensor.tex}
\end{equation*} 
\begin{itemize}
	\item[] \textbf{Notation.} While we use the same symbol for consistency with the existing literature, the tensor of two 1-cells, as defined here, must not be confused with the product of two 1-cells in the tensor product of two polygraphs, which is a 2-cell.
\end{itemize}

A \emph{right hom} from $a$ to $b$ is a 1-cell $\rimp{a}{b}$, together with a 2-cell $e^R_{a,b}: (a, \rimp{a}{b}) \to b$ that is divisible for compositions of type $(r_*)$; that is, for all
\begin{equation*} 
p: (a, \Gamma) \to (b, \Delta),
\end{equation*}
there exists a unique 
\begin{equation*} 
\tilde{p}: (\Gamma) \to (\rimp{a}{b}, \Delta),
\end{equation*}
such that
\begin{equation*}
\input{img/c3_3_righthom.tex}
\end{equation*} 

A \emph{left hom} from $a$ to $b$ is a 1-cell $\limp{a}{b}$, together with a 2-cell $e^L_{a,b}: (\limp{a}{b}, a) \to b$ that is divisible for compositions of type $(l_*)$. \index{hom} \index{left hom|see {hom}} \index{right hom|see {hom}}

We say that $X$ is \emph{tensor 1-representable}, \emph{right closed}, and \emph{left closed}, respectively, if it has tensors, right homs, and left homs, respectively, for all pairs of 1-cells in the appropriate configuration. \index{poly-bicategory!1-representable!tensor} \index{poly-bicategory!closed}
\end{dfn}
A right hom in $X$ is the same as a left hom in $\opp{X}$. By usual reasoning with universal properties, the 2-cells exhibiting binary tensors, homs, and their duals can be composed to produce $n$-ary tensors, homs, and their duals. For example,
\begin{equation*}
\input{img/c3_3_composites.tex}
\end{equation*} 
exhibit $(a \otimes b) \otimes c$ as a tensor of $(a, b, c)$, and $\rimp{a}{(\rimp{b}{c})}$ as a right hom from $(b,a)$ to $c$, respectively. 

All the definitions can also be dualised to $\coo{X}$; given a 2-cell $p: c \to (a,b)$, there are three ways in which it can be composed with another 2-cell, dual to those in diagram (\ref{eq:2divisibility}):
\begin{equation*}
\input{img/c3_3_co2divisibility.tex}
\end{equation*} 
\begin{dfn}
Let $a$, $b$ be two 1-cells in a regular poly-bicategory $X$, with compatible boundaries as required separately by each definition. \index{divisible!2-cell!$(c_*), (l^*), (r^*)$}

A \emph{par} of $a$ and $b$ is a 1-cell $a \parr b$, together with a 2-cell $p_{a,b}: a \parr b \to (a, b)$ that is divisible for compositions of type $(c_*)$. A \emph{right cohom} from $a$ to $b$ is a 1-cell $\rcimp{a}{b}$, together with a 2-cell $c^R_{a,b}: b \to (a, \rcimp{a}{b})$ that is divisible for compositions of type $(r^*)$. A \emph{left cohom} from $a$ to $b$ is a 1-cell $\lcimp{a}{b}$, together with a 2-cell $c^L_{a,b}: b \to (\lcimp{a}{b}, a)$ that is divisible for compositions of type $(l^*)$. \index{par} \index{cohom} \index{left cohom|see {cohom}} \index{right cohom|see {cohom}}

We say that $X$ is \emph{par 1-representable}, \emph{right coclosed}, and \emph{left coclosed}, respectively, if it has pars, right cohoms, and left cohoms, respectively, for all pairs of 1-cells in the appropriate configuration. \index{poly-bicategory!1-representable!par}
\end{dfn}
Clearly, $p: (a,b) \to c$ exhibits a tensor, a right hom, a left hom in $X$, respectively, if and only if $\coo{p}: c \to (a,b)$ exhibits a par, a right cohom, a left cohom in $\coo{X}$, respectively.

For all compatible pairs of 1-cells, we expect these constructions to define a unique 1-cell ``up to isomorphism'' --- except, we have not shown how to define isomorphisms in the absence of units for composition. The idea is to use, again, divisibility properties.
\begin{dfn}
A 2-cell $id_a: a \to a$ in a regular poly-bicategory $X$ is a \emph{unit} on $a$ if the result of any composition of $id_a$ with another 2-cell $p$ is $p$. 

A 2-cell $p: a \to a'$ with a single input and output is \emph{divisible} if it is both $(c^*)$-divisible and $(c_*)$-divisible. \index{divisible!2-cell}
\end{dfn}

\begin{lem} \label{lem:2divunit}
Let $p: a \to a'$ be a divisible 2-cell in a regular poly-bicategory $X$, and let $id: a \to a$, $id': a' \to a'$ be the unique 2-cells such that 
\begin{equation*}
	\cutt{1,1}(id, p) = p, \quad\quad \cutt{1,1}(p, id') = p.
\end{equation*}
Then $id$ is a unit on $a$, and $id'$ is a unit on $a'$.
\end{lem}
\begin{proof}
Consider a 2-cell $q: (\Gamma_1, a, \Gamma_2) \to (\Delta)$ in $X$; by $(c^*)$-divisibility of $p$, there exists a unique $q': (\Gamma_1, a', \Gamma_2) \to (\Delta)$ such that $\cutt{1,i}(p,q') = q$. Then,
\begin{align*}
	\cutt{1,i}(id, q) & = \cutt{1,i}(id, \cutt{1,i}(p,q')) = \qquad \qquad \qquad \qquad \text{(associativity)} \\
	& = \cutt{1,i}(\cutt{1,1}(id,p),q') = \cutt{1,i}(p,q') = q.
\end{align*}
Next, consider a 2-cell $r: (\Gamma) \to (\Delta_1, a, \Delta_2)$. Composing it with $p$, we find
\begin{equation*}
	\cutt{j,1}(r, p) = \cutt{j,1}(r, \cutt{1,1}(id,p)) = \cutt{j,1}(\cutt{j,1}(r,id),p);
\end{equation*}
therefore, by $(c_*)$-divisibility of $p$, we obtain $r = \cutt{j,1}(r,id)$. This proves that $id$ is a unit on $a$; the same proof applied to $\coo{X}$ proves that $id'$ is a unit on $a'$.
\end{proof}

\begin{prop} \label{prop:2units}
Let $X$ be a regular poly-bicategory. Then, the following are equivalent:
\begin{enumerate}
	\item for all 1-cells $a$ of $X$, there exist a 1-cell $\overline{a}$ and a divisible 2-cell $p: a \to \overline{a}$;
	\item for all 1-cells $a$ of $X$, there exist a 1-cell $\underline{a}$ and a divisible 2-cell $p': \underline{a} \to a$;
	\item for all 1-cells $a$ of $X$, there exists a (necessarily unique) unit $id_a$ on $a$.
\end{enumerate}
\end{prop}
\begin{proof}
Let $a$ be a 1-cell of $X$. If there exists a unit $id_a: a \to a$, it is clearly divisible, and fulfills the other two conditions for $a$.

Conversely, by Lemma \ref{lem:2divunit}, from any divisible 2-cell $e: a \to \overline{a}$, and from any divisible 2-cell $e': \underline{a} \to a$, we can construct a unit on $a$.
\end{proof}

\begin{dfn}
We say that a regular poly-bicategory $X$ is \emph{1-representable} if it satisfies any of the equivalent conditions of Proposition \ref{prop:2units}. \index{poly-bicategory!1-representable}
\end{dfn} 

\begin{cor} \label{cor:onedivisible}
Let $X$ be a 1-representable regular poly-bicategory, and $p: a \to a'$ a 2-cell of $X$ with a single input and output. Then, if $p$ satisfies any divisibility property, it satisfies all of them.
\end{cor}
\begin{proof}
Obviously, $(c^*)$, $(l^*)$, $(r^*)$-divisibility collapse to the same property in this case, and the same holds for $(c_*)$, $(l_*)$, $(r_*)$. 

Suppose $p$ is $(c^*)$-divisible, and let $id_a, id_{a'}$ be the units on $a$, $a'$. Dividing $id_a$ by $p$, we obtain a unique $\bar{p}: a' \to a$ such that $\cutt{1,1}(p,\bar{p}) = id_a$. Since 
\begin{equation*}
	\cutt{1,1}(p, \cutt{1,1}(\bar{p},p)) = \cutt{1,1}(id_a,p) = p = \cutt{1,1}(p, id_{a'}),
\end{equation*}
it follows by uniqueness that $\cutt{1,1}(\bar{p},p) = id_{a'}$. The other direction is analogous.
\end{proof}
Thus, assuming 1-representability, we retrieve the familiar notion of isomorphism: a 2-cell $p: a \to a'$ that has a two-sided inverse with respect to the unique units $id_a, id_{a'}$, or that, equivalently, satisfies any divisibility property. We can then apply standard reasoning to show that $a \otimes b, \rimp{a}{b}, \limp{a}{b}$, and their duals, when they exist, are unique up to a unique isomorphism. We will write $a \simeq a'$ if an isomorphism $p: a \to a'$ exists.

The following is a useful technical lemma.
\begin{lem} \label{lem:twodivisible}
Let $t: (a, b) \to c$ be a 2-cell in a 1-representable regular poly-bicategory, and suppose $t$ satisfies two different divisibility properties. Then:
\begin{itemize}
	\item if $t$ is $(c^*)$ and $(r_*)$-divisible, then any $(c^*)$-divisible $u: (a, b) \to c'$ is also $(r_*)$-divisible, and any $(r_*)$-divisible $u: (a, b') \to c$ is also $(c^*)$-divisible;
	\item if $t$ is $(c^*)$ and $(l_*)$-divisible, then any $(c^*)$-divisible $u: (a, b) \to c'$ is also $(l_*)$-divisible, and any $(l_*)$-divisible $u: (a', b) \to c$ is also $(c^*)$-divisible;
	\item if $t$ is $(r_*)$ and $(l_*)$-divisible, then any $(r_*)$-divisible $u: (a, b') \to c$ is also $(l_*)$-divisible, and any $(l_*)$-divisible $u: (a', b) \to c$ is also $(r_*)$-divisible.
\end{itemize}
\end{lem}
\begin{proof}
We only consider the case in which $t$ is $(c^*)$ and $(r_*)$-divisible, and $u: (a, b') \to c$ is $(r_*)$-divisible; the others are completely analogous. By factoring $u$ through $t$, we obtain
\begin{equation*}
\input{img/c3_3_twodivisible_fix.tex}
\end{equation*} 
for a unique $p: b' \to b$, which due to 1-representability must be an isomorphism. As composition with isomorphisms does not affect divisibility properties, and $t$ is $(c^*)$-divisible, $u$ is also $(c^*)$-divisible.
\end{proof}

We obtain an analogous result for cells $t': c \to (a,b)$, by dualising to $\coo{X}$.

\subsection{Weak units and divisible 1-cells}
In the proof of Proposition \ref{prop:2units}, we assumed a representability condition, and constructed units for a strictly associative composition. That proof is going to be our blueprint for the construction of weak unit 1-cells, relative to the two internal, weakly associative notions of composition given by tensors and pars. We will proceed as follows:
\begin{enumerate}
	\item first, we will define a notion of weak unit 1-cell appropriate for regular poly-bicategories;
	\item then, we will show that the existence of weak units is equivalent to a lower-dimensional representability condition;
	\item finally, we will show that our weak units induce the same coherent structure on $GX$ as the ones based on degenerate boundaries in \cite{hermida2000representable, cockett2003morphisms}.
\end{enumerate}
We restrict our attention to tensor units first. 

\begin{dfn} \label{dfn:unit} \index{unit!tensor} \index{weak unit|see {unit}}
Let $x$ be a 0-cell in a regular poly-bicategory $X$. A 1-cell $1_x: x \to x$ is a \emph{tensor unit} on $x$ if, for all $a: x \to y$, $b: z \to x$, there exist 2-cells
\begin{equation*}
\input{img/c3_3_tensorunit.tex}
\end{equation*}
that are, respectively, $(c^*)$, $(r_*)$-divisible, and $(c^*), (l_*)$-divisible, that is, they exhibit $a$ as both $1_x \otimes a$ and $\rimp{1_x}{a}$, and $b$ as both $b \otimes 1_x$ and $\limp{1_x}{b}$.
\end{dfn}

\begin{remark} \label{remark:units}
By Lemma \ref{lem:twodivisible}, if $X$ is 1-representable, and $1_x: x \to x$ is a tensor unit on $x$, a seemingly stronger claim can be made: all $(c^*)$-divisible 2-cells of the form $t: (1_x, a) \to a'$ are also $(r_*)$-divisible, all $(r_*)$-divisible 2-cells of the form $u: (1_x, a') \to a$ are also $(c^*)$-divisible, and so on. The same reasoning can be applied to par units.
\end{remark}

In the following definition, we introduce a notion of equivalence 1-cell which does not depend on any specified units, nor, in fact, on the existence of units.

\begin{dfn} \index{divisible!1-cell!tensor} \index{tensor divisible|see {divisible}}
Let $e: x \to x'$ be a 1-cell in a regular poly-bicategory $X$. We say that $e$ is \emph{tensor left divisible} if, for all $a: x \to y$, $a': x' \to y$, a right hom and tensor
\begin{equation*}
\input{img/c3_3_tensordiv.tex}
\end{equation*}
exist and are both $(c^*)$ and $(r_*)$-divisible, that is, they exhibit $a$ as $e \otimes (\rimp{e}{a})$ and $a'$ as $\rimp{e}{(e \otimes a')}$.

Dually, we say that $e$ is \emph{tensor right divisible} if, for all $b: z \to x$, $b': z \to x'$, a left hom and tensor
\begin{equation*}
\input{img/c3_3_tensordivright.tex}
\end{equation*}
exist and are both $(c^*)$ and $(l_*)$-divisible, that is, they exhibit $b'$ as $(\limp{e}{b'}) \otimes e$ and $b$ as $\limp{e}{(b \otimes e)}$. A 1-cell $e$ is \emph{tensor divisible} if it is both tensor left and tensor right divisible.
\end{dfn}

\begin{remark}
A \emph{quasigroup}, in the equational formulation \cite[Section 1.2]{smith2006introduction}, is a set $Q$ together with three binary operations $\cdot, \diagup, \diagdown$ satisfying the axioms
\begin{align*}
	& x \cdot (\rcimp{x}{y}) = x, & \rcimp{x}{(x \cdot y)} = x, \\
	& (\lcimp{x}{y}) \cdot x = x, & \lcimp{x}{(y \cdot x)} = x.
\end{align*}
The isomorphisms enforced by divisibility can be seen as a categorified version of these equations.
\end{remark}

The following shows that tensor divisible 1-cells satisfy a ``two-out-of-three'' property, a common requirement for classes of weak equivalences.

\begin{thm} \label{thm:2outof3}
Let $e: x \to x'$, $e': x' \to x''$, $e'': x \to x''$ be 1-cells in a 1-representable regular poly-bicategory $X$, and suppose that
\begin{equation*}
\input{img/c3_3_2outof3.tex}
\end{equation*}
satisfies two different divisibility properties. If two of the three 1-cells $e$, $e'$, and $e''$ are tensor divisible, then the third is also tensor divisible, and $p$ is at once $(c^*)$, $(l_*)$, and $(r_*)$-divisible.
\end{thm}
\begin{proof}
Suppose first that $e, e'$ are tensor divisible. Because $p$ satisfies two different divisibility properties, it satisfies at least one of $(l_*)$ and $(r_*)$-divisibility. Then, factoring it through $e^L_{e',e''}$, in the first case, or through $e^R_{e,e''}$, in the second case, we find that it is also $(c^*)$-divisible, that is, it exhibits $e''$ as $e \otimes e'$. Finally, factoring through $t_{e,e'}$, which, by tensor divisibility of both $e$ and $e'$ and an application of Lemma \ref{lem:twodivisible}, is at once $(c^*)$, $(l_*)$, and $(r_*)$-divisible, we obtain that $p$ has the same property. 

To check that $e''$ is tensor left divisible, consider an arbitrary 1-cell $a : x \to y$; then, the unique 2-cell $q$ obtained in the factorisation
\begin{equation} \label{eq:2outof3_1}
\input{img/c3_3_2outof3_tensor.tex}
\end{equation} 
using the $(c^*)$-divisibility of $p$ is both $(c^*)$ and $(r_*)$-divisible, as we will now show. To prove $(c^*)$-divisibility of $q$, given a 2-cell 
\begin{equation*}
	p': (\Gamma_1, e'', \rimp{e'}{(\rimp{e}{a})}, \Gamma_2) \to (\Delta), 
\end{equation*}
we precompose it with $p$ and factorise as
\begin{equation*} 
\input{img/c3_3_2outof3_factor3.tex}
\end{equation*} 
for a unique 2-cell $p''$, using the $(c^*)$-divisibility of both sides of equation (\ref{eq:2outof3_1}). By $(c^*)$-divisibility of $p$, we obtain a necessarily unique factorisation of $p'$ through $q$. The $(r_*)$-divisibility of $q$ is obtained by an analogous argument. 

Similarly, given an arbitrary 1-cell $a': x'' \to y$, the unique 2-cell $q'$ obtained in the factorisation
\begin{equation*} 
\input{img/c3_3_2outof3_tensor2.tex}
\end{equation*} 
using the $(c^*)$-divisibility of $p$ is shown to be both $(c^*)$ and $(r_*)$-divisible. This proves that $e''$ is tensor left divisible; a symmetrical argument shows that it is tensor right divisible.

Next, we consider the case in which $e$ and $e''$ are tensor divisible. Because $p$ satisfies two different divisibility properties, it satisfies at least one of $(c^*)$ and $(r_*)$-divisibility; factoring it through $t_{e,e'}$, in the first case, or through $e^R_{e,e''}$, in the second case, we find that it must actually satisfy both of them. 

To check that $e'$ is tensor left divisible, consider a 1-cell $b : x' \to y$. The unique 2-cell $r$ obtained in the factorisation
\begin{equation} \label{eq:2outof3_2}
\input{img/c3_3_2outof3_righthom.tex}
\end{equation} 
using the $(r_*)$-divisibility of $t_{e,b}$ is easily shown to be $(r_*)$-divisible. To prove that it is also $(c^*)$-divisible, consider a 2-cell
\begin{equation*}
	p': (\Gamma_1, e', \rimp{e''}{(e \otimes b)}, \Gamma_2) \to (\Delta).
\end{equation*}
Suppose that $(\Gamma_1) = (\Gamma_1', c)$ for some $c: z \to x'$; then, we can perform the following sequence of factorisations (labels of 0-cells are omitted):
\begin{equation*} 
\input{img/c3_3_2outof3_factor.tex}
\end{equation*} 
\begin{equation*} 
\input{img/c3_3_2outof3_factor2.tex}
\end{equation*} 
for unique 2-cells $\tilde{p}$ and $p''$, where we used first the $(c^*)$-divisibility of the two sides of equation (\ref{eq:2outof3_2}), then the $(c^*)$-divisibility of $e^L_{e,c}$. Cancelling the latter, we obtain a necessarily unique factorisation of $p'$ through $r$. In case $\Gamma_1$ is empty, let $(\Delta) = (c', \Delta')$ for some $c': x' \to z$, and apply the same reasoning to the postcomposition of $p'$ with $t_{e,c}: (e, c) \to e \otimes c$.

Next, consider a 1-cell $b': x'' \to y$. The unique 2-cell $r'$ obtained in the factorisation
\begin{equation*} 
\input{img/c3_3_2outof3_righthom2.tex}
\end{equation*} 
using the $(r_*)$-divisibility of $e^R_{e,e'' \otimes b'}$ is shown to be both $(c^*)$ and $(r_*)$-divisible by a similar argument. This proves that $e'$ is tensor left divisible.

The proof that $e'$ is tensor right divisible is similar, and involves the 2-cells $s, s'$ obtained from factorisations
\begin{equation*} 
\input{img/c3_3_2outof3_lefthom.tex}
\end{equation*} 
\begin{equation*} 
\input{img/c3_3_2outof3_lefthom2.tex}
\end{equation*} 
Once we have acquired that $e'$ is tensor divisible, the fact that $p$ is also $(c_*)$-divisible is a consequence of Lemma \ref{lem:twodivisible}.

Finally, the statement in the case where $e', e''$ are tensor divisible follows from the previous case applied to $\opp{X}$.
\end{proof}

\begin{cor} \label{cor:div_closure}
The class of tensor divisible 1-cells in a 1-representable regular poly-bicategory is closed under tensors, left homs, and right homs.
\end{cor}
\begin{proof}
If $e''$ is obtained as a tensor, right hom, or left hom of tensor divisible 1-cells, then the 2-cell that exhibits it falls under the hypotheses of Theorem \ref{thm:2outof3}. It follows that $e''$ is also tensor divisible.
\end{proof}

\begin{lem} \label{lem:1sideunit}
Let $e: x \to x'$ be a tensor divisible 1-cell in a 1-representable regular poly-bicategory. Then, $\limp{e}{e}: x \to x$ is a tensor unit on $x$, and $\rimp{e}{e}: x' \to x'$ is a tensor unit on $x'$.
\end{lem}
\begin{proof}
Let $a: x \to y$ be a 1-cell. Then, the unique 2-cell $p$ obtained in the factorisation
\begin{equation*}
\input{img/c3_3_leftunit.tex}
\end{equation*} 
by using the $(c^*)$-divisibility of $e^R_{e,a}$ is easily determined to be $(c^*)$-divisible. By Corollary \ref{cor:div_closure}, $\limp{e}{e}$ is tensor divisible, so by Lemma \ref{lem:twodivisible} $p$ is also $(r_*)$-divisible.

Similarly, take any $a': x' \to y$. The unique 2-cell $p'$ obtained in the factorisation
\begin{equation*}
\input{img/c3_3_leftunitimpl.tex}
\end{equation*} 
by using the $(r_*)$-divisibility of $t_{e,a'}$ is itself $(r_*)$-divisible. Since $\rimp{e}{e}$ is tensor divisible, $p'$ is also $(c^*)$-divisible. 

This proves the left-handed condition of tensor units for both $\limp{e}{e}$ and $\rimp{e}{e}$; a dual argument in $\opp{X}$ leads to the right-handed condition.
\end{proof}

\begin{thm} \label{thm:0representable}
Let $X$ be a 1-representable regular poly-bicategory. Then, the following conditions are equivalent:
\begin{enumerate}
	\item for all 0-cells $x$ of $X$, there exist a 0-cell $\overline{x}$ and a tensor divisible 1-cell $e: x \to \overline{x}$;
	\item for all 0-cells $x$ of $X$, there exist a 0-cell $\underline{x}$ and a tensor divisible 1-cell $e': \underline{x} \to x$;
	\item for all 0-cells $x$ of $X$, there exists a tensor unit $1_x$ on $x$.
\end{enumerate}
\end{thm}
\begin{proof}
Let $x$ be a 0-cell of $X$. If there exists a tensor unit $1_x: x \to x$, then it is clearly tensor divisible, and fulfills the condition for $x$ on both sides.

Conversely, suppose $e: x \to \overline{x}$ is a tensor divisible 1-cell, and define $1_x := \limp{e}{e}$. By Lemma \ref{lem:1sideunit}, $1_x$ is a tensor unit on $x$. Similarly, if $e': \underline{x} \to x$ is tensor divisible, then ${1_x}': \rimp{e'}{e'}$ is a tensor unit on $x$. This completes the proof.
\end{proof}

\begin{dfn} \index{poly-bicategory!representable!tensor}
We say that a 1-representable regular poly-bicategory $X$ is \emph{tensor 0-representable} if it satisfies any of the equivalent conditions of Theorem \ref{thm:0representable}. \index{poly-bicategory!0-representable!tensor}

We say that $X$ is \emph{tensor representable} if it is 1-representable, tensor 0-representable, and tensor 1-representable. 
\end{dfn}

When $X$ is tensor representable, we have isomorphisms $1_x \otimes a \simeq a$, $b \otimes 1_x \simeq b$ in $GX$ whenever the left-hand side is defined; it only remains to be show that a coherent set of units and isomorphisms --- that is, one that satisfies Mac Lane's triangle axioms --- can be chosen. The point is that our definition of tensor unit subsumes the notion called \emph{Saavedra unit} in \cite{kock2008elementary}, on which the Joyal-Kock weak units defined in Section \ref{sec:dimension} are also based: \index{unit!Saavedra}
\begin{enumerate}
	\item the $l_{1_x}$ or $r_{1_x}$ are witnesses of $1_x \otimes 1_x \simeq 1_x$, the \emph{idempotency} of the $1_x$;
	\item the divisibility conditions on the $l_a, r_a$ imply the \emph{cancellability} of the $1_x$.
\end{enumerate}
The following proof is based on \cite[Proposition 2.6]{kock2008elementary}. It is due remarking that this, like all theorems of ``coherence via universality'', relies on a choice of representatives, so foundational \emph{caveat}s apply, relative to the axiom of choice.

\begin{prop} \label{thm:bicategory}
Let $X$ be a regular poly-bicategory, and suppose $X$ is tensor representable. Then $GX$ admits the structure of a bicategory, with tensors as composites of 1-cells.

In particular, if $X$ has a single 0-cell, $GX$ admits the structure of a monoidal category. 
\end{prop}
\begin{proof}
The part relative to the coherent weak associativity of the composition is handled as in \cite[Definition 9.6]{hermida2000representable}, so we only need to show that a coherent set of units and unitors can be chosen. For all composable 1-cells $a, b$, let $t_{a,b}: (a, b) \to a \otimes b$ be the $(c^*)$-divisible cell chosen as a witness of the composition.

By Remark \ref{remark:units}, given a tensor unit $1_x$, we know that the $t_{1_x,a}: (1_x, a) \to 1_x \otimes a$ are also $(r_*)$-divisible. We define, for all $a$, an isomorphism $\lambda_a: 1_x \otimes a \to a$ as the unique 2-cell induced by the factorisation
\begin{equation*}
\input{img/c3_3_unitor_left.tex}
\end{equation*} 
Let $p;q := \cutt{1,1}(p,q)$ for two composable 2-cells $p,q$ in $GX$. Naturality of the $\lambda_a$ in $a$ follows from the fact that, for all $p: a \to b$, the two sequences of factorisations of
\begin{equation*}
\input{img/c3_3_unitor_natural.tex}
\end{equation*} 
respectively leading to $(1_x\otimes p);\lambda_b$ and to $\lambda_a;p$ must yield the same result. Similarly, we define isomorphisms $\rho_b: b \otimes 1_x \to b$, natural in $b$, as the unique 2-cells induced by the factorisations
\begin{equation} \label{eq:unitor_right}
\input{img/c3_3_unitor_right.tex}
\end{equation} 
Observe that, by Lemma \ref{lem:twodivisible}, we can always pick $l_{1_x} = r_{1_x}$, so there is no loss of generality in starting from one rather than the other.

The proof that the $\lambda_a$ and $\rho_a$ satisfy Mac Lane's triangle axioms then proceeds as in \cite[Proposition 2.6]{kock2008elementary}: omitting associators, which can be inserted where appropriate, for all compatible $a, b$, the invertible 2-cells $\rho_a \otimes 1_x \otimes b$ and $a \otimes 1_x \otimes \lambda_b$ are equal by definition; since $(\rho_a \otimes 1_x \otimes b);(a \otimes \lambda_b) = (a \otimes 1_x \otimes \lambda_b);(\rho_a \otimes b)$ up to associators, it follows that $\rho_a \otimes b = a \otimes \lambda_b$ up to associators.
\end{proof}

\begin{remark} \label{remark:unitor_preserve}
We can modify the first part of the previous proof to show the following. Let $l_a: (1_x, a) \to a$, $r_a: (b, 1_x) \to b$ be arbitrary witnesses of unitality in a 1-representable, tensor 0-representable regular poly-bicategory, and assume without loss of generality $r_{1_x} = l_{1_x}$. Then, we can define new families $\tilde{l}_a, \tilde{r}_a$ as follows: let $\tilde{l}_a := \cutt{1,1}(l_a, e_a)$, where $e_a: a \to a$ is obtained by the factorisation 
\begin{equation*}
\input{img/c3_3_unitor_left_bis.tex}
\end{equation*}
and $\tilde{r}_a := \cutt{1,1}(r_a, e'_a)$, where $e'_a$ comes from the dual factorisation, based on the $r_a$ and diagram (\ref{eq:unitor_right}). Then, the $\{\tilde{r}_a, \tilde{l}_a\}$ satisfy poly-bicategorical versions of naturality and the triangle axioms: that is, for all $p, p', q$ as pictured,
\begin{equation} \label{eq:preserve}
\input{img/c3_3_unitor_preserve.tex}
\end{equation}
\begin{equation} \label{eq:preserve2}
\input{img/c3_3_unitor_preserve2.tex}
\end{equation}
\begin{equation} \label{eq:triangle}
\input{img/c3_3_unitor_triangle.tex}
\end{equation}
To prove the first one, precompose the left-hand side with $l_{1_x}$; then, we have two possible factorisations (labels of 1-cells are omitted to avoid clutter):
\begin{equation*}
\input{img/c3_3_unitor_proof.tex}
\end{equation*}
by definition of $\tilde{l}_b$; and, for some unique $\tilde{p}$ given by $(c^*)$-divisibility of $l_a$,
\begin{equation*}
\input{img/c3_3_unitor_proof2.tex}
\end{equation*}
\begin{equation*}
\input{img/c3_3_unitor_proof3.tex}
\end{equation*}
Comparing the two, and using the $(r^*)$-divisibility of $\tilde{l}_b = \cutt{1,1}(l_b,e_b)$, we obtain equation (\ref{eq:preserve}). Equation (\ref{eq:preserve2}) is proved similarly.

The same reasoning used in the last part of the proof of Proposition \ref{thm:bicategory} can then be applied in order to prove equation (\ref{eq:triangle}).
\end{remark}
The equivalence of our definition of units with the usual construction, at the level of $GX$, can also be extended to combinations with other representability properties. Closed categories, where a notion of hom exists independently of any monoidal structure, were first introduced by Eilenberg and Kelly in \cite{eilenberg1966closed}; in \cite{manzyuk2012closed}, Manzyuk proved a coherence-via-universality theorem relating them to closed multicategories (with possibly degenerate boundaries). That result can be amended as follows.

\begin{prop}
Let $X$ be a regular polycategory, and suppose that $X$ is 1-representable and tensor 0-representable. Then:
\begin{enumerate}
	\item if $X$ is right, respectively, left closed, then $GX$ admits the structure of a right, respectively, left closed category;
	\item if $X$ is tensor 1-representable, right and left closed, then $GX$ admits the structure of a monoidal biclosed category.
\end{enumerate}
\end{prop}
\begin{proof}
Combine the proof of \cite[Proposition 4.3]{manzyuk2012closed} with the arguments of the proof of Proposition \ref{thm:bicategory}.
\end{proof}

Moreover, we can translate everything to $\coo{X}$, defining notions of \emph{par units} $\bot_x: x \to x$, \emph{par divisible 1-cells}, and \emph{par 0-representability}. The dual of Theorem \ref{thm:0representable}, that the existence of par units is equivalent to the existence of enough par divisible cells, holds, and so do all the ``coherence via universality'' theorems, with pars replacing tensors as the chosen composition of 1-cells, and cohoms replacing homs as part of a closed structure. \index{unit!par} \index{poly-bicategory!0-representable!par} \index{divisible!1-cell!par} \index{poly-bicategory!representable!par} \index{1-representable|see {poly-bicategory}} \index{0-representable|see {poly-bicategory}} \index{par divisible|see {divisible}}

\begin{dfn}
A regular poly-bicategory is \emph{representable} if it is tensor representable and par representable.
\end{dfn} \index{poly-bicategory!representable}

In \cite{cockett2000introduction}, Cockett, Koslowski and Seely introduced the notion of \emph{linear bicategory}: this is a 2-globular set together with two different structures of bicategory, with the same composition and units for 2-cells, but two different compositions and units for 1-cells --- tensor and par --- related by natural transformations 
\begin{align*}
	\delta^L_{a,b,c} : \; & \; a \otimes (b \parr c) \to (a \otimes b) \parr c, \\
	\delta^R_{a,b,c} : \; & \; (a \parr b) \otimes c \to a \parr (b \otimes c),
\end{align*}
called linear distributors, and satisfying certain coherence conditions; see \cite[Definition 2.1]{cockett2000introduction} for details. Linear bicategories are the multi-coloured version of linearly distributive categories, defined in \cite{cockett1997weakly}.

\begin{prop} \label{prop:linearbicat}
Let $X$ be a representable regular poly-bicategory. Then $GX$ admits the structure of a linear bicategory.

In particular, if $X$ has a single 0-cell, $GX$ admits the structure of a linearly distributive category. 
\end{prop}
\begin{proof}
By Proposition \ref{thm:bicategory} and its dual, $GX$ admits two separate structures of bicategory. The construction of linear distributors relating them, and the proof of the relative coherences are the same as in the proof of \cite[Theorem 2.1]{cockett1997weakly}.
\end{proof}

In \cite[Definition 3.1]{cockett2000introduction}, a notion of linear adjunction was defined for 1-cells in a linear bicategory. This was proven to be equivalent to an algebraic definition of adjunction in a poly-bicategory given by \cite[Equation 4]{cockett2003morphisms}; however, we cannot utilise it, because it is based on cells with degenerate boundaries. On the other hand, an equivalent definition based on homs and cohoms still works, relativised to a choice of weak units.
\begin{prop} \label{prop:adjoints} \index{linear adjoint}
Let $X$ be a representable regular poly-bicategory, $f: x \to y$, $g: y \to x$ be 1-cells, and $1_x: x \to x$, $\bot_y: y \to y$ a tensor and a par unit in $X$. Then, $f$ is left linear adjoint to $g$ in $GX$ if any of the following, equivalent conditions holds: 
\begin{enumerate}
	\item there exists a $(r_*)$-divisible cell $\varepsilon: (g, f) \to \bot_y$;
	\item there exists a $(l_*)$-divisible cell $\varepsilon: (g, f) \to \bot_y$;
	\item there exists a $(r^*)$-divisible cell $\eta: 1_x \to (f,g)$;
	\item there exists a $(l^*)$-divisible cell $\eta: 1_x \to (f,g)$.
\end{enumerate}
\end{prop}
\begin{proof}
By Remark \ref{remark:unitor_preserve}, and its dual version, we can fix witnesses $l^\otimes_f$, $r^\otimes_g$, $l^\parr_g$, $r^\parr_f$ of the unitality of $1_x$ and $\bot_y$ that can be ``moved around'' freely, induce natural unitor isomorphisms on $GX$, and satisfy the poly-bicategorical triangle axiom. The proof of \cite[Proposition 1.7]{cockett2003morphisms}, based on a characterisation of ordinary adjunctions by Street and Walters \cite[Proposition 2]{street1978yoneda}, can then be made relative to this choice. For example, starting from an $(r_*)$-divisible $\varepsilon: (g, f) \to \bot_y$, begin with the factorisation
\begin{equation} \label{eq:adjoint1}
\input{img/c3_3_adjoint.tex}
\end{equation}
which yields in $GX$ one of the equations of linear adjunctions. For the other, observe that
\begin{equation*}
\input{img/c3_3_adjoint2.tex}
\end{equation*}
by equation (\ref{eq:adjoint1}) and our choice of $l^\otimes_f$, $r^\otimes_g$, $l^\parr_g$, $r^\parr_f$. Since $\varepsilon$ is $(r_*)$-divisible, it follows that
\begin{equation*}
\input{img/c3_3_adjoint3.tex}
\end{equation*}
which yields in $GX$ the other equation of linear adjunctions.
\end{proof}

If $X$ is left and right closed, and has both tensor and par units, then every 1-cell has both a left and a right linear adjoint. A choice of adjoints induces an equivalence between $X$ and $X^\mathrm{op\,co}$; in particular, $X$ is automatically left and right coclosed, and if it is tensor representable, then it is also par representable.

As discussed in \cite[Section 3]{cockett2000introduction}, a linearly distributive category where every 1-cell has a left and a right linear adjoint is the same as a nonsymmetric $*$-autonomous category, in the sense of \cite{barr1995nonsymmetric}. Therefore, we immediately obtain the following.
\begin{cor}
Let $X$ be a representable regular polycategory. If $X$ is also closed, then $GX$ admits the structure of a nonsymmetric $*$-autonomous category.
\end{cor}

\subsection{Special poly-bicategories. Towards higher dimensions} \label{sec:special}

By Proposition \ref{thm:bicategory}, if our goal were simply to reconstruct the structure of bicategories in the language of regular polygraphs, we could restrict our attention to multi-bicategories, or polygraphs whose globes always have a single output cell; which is what is done by Hermida, or more generally in the opetopic approach. 

On the other hand, one of our \emph{desiderata} concerns string diagrams, and string diagrams are commonly used to reason about non-strict 2-categories or monoidal categories. This is sometimes justified with an appeal to Mac Lane's coherence theorem, that is, by claiming that a strictification is always implicitly assumed. 

A more convincing formalisation, in our opinion, is that string diagrams are always interpreted in an ambient regular polygraph $X$, and can only be used to infer results about weak algebraic structure by passing to representatives in $GX$, whenever these exist. This provides a general heuristic as to what can or cannot be proved diagrammatically in certain contexts.

For example, when reasoning about monoidal categories with string diagrams, the following ``equation'' is sometimes postulated:
\begin{equation} \label{eq:two_wires_one}
\input{img/c3_3_two_wires_one.tex}
\end{equation}
yet this becomes forbidden when reasoning about $*$-autonomous or linearly distributive categories, where one is told that the interpretation of a string diagram with inputs $(a_1, \ldots, a_n)$ and outputs $(b_1, \ldots, b_m)$ is a morphism $a_1 \otimes \ldots \otimes a_n \to b_1 \parr \ldots \parr b_m$. In general, there is not even a cell $a \otimes b \to a \parr b$ to replace the unit in the second diagram.

Of course, the ``equation'' (\ref{eq:two_wires_one}) itself is nonsensical in a regular poly-bicategory, unless either $a$ or $b$ is a weak unit that can be merged with the background (however, it can make sense in a 2-category whose 1-skeleton is not a polygraph). On the other hand, both sides can be interpreted in a regular 2-polygraph or poly-bicategory $X$, and if the poly-bicategory is representable, the following also has an interpretation as a 2-globe in $X$:
\begin{equation} \label{eq:parthentensor}
\input{img/c3_3_parthentensor.tex}
\end{equation}
where the two 2-cells exhibit $a \parr b$, $a \otimes b$ as a par and a tensor of $a, b$, respectively. 

Such nodes commonly appear in proof nets for linear logic and its fragments \cite{girard1996proof}, yet the particular circuit of diagram (\ref{eq:parthentensor}) is invariably forbidden: that is, it does not represent a ``valid'' proof net.

In the literature on proof nets, a lot of care is devoted to the conditions of validity of a proof net, see for example \cite{danos1995proof, blute1996natural}. The non-validity of diagram (\ref{eq:parthentensor}) is easily understood in the context of a representable regular poly-bicategory $X$: it is the image of a 2-globe in $X$, but not of a simple 2-globe, so it does not have a composite, that is, it does not ``equal'' any 2-cell in $GX$ after composition. 

The analysis is complicated by the fact that, in proof nets, the node $(a, b) \to a \otimes b$ is often assumed to have a dual $a \otimes b \to (a, b)$, and the following equation is considered to be valid:
\begin{equation} \label{eq:tensorthentensor}
\input{img/c3_3_tensorthentensor.tex}
\end{equation}
This is still not a composition in a regular poly-bicategory, but at least it is a well-formed equality in a 2-truncated regular polygraph; and, in fact, in monoidal categories, ``string diagrams with loops'' such as (\ref{eq:parthentensor}) are always assumed to be well-formed, that is, to have a composite.

Assuming that equation (\ref{eq:tensorthentensor}) is, in fact, the formally correct version of equation (\ref{eq:two_wires_one}), in order for it to be a theorem about regular polygraphs $X$ such that $GX$ admits the structure of a monoidal category, we need $X$ to have
\begin{enumerate}
	\item universal 2-cells $a \otimes b \to (a, b)$, and
	\item composites of ``string diagrams with loops''.
\end{enumerate}
The first we can already have in a representable regular poly-bicategory, in the case that $a \otimes b \simeq a \parr b$ for all $a, b$. For the second, however, we have to extend the algebraic composition from simple to generic 2-globes.

\begin{dfn} \index{poly-bicategory!special}
A \emph{special poly-$n$-category} is a regular $n$-polygraph $X$, together with
\begin{enumerate}
	\item for all $n$-globes $G$ with two $n$-cells, reducing to the atomic $n$-globe $G'$ with $\bord{}{}G' = \bord{}{}G$ after a sequence of simple mergers, a composition operation 
	\begin{equation*}
		\cutt{G}: \mathrm{Hom}(G,X) \to X(G'), 
	\end{equation*}
	restricting to the identity on $\mathrm{Hom}(\bord{}{\alpha}G, X) = \mathrm{Hom}(\bord{}{\alpha}G', X)$, and such that
	\item for all $n$-globes $G_{1,2,3}$ with three $n$-cells $x_1, x_2, x_3$, reducing to the atomic $n$-globe $G'$ with $\bord{}{}G_{1,2,3} = \bord{}{}G'$ in two stages, through either a sequence of simple mergers involving $U_{x_1}$ and $U_{x_2}$, or a sequence of simple mergers involving $U_{x_2}$ and $U_{x_3}$, the following diagram commutes:
\begin{equation*}
\begin{tikzpicture}[baseline={([yshift=-.5ex]current bounding box.center)}]
	\node[scale=1.25] (0) at (-2,2) {$\mathrm{Hom}(G_{1,2,3},X)$};
	\node[scale=1.25] (1) at (2,2) {$\mathrm{Hom}(G_{12,3},X)$};
	\node[scale=1.25] (2) at (-2,0) {$\mathrm{Hom}(G_{1,23},X)$};
	\node[scale=1.25] (3) at (2,0) {$X(G').$};
	\draw[1c] (0) to node[auto] {$\cutt{G_{1,2}}$} (1);
	\draw[1c] (2) to node[auto,swap] {$\cutt{G_{1,23}}$} (3);
	\draw[1c] (0) to node[auto,swap] {$\cutt{G_{2,3}}$} (2);
	\draw[1c] (1) to node[auto] {$\cutt{G_{12,3}}$} (3);
\end{tikzpicture}
\end{equation*}
\end{enumerate}
A morphism $f: X \to Y$ of special poly-$n$-categories is a map of the underlying regular polygraphs such that
\begin{equation*}
\begin{tikzpicture}[baseline={([yshift=-.5ex]current bounding box.center)}]
	\node[scale=1.25] (0) at (-1,1.5) {$\mathrm{Hom}(G,X)$};
	\node[scale=1.25] (1) at (2,1.5) {$X(G')$};
	\node[scale=1.25] (2) at (-1,0) {$\mathrm{Hom}(G,Y)$};
	\node[scale=1.25] (3) at (2,0) {$Y(G').$};
	\draw[1c] (0) to node[auto] {$\cutt{G}$} (1);
	\draw[1c] (2) to node[auto,swap] {$\cutt{G}$} (3);
	\draw[1c] (0) to node[auto,swap] {$\mathrm{Hom}(G,f)$} (2);
	\draw[1c] (1) to node[auto] {$f_{G'}$} (3);
\end{tikzpicture}
\end{equation*}
commutes for all suitable $G \leadsto G'$.
\end{dfn}
The terminology, here, is not meant to imply that special poly-$n$-categories are not regular as polygraphs, but rather be reminiscent of the distinction between Frobenius algebras and special Frobenius algebras: special poly-bicategories admit compositions ``with loops'' as in diagram (\ref{eq:spider_special}).

While in a regular poly-bicategory, only a single-input, single-output 2-cell could be both $(c^*)$ and $(c_*)$-divisible, in a special poly-bicategory this makes sense for arbitrary 2-cells.
\begin{dfn} \index{divisible!2-cell}
Let $p: (\Gamma) \to (\Delta)$ be a 2-cell in a special poly-bicategory $X$. We say that $p$ is \emph{divisible} if it is $(c^*)$-divisible and $(c_*)$-divisible.
\end{dfn}
\begin{remark}
Because binary tensors and pars can be composed to construct $n$-ary ones, $n > 2$, it suffices to formulate the 1-representability condition for 1-globes $\Gamma$ with one or two 1-cells.
\end{remark}

The proofs of Proposition \ref{prop:2units} and Corollary \ref{cor:onedivisible} can then be adapted to show the following.
\begin{prop} \label{prop:represpecial}
Let $X$ be a special poly-bicategory. Then, the following are equivalent:
\begin{enumerate}
	\item for all 1-globes $\Gamma$ of $X$, there exist a 1-cell $\overline{a}$ and a divisible 2-cell $p: (\Gamma) \to \overline{a}$;
	\item for all 1-globes $\Gamma$ of $X$, there exist a 1-cell $\underline{a}$ and a divisible 2-cell $p': \underline{a} \to (\Gamma)$;
	\item for all 1-globes $\Gamma$ of $X$, there exists a (necessarily unique) unit $id_\Gamma$ on $\Gamma$.
\end{enumerate}
\end{prop}

\begin{dfn}
We say that a special poly-bicategory $X$ is \emph{1-representable} if it satisfies any of the equivalent conditions of Proposition \ref{prop:represpecial}. \index{poly-bicategory!1-representable}
\end{dfn}

\begin{cor}
Let $X$ be a 1-representable special poly-bicategory, and $p: (\Gamma) \to (\Delta)$ a 2-cell of $X$. Then, $p$ is $(c^*)$-divisible if and only if it is $(c_*)$-divisible.
\end{cor}

In a 1-representable special poly-bicategory, tensors and pars of 1-cells collapse to the same notion: a $(c^*)$-divisible $t_{a,b}: (a, b) \to a \otimes b$ has a $(c_*)$-divisible inverse $\invrs{t}_{a,b}: a \otimes b \to (a,b)$. Moreover, if $X$ is also tensor and par 0-representable, we can show that tensor and par units coincide in $X$, and so do tensor and par divisibility.

\begin{dfn} \index{poly-bicategory!representable}
A special poly-bicategory is \emph{representable} if it is 1-representable and tensor and par 0-representable.
\end{dfn}

\begin{prop}
Let $X$ be a representable special poly-bicategory. Then, the tensor units in $X$ are also par units, and if $\{l_a, r_a\}$ are 2-cells exhibiting $\{1_x\}$ as tensor units, $\{\invrs{l}_a, \invrs{r}_a\}$ exhibit them as par units.
\end{prop}
\begin{proof}
Given a tensor unit $1_x$ and a par unit $\bot_x$ on $x$, since tensors and pars coincide in $X$, we have $1_x \simeq 1_x \parr \bot_x = 1_x \otimes \bot_x \simeq \bot_x$. Then, we can precompose a family of 2-cells $\{l^\parr_a, r^\parr_a\}$ exhibiting $\bot_x$ as a par unit with isomorphisms $1_x \to \bot_x$, and apply Lemma \ref{lem:twodivisible} to these 2-cells and to the $\{\invrs{l}_a, \invrs{r}_a\}$, which are known to be $(c_*)$-divisible, to conclude that they are also $(r^*)$ or $(l^*)$-divisible.
\end{proof}
Of course, the dual theorem, with the roles of tensor and par units switched, also holds. Therefore, we can speak simply of \emph{units} in a representable special poly-bicategory.
\begin{cor}
Let $X$ be a representable special poly-bicategory. Then a 1-cell in $X$ is tensor divisible if and only if it is par divisible.
\end{cor}
\begin{proof}
Let $e: x \to y$ be a tensor divisible 1-cell. By the degenerate case of Proposition \ref{prop:adjoints}, dividing a unit $1_x$ or $1_y$ by $e$ produces a cell $\bar{e}: y \to x$ which is part of an adjoint equivalence with $e$ (relative to a coherent choice of units). By standard reasoning, it can then be shown that, for all $a: x \to z$, the 1-cell $\bar{e} \otimes a$ is both a right hom and a right cohom from $e$ to $a$, such that $e \otimes (\bar{e} \otimes a) \simeq a$; similarly, for all $b: z \to y$, the 1-cell $b \otimes \bar{e}$ is both a left hom and a left cohom from $e$ to $b$, such that $b \simeq (b \otimes \bar{e}) \otimes e$.
\end{proof}

It follows that, in the definition of a representable special poly-bicategory, we can demand, without loss of generality, the existence of 1-cells that are simultaneously tensor and par divisible.

\begin{remark}
It is possible that some divisibility properties are collapsed when $X$ is 1-representable even in the absence of 0-representability; we have not attempted a complete classification yet.
\end{remark}

\begin{prop}
Let $X$ be a representable special poly-bicategory. Then $GX$ admits the structure of a bicategory.
\end{prop}
\begin{proof}
Follows from the previous considerations, as a degenerate case of Proposition \ref{prop:linearbicat} where tensor and par coincide.
\end{proof}

\begin{remark}
In the light of this result, we can attempt a broad classification of the main uses of string diagrams in reasoning about weak 2-dimensional structures:
\begin{itemize}
	\item string diagrams with loops, and with two-sided tensor nodes $(a, b) \to a\otimes b$ and $a\otimes b \to (a,b)$, inverse to each other, can be interpreted in all representable special poly-bicategories $X$ to infer results about the bicategory $GX$;
	\item string diagrams without loops, and with one-sided tensor and par nodes $(a, b) \to a \otimes b$ and $a \parr b \to (a,b)$ can be interpreted in all representable regular poly-bicategories $X$ to infer results about the linear bicategory $GX$.
\end{itemize}
Proof nets with two-sided tensor and par nodes, as in \cite{blute1996natural}, seem to belong in a sort of hybrid: a 2-polygraph with separate tensor and par sub-polygraphs, separately forming representable special poly-bicategories; and a ``mixed'' sub-polygraph, forming a representable regular poly-bicategory. What they tell about $GX$, if anything, will depend on the particular mixture of partial composition and partial representability. 

If all compositions can be subsumed by representability properties, regular polygraphs remain as a potential unifying framework for all such diagrammatic languages.
\end{remark}

Something that we have not addressed so far is functoriality of the constructions. Consider first the case of regular poly-bicategories. If $f: X \to Y$ is a morphism of tensor or par representable regular poly-bicategories, its restriction $f: GX \to GY$ is in general only lax for tensor compositions, and colax for par compositions: this is because the images $f(t_{a,b}): (fa, fb) \to f(a \otimes b)$ of $(c^*)$-representable 2-cells of $X$ will factor through $t_{fa,fb}: (fa, fb) \to fa \otimes fb$ in $Y$, and dually for $(c_*)$-representable 2-cells. On the other hand, nothing about units, as we defined them, is necessarily preserved.

While this makes the connection with bicategories somewhat less strong, it is worth considering that there might be something deeper to the mismatch. One of the fundamental slogans of higher category theory is that (weak) equivalence, rather then identity, witnessed by unit cells, is the ``right'' notion of equality, one that preserves all the relevant structure, and nothing more. The need for making this intuition precise has been one of the main drives of the Univalent Foundations programme \cite{hottbook}. 

Our definition of divisible 1-cell can be seen as an elementary definition of equivalence, which does not depend \emph{a priori} on any notion of identity. If we accept that representability is the fundamental concept connecting the combinatorics of regular polygraphs with the weak coherent structure of higher categories, we are led to require that enough equivalences, rather than units, exist; and yet, by Theorem \ref{thm:0representable} and its dual, the two conditions coincide. From this perspective, units have no special status: they just happen to exist among other equivalences.  

In ordinary (higher) category theory, preservation of isomorphisms or equivalences by functors is usually derived from the requirement that they preserve, at least weakly, units and composition. By Corollary \ref{cor:onedivisible}, we can treat the different divisibility properties of 2-cells as a variety of generalisations of the one notion of isomorphism available for single input, single output cells --- generalisations that cannot be captured by the ``existence of inverses''. 

The most natural thing to do, perhaps, is to ask that morphisms preserve all divisibility properties. \index{equivalence!in regular polygraphs}
\begin{dfn}
Let $f: X \to Y$ be a morphism of regular or special poly-bicategories. We say that $f$ is \emph{very strong} if it preserves all divisibility properties of 1-cells and 2-cells.
\end{dfn}

\begin{prop} \label{prop:verystrong}
Let $f: X \to Y$ be a very strong morphism between tensor representable regular poly-bicategories. Then its restriction $f: GX \to GY$ is a pseudofunctor of bicategories.
\end{prop}
\begin{proof}
Since $f$ preserves $(c^*)$-divisibility, for all composable 1-cells $a, b$ of $X$, clearly $f(a \otimes b) \simeq fa \otimes fb$ in $Y$. Let $1_x$ be a tensor unit in $X$. Then $f(1_x)$ is tensor divisible in $Y$; moreover, 
\begin{equation*}
	f(1_x) \simeq f(\rimp{1_x}{1_x}) \simeq \rimp{f(1_x)}{f(1_x)},
\end{equation*}
because $f$ preserves $(r_*)$-divisibility. From the proof of Theorem \ref{thm:0representable}, we know that, in a tensor 0-representable regular poly-bicategory, $\rimp{e}{e}$ is a tensor unit for all tensor divisible $e$. It follows that $f(1_x) \simeq 1_{fx}$. Coherence comes automatically from universality.

The proof that units for 2-composition are strictly preserved is similar, and simpler; by definition, 2-composition is preserved on the nose.
\end{proof}
The dual for par representable regular poly-bicategories also holds, and everything applies to special poly-bicategories as well. The problem is that $f$ is not just any pseudofunctor: it is one that also preserves all homs and cohoms, or absolute Kan extensions and lifts in bicategorical terminology.

At least in the case of special poly-bicategories, then, we might consider only requiring the preservation of divisible 1-cells and divisible 2-cells. This turns out to be sufficient for the preservation of units, by the following argument.

\begin{dfn}
Let $f: X \to Y$ be a morphism of special poly-bicategories. We say that $f$ is \emph{strong} if it preserves divisible 1-cells and 2-cells.
\end{dfn}

\begin{prop}
Let $f: X \to Y$ be a strong morphism between representable special poly-bicategories. Then its restriction $f: GX \to GY$ is a pseudofunctor of bicategories.
\end{prop}
\begin{proof}
Preservation of composition is proved as in Proposition \ref{prop:verystrong}. If $1_x: x \to x$ is a unit in $X$, then $f(1_x): fx \to fx$ is divisible in $Y$, so there exists some $e: fx \to fx$ such that $e \otimes f(1_x) \simeq 1_{fx}$. Then, because $1_x \simeq 1_x \otimes 1_x$, we obtain the following chain of isomorphisms:
\begin{equation*}
	1_{fx} \simeq e \otimes f(1_x) \simeq e \otimes (f(1_x) \otimes f(1_x)) \simeq 1_{fx} \otimes f(1_x) \simeq f(1_x).
\end{equation*}
The proof that units for 2-composition are strictly preserved is similar.
\end{proof}

To conclude, we outline a plan for the extension of these constructions to higher dimensions. In the 2-dimensional case, we opted to focus on regular poly-bicategories first, because of the elegant shape that the combinatorics of divisibility assume: the six types of 2-cell divisibility are all on an equal footing, each corresponding to a different connective of non-commutative linear logic. 

In this sense, special poly-bicategories look like a poorer, degenerate relative of regular poly-bicategories; however, it is the latter that we consider to be a stepping stone towards higher-dimensional directed spaces. This is due mainly to the fact that $n$-globes have better closure properties than simple $n$-globes: a simple $n$-globe may not have simple boundaries, and even if $X$ and $Y$ are simple globes ``all the way down'', $X \otimes Y$ may not have this property; see for example the cylinder of diagram (\ref{eq:multiply_cylinder}), whose output boundary is not simple, even though it is the tensor product of simple globes with simple boundaries. The fact that we cannot truncate a cylinder on a regular poly-bicategory, and obtain another regular poly-bicategory, is also the reason why their morphisms admit only quite degenerate forms of natural transformations; see \cite{cockett2003morphisms} for more details on this subject.

So, moving on to special poly-3-categories, we would
\begin{enumerate}
	\item reformulate 2-composition as a representability property: for all 2-globes $G$ in $X$, there exists a 2-cell $p$ and 3-cells $G \to p$, $p \to G$, divisible for the algebraic 3-composition;
	\item reformulate divisibility for 2-cells relative to the weak composition: instead of requiring uniqueness of the result of division, ask for a higher-dimensional analogue of $(r^*)$-divisibility and its duals, as appropriate, for the 3-cells exhibiting composites.
\end{enumerate}
Tricategories are the highest-dimensional example of weak higher category for which a workable, complete algebraic definition exist, so ``coherence via universality'' could be worked out explicitly in this case. We would expect, in particular, unit 1-cells to subsume the notion of Joyal-Kock weak unit in this case. One dimension higher, it would be interesting to compare the result with the quasistrict 4-categories of \cite{bar2016data}, an algebraic definition which, like our approach, is strongly influenced by string-diagrammatic algebra.

The possibility of an effective generalisation relies on a better combinatorial understanding of higher globes, and associated types of divisibility. If that were in place, however, nothing would prevent us from generalising inductively to arbitrary dimensions. Moreover, since divisibility of lower-dimensional cells is only defined by the \emph{existence} of higher-dimensional divisible cells, the possibility of a coinductive definition remains open; this is barred in the opetopic definition, as given in \cite[Chapter 4]{cheng2004higher}, by the presence of a universal quantification over higher universal cells. Ultimately, we would like to be able to state something of the form:
\begin{center}
	\emph{If $X$ is a representable regular polygraph, $GX$ admits the structure of a weak $\omega$-category}.
\end{center} \index{polygraph!representable}

There is, then, the question of compositionality. Even when $X$ and $Y$ are representable, we cannot expect $X \otimes Y$ to be representable --- but we can hope to be able to universally complete it. Here, the paragon is \cite{hermida2000representable}, where a pseudomonad $T$ on the bicategory of multicategories is constructed, such that $TM$ is representable for all objects $M$, and in a certain sense equivalent to $M$ whenever $M$ is already representable. 

This completion operation is then used in an abstract proof of the coherence theorem for monoidal categories. There is no reason why the basic structure of this proof could not be generalised, yet coherence results for higher categories are notoriously elusive: Simpson's conjecture, that higher homotopy types are equivalent to homotopy types with everything strict but the units, has only been settled in dimension 3 \cite{simpson2009homotopy, joyal2007weak}.

Our hope is that the relevant traits of those proofs that do exist will be better understood, and extrapolated for further use, by keeping higher-dimensional spaces grounded in combinatorics and in the intuition of string diagrams.

With this, we conclude the first part of the thesis, the one strictly devoted to the ``geometry of composition'', and move on to a more specific application of diagrammatic algebra.

  \thispagestyle{empty} 
\chapter{Categorical quantum theory} \label{chap:quantum}
\thispagestyle{plain}

\noindent\emph{In this chapter:}
\begin{itemize}
	\item[$\triangleright$] We define some of the structures relevant to the formulation of quantum theory within category theory. We discuss some issues relative to the ``category of Hilbert spaces'', and prove a new reconstruction result, guided by categorical universal algebra and the principle of representability, which exhibits the category of Hilbert spaces and short linear maps as the truncation of a bicategory of cospans of isometries. --- \emph{Section \ref{sec:abelian}}
	\item[$\triangleright$] Frobenius algebras were introduced to categorical quantum mechanics as an algebraic counterpart to orthonormal bases. We review how this led to ZX calculi, and the search for complete diagrammatic axiomatisations of fragments of quantum theory.  --- \emph{Section \ref{sec:frobenius}}
	\item[$\triangleright$] We show how the classification of Frobenius algebras on Hilbert spaces, combined with the experience of string-diagrammatic calculi for quantum theory, led to an original approach to the operational classification of entangled quantum states, and we outline a plan for its development. --- \emph{Section \ref{sec:entanglement}}
\end{itemize}

\section{From abelian groups to Hilbert spaces} \label{sec:abelian}
This chapter introduces and elaborates on some of the ideas of categorical quantum mechanics, a research programme initiated by Abramsky and Coecke in \cite{abramsky2004categorical, abramsky2005abstract}. The premise of categorical quantum mechanics is that, with the rise of quantum information and computation theory \cite{nielsen2009quantum}, it is sensible to study quantum theory as an abstract model of computation, or theory of processes, with a neutral stance on its physical meaning or implementation; in other words, look at the abstract structure of quantum algorithms, rather than the physical supports which may run them.

In the light of the Curry-Howard-Lambek correspondence between programs and morphisms of a category \cite{sorensen2006}, a model of computation can be seen as equivalent to the class of categories that support an interpretation of its characteristic processes, in the sense of categorical universal algebra. The original example is the correspondence between the typed $\lambda$-calculus and cartesian closed categories; see for example \cite[Chapter 5]{barendregt1984lambda}.

For quantum computation, a process-theoretic understanding seems particularly desirable, as it is still unclear what processes deserve to be called ``fundamental'', or will be the basic tools of future quantum programmers. In fact, many well-known results on quantum theory as a model of computation come in the negative form of ``no-go'' theorems: it is \emph{not} possible, for instance, to copy arbitrary states \cite{dieks1982communication, wootters1982single, barnum1996noncommuting}, nor to delete them \cite{pati2000impossibility}; the former has been reinterpreted by Abramsky as a general result about monoidal categories \cite{abramsky2009no}. Both operations are, of course, staples of classical computation.

On the constructive side, we know a number of quantum algorithms that achieve a speedup with respect to their optimal classical counterparts, yet we still lack a precise understanding of the class of problems for which such a speedup is achievable \cite{aaronson2014need}, or of the irreducible ``quantum core'' of these algorithms, which separates them from classical ones. 

This failure has been attributed to the shortcomings of the traditional formalism for quantum theory, based on Hilbert spaces, which has remained mostly invaried since it was codified by von Neumann \cite{neumann2013mathematische} --- even though its inventor himself expressed doubts, and started the field of quantum logic as an alternative \cite{birkhoff1958neumann}. The verification of quantum algorithms by matrix or vector calculations has been compared to the verification of computer programs by the examination of machine code \cite{coecke2010picturalism}: an extremely ``low-level'' language, in programming jargon.

If the study of quantum computation is to be grounded in categorical universal algebra, string diagrams are an obvious candidate for a higher-level language, and indeed they have been core to categorical quantum mechanics since its inception. The recent \emph{summa} of the subject so far \cite{coecke2017picturing}, written by Coecke and Kissinger, develops the theory starting from abstract diagrammatic algebra and process theory, rather than category theory.

Diagrammatic algebra has helped revealing structural similarities between quantum algorithms \cite{vicary2013topological, zeng2015abstract}, and the equivalence of certain protocols in terms of the abstract operations needed to implement them \cite{vicary2012higher, reutter2017shaded}. Moreover, it has naturally led to the problem of axiomatising fragments of quantum theory as independent higher algebraic theories: once we have packaged such a fragment into a theory, we can study it and tweak it on its own, forgetting the context from which it came, and gaining new perspective from its interpretation in different contexts. Examples include the theory of biunitaries, which underlies several constructions in quantum information theory \cite{reutter2016biunitary}, but also classical cryptography with ``one time pads'' \cite{stay2013bicategorical}; the various theories called ``ZX calculi'', axiomatisations of fragments of the theory of qubits, including stabiliser quantum mechanics \cite{backens2014zx}; and our own ZW calculus \cite{hadzihasanovic2015diagrammatic}, the subject of Chapter \ref{chap:zwcalculus}, towards which we will be working in this chapter.

The ``vanilla'' version of the ZW calculus axiomatises a full subcategory of the category $\cat{Ab}$ of abelian groups and homomorphisms; its extensions are not significantly different. Our first aim is to show precisely why this tells anything about Hilbert spaces. In the way, we will try to clarify some subtleties about ``the category of Hilbert spaces'', some of which are folklore, yet rarely discussed in print. 

We will not, however, devote much space to physical interpretations, for which any introductory textbook on quantum theory can be consulted, including the ``pictorial'' \cite{coecke2017picturing}; nor will we discuss cognates of quantum theory with regard to abstract categorical properties rather than concrete algebraic models, for which good starting points are Heunen's PhD thesis \cite{heunen2009categorical}, and the lecture notes by Heunen and Vicary for the Oxford course on categorical quantum mechanics \cite{heunen2017lectures}. Although we will briefly recall the basic definitions, a previous knowledge of the theory of Hilbert spaces will be helpful; \cite{conway2007course} is a possible source.

In line with the two previous chapters, we adopt two guiding principles in defining (higher) categories:
\begin{enumerate}
	\item whenever objects are algebras of a theory, the morphisms should strictly preserve all the algebraic structure;
	\item natural notions of composition come from representability properties.
\end{enumerate}
With regard to the second point, we will see that both the important compositions of Hilbert spaces --- tensor products and direct sums --- can be handled conceptually at the level of abelian groups and homomorphisms, with the constructions carrying through at each successive addition of structure.

\begin{remark}
In what follows, we need to consider models of algebraic theories, which we defined as PROs using the $\omega$-categorical formalism in Chapter \ref{chap:interacting}, in regular and special multicategories and polycategories, for which we used a different formalism in Chapter \ref{chap:directed}. The two are quite easily intertranslatable, and we will not be pedantic about the distinction. The only non-trivial adjustment that we need to make is that, whenever we need to interpret a PRO containing 0-ary operations $f: [0] \to [n]$ or $[n] \to [0]$ in a regular or special poly-bicategory, we assume that $[0]$ has been replaced with a new, non-degenerate 1-cell, and that cell is mapped to a weak unit.
\end{remark}

\begin{cons} \index{abelian group} \index{Ab@$\cat{Ab}, \textit{Ab}_\otimes, \textit{Ab}_\oplus$}
Let $\textit{Ab}_\otimes$ be the regular multicategory whose 1-cells are abelian groups, and 2-cells $f: (G_1, \ldots, G_n) \to G$ are multilinear maps, that is, functions $f: G_1 \times \ldots \times G_n \to G$ such that, for all $i = 1, \ldots, n$, and $g_i, h_i \in G_i$,
\begin{equation*}
	f(g_1, \ldots, g_i + h_i, \ldots, g_n) = f(g_1, \ldots, g_i, \ldots,  g_n) + f(g_1, \ldots, h_i, \ldots, g_n).
\end{equation*}
Given maps $f': (G_{i1}, \ldots, G_{in}) \to G_i$ and $f: (G_1, \ldots, G_m) \to G$, for some $1 \leq i \leq m$, the composite $\mathrm{cut}_{1,i}(f',f)$ is the multilinear map 
\begin{equation*}
	(g_1, \ldots, g_{i-1}, g_{i1}, \ldots, g_{in}, g_{i+1}, \ldots, g_m) \mapsto f(g_1, \ldots, g_{i-1}, f'(g_{i1}, \ldots, g_{in}), g_{i+1}, \ldots, g_m).
\end{equation*}

Then, $\textit{Ab}_\otimes$ is tensor representable: the tensor unit is the group $\mathbb{Z}$ of integers, with witnessing cells $l_G: (\mathbb{Z}, G) \to G$ and $r_G: (G, \mathbb{Z}) \to G$ uniquely determined by \index{tensor product!of abelian groups}
\begin{equation*}
	l_G(1, g) = r_G(g, 1) = g,
\end{equation*}
and the tensor of two 1-cells $G, H$ is the tensor product $G \otimes H$, that is, the abelian group generated under finite sums by $\{g \otimes h \,|\, g \in G, h \in H\}$ with relations $(g + g') \otimes h \sim g \otimes h + g' \otimes h$, and $g \otimes (h + h') \sim g \otimes h + g \otimes h'$, together with the witness $(G, H) \to G \otimes H$ defined by $(g, h) \mapsto g \otimes h$. 
\end{cons}

\begin{cons}
There is another multicategory with the same 1-cells, inducing a different monoidal structure. Let $\textit{Ab}_\oplus$ have the same 1-cells as $\textit{Ab}_\otimes$, while 2-cells $f: (G_1, \ldots, G_n) \to G$ are multi-cospans of group homomorphisms, that is, families $\{f_i: G_i \to G\}_{i=1}^n$ of homomorphisms with the same codomain. Given multi-cospans $f': (G_{i1}, \ldots, G_{in}) \to G_i$ and $f: (G_1, \ldots, G_m) \to G$, where $1 \leq i \leq m$, the composite $\mathrm{cut}_{1,i}(f',f)$ is the multi-cospan obtained as the disjoint union
\begin{equation*}
	\{f_j: G_j \to G\}_{j=1}^{i-1} + \{f_i f'_j: G_{ij} \to G\}_{j=1}^n + \{f_j: G_j \to G\}_{j=i+1}^{m}.
\end{equation*}

Then, $\textit{Ab}_\oplus$ is also tensor representable: the tensor unit is the trivial group $0$, and the tensor of two 1-cells $G, H$ is the direct sum $G \oplus H$, generated under finite sums by $\{g \oplus h \,|\, g \in G, h \in H\}$ with relations $(g + g') \oplus (h + h') \sim g \oplus h + g' \oplus h'$, together with the witness $(G, H) \to G \oplus H$, given by the cospan $g \mapsto g \oplus 0$, $h \mapsto 0 \oplus h$. \index{direct sum!of abelian groups}

Notice that 2-cells with a single input are just group homomorphisms, both in $\textit{Ab}_\otimes$ and $\textit{Ab}_\oplus$.
\end{cons}

\begin{dfn}
We write $\cat{Ab}$ for the category $G(\textit{Ab}_\otimes) = G(\textit{Ab}_\oplus)$, with the two monoidal structures induced by tensor representability of $\textit{Ab}_\otimes$ and $\textit{Ab}_\oplus$.
\end{dfn}

\begin{remark}
There are canonical isomorphisms $G \otimes (H \oplus K) \simeq (G \otimes H) \oplus (G \otimes K)$, making $\cat{Ab}$ a distributive monoidal category.
 
Both the monoidal structures are, in fact, symmetric, with the obvious swap morphisms defined by $g \otimes h \mapsto h \otimes g$, and by $g \oplus h \mapsto h \oplus g$ on the generators. This can be traced back to $\textit{Ab}_\otimes$ and $\textit{Ab}_\oplus$ being \emph{symmetric} multicategories, in the sense that, for all $n$, the permutation group $S_n$ acts on the set of 2-cells with $n$ inputs. With a slight abuse of notation (that can be made formally precise with a slight effort), we depict permutations of the inputs of a 2-cell $f$ with string diagrams like
\begin{equation*}
\input{img/c4_1_permutation.tex}
\end{equation*} 
This also allows us to make sense of ``commutative monoids in $\textit{Ab}_\otimes$'', that is, maps from the PROP $\textit{CMon}$ into $\textit{Ab}_\otimes$; these correspond precisely to commutative rings.
\end{remark}

\begin{remark} By defining, for all sequences $G_1, \ldots, G_n$ of abelian groups, the $n$-ary tensor product $G_1 \otimes \ldots \otimes G_n$ to be the abelian group generated under finite sums by $\{g_1 \otimes \ldots \otimes g_n \,|\, g_i \in G_i, i \leq n\}$ with multilinearity relations, we can construct a special polycategory $\textit{Ab}$ as follows: a 2-cell $f: (G_1, \ldots, G_n) \to (H_1, \ldots, H_m)$ of $\textit{Ab}$ is just a group homomorphism $f: G_1 \otimes \ldots \otimes G_n \to H_1 \otimes \ldots \otimes H_m$. 

As an example of composition, let $f: G \to (H, K)$ and $f': (K, L) \to M$ be 2-cells; for all $g \in G$, $k \in K$, $l \in L$, there are expressions
\begin{equation*}
	f(g) = \sum_i h_i \otimes k_i, \quad f'(k \otimes l) = \sum_j m_j.
\end{equation*}
Then, the composite of $f$ and $f'$ along $K$ is the 2-cell $(G, L) \to (H, M)$ defined by
\begin{equation*}
	g \otimes l \mapsto f(g) \otimes l = \sum_i h_i \otimes k_i \otimes l \mapsto \sum_i h_i \otimes f'(k_i \otimes l) = \sum_{i,j} h_i \otimes m_{ij};
\end{equation*}
this is well-defined by the definition of tensor products by generators and relations. $\textit{Ab}$ is representable by construction, with witnesses of compositions $(A, B) \to A \otimes B$ and $A \otimes B \to (A,B)$ given by identities in $\cat{Ab}$, and $G(\textit{Ab})$ is clearly just $\cat{Ab}$ again, with a strictified monoidal structure. 

The special polycategory $\textit{Ab}$ can be seen as a strictly associative version of the monoidal category $(\cat{Ab}, \otimes, \mathbb{Z})$. As a consequence of Mac Lane's coherence theorem, such a construction exists for any monoidal category, which is what allows us to interpret string diagrams with multiple inputs and outputs ``in $\cat{Ab}$'', and in all subsequent examples of monoidal categories.
\end{remark}

\begin{cons} \index{module} \index{RMod@$R\cat{Mod},\textit{RMod}_\otimes,\textit{RMod}_\oplus$} \index{Rmodule@$R$-module|see {module}}
Given a commutative ring $R$, that is, a commutative monoid in $\textit{Ab}_\otimes$, we can form a regular multicategory $\textit{RMod}_\otimes$, whose 1-cells are left $R$-modules, that is, abelian groups $G$ together with a left action of $R$, and 2-cells $f: (G_1, \ldots, G_n) \to G$ are multilinear maps of the underlying abelian groups, such that, for all $r \in R$, all $i = 1,\ldots,n$, and $g_i \in G_i$,
\begin{equation*}
	f(g_1, \ldots, r \cdot g_i, \ldots, g_n) = r \cdot f(g_1, \ldots, g_i, \ldots g_n),
\end{equation*}
or, graphically in $\textit{Ab}_\otimes$,
\begin{equation*}
\input{img/c4_1_equiaction.tex}
\end{equation*} 
Then $\textit{RMod}_\otimes$ is tensor representable: the tensor unit is the underlying group of $R$, with its action by multiplication on itself, and the tensor of two 1-cells $G, H$ is the quotient $G \otimes_R H$ of $G \otimes H$ by the relation $(r\cdot g) \otimes h \sim g \otimes (r \cdot h)$, with $r \cdot (g \otimes h)$ defined to be the result of the identification of $(r\cdot g) \otimes h$ and $g \otimes (r \cdot h)$. \index{tensor product!of $R$-modules}
\end{cons}

\begin{cons}
Similarly to the case of abelian groups, we can also define a regular multicategory $\textit{RMod}_\oplus$, with the same 1-cells, whose 2-cells $f: (G_1, \ldots, G_n) \to G$ are multi-cospans $\{f_i: G_i \to G\}_{i=1}^n$ of group homomorphisms such that $f_i(r\cdot g) = r \cdot f_i(g)$. This is also tensor representable: the tensor unit is the trivial group with the trivial action, and the tensor of two 1-cells $G, H$ is the direct sum $G \oplus H$, with the action of $R$ defined by $r\cdot(g \oplus h) := (r\cdot g) \oplus (r \cdot h)$. \index{direct sum!of $R$-modules}
\end{cons}

\begin{remark}
The construction of $\textit{RMod}_\otimes$ from $\textit{Ab}$ is an instance of the more general ``bimodule construction'', which can be used to define a multi-bicategory from any multicategory, see for example \cite{cruttwell2010unified}: it corresponds to the special case where a left $R$-module is considered to be an $R$-$R$-bimodule with the right action
\begin{equation*}
\input{img/c4_1_right_action.tex}
\end{equation*} 
well-defined due to the commutativity of $R$.
\end{remark}

\begin{dfn}
We write $R\cat{Mod}$ for the category $G(\textit{RMod}_\otimes) = G(\textit{RMod}_\oplus)$, with the two monoidal structures induced by tensor representability of $\textit{RMod}_\otimes$ and $\textit{RMod}_\oplus$, making it a distributive monoidal category.

When $R$ is a field $k$, $R$-modules are called $k$-vector spaces, and we write $\cat{Vec}_k$ instead of $k\cat{Mod}$. We also write $\cat{FVec}_k$ for the full subcategory of $\cat{Vec}_k$ on finite-dimensional vector spaces.
\end{dfn} \index{vector space} \index{Vec@$\cat{Vec}_k, \cat{FVec}_k$}

\begin{cons} For all endomorphisms $a: R \to R$ of a ring $R$, there is a morphism of regular multicategories $a^*: \textit{RMod}_\otimes \to \textit{RMod}_\otimes$, which maps the group $G$ with the action $(r,g) \mapsto r \cdot g$ to the same group with the action $(r, g) \mapsto a(r) \cdot g$, without affecting 2-cells. This morphism preserves tensors, but not always tensor units, and induces a ``semi-monoidal'' endofunctor on $R\cat{Mod}$.

In particular, when $R$ is the field $\mathbb{C}$ of complex numbers, there is an involutive endofunctor $\overline{(-)}$ on $\cat{Vec}_\mathbb{C}$, induced by complex conjugation: given a complex vector space $V$, the vector space $\overline{V}$ has the same underlying abelian group, with the conjugate action $(\lambda, v) \mapsto \overline{\lambda}\cdot v$ on vectors $v \in V$.
\end{cons}

\begin{dfn}
A \emph{pre-Hilbert space} is a complex vector space $H$ together with an inner product, that is, a bilinear map $(\overline{H}, H) \to \mathbb{C}$, denoted by $(v,w) \mapsto \braket{v}{w}$, such that
\begin{enumerate}
	\item $\braket{v}{w} = \overline{\braket{w}{v}}$, implying $\braket{v}{v} \in \mathbb{R}$ for all $v \in H$ (conjugate symmetry);
	\item $\braket{v}{v} \geq 0$, and $\braket{v}{v} = 0$ implies $v = 0$ (positive definiteness).
\end{enumerate}
A pre-Hilbert space $H$ becomes a normed vector space with the norm $\| v \| := \braket{v}{v}^{\frac{1}{2}}$; then $H$ is a \emph{Hilbert space} if it is complete as a metric space with respect to this norm.
\end{dfn} \index{pre-Hilbert space} \index{Hilbert space} \index{inner product} 

While the definition of Hilbert space is standard, that of a morphism of Hilbert spaces is not quite as straightforward. There are, in fact, two common choices, defined as follows.

\begin{dfn}
Let $H, K$ be normed vector spaces. A \emph{bounded} linear map $f: H \to K$ is a linear map such that, for some $c \geq 0$, it holds that $\|fv\|_K \leq c\|v\|_H$ for all $v \in H$.

A \emph{short} linear map $f: H \to K$ is a linear map such that $\|fv\|_K \leq \|v\|_H$ for all $v \in H$. \index{bounded map} \index{short map}

We write $\cat{Hilb}$ for the category of Hilbert spaces and bounded linear maps, and $\cat{Hilb}_{\leq 1}$ for the category of Hilbert spaces and short linear maps. We also write $\cat{FHilb}$ and $\cat{FHilb}_{\leq 1}$, respectively, for their full subcategories on finite-dimensional Hilbert spaces. \index{Hilb@$\cat{Hilb}, \cat{FHilb}$} \index{Hilbshort@$\cat{Hilb}_{\leq 1}, \cat{FHilb}_{\leq 1}$}
\end{dfn}

As discussed in \cite{baez2006quantum}, the choice of bounded linear maps as morphisms has the unwanted consequence of making $\cat{FHilb}$ equivalent, as a category, to $\cat{FVec}_\mathbb{C}$. The missing information about the structure of Hilbert spaces, not provided by objects and morphisms, must be supplemented as additional structure: namely, the involutive functor $\dagg{(-)}: \opp{\cat{Hilb}} \to \cat{Hilb}$, sending a map $f: H \to K$ to its \emph{adjoint}, that is, the unique map $\dagg{f}: K \to H$ such that, for all $v \in K$ and $w \in H$, \index{adjoint|see {dagger}} \index{dagger}
\begin{equation*}
	\braket{v}{fw}_K = \braket{\dagg{f}v}{w}_V.
\end{equation*}
The existence and uniqueness of $\dagg{f}$ is a basic result of the theory of Hilbert spaces. Identifying a vector $v \in H$ with the evaluation of a map $v: \mathbb{C} \to H$ at $1$, we can reconstruct the inner product on $H$ from the functor $\dagg{(-)}$, by the identity
\begin{equation*}
	\braket{v}{w} = \dagg{v}w.
\end{equation*}
This equality is what motivates Dirac's bra-ket notation, where a vector $v \in H$, identified with the map $v: \mathbb{C} \to H$, is written $\ket{v}$, and its adjoint $\dagg{v}: H \to \mathbb{C}$ is written $\bra{v}$. A map $\mathbb{C} \to H$ is commonly called a \emph{state}, and a map $H \to \mathbb{C}$ an \emph{effect}, in the context of quantum theory. \index{state} \index{effect}

Formally, this makes $\cat{Hilb}$ a kind of ``category with structure'', as follows.
\begin{dfn} \index{dagger!category}
A \emph{dagger category} is a category $\cat{C}$ together with an involutive, identity-on-objects functor $(-)^\dagger: \opp{\cat{C}} \to \cat{C}$, called the \emph{dagger}. A functor $F: \cat{C} \to \cat{D}$ of dagger categories is a functor satisfying $F(f^\dagger) = (Ff)^\dagger$ for all morphisms $f$ in $\cat{C}$.
\end{dfn}
This terminology was first introduced in \cite{selinger2007dagger}, and has since become the standard in categorical quantum mechanics, where ``the category of Hilbert spaces'' is usually synonymous with the dagger category $\cat{Hilb}$. On the other hand, this choice carries a number of subtle technical and conceptual issues.
\begin{enumerate}
	\item As discussed on \cite{nlabequivalence}, the structure of dagger category on a category is not equivalence-invariant; this is related to the fact that the intended notion of isomorphism in a dagger category is a \emph{unitary} isomorphism, that is, a morphism $f$ whose inverse is $\dagg{f}$, yet a dagger category may contain non-unitary isomorphisms (as does $\cat{Hilb}$). Thus, rather than being higher algebraic structure on a category, which would make it equivalence-invariant, the dagger structure has to be taken as a primitive, akin to composition. While there is nothing wrong with this \emph{per se}, it does not seem ideal from the standpoint of a principle of economy.
	
	\item In the algebraic formulation, one proceeds ``locally'' from complex vector spaces to Hilbert spaces, by endowing a vector space with some algebraic structure, and then verifying it satisfies a certain property. If we believe in dagger categories as a foundation for quantum theory, we have to define Hilbert spaces through the ``global'' information of a dagger functor. And what category should this be defined on? 
	
	In the finite-dimensional case, we can say it is $\cat{FVec}_\mathbb{C}$; in \cite[Theorem 5.2]{vicary2011completeness}, Vicary lists some reasonable categorical properties that characterise ``the'' dagger functor turning $\cat{FVec}_\mathbb{C}$ into $\cat{FHilb}$ among all possibilities. On the other hand, not all infinite-dimensional vector spaces even admit the structure of a Hilbert space \cite{kruse1964badly}, and not all linear maps between those that do are bounded, so $\cat{Vec}_\mathbb{C}$ cannot be taken as a starting point for an abstract definition of the dagger category $\cat{Hilb}$. \index{unitary}
	\item There are notions of tensor product and direct sum of Hilbert spaces, where the first, in particular, has a fundamental physical significance, describing how quantum systems compose. If we wanted to generalise their construction from $R$-modules to Hilbert spaces as a dagger category, in the case of direct sums we could define a ``dagger functor'' of polycategories taking a multi-cospan to its adjoint multi-span, and indeed this would preserve universal 2-cells: direct sums in $\cat{Hilb}$ are dagger biproducts \cite[Section 3.2]{heunen2009categorical}, that is, they are both categorical products and coproducts, with universal diagrams adjoint to one another.
	
	On the other hand, there is no sensible notion of adjoint of a multilinear map, so there does not seem to be any clear way of factoring the dagger structure into a construction of tensor products via representability.
	
	\item Finally, in a certain sense, being a Hilbert space can be seen as a property of a complex Banach space, that is, a complete normed vector space: Hilbert spaces are precisely those Banach spaces whose norm satisfies the parallelogram law
	\begin{equation*}
		\|v\|^2 + \|w\|^2 = \frac{1}{2}(\|v + w\|^2 + \|v - w\|^2).
	\end{equation*}
	So, we would expect ``the category of Hilbert spaces'' to be a full subcategory of the category of complex Banach spaces; yet the latter is not a dagger category. This leads to the odd situation where what is taken to be the fundamental structure of the category of Hilbert spaces makes it incomparable with its most obvious generalisation.
\end{enumerate}

The category $\cat{Hilb}_{\leq 1}$ of Hilbert spaces and short linear maps is immune from some, but not all of these problems. All the properties of Hilbert spaces that only depend on their structure as normed vector spaces are now captured by the morphisms: in particular, isomorphisms in $\cat{Hilb}_{\leq 1}$ are isometric isomorphisms, that is, unitaries. In this sense, $\cat{Hilb}_{\leq 1}$ fits naturally into the category $\cat{Ban}_{\leq 1}$ of Banach spaces and short linear maps, the standard setting of functional analysis. \index{Banach space} \index{Ban@$\cat{Ban}_{\leq 1}$} 

\begin{remark} This choice also has a physical motivation, by the interpretation of inner products as probabilistic amplitudes of observations on quantum systems (Born rule): the short linear maps are the ones that do not lead to ``probabilities'' larger than 1. Nevertheless, specially in the finite-dimensional case, where all linear maps are bounded, it is customary to work in $\cat{Hilb}$, even when ultimately the physical interpretation belongs in $\cat{Hilb}_{\leq 1}$; it is assumed that the normalisation of maps can always be postponed to a later stage.
\end{remark}

On the other hand, all notions that depend explicitly on the inner product --- such as positive or self-adjoint operators --- do not have a clear categorical description in terms of morphisms only; thus, endowing $\cat{Hilb}_{\leq 1}$ with the restriction of the dagger of $\cat{Hilb}$ is still a non-trivial addition. Moreover, the natural notion of short multilinear map of normed vector spaces, that is, a multilinear map $f: (H_1, \ldots, H_n) \to H$ such that, for all $i = 1,\ldots, n$, and $h_i \in H_i$,
\begin{equation*}
	\|f(h_1, \ldots, h_n)\|_H \leq \|h_1\|_{H_1}\ldots\|h_n\|_{H_n},
\end{equation*}
does define a tensor representable multicategory, but the tensor of two Hilbert spaces receives a ``wrong'' norm, that does not satisfy the parallelogram law. In general, the norms arising from universal constructions of normed vector spaces are $1$-norms and $\infty$-norms, which do not come from inner products; see \cite{egger2014notes}.

We propose here an alternative pathway to $\cat{Hilb}_{\leq 1}$ and its dagger, in accordance with our general guidelines. First of all, we observe that the constructions of tensor products and direct sums easily generalise when morphisms of Hilbert spaces are restricted to isometries, which amounts to taking the inner product literally as algebraic structure that needs to be strictly preserved. \index{isometry} 
\begin{cons} \index{Hilbiso@$\cat{Hilb}_1, Hilb_{1, \otimes}, Hilb_{1, \oplus}$}
Let $\textit{Hilb}_{1,\otimes}$ be the regular multicategory whose 1-cells are Hilbert spaces, and 2-cells $f: (H_1, \ldots, H_n) \to H$ are multilinear isometries, that is, multilinear maps such that 
\begin{equation*}
	\braket{f(v_1,\ldots,v_n)}{f(w_1,\ldots,w_n)}_H = \braket{v_1}{w_1}_{H_1} \ldots \braket{v_n}{w_n}_{H_n},
\end{equation*}
or, graphically in $\mathbb{C}\textit{Mod}_\otimes$,
\begin{equation*}
\input{img/c4_1_multihilb.tex}
\end{equation*} \index{tensor product!of Hilbert spaces} 
Then, $\textit{Hilb}_{1,\otimes}$ is tensor representable: the tensor unit is $\mathbb{C}$ with the inner product $(\lambda, \mu) \mapsto \overline{\lambda}\mu$, and the tensor of two 1-cells $H, K$ is their tensor product as Hilbert spaces, that is, the completion $H \otimes K$ of the vector space tensor product of $H, K$ with respect to the norm induced by the inner product
\begin{equation*}
	\braket{v \otimes w}{v' \otimes w'} := \braket{v}{v'}_H \braket{w}{w'}_K.
\end{equation*}
\end{cons}

An isometry $f: H \to K$ identifies $H$ with the closed subspace $f(H)$ of $K$, and closed subspaces of a Hilbert space form a lattice; in particular, given a multi-cospan of isometries $\{f_i: H_i \to H\}_{i=1}^n$, it makes sense to speak of the \emph{intersection} of the subspaces $H_i$ in $H$. 

\begin{cons}
Let $\textit{Hilb}_{1,\oplus}$ be the regular multicategory whose 1-cells are Hilbert spaces, and 2-cells $f: (H_1, \ldots, H_n) \to H$ are mutually independent multi-cospans of isometries, that is, families $\{f_i: H_i \to H\}_{i=1}^n$ of isometries such that $H_i$ and $H_j$ have trivial intersection in $H$ whenever $i \neq j$. Then, $\textit{Hilb}_{1,\oplus}$ is tensor representable: the tensor unit is the trivial Hilbert space $0$, and the tensor of two 1-cells $H, K$ is their direct sum $H \oplus K$ as vector spaces, with the inner product defined by \index{direct sum!of Hilbert spaces}
\begin{equation*}
	\braket{v \oplus w}{v' \oplus w'} := \braket{v}{v'}_H + \braket{w}{w'}_K.
\end{equation*}
\end{cons}

\begin{remark}
Equivalently, by the positive definiteness of inner products, the 2-cells of $\textit{Hilb}_{1,\oplus}$ can be characterised as the \emph{orthogonal} multi-cospans, that is, the ones satisfying
\begin{equation*}
	\braket{f_i (v)}{f_j (w)} =  \begin{cases}
		\braket{v}{w}, & i = j, \\
		0, & i \neq j.
	\end{cases}
\end{equation*}
By a basic result in the theory of Hilbert spaces \cite[Section I.2]{conway2007course}, for all closed subspaces $H$ of a Hilbert space $K$, 
\begin{equation*}
	H^\bot := \{v \in K \,|\, \text{for all $w \in H$}, \braket{v}{w} = 0\}
\end{equation*} 
is a closed subspace of $K$, and the inclusions $\{H, H^\bot \hookrightarrow K\}$ form an orthogonal cospan exhibiting $K$ as $H \oplus H^\bot$.
\end{remark}

\begin{dfn} 
We write $\cat{Hilb}_1$ for the category $G(\textit{Hilb}_{1,\otimes}) = G(\textit{Hilb}_{1,\oplus})$, with the two monoidal structures induced by tensor representability of $\textit{Hilb}_{1,\otimes}$ and $\textit{Hilb}_{1,\oplus}$, making it a distributive monoidal category.
\end{dfn}

Thus, the construction of tensor products and disjoint sums seems to best generalise to Hilbert spaces when morphisms are restricted to isometries, with a slight restriction on the multi-cospans involved in the definition of disjoint sums, which, nevertheless, relies on the natural concept of orthogonality. 

However, the restriction to isometries is a very strong one; in particular, the only isometries whose adjoint is an isometry are the unitaries. Yet, in a certain sense, isometries are exactly ``half'' the morphisms one needs to obtain all short linear maps.

\begin{cons} \index{Cosp@$\cat{Cosp}(\cat{Hilb}_1), \textit{Cosp}(\textit{Hilb}_1)$}
We define $\textit{Cosp}(\textit{Hilb}_1)$ to be the following regular multi-bicategory: 0-cells are Hilbert spaces, 1-cells $(K, i, j): H \to H'$ are cospans 
\begin{equation*}
\begin{tikzpicture}[baseline={([yshift=-.5ex]current bounding box.center)}]
	\node[scale=1.25] (0) at (0,0) {$H$};
	\node[scale=1.25] (1) at (2,0) {$K$};
	\node[scale=1.25] (2) at (4,0) {$H'$};
	\draw[1c] (0) to node[auto] {$i$} (1);
	\draw[1c] (2) to node[auto,swap] {$j$} (1);
\end{tikzpicture}
\end{equation*}
of isometries, and 2-cells $f: ((K_1, i_1,j_1), \ldots, (K_n, i_n,j_n)) \to (K, i,j)$ are commutative diagrams
\begin{equation*}
\begin{tikzpicture}[baseline={([yshift=-.5ex]current bounding box.center)}]
	\node[scale=1.25] (0) at (-.5,0) {$H_1$};
	\node[scale=1.25] (1) at (.9,-.7) {$K_1$};
	\node[scale=1.25] (2) at (2.5,-1) {$H_2$};
	\node[scale=1.25] (3) at (4,-1) {$\ldots$};
	\node[scale=1.25] (4) at (5.5,-1) {$H_{n}$};
	\node[scale=1.25] (5) at (7.1,-.7) {$K_{n}$};
	\node[scale=1.25] (6) at (8.5,0) {$H_{n+1}$};
	\node[scale=1.25] (7) at (4,1) {$K$};
	\draw[1c] (0) to node[auto] {$i_1$} (1);
	\draw[1c] (2) to node[auto,swap] {$j_1$} (1);
	\draw[1c] (2) to node[auto] {$i_2$} (3);
	\draw[1c] (4) to node[auto,swap] {$j_{n-1}$} (3);
	\draw[1c] (4) to node[auto] {$i_n$} (5);
	\draw[1c] (6) to node[auto,swap] {$j_n$} (5);
	\draw[1c, out=22, in=180] (0) to node[auto] {$i$} (7);
	\draw[1c, out=158, in=0] (6) to node[auto,swap] {$j$} (7);
	\draw[1c, out=45, in=-165] (1) to node[auto,swap] {$f_1$} (7);
	\draw[1c, out=135, in=-15] (5) to node[auto] {$f_n$} (7);
\end{tikzpicture}
\end{equation*}
of isometries, with the following property: for all $k = 1,\ldots, n-1$, the intersection of $K_k$ and $K_{k+1}$ as closed subspaces of $K$ is $H_{k+1}$, that is, the squares
\begin{equation*}
\begin{tikzpicture}[baseline={([yshift=-.5ex]current bounding box.center)},scale=.7]
	\node[scale=1.25] (0) at (0,2) {$K$};
	\node[scale=1.25] (1) at (-2,0) {$K_k$};
	\node[scale=1.25] (2) at (2,0) {$K_{k+1}$};
	\node[scale=1.25] (3) at (0,-2) {$H_{k+1}$};
	\draw[1c] (1) to node[auto] {$f_k$} (0);
	\draw[1c] (2) to node[auto,swap] {$f_{k+1}$} (0);
	\draw[1c] (3) to node[auto] {$i_k$} (1);
	\draw[1c] (3) to node[auto,swap] {$j_{k+1}$} (2);
	\draw[edge] (-.6,-.9) to (0,-.3) to (.6,-.9);
\end{tikzpicture}
\end{equation*}
are pullback diagrams in $\cat{Hilb}_1$. Composition is pasting of commutative diagrams.

Then, $\textit{Cosp}(\textit{Hilb}_1)$ is tensor representable. A tensor unit on the Hilbert space $H$ is the cospan 
\begin{equation*}
\begin{tikzpicture}[baseline={([yshift=-.5ex]current bounding box.center)}]
	\node[scale=1.25] (0) at (0,0) {$H$};
	\node[scale=1.25] (1) at (2,0) {$H$};
	\node[scale=1.25] (2) at (4,0) {$H$,};
	\draw[1c] (0) to node[auto] {$\mathrm{id}_H$} (1);
	\draw[1c] (2) to node[auto,swap] {$\mathrm{id}_H$} (1);
\end{tikzpicture}
\end{equation*}
as witnessed, for all 1-cells $(K, i,j): H \to H'$, by the 2-cells
\begin{equation*}
\begin{tikzpicture}[baseline={([yshift=-.5ex]current bounding box.center)}]
	\node[scale=1.25] (0) at (-.5,0) {$H$};
	\node[scale=1.25] (1) at (.9,-.7) {$H$};
	\node[scale=1.25] (2) at (2.5,-1) {$H$};
	\node[scale=1.25] (5) at (4.1,-.7) {$K$};
	\node[scale=1.25] (6) at (5.5,0) {$H'$};
	\node[scale=1.25] (7) at (2.5,1) {$K$};
	\draw[1c] (0) to node[auto] {$\mathrm{id}_H$} (1);
	\draw[1c] (2) to node[auto,swap] {$\mathrm{id}_H$} (1);
	\draw[1c] (2) to node[auto] {$i$} (5);
	\draw[1c] (6) to node[auto,swap] {$j$} (5);
	\draw[1c, out=30, in=180] (0) to node[auto] {$i$} (7);
	\draw[1c, out=150, in=0] (6) to node[auto,swap] {$j$} (7);
	\draw[1c, out=75, in=-165] (1) to node[auto,swap] {$i$} (7);
	\draw[1c, out=105, in=-15] (5) to node[auto] {$\mathrm{id}_K$} (7);
\end{tikzpicture}
\end{equation*}
on the left, and similar ones on the right. 

Let $(K, i, j): H \to H', (K', i', j'): H' \to H''$ be two 1-cells, and decompose $K$ as $(H')_K^\bot \oplus H'$, and $K'$ as $H' \oplus (H')_{K'}^\bot$, so that $j = 0 \oplus \mathrm{id}_{H'}$, and $i' = \mathrm{id}_{H'} \oplus 0$. Let $K'' := (H')_K^\bot \oplus H' \oplus (H')_{K'}^\bot$; then, the tensor of the two 1-cells is the composite cospan
\begin{equation*}
\begin{tikzpicture}[baseline={([yshift=-.5ex]current bounding box.center)}]
	\node[scale=1.25] (0) at (0,0) {$H$};
	\node[scale=1.25] (1) at (2,0) {$K$};
	\node[scale=1.25] (1b) at (4.5,0) {$K''$};
	\node[scale=1.25] (2b) at (7,0) {$K'$};
	\node[scale=1.25] (2) at (9,0) {$H'',$};
	\draw[1c] (0) to node[auto] {$i$} (1);
	\draw[1c] (1) to node[auto] {$\mathrm{id}_K \oplus 0$} (1b);
	\draw[1c] (2b) to node[auto,swap] {$0 \oplus \mathrm{id}_{K'}$} (1b);
	\draw[1c] (2) to node[auto,swap] {$j'$} (2b);
\end{tikzpicture}
\end{equation*}
together with the 2-cell 
\begin{equation*}
\begin{tikzpicture}[baseline={([yshift=-.5ex]current bounding box.center)}]
	\node[scale=1.25] (0) at (-.5,0) {$H$};
	\node[scale=1.25] (1) at (.9,-.7) {$K$};
	\node[scale=1.25] (2) at (2.5,-1) {$H'$};
	\node[scale=1.25] (5) at (4.1,-.7) {$K'$};
	\node[scale=1.25] (6) at (5.5,0) {$H''$};
	\node[scale=1.25] (7) at (2.5,1) {$K''$};
	\draw[1c] (0) to node[auto] {$i$} (1);
	\draw[1c] (2) to node[auto,swap] {$j$} (1);
	\draw[1c] (2) to node[auto] {$i'$} (5);
	\draw[1c] (6) to node[auto,swap] {$j'$} (5);
	\draw[1c, out=30, in=180] (0) to node[auto] {$(\mathrm{id}_K \oplus 0)i$} (7);
	\draw[1c, out=150, in=0] (6) to node[auto,swap] {$(0 \oplus \mathrm{id}_{K'})j'$} (7);
	\draw[1c, out=75, in=-165] (1) to node[auto,swap, pos=0.3] {$\mathrm{id}_K \oplus 0$} (7);
	\draw[1c, out=105, in=-15] (5) to node[auto, pos=0.6] {$0 \oplus \mathrm{id}_{K'} \!\!\!$} (7);
	\node[scale=1.25] at (6,-.7) {,};
\end{tikzpicture}
\end{equation*}
which is $(c^*)$-divisible because the square formed by $j, i', \mathrm{id}_K \oplus 0,$ and $0 \oplus \mathrm{id}_{K'}$ is a pushout in $\cat{Hilb}_1$.
\end{cons}

\begin{dfn}
We write $\cat{Cosp}(\cat{Hilb}_1)$ for the bicategory $G(\textit{Cosp}(\textit{Hilb}_1))$, with the composition and units induced by tensor representability.
\end{dfn}

Any strictification of $(\cat{Hilb}_1, \oplus, 0)$ induces a strictification of $\cat{Cosp}(\cat{Hilb}_1)$, by the existence of canonical composites of cospans provided by orthogonal decompositions. In what follows, we will treat $\cat{Cosp}(\cat{Hilb}_1)$ as a strict 2-category.

Using the two monoidal structures on $\cat{Hilb}_1$, we can tensor together 2-cells of $G(\textit{Cosp}(\textit{Hilb}_1))$, producing other 2-cells, both with tensor products and direct sums; the only thing to verify is preservation of pullback squares of isometries, which is easily checked. These operations also preserve tensors and tensor units, by distributivity of tensor products over direct sums in one case, and by symmetry of direct sums in the other, so they induce the structure of a distributive monoidal 2-category on $\cat{Cosp}(\cat{Hilb}_1)$.

Finally, like all (higher) categories whose 1-cells are spans or cospans of morphisms, $\cat{Cosp}(\cat{Hilb}_1)$ admits a canonical dagger: for all $(K, i, j): H \to H'$, define $(K, i, j)^\dagger := (K, j, i): H' \to H$.

\begin{lem}
Let $i: H \to K$ be an isometry of Hilbert spaces. Then:
\begin{enumerate}[label=(\alph*)]
	\item $\dagg{i} i = \mathrm{id}_H$, and
	\item if $j: H' \to K$ is another isometry, then $\dagg{j} i: H \to H'$ is a short linear map.
\end{enumerate}
\end{lem}
\begin{proof}
The first point is \cite[Proposition II.2.17]{conway2007course}. For the second point, by \cite[Proposition II.2.7]{conway2007course}, if a linear map $f$ is bounded with constant $c \leq 0$, then $\dagg{f}$ is bounded with the same constant; in particular, if $j$ is short, $\dagg{j}$ is also short. Because isometries are short, and composites of short linear maps are short, we conclude.
\end{proof}

\begin{lem} \label{lem:hilbconnect}
Let $t: H \to H'$ be a short linear map of Hilbert spaces. Then there exist a Hilbert space $K$ and a cospan $(K,i,j): H \to H'$ of isometries such that
\begin{enumerate}
	\item $t = \dagg{j} i$, and
	\item if $(K', i', j'): H \to H'$ is another cospan of isometries such that $t = \dagg{(j')}i'$, then there is a 2-cell $f: (K, i, j) \to (K', i', j')$ in $\cat{Cosp}(\cat{Hilb}_1)$.
\end{enumerate}
\end{lem}
\begin{proof}
Given any short linear map $t: H \to H'$, the operators $(\mathrm{id}_H - t^\dagger t): H \to H$ and $(\mathrm{id}_{H'} - tt^\dagger): H' \to H'$ are positive, therefore they have unique non-negative square roots
\begin{equation*}
	d_t := (\mathrm{id}_H - t^\dagger t)^\frac{1}{2}, \quad \quad d_{t^\dagger} := (\mathrm{id}_{H'} - t t^\dagger)^\frac{1}{2}.
\end{equation*}
As shown, for instance, in \cite[Section I.3]{sznagy2010harmonic}, these operators satisfy
\begin{equation} \label{eq:defect}
	t d_t = d_{t^\dagger} t, \quad \quad t^\dagger d_{t^\dagger} = d_t t^\dagger,
\end{equation}
so letting $D_t$ be the closure of $d_t(H)$ in $H$ and $D_{t^\dagger}$ the closure of $d_{t^\dagger}(H')$ in $H'$, we have that $t(D_t)$ is contained in $D_{t^\dagger}$, and $t^\dagger(D_{t^\dagger})$ in $D_t$.

Then, we can define $u: H \oplus D_{t^\dagger} \to H' \oplus D_t$ to be, in matrix notation,
\begin{equation} \label{eq:dilation}
\begin{pmatrix} t & d_{t^\dagger} \\ d_t & -t^\dagger \end{pmatrix};
\end{equation}
using equation $(\ref{eq:defect})$, this is shown to be unitary. Letting $K := H' \oplus D_t$, $i := u(\mathrm{id}_H \oplus 0): H \to K$, and $j := \mathrm{id}_{H'} \oplus 0: H' \to K$, we obtain a cospan $(K, i, j)$ of isometries such that $t = j^\dagger i$.

To show that $(K, i, j)$ is initial among the cospans of isometries with this property, we adapt the argument of \cite[Proposition 2.1]{levy2014dilation}. Given a cospan $(K', i', j')$ such that $t = \dagg{(j')}i'$, define $f: K \to K'$ on elements of the form $w \oplus d_t(v)$ by
\begin{equation*}
	f(w \oplus d_t(v)) := i'(v) + j'(w - t(v)),
\end{equation*}
and extend uniquely along limits to their closure. Then, $i' = fi$ and $j' = fj$, and $f$ is an isometry; this can be verified using the equations
\begin{align*}
	\braket{d_t(v)}{d_t(v')} & = \braket{v}{d_t^2(v')} = \braket{v}{v'} - \braket{v}{t^\dagger t(v')} = \\
	& = \braket{v}{v'} - \braket{t(v)}{t(v')},
\end{align*}
and
\begin{equation*}
	\braket{i'(v)}{j'(w)} = \braket{(j')^\dagger i'(v)}{w} = \braket{t(v)}{w},
\end{equation*}
for all $v, v' \in H$, $w \in H'$.
\end{proof}

\begin{prop} \label{prop:hilbtrunc}
The category $\cat{Hilb}_{\leq 1}$ is the 1-truncation of $\cat{Cosp}(\cat{Hilb}_1)$. The truncation map preserves the dagger structure, tensor products and direct sums.
\end{prop}
\begin{proof}
Let $\tau: \cat{Cosp}(\cat{Hilb}_1) \to \cat{Hilb}_{\leq 1}$ be the identity on 0-cells, send a 1-cell $(K, i, j): H \to H'$ to the short linear map $j^\dagger i: H \to H'$, and 2-cells to unit 2-cells. We will show that $\tau$ is a 1-truncation map.

First of all, $\tau$ is well-defined. Unit cospans $(H, \mathrm{id}_H, \mathrm{id}_H)$ on $H$ are mapped to units $\mathrm{id}_H^\dagger \mathrm{id}_H = \mathrm{id}_H$ in $\cat{Hilb}_{\leq 1}$. To show that composition of 1-cells is preserved, it suffices to observe that squares of the form
\begin{equation*}
\begin{tikzpicture}[baseline={([yshift=-.5ex]current bounding box.center)},scale=.8]
	\node[scale=1.25] (0) at (0,2) {$H' \oplus H \oplus H''$};
	\node[scale=1.25] (1) at (-2,0) {$H' \oplus H$};
	\node[scale=1.25] (2) at (2,0) {$H \oplus H''$};
	\node[scale=1.25] (3) at (0,-2) {$H$};
	\draw[1c] (1) to node[auto] {$\mathrm{id}_{H' \oplus H} \oplus 0$} (0);
	\draw[1c] (0) to node[auto] {$(0 \oplus \mathrm{id}_{H \oplus H''})^\dagger$} (2);
	\draw[1c] (1) to node[auto,swap] {$(0 \oplus \mathrm{id}_H)^\dagger$} (3);
	\draw[1c] (3) to node[auto,swap] {$\mathrm{id}_H \oplus 0$} (2);
\end{tikzpicture}
\end{equation*}
commute in $\cat{Hilb}_{\leq 1}$. Finally, for all 2-cells 
\begin{equation*}
\begin{tikzpicture}[baseline={([yshift=-.5ex]current bounding box.center)}]
	\node[scale=1.25] (0) at (0,0) {$H$};
	\node[scale=1.25] (1) at (2,-1) {$K$};
	\node[scale=1.25] (2) at (2,1) {$K'$};
	\node[scale=1.25] (3) at (4,0) {$H'$};
	\draw[1c, out=45, in=180] (0) to node[auto] {$i'$} (2);
	\draw[1c, out=135, in=0] (3) to node[auto,swap] {$j'$} (2);
	\draw[1c, out=-45, in=180] (0) to node[auto,swap] {$i$} (1);
	\draw[1c, out=-135, in=0] (3) to node[auto] {$j$} (1);
	\draw[1c] (1) to node[auto] {$f$} (2);
\end{tikzpicture}
\end{equation*}
in $\cat{Cosp}(\cat{Hilb}_1)$, we have $t:= (j')^\dagger i' = (fj)^\dagger (fi) = j^\dagger(f^\dagger f)i = j^\dagger i$ because $f$ is an isometry, so it makes sense to map $f$ to the unit $\idd{}t: t \to t$.

To prove that $\tau$ is a 1-truncation map, it suffices to show that:
\begin{enumerate}
	\item all short linear maps $t$ of Hilbert spaces are of the form $j^\dagger i$ for some cospan $(K, i, j)$ of isometries, and
	\item whenever two 1-cells $(K_1, i_1, j_1)$ and $(K_2, i_2, j_2)$ are mapped to the same linear map $t$, they are connected in $\cat{Cosp}(\cat{Hilb}_1)$, in the sense that there exists a third cospan $(K, i, j)$ and a pair of 2-cells $f_k$, $k = 1, 2$, that have $(K, i, j)$ on one boundary, and $(K_k, i_k, j_k)$ on the other.
\end{enumerate}
This follows from Lemma \ref{lem:hilbconnect}, since the cospan $(K, i, j)$ constructed there is connected to all cospans mapped to $t$.

To conclude, we have that $\tau((K,i,j)^\dagger) = \tau(K,j,i) = i^\dagger j = (j^\dagger i)^\dagger = (\tau(K,i,j))^\dagger$, so $\tau$ is compatible with the dagger; preservation of tensor products and direct sums follows from their compatibility with adjoints. This completes the proof.
\end{proof}

By Proposition \ref{prop:hilbtrunc}, we can see $\cat{Cosp}(\cat{Hilb}_1)$ as a natural higher-dimensional refinement of $\cat{Hilb}_{\leq 1}$, meeting all our requirements: the ``primitive'' morphisms are morphisms of algebras, all compositions are defined by universal properties, and the dagger structure is canonical. Moreover, we can define a multi-bicategory $\textit{Cosp}(\textit{Ban}_1)$ whose 1-cells are cospans of isometries of Banach spaces, into which $\textit{Cosp}(\textit{Hilb}_1)$ fully embeds; what will be missing is a canonical composition of cospans, due to the lack of orthogonal decompositions of arbitrary Banach spaces. 

The important class of self-adjoint operators on Hilbert spaces can be lifted to $\cat{Cosp}(\cat{Hilb}_1)$, with an appropriate weakening: the corresponding cospans are those that are ``self-adjoint up to isomorphism''.
\begin{prop} \label{prop:selfadjoint}
Let $(K,i,j): H \to H$ be a 1-cell in $\cat{Cosp}(\cat{Hilb}_1)$, and $t = \dagg{j}i$. Then $t = \dagg{t}$ if and only if $(K,i,j)$ and $(K,j,i)$ are isomorphic 1-cells.
\end{prop}
\begin{proof}
Suppose $u: (K,i,j) \to (K,j,i)$ is an isomorphism of 1-cells, that is, $u: K \to K$ is a unitary operator such that $ui = j$, $uj = i$. Then also $j = \dagg{u}i$, so $\dagg{j}i = \dagg{i}\dagg{u}i = \dagg{i}j$, that is, $t = \dagg{t}$.

Conversely, suppose $t = \dagg{t}$. In this case, the minimal unitary dilation (\ref{eq:dilation}) has the form
\begin{equation*}
u := \begin{pmatrix} t & d_t \\ d_t & -t \end{pmatrix},
\end{equation*}
which is also self-adjoint, that is, it is an involution on $H \oplus D_t$. The corresponding cospan of isometries is $(K', i', j') := (H \oplus D_t, u(\mathrm{id}_H \oplus 0), \mathrm{id}_H \oplus 0)$, and $u$ itself induces an isomorphism $u: (K', i', j') \to (K', j', i')$.

In general, given an isometry $f: (K', i', j') \to (K, i, j)$ constructed as in the proof of Proposition \ref{prop:hilbtrunc}, we can decompose $K$ as $f(K') \oplus f(K')^\bot$, so that $\dagg{f}$ restricts to a unitary $\invrs{f}: f(K') \to K'$. Then, $v := fu\invrs{f} \oplus \mathrm{id}_{f(K')^\bot}$ is an involution on $K$ satisfying $vi = vfi' = fj' = j$, hence also $vj = i$.
\end{proof}

\begin{remark} In fact, as the proof shows, Proposition \ref{prop:selfadjoint} can be strengthened by requiring that the isomorphism $(K, i, j) \to (K, j, i)$ be an involution. Note that the ``strictly self-adjoint'' cospans of isometries, those of the form $(K, i, i)$, are all mapped to identities by the truncation map. 
\end{remark}

We note that categories of cobordisms, used in \cite{baez2006quantum} as a motivating analogy for the study of Hilbert spaces as a dagger category, also admit a description of morphisms as cospans of embeddings of manifolds. Thus, it is worth asking whether certain naturally arising functors of dagger categories may lift to functors of bicategories $\cat{Cosp}(\cat{C}) \to \cat{Cosp}(\cat{D})$ induced by functors $\cat{C} \to \cat{D}$, or similarly with spans replacing cospans, and in these cases one could dispense entirely with the dagger. 

A possible case study is the CP construction \cite{coecke2016pictures}, an infinite-dimensional analogue of Selinger's CPM construction \cite{selinger2007dagger}, which, applied to $\cat{Hilb}_{\leq 1}$, produces the category of quantum operations, a broader family of physical processes including the probabilistic mixing of quantum states. We leave this question to future work.

\section{Frobenius algebras and ZX calculi} \label{sec:frobenius}

Hilbert spaces are very homogeneous spaces: to be able to perform computations with specific quantum states, one needs a reference frame, in the form of an orthonormal basis. Moreover, because maps that send basis vectors to basis vectors are exactly as expressive as functions of sets, orthonormal bases are commonly used as a stand-in for classical information in a quantum context. Thus, one of the early challenges of categorical quantum mechanics was to find a suitably abstract formalism for their handling. \index{orthonormal basis}

Depending on the way that a vector space or a Hilbert space is defined, it may come with a canonical choice of basis. One particular case is that of vector spaces (or $R$-modules, or abelian groups) that are free on a set $X$. The categories that we defined in the previous section come with a sequence of forgetful functors
\begin{equation*}
\begin{tikzpicture}[baseline={([yshift=-.5ex]current bounding box.center)}]
	\node[scale=1.25] (0) at (0,0) {$R\cat{Mod}$};
	\node[scale=1.25] (1) at (2,0) {$\cat{Ab}$};
	\node[scale=1.25] (2) at (4,0) {$\cat{Set}$};
	\draw[1c] (0) to (1);
	\draw[1c] (1) to (2);
\end{tikzpicture}
\end{equation*}
which have left adjoints: the free abelian group on the set $X$ is the group $\mathbb{Z}X$ of formal finite sums of elements of $X$, and the free $R$-module on an abelian group $G$ is the tensor product $R \otimes G$, with the action $r \cdot (r'\otimes g) := rr' \otimes g$. The ``free algebra'' functors also preserve monoidal structures, in the sense that disjoint unions and cartesian products of sets are mapped to direct sums and tensor products of abelian groups and $R$-modules. \index{module!free} \index{abelian group!free} 

We will write, in particular, $\mathbb{C}X$ for the free complex vector space on a set $X$, which can be seen as the space of formal finite sums of elements of $x$ with complex coefficients. \index{vector space!free}

We could try to extend the sequence with a forgetful functor from some version of the category of Hilbert spaces to $\cat{Vec}_\mathbb{C}$, yet this would have no left adjoint. Nevertheless, there is a construction that comes close to the intuition of a ``free Hilbert space on a set $X$'': this is the Hilbert space
\begin{equation*}
	\ell^2(X) := \Big\{f: X \to \mathbb{C} \,\Big|\, \sum_{x \in X} |f(x)|^2 < +\infty\Big\},
\end{equation*}
with the inner product $\braket{f}{g} = \sum_{x \in X} \overline{f(x)}g(x)$. It comes with a canonical orthonormal basis, indexed by elements $x \in X$, given by the functions $\delta_x : X \to \mathbb{C}$ such that, for all $y \in X$,
\begin{equation*}
	\delta_x(y) =  \begin{cases}
		1, & y = x, \\
		0, & y \neq x.
	\end{cases}
\end{equation*}
As explained in \cite{heunen2013functor}, this construction extends to a functor not on $\cat{Set}$, but on the category $\cat{PInj}$ of sets and partial injections. 

\begin{remark}
In fact, $\ell^2$ can be first seen as a functor $\cat{Inj} \to \cat{Hilb}_1$, mapping injections of sets to isometries of Hilbert spaces. This induces a map $\cat{Cosp}(\cat{Inj}) \to \cat{Cosp}(\cat{Hilb}_1)$ of 2-categories, whose domain is defined similarly to the codomain, with disjoint unions and complements of subsets replacing direct sums and orthogonal complements. Passing to 1-truncations, we recover Heunen's functor $\ell^2: \cat{PInj} \to \cat{Hilb}_{\leq 1}$. \index{l2@$\ell^2(-)$}
\end{remark}

In quantum computation, one is often only interested in finite-dimensional Hilbert spaces, which, at least, makes certain distinctions less meaningful. As we already discussed, the underlying category of $\cat{FHilb}$ is equivalent to $\cat{FVec}_\mathbb{C}$, and if $X$ is a finite set, $\ell^2(X)$ has $\mathbb{C}X$ as its underlying vector space.

The ZW calculus of Chapter \ref{chap:zwcalculus}, in its vanilla version, is an axiomatisation of a full subcategory of $\cat{Ab}$ whose objects are all free on a finite set. Through the ``free $\mathbb{C}$-module'' construction, this becomes a subcategory of $\cat{FVec}_\mathbb{C}$, whose objects are also free on a finite set; and through the identification of $\mathbb{C}X$ with $\ell^2(X)$, a subcategory of $\cat{FHilb}$. Formally, this is how an algebraic theory which, \emph{per se}, describes abelian groups can be interpreted as a theory of quantum systems.

When $X$ is finite, $\mathbb{C}X$ also comes with a canonical self-duality, that is, a self-adjunction in $(\cat{Vec}_\mathbb{C}, \otimes, \mathbb{C})$, seen as a 2-category with a single 0-cell. This is given by the pair of maps
\begin{align*}
	\eta: \mathbb{C} \to \mathbb{C}X \otimes \mathbb{C}X, \quad \quad \quad & 1 \mapsto \sum_{x \in X} x \otimes x, \\
	\varepsilon: \mathbb{C}X \otimes \mathbb{C}X \to \mathbb{C}, \quad \quad \quad  & x \otimes y \mapsto \begin{cases}
		1, & x = y, \\
		0, & x \neq y,
	\end{cases}
\end{align*}
for all $x, y \in X$. In string diagrams, they can be portrayed as
\begin{equation*} 
\input{img/c4_2_selfdual.tex}
\end{equation*} 
and satisfy the equations
\begin{equation} \label{eq:selfdual}
\input{img/c4_2_snake.tex}
\end{equation} 
\begin{equation} \label{eq:selfdual2}
\input{img/c4_2_duality_comm.tex}
\end{equation} 
Moreover, under the identification of $\mathbb{C}X$ and $\ell^2(X)$, the two maps are each other's adjoint in $\cat{FHilb}$: $\eta^\dagger = \varepsilon$. 

This, coupled with the fact that every finite-dimensional vector space is isomorphic to $\mathbb{C}X$, and every finite-dimensional Hilbert space is unitarily isomorphic to $\ell^2(X)$ for some finite set $X$, makes $\cat{FVec}_\mathbb{C}$ a compact closed category, and $\cat{FHilb}$ a \emph{dagger} compact closed category \cite{selinger2011finite}. 

\begin{remark}
In fact, Selinger proved that $\cat{FHilb}$ is complete as a model of the algebraic theory of dagger compact closed categories: if some equation in the language of dagger compact closed categories holds for all morphisms in $\cat{FHilb}$, then it holds in the theory.
\end{remark}

If we are content with fixing one orthonormal basis and consistently working with it, we can restrict ourselves to Hilbert spaces of the form $\ell^2(X)$. However, many fundamental features of quantum theory, such as complementarity in the sense of Bohr, rely on the interplay between different observables: formally, different self-adjoint operators, which come with different orthonormal bases of eigenvectors. It is desirable, then, to have a categorical formalism for describing different orthonormal bases on the same Hilbert space.

First of all, we can see the equations (\ref{eq:selfdual}) and (\ref{eq:selfdual2}) as the defining axioms of a PRO.
\begin{dfn}
The theory of \emph{self-dualities} is the PRO $\textit{SDual}$ presented by the 2-generators $\eta: [0] \to [2]$, $\varepsilon: [2] \to [0]$ and the axioms $adj_L, adj_R$. \index{self-duality} \index{SDual@$\textit{SDual}, \textit{CSDual}$}

The theory $\textit{CSDual}$ is the PROP containing $\textit{SDual}$, satisfying the additional axioms $com$, $\coo{com}$, together with the usual axioms of PROPs.
\end{dfn}
Then, we observe that $\textit{SDual}$ is a sub-PRO of $\textit{Frob}$, and $\textit{CSDual}$ a sub-PROP of $\textit{CFrob}$, through the assignments 
\begin{equation*} 
\input{img/c4_2_dual_frobenius.tex}
\end{equation*} 
that is, assigning a Frobenius algebra structure to an object of a monoidal category determines a self-duality on the same object.

In fact, the two maps $\eta: \mathbb{C} \to \mathbb{C}X \otimes \mathbb{C}X$, $\varepsilon: \mathbb{C}X \otimes \mathbb{C}X \to \mathbb{C}$ can be decomposed into the maps
\begin{align*}
	c: \mathbb{C}X \to \mathbb{C}X \otimes \mathbb{C}X, \quad \quad \quad & x \mapsto x \otimes x, \\
	m: \mathbb{C}X \otimes \mathbb{C}X \to \mathbb{C}X, \quad \quad \quad & x \otimes y \mapsto \begin{cases}
		x, & x = y, \\
		0, & x \neq y, 	\end{cases}
\end{align*}
\begin{align*}
	d: \mathbb{C}X \to \mathbb{C}, \quad \quad \quad & x \mapsto 1, \\
	u: \mathbb{C} \to \mathbb{C}X, \quad \quad \quad & 1 \mapsto \sum_{x \in X} x,
\end{align*}
for all $x, y \in X$, defining a \emph{special} commutative Frobenius algebra structure on $\mathbb{C}X$. Furthermore, on $\ell^2(X)$, this is a \emph{dagger} special commutative Frobenius algebra, in the sense that $m = c^\dagger$, and $u = d^\dagger$. 

We will refer to dagger special commutative Frobenius algebras shortly as $\dagger$-SCFAs; note that a $\dagger$-SCFA on a Hilbert space is completely specified by either its monoid part $(m, u)$, or its comonoid part $(c,d)$. \index{daggerSCFA@$\dagger$-SCFA}

In \cite{coecke2012new}, this was taken as the starting point of a characterisation of orthonormal bases in the spirit of categorical universal algebra.
\begin{thm} \label{thm:bases}
Orthonormal bases on a finite-dimensional Hilbert space $H$ are in a biunivocal correspondence with $\dagger$-SCFAs on $H$ in $\cat{FHilb}$.
\end{thm}
\begin{proof}
This is \cite[Theorem 5.1]{coecke2012new}.
\end{proof}

In the infinite-dimensional case, the ``copy'' and ``multiplication'' maps are still definable, and are in fact short, that is, they are morphisms in $\cat{Hilb}_{\leq 1}$. On the other hand, the unit and discard maps are not even bounded, so the characterisation cannot be immediately generalised. Nevertheless, a refinement of this result was proved in \cite{abramsky2012algebras}, as we proceed to explain.

\begin{dfn}
The theory of \emph{non-unital Frobenius algebras} is the sub-PRO $\textit{Frob}_*$ of $\textit{Frob}$ whose 2-generators are restricted to $c: [1] \to [2]$ and $m: [2] \to [1]$. \index{Frobenius algebra!non-unital} \index{Frobstar@$\textit{Frob}_*$}

There are corresponding PROs and PROPs of non-unital special, commutative, and special commutative Frobenius algebras $\textit{SFrob}_*$, $\textit{CFrob}_*$, and $\textit{SCFrob}_*$, with the same restriction on the generators.
\end{dfn}
In a compact closed category, non-unital Frobenius algebras are automatically unital \cite{carboni1991matrices}, so considering models of $\textit{Frob}_*$ instead of $\textit{Frob}$ in $\cat{FHilb}$ is no restriction. 

\begin{thm} \label{thm:bases_infinite}
Orthonormal bases on a Hilbert space $H$ are in a biunivocal correspondence with non-unital $\dagger$-SCFAs on $H$ in $\cat{Hilb}_{\leq 1}$, which furthermore satisfy the following condition: for all $v: \mathbb{C} \to H$, there exists a $v_*: \mathbb{C} \to H$ such that
\begin{equation*} 
\input{img/c4_2_axiom_h.tex}
\end{equation*} 
\end{thm}
\begin{proof}
This is part of \cite[Theorem 22]{abramsky2012algebras}.
\end{proof}
The property holds automatically for unital Frobenius algebras: Theorem \ref{thm:bases_infinite} is a proper generalisation of Theorem \ref{thm:bases}. 

\begin{remark} In \cite{gogioso2017infinite}, using ideas of non-standard analysis, the authors showed that, in some situations, it is possible to reason with orthonormal bases in infinite dimension ``as if'' they were unital Frobenius algebras. 

Theorem \ref{thm:bases} was also generalised in a different direction: in \cite{vicary2010categorical, coecke2014categories}, certain \emph{non-commutative} dagger Frobenius algebras in $\cat{FHilb}$ are shown to correspond to finite-dimensional C$^*$-algebras, which can be used to model quantum information in the same way as orthonormal bases model classical information.
\end{remark}

Both in the finite and in the infinite-dimensional case, the elements of the orthonormal basis are recovered as the ``copyables'' of the copy operation of the Frobenius algebra: that is, $v: \mathbb{C} \to H$ is an element of the orthonormal basis if and only if
\begin{equation*} 
\input{img/c4_2_copyables.tex}
\end{equation*}

Besides the copyables, a $\dagger$-SCFA comes with another distinguished family of states. 
\begin{dfn} \index{phase} \index{phase!shift}
Let $(c, d)$ be a $\dagger$-SCFA on a finite-dimensional Hilbert space $H$. A state $\varphi: \mathbb{C} \to H$ is a \emph{phase} for $(c,d)$ if
\begin{equation*} 
\input{img/c4_2_phase_state.tex}
\end{equation*} 
If $\varphi$ is a phase for $(c,d)$, the corresponding \emph{phase shift} on $H$ is the unitary operator on $H$ defined by
\begin{equation*} 
\input{img/c4_2_phase_shift.tex}
\end{equation*} 
\end{dfn}
If a $\dagger$-SCFA on $H$ corresponds to the orthonormal basis $\{v_i\}_{i=1}^n$, its phases are precisely the states of the form
\begin{equation*}
	\varphi = e^{i\varphi_1}v_1 + \ldots + e^{i\varphi_n}v_n,
\end{equation*}
for some $\varphi_i \in [0, 2\pi)$, and its phase shifts are the unitary operators defined by $v_k \mapsto e^{i\varphi_k}v_k$, $k = 1, \ldots, n$. 

In \cite{coecke2011phase}, the authors showed abstractly, for $\dagger$-SCFAs in arbitrary dagger monoidal categories, that phase shifts form an abelian group under composition; for a $\dagger$-SCFA on an $n$-dimensional Hilbert space, this group is an $n$-fold product of the circle group $U(1)$. 

By construction, phase shifts also satisfy 
\begin{equation*} 
\input{img/c4_2_spider_phase.tex}
\end{equation*} 
for algebra operations of any arity. From the point of view of diagrammatic algebra, this enables us to extend the ``spider'' presentation of $\textit{SCFrob}$ to spiders labelled with elements of the phase group. The $cut$ axioms are amended as follows:
\begin{equation} \label{eq:spider_phase}
\input{img/c4_2_spider_merge.tex}
\end{equation} 
where one or more wires connect the nodes on the left hand side, and $\varphi + \psi$ is composition in the phase group. \index{spider presentation!labelled}

\begin{dfn} \index{Frobenius algebra!labelled} \index{FrobG@$\textit{Frob}_G$}
Let $G$ be an abelian group. The theory of \emph{$G$-labelled Frobenius algebras} is the PRO $\textit{Frob}_G$, with generators $s_n^m[g]: [n] \to [m]$, for $m, n \in \mathbb{N}$, $g \in G$, satisfying $s_1^1[e_G] = \idd{}a$, where $e_G$ is the unit of $G$, and the $cut$ axioms of the spider presentation of $\textit{Frob}$, modified as in equation (\ref{eq:spider_phase}) to include the composition of elements of $G$.

There are corresponding PROs and PROPs of $G$-labelled special, commutative, and special commutative Frobenius algebras $\textit{SFrob}_G$, $\textit{CFrob}_G$, $\textit{SCFrob}_G$, similarly defined.
\end{dfn}

\begin{remark}
Any homomorphism $f: G \to H$ of abelian groups induces a map of PROs $\textit{Frob}_G \to \textit{Frob}_H$, by the assignment $s_n^m[g] \mapsto s_n^m[f(g)]$. In particular, since $\textit{Frob}$ is equivalent to $\textit{Frob}_1$, where $1$ is the trivial group, it is a sub-PRO of all the $\textit{Frob}_G$, with the inclusion $s_n^m \mapsto s_n^m[e_G]$.
\end{remark}

Let $H$ be a Hilbert space with two (possibly non-unital) $\dagger$-SCFAs on it; this is the same as $H$ having two fixed orthonormal bases. We draw the Frobenius algebra operations in spider notation, with dots of different colours to distinguish between the two. In \cite{coecke2008interacting}, Coecke and Duncan proved that, in the finite-dimensional case, complementarity of the two bases is equivalent to the following relation of the corresponding Frobenius algebras: \index{orthonormal basis!complementary}
\begin{equation} \label{eq:complementary}
\input{img/c4_2_complementary.tex}
\end{equation} 
where 
\begin{equation*} 
\input{img/c4_2_dimension.tex}
\end{equation*} 
evaluates on $1$ to the dimension of $H$. 
\begin{remark}
Because the two Frobenius algebras are dagger, any equation of their operations is equivalent to its 2-opposite, that is, its vertical reflection in string diagrams. In particular, equation \ref{eq:complementary} also holds with inverted colours. Moreover, because of commutativity, any equation is also equivalent to its 1-opposite, that is, its horizontal reflection in string diagrams.
\end{remark}

Equation (\ref{eq:complementary}) has the same form as the axiom $hopf_R$ of Hopf algebras, which is equivalent to $hopf_L$ with commutative monoids or cocommutative comonoids, with the assignment
\begin{equation*} 
\input{img/c4_2_hopf_antipode.tex}
\end{equation*} 
In fact, Coecke and Duncan showed that equation (\ref{eq:complementary}) is derivable whenever the monoid part of one Frobenius algebra and the comonoid part of the other form a bialgebra ``up to an invertible scalar'', that is, up to a non-zero map $\mathbb{C} \to \mathbb{C}$.

\begin{dfn}
Two $\dagger$-SCFAs on a Hilbert space $H$ are \emph{strongly complementary} when they satisfy the following equations: \index{daggerSCFA@$\dagger$-SCFA!strongly complementary} \index{complementary|see {dagger SCFA@$\dagger$-SCFA}}
\begin{equation*} 
\input{img/c4_2_bialgebra_scaled.tex}
\end{equation*} 
\begin{equation*} 
\input{img/c4_2_bialgebra2_scaled.tex}
\end{equation*} 
\end{dfn}
The derivation of equation (\ref{eq:complementary}) from the strong complementarity equations proceeds as follows. First of all, the following hold for all strongly complementary $\dagger$-SCFAs in $\cat{FHilb}$:
\begin{equation*} 
\input{img/c4_2_transpose.tex}
\end{equation*} 
in general, we may want to impose them as additional axioms, in which case equation $sc_d$ becomes derivable.
Then,
\begin{equation*} 
\input{img/c4_2_dimension_split.tex}
\end{equation*} 
moreover, 
\begin{equation*} 
\input{img/c4_2_copy_transpose.tex}
\end{equation*} 
and similarly with inverted colours. It follows that
\begin{equation*} 
\input{img/c4_2_antipode_proof.tex}
\end{equation*} 
\begin{equation*} 
\input{img/c4_2_antipode_proof2.tex}
\end{equation*} 

Thus, strong complementarity of $\dagger$-SCFAs is indeed stronger than complementarity, and strictly so: while the classification of complementary pairs is still an open problem, all pairs of strongly complementary $\dagger$-SCFAs in $\cat{FHilb}$ have been classified.
\begin{thm} \label{thm:sc_group}
Let $(c_\gbullet, d_\gbullet)$, $(c_\wbullet, d_\wbullet)$ be two $\dagger$-SCFAs on a Hilbert space $H$ of finite dimension $n$. Then, there exist an abelian group $G$ and an isometric isomorphism $H \simeq \ell^2(G)$, such that $(c_\gbullet, d_\gbullet)$ is defined, in the canonical basis of $\ell^2(G)$, by
\begin{equation*}
	c^\dagger_\gbullet: g \otimes h \; \mapsto \; \frac{1}{\sqrt{n}} \, gh, \quad \quad \quad d^\dagger_\gbullet: 1 \; \mapsto \;\sqrt{n}\, e_G,
\end{equation*}
and $(c_\wbullet, d_\wbullet)$ by
\begin{equation*}
	c_\wbullet: g \; \mapsto \; g \otimes g, \quad \quad \quad d_\wbullet: g \; \mapsto \; 1,
\end{equation*}
for all $g, h \in G$, where $e_G$ is the unit of $G$.
\end{thm}
\begin{proof}
This is \cite[Corollary 3.10]{coecke2012strong}. The proof relies on the fact that copyables of one $\dagger$-SCFA form a subgroup of the phase group of the other.
\end{proof}

The theory of quantum computation often has a special regard for \emph{qubits}, that is, 2-dimensional Hilbert spaces, as a quantum version of classical bits. In this case, by Theorem \ref{thm:sc_group}, up to isomorphism there is a unique pair of $\dagger$-SCFAs, induced by the group $\mathbb{Z}_2$. The canonical basis elements $\ket{0}, \ket{1}$ of $\ell^2(\mathbb{Z}_2)$ are the copyables of $(c_\wbullet, d_\wbullet)$, while the copyables of $(c_\gbullet, d_\gbullet)$ are \index{qubit}
\begin{equation*}
	\ket{+} := \frac{1}{\sqrt{2}}(\ket{0} + \ket{1}), \quad \quad \ket{-} := \frac{1}{\sqrt{2}}(\ket{0} - \ket{1}).
\end{equation*}
The orthonormal bases $\{\ket{0},\ket{1}\}$ and $\{\ket{+}, \ket{-}\}$ are commonly called the \emph{Z basis} and the \emph{X basis}, respectively, of the qubit. \index{Z basis@$Z$ basis} \index{X basis@$X$ basis}

In \cite{coecke2008interacting}, Coecke and Duncan observed that these two $\dagger$-SCFAs, together with phase shifts in $U(1) \hookrightarrow U(1) \times U(1)$, $\varphi \mapsto (0,\varphi)$, form a ``universal set'' for quantum computation. More formally, let $\cat{Qubit}$ be the full monoidal subcategory of $\cat{FHilb}$ whose objects are tensor products of a finite number of copies of $\ell^2(\mathbb{Z}_2)$; this is, by definition, a PROP. Then, the map of PROPs \index{Qubitcat@$\cat{Qubit}, \cat{Qubit}_{G,H}$} \index{universality}
\begin{equation*}
	zx: \textit{SCFrob}_{U(1)} \uplus \textit{SCFrob}_{U(1)} \to \cat{Qubit},
\end{equation*}
sending one copy of $\textit{SCFrob}_{U(1)}$ to the $Z$ basis $\dagger$-SCFA, and the other to the $X$ basis $\dagger$-SCFA, with their respective phase shifts, is surjective and full, that is, surjective on both 1-cells and 2-cells.

This led naturally to the question: is it possible, by adding axioms to the spider presentation of $\textit{SCFrob}_{U(1)} \uplus \textit{SCFrob}_{U(1)}$, to reach a complete axiomatisation of $\cat{Qubit}$, that is, use a pair of Frobenius algebras as a foundation for an algebraic theory of qubits?

More in general, it may be good enough to have a finite set of basic operations which is \emph{approximately} universal; that is, for all $\varepsilon > 0$, and unitary operators $u: H \to H$, be able to construct a unitary $u_\varepsilon: H \to H$ such that $\|u - u_\varepsilon\| < \varepsilon$ in the norm of $\rimp{H}{H} \simeq H \otimes H$. Even non-universal sets of operations which can be efficiently simulated by classical means, such as those in the stabiliser fragment of quantum mechanics \cite{gottesman1997stabilizer}, may be interesting as a showcase of different features of quantum theory.

Picking one such fragment corresponds to picking a sub-PROP of $\cat{Qubit}$. Often, these are obtained by restricting the $Z$ and $X$ phase groups to some finite pair of subgroups $G, H \hookrightarrow U(1)$; for instance, the stabiliser fragment corresponds to taking $G = H = \mathbb{Z}_4$, with the inclusion $k \mapsto k\frac{\pi}{2}$, $k = 0,\ldots,3$. 

Write $\cat{Qubit}_{G,H}$ for the sub-PROP so obtained. In these cases, we have a surjective and full map
\begin{equation*}
	zx_{G,H}: \textit{SCFrob}_{G} \uplus \textit{SCFrob}_{H} \to \cat{Qubit}_{G,H},
\end{equation*}
and, again, we may ask what additional axioms are needed to obtain an equivalence.

A clear choice are the axioms $sc_a$, $sc_b$, $sc_c$ of strong complementarity. The pair of $\dagger$-SCFAs corresponding to the $Z$ and $X$ bases satisfies the additional equation
\begin{equation} \label{eq:zx_duality}
\input{img/c4_2_zx_duality.tex}
\end{equation} 
from which $sc'_a$, $sc'_b$, and $sc_d$ are derivable. 

\begin{dfn} \index{ZX calculus} \index{completeness}
A \emph{ZX calculus} is any presentation $ZX'$ of a PROP $ZX$, such that $ZX'$ contains the spider presentation of $\textit{SCFrob}_{G} \uplus \textit{SCFrob}_{H}$, for some pair of abelian groups $G, H$, together with the axioms $sc_a$, $sc_b$, $sc_c$ and $cup$.

We say that $ZX'$ is \emph{complete} for the sub-PROP $\cat{Qubit}_{G,H}$ of $\cat{Qubit}$ if there is an isomorphism of PROs $f: ZX \to \cat{Qubit}_{G,H}$ such that the following diagram commutes:
	\begin{equation*}
\begin{tikzpicture}[baseline={([yshift=-.5ex]current bounding box.center)}]
	\node[scale=1.25] (1) at (2,1.8) {$ZX$};
	\node[scale=1.25] (2) at (-3,0) {$\textit{SCFrob}_{G} \uplus \textit{SCFrob}_{H}$};
	\node[scale=1.25] (3) at (2,0) {$\cat{Qubit}_{G,H}$.};
	\draw[1c] (2) to (1);
	\draw[1c] (2) to node[auto,swap] {$zx_{G,H}$} (3);
	\draw[1c] (1) to node[auto] {$f$} (3);
\end{tikzpicture}
\end{equation*}
\end{dfn}

Since \cite{coecke2008interacting}, a number of ZX calculi have been produced, with both negative and positive results about completeness. In \cite{bonchi2014interacting}, a basic ZX calculus with $G = H = 1$ was reconstructed, as a theory of ``interacting bialgebras'', in Lack's framework of composing PROPs. In \cite{backens2014zx}, Backens proved completeness of a ZX calculus for the stabiliser fragment $\cat{Qubit}_{\mathbb{Z}_4,\mathbb{Z}_4}$, and then in \cite{backens2014single} extended the result to the approximately universal Clifford+T fragment ($G = \mathbb{Z}_8$, $H = \mathbb{Z}_4$) only for single-qubit maps. A simplified, equivalent axiomatisation was given in \cite{backens2017simplified}. 

The same ZX calculus, with phase groups extended to the whole of $U(1)$, was proven incomplete for $\cat{Qubit}$ in \cite{schroder2014incomplete}. A more expressive ZX calculus was presented and shown to be complete for the Clifford+T fragment in \cite{jeandel2017complete}, and recently a complete axiomatisation for the full theory of qubits has been produced \cite{ng2017universal}; both proofs rely on our own completeness proof for the ZW calculus, through a translation of ZX into ZW.

\section{The entanglement classification problem} \label{sec:entanglement}

Combined with the axioms of commutative Frobenius algebras, the axiom $cup$ has the remarkable effect that, in the spider presentation, string diagrams can be manipulated as undirected graphs: any shuffling of the wires attached to the same node, including the exchange of inputs and outputs using self-duality, leaves the corresponding operation unmodified. In this sense, ZX calculi can be seen as graph rewriting systems as much as presentations of algebraic theories, and their development has gone hand in hand with the development of the automated graph rewriting software \emph{Quantomatic} \cite{kissinger2015quantomatic}.

When working with ZX calculus diagrams, we can forget the distinction between inputs and outputs, between states, maps and effects, and focus entirely on the network structure of the diagram. This is sometimes phrased as an ``meta-axiom'' of ZX calculi: \emph{only the topology matters}. For example, the ternary spiders
\begin{equation*} 
\input{img/c4_3_ternary_Z.tex}
\end{equation*} 
interpreted, respectively, as the linear map $\ketbra{00}{0} + \ketbra{11}{1}$, the state $\ket{000} + \ket{111}$, and the effect $\bra{000} + \bra{111}$, can all be seen as manifestations of the same abstract network
\begin{equation} \label{eq:ternary_Z}
\begin{tikzpicture}[baseline={([yshift=-.5ex]current bounding box.center)}]
	\draw[edge, out=180, in=-45] (0,0) to node[auto, pos=1] {$A$} (-1, .5);
	\draw[edge, out=30, in=-90] (0,0) to node[auto, pos=1] {$B$} (.7, .9);
	\draw[edge, out=-60, in=90] (0,0) to node[auto, pos=1] {$C$} (.25,-1);
	\node[dotwhite] at (0,0) {};
	\node[scale=1.25] at (1.5,-1) {.};
\end{tikzpicture}

\end{equation} 

Here, we can view the white node as a ``hub'', that can be accessed through the ``ports'' $A$, $B$, $C$. In fact, a similar point of view is implicit in much of quantum information theory, where quantum states are shared resources between agents of a multi-party protocol; and it is perhaps most clear in so-called entanglement-based key exchange protocols, pioneered by Ekert's E91 \cite{ekert1991quantum}, where a shared quantum state is the only device for the sharing of information.

While the actual transmission of information from one agent to the other is not possible solely through a shared state --- the ``network'' does not exist in physical space, but in a kind of abstract, informational space --- what is possible is the random generation of shared bits of knowledge, through the correlations between results of local measurements. For example, in the network (\ref{eq:ternary_Z}), corresponding, after normalisation, to the sharing of a state \index{GHZ state@\emph{GHZ} state}
\begin{equation*}
	\ghz := \frac{1}{\sqrt{2}}(\ket{000} + \ket{111}),
\end{equation*}
if all three agents perform measurements in the $X$ basis, they will obtain bits $a, b, c \in \{0,1\}$ which satisfy $a + b + c \equiv 0\,\mathrm{mod}\,2$. This is the basis of a protocol where agents at $B$ and $C$ can share their bits $b, c$ in order to reconstruct bit $a$. Incidentally, the state $\ghz$ was the basis of a celebrated logical proof of quantum contextuality \cite{greenberger1990bell}.

The possibility of sharing information relies on the qubit at $A$ being \emph{entangled} with the qubits at $B$ and $C$: there is no way of writing $\ghz$ as a tensor product $v_A \otimes v'_{BC}$. In string diagrams, this means that there is no way of picturing network (\ref{eq:ternary_Z}) as a pair of disconnected diagrams; on the contrary, in the ZX network \index{entanglement}
\begin{equation*} 
\input{img/c4_3_disconnect_network.tex}
\end{equation*} 
corresponding, after normalisation, to the sharing of a state
\begin{equation*}
	\frac{1}{\sqrt{2}}(\ket{000} + \ket{011}) = \ket{0} \otimes \frac{1}{\sqrt{2}}(\ket{00} + \ket{11}),
\end{equation*}
no matter what measurement is performed at $A$, its results will not be correlated with the results of measurements at $B$ and $C$.

Thus, an information-theoretic property --- the absence of correlation between certain measurements --- is reflected in a topological property: the disconnectedness of string diagrams. In comparison, the vector notation is very uninformative.

On the other hand, the connectedness of a string diagram does not necessarily imply the presence of entanglement: for example, the connected ZX network
\begin{equation*} 
\input{img/c4_3_fake_connect.tex}
\end{equation*} 
corresponds to the sharing of a totally disconnected state $\ket{000}$. Nevertheless, the network can be disconnected using graph versions of the Frobenius $cut$ axioms and the $sc_b$ axiom:
\begin{equation*} 
\input{img/c4_3_fake_connect_2.tex}
\end{equation*} 
where $\diamond$ denotes equality up to an invertible scalar. One aim of completeness for ZX and related calculi is, precisely, the ability to design such networks and check their properties purely by diagrammatic reasoning, without falling back on linear algebra.

The distinction between connected and disconnected networks is quite a crude one, and one for which simple algebraic criteria exist (in the case of pure quantum states). More in general, we may be interested in how a connected network, that is, an entangled state responds to certain actions of the individual agents. 

For example, in the protocol based on $\ghz$ described earlier, suppose that the agent at $A$ decides not to cooperate, which can be described as her measuring in the complementary $Z$ basis, instead of the $X$ basis. If she obtains the result $a \in \{0,1\}$, the network becomes, up to a scalar,
\begin{equation*} 
\input{img/c4_3_uncoop_ghz.tex}
\end{equation*} 
and $B$ and $C$ are disconnected either way. This can be phrased as the network not being ``robust'' to the loss of one party.

This is not necessarily the case for all tripartite states: for example, there is no measurement of one agent on the state \index{W state@$W$ state}
\begin{equation*}
	\wstate := \frac{1}{\sqrt{3}}(\ket{001} + \ket{010} + \ket{100})
\end{equation*}
that will disconnect the other two with certainty \cite{briegel2001persistent}. This enables the state $\wstate$ to support protocols for which $\ghz$ is inadequate, see for example \cite{dhondt2006computational}.

In general, for all protocols that involve some degree of cooperation, it is important to know what an individual agent can achieve on her own. This may range from the destructive, like the $Z$ basis measurement of $\ghz$ disconnecting the network, to the completely innocuous, like performing a local unitary; in the latter case, the resulting network is equivalent to the original one for all information-theoretic purposes.

\begin{dfn}
Let $v, w: \mathbb{C} \to H \otimes \ldots \otimes H$ be two states of the $n$-fold tensor product of a Hilbert space $H$. We say that $v$ is \emph{SLOCC-reducible} to $w$, and write $w \leq_{\mathrm{SL}} v$, if there exist bounded operators $f_1, \ldots, f_n: H \to H$ such that 
\begin{equation*}
	w = (f_1 \otimes \ldots \otimes f_n) v.
\end{equation*}
Two states $v, w$ are \emph{SLOCC-equivalent} if $v \leq_{\mathrm{SL}} w$ and $w \leq_{\mathrm{SL}} v$.
\end{dfn}

\begin{remark} \index{SLOCC class}
The acronym SLOCC stands for \emph{stochastic local operations and classical communication}. In fact, SLOCC-reducibility was first defined information-theoretically, in terms of multi-party protocols converting one state into the other with a non-zero probability, but the two formulations are proven equivalent in \cite{dur2000three}.
\end{remark}

For states of two qubits, there are only two equivalence classes in the poset reflection of $\leq_{\mathrm{SL}}$, corresponding to entangled and product states. For states of three qubits, on the other hand, it was shown in \cite{dur2000three} that there are exactly two inequivalent classes of states that are maximal in this poset: one containing $\ghz$, and the other containing $\wstate$. 

For $n \geq 4$ qubits, the degrees of freedom of a state overcome those of the available local operations, which produces an infinity of SLOCC classes. These can be reduced to a finite number of ``super-classes'' by choosing a suitable parametrisation of states, which was the route taken in \cite{lamata2006inductive} to define an inductive SLOCC classification for arbitrary numbers of qubits, and then explicitly enumerate the classes of states of four qubits in \cite{lamata2007inductive}; a few cases originally missed were added by Backens in \cite{backens2017number}. Arguably, though, the point of SLOCC classification is to provide some operational insight on the uses of a quantum state in information-theoretic protocols, and such inductive, heavily algebraic approaches fall short on the subject.

Because the $\ghz$ state and its $n$-qubit generalisations
\begin{equation*}
	\ghz_n = \frac{1}{\sqrt{2}}(\ket{0\ldots 0} + \ket{1 \ldots 1})
\end{equation*}
correspond, up to a scalar, to the $Z$ spiders, there was some hope, initially, that ZX calculi could offer some insight on the problem. Indeed, by universality, the $\wstate$ state is expressible in a ZX calculus; up to a scalar, it can be pictured as the network
\begin{equation*} 
\input{img/c4_3_w_in_zxcalculus.tex}
\end{equation*} 
This, however, turned out not to be particularly fruitful, as the phases that make a network of the general form
\begin{equation*} 
\input{img/c4_3_winzx_phases.tex}
\end{equation*} 
SLOCC-equivalent to $\wstate$ do not seem to have, in general, any particular algebraic relation to the $Z$ and $X$ Frobenius algebras.

More fundamentally, we have on one side a strict dichotomy between two classes of 3-qubit states, given by SLOCC inequivalence, and on the other a calculus which is completely symmetric with respect to the exchange of its two building blocks: the ternary spider of any $\dagger$-SCFA on the qubit is SLOCC-equivalent to $\ghz$, since $\dagger$-SCFAs correspond to orthonormal bases, and any two bases are related by a unitary operator. 

In \cite{coecke2010compositional}, Coecke and Kissinger observed that there exists, in fact, a commutative Frobenius algebra $(c_\bbullet, d_\bbullet, m_\bbullet, u_\bbullet)$ on the qubit, whose spider with three outputs and no inputs is proportional to the $\wstate$ state. This is defined by \index{W state@$W$ state!Frobenius algebra}
\begin{equation*}
	c_\bbullet := \ketbra{00}{0} + \ketbra{01}{1} + \ketbra{10}{1}, \quad \quad \quad d_\bbullet := \bra{0},
\end{equation*}
\begin{equation*}
	m_\bbullet := \ketbra{1}{11} + \ketbra{0}{01} + \ketbra{0}{10}, \quad \quad \quad u_\bbullet := \ket{1}.
\end{equation*}
Anticipating a later decomposition of these operations, we are going to depict them as
\begin{equation} \label{eq:w_algebra}
\input{img/c4_3_w_algebra.tex}
\end{equation} 
respectively; general spiders should have an additional dot on each input wire. For now, however, these should be seen as ``atomic'' string diagrams.

This is not a special Frobenius algebra: it satisfies, instead, the equation
\begin{equation*} 
\input{img/c4_3_antispecial.tex}
\end{equation*} 
called \emph{anti-specialness} in \cite{coecke2010compositional}, and shown in \cite{herrmann2010models} to be canonical for Frobenius algebras whose ``loop'' operation $mc: [1] \to [1]$ disconnects. \index{Frobenius algebra!anti-special}

Coecke and Kissinger identified a class of tripartite states, including both $\ghz$ and $\wstate$, from which one can build commutative Frobenius algebras.
\begin{dfn} \index{FrobSt@$\textit{FrobSt}$}
We define $\textit{FrobSt}$ to be the PROP presented by 2-generators $v: [0] \to [3]$, $d: [1] \to [0]$, and $\varepsilon: [2] \to [0]$, satisfying the axioms 
\begin{equation*} 
\input{img/c4_3_ternary_symmetric.tex}
\end{equation*} 
\begin{equation*} 
\input{img/c4_3_frobenius_state.tex}
\end{equation*} 
together with the usual axioms of PROPs.
\end{dfn}

\begin{prop}
$\textit{FrobSt}$ and $\textit{CFrob}$ are isomorphic PROPs.
\end{prop}
\begin{proof}
These are \cite[Theorem 5 and Theorem 6]{coecke2010compositional}.
\end{proof}

As could be expected, the generators of $\textit{FrobSt}$ correspond to the spiders with the same arity in $\textit{CFrob}$. It follows that a commutative Frobenius algebra in $\cat{FHilb}$ can be specified as a triple $(v,d,\varepsilon)$ of a state and two effects.

\begin{dfn} \index{state!symmetric} \index{state!Frobenius}
Let $H$ be a Hilbert space. A state $v: \mathbb{C} \to H \otimes H \otimes H$ is \emph{symmetric} if it satisfies the equations $sym_L$, $sym_R$. It is a \emph{Frobenius state} if there exist effects $d: H \to \mathbb{C}$, $\varepsilon: H \otimes H \to \mathbb{C}$ such that $(v,d,\varepsilon)$ is a commutative Frobenius algebra on $H$.
\end{dfn}

\begin{thm} \label{thm:slocc_classification}
Let $v$ be a SLOCC-maximal, symmetric 3-qubit state. Then, $v$ is a Frobenius state. Moreover, if $(v,d,\varepsilon)$ is the associated commutative Frobenius algebra, then
\begin{itemize}
	\item $v$ is SLOCC-equivalent to $\ghz$ if and only if $(v,d,\varepsilon)$ is special, and
	\item $v$ is SLOCC-equivalent to $\wstate$ if and only if $(v,d,\varepsilon)$ is anti-special.
\end{itemize}
\end{thm}
\begin{proof}
These are \cite[Theorem 9 and Theorem 10]{coecke2010compositional}.
\end{proof}

This was partially extended to the case of qutrits, that is, 3-dimensional Hilbert spaces, where an intermediate class between special and anti-special appears, in \cite{honda2012graphical}.

So, this result traced 3-qubit SLOCC classification back to a classification of Frobenius algebras on the qubit. Moreover, Coecke and Kissinger showed that the Frobenius algebras associated with the $\ghz$ and $\wstate$ states, together with arbitrary single-qubit states, are universal for quantum computation, and listed several equations satisfied by the pair, some of them reminiscent of the axioms of ZX calculi. 

This led them to conjecturing that the two algebras could be building blocks for a new algebraic-diagrammatic theory, that could give access to a \emph{compositional} SLOCC classification for states of an arbitrary number of qubits; that is, one where representatives of higher-dimensional SLOCC classes are assembled through the composition of $Z$ spiders and $W$ spiders. 

Some additional equations, and a few preliminary results, mostly pertaining to universality, were given in \cite{coecke2011rational} and in Kissinger's PhD thesis \cite{kissinger2011pictures}. These were all based on the idea that, similarly to ZX calculi, the starting point of an axiomatisation should be $\textit{CFrob} \uplus \textit{CFrob}$, the theory of a pair of commutative Frobenius algebras.

Because the axiom $cup$ does not hold in this case, one of the most appealing features of ZX calculi --- the ability to treat string diagrams as undirected graphs --- seems to be lost. On the other hand, if we take the self-duality of ZX calculi as a primitive, or, in other words, we work with a fixed orthonormal basis, the $\wstate$ state 
does become ``fully symmetric'', in the way that the ternary $Z$ and $X$ spiders are; we can define operations
\begin{equation*} 
\input{img/c4_3_ternary_W.tex}
\end{equation*} 
and so on, interpreted as the linear maps $\ketbra{00}{1} + \ketbra{01}{0} + \ketbra{10}{0}$, the state $\ket{001} + \ket{010} + \ket{100}$, and the effect $\bra{001} + \bra{010} + \bra{100}$, each a transpose of another with respect to the fixed self-duality, and see them as manifestations of the same network
\begin{equation*} 
\begin{tikzpicture}[baseline={([yshift=-.5ex]current bounding box.center)}]
	\draw[edge, out=180, in=-45] (0,0) to node[auto, pos=1] {$A$} (-1, .5);
	\draw[edge, out=30, in=-90] (0,0) to node[auto, pos=1] {$B$} (.7, .9);
	\draw[edge, out=-60, in=90] (0,0) to node[auto, pos=1] {$C$} (.25,-1);
	\node[dotdark] at (0,0) {};
	\node[scale=1.25] at (1.5,-1) {.};
\end{tikzpicture}

\end{equation*} 
If we introduce an additional, equally symmetric operation
\begin{equation} \label{eq:binary_w}
\begin{tikzpicture}[baseline={([yshift=-.5ex]current bounding box.center)}]
\begin{pgfonlayer}{bg}
	\path[fill, color=gray!10] (-1,-1) rectangle (1,1);
\end{pgfonlayer}
\begin{pgfonlayer}{mid}
	\draw[edge] (0,-1) to (0,1);
	\node[dotdark] at (0,0) {};
	\node[scale=1.25] at (1.5,-1) {,};
\end{pgfonlayer}
\end{tikzpicture}

\end{equation} 
interpreted as the unitary operator $\ketbra{0}{1} + \ketbra{1}{0}$, we recover the generators of the $W$ Frobenius algebra as compositions of the new generators, as in the diagrams (\ref{eq:w_algebra}). Note that the operation (\ref{eq:binary_w}) is definable in the theory of a pair of commutative Frobenius algebras, interpreted as the $Z$ and the $W$ Frobenius algebra, as the composite map
\begin{equation*} 
\input{img/c4_3_w_tick.tex}
\end{equation*} 

It could be objected that it is undue to favour simpler diagrammatics over a ``hard'' result such as the classification of Frobenius algebras, and treat the $W$ Frobenius algebra as a derived object. On the other hand, one should be cautious not to over-interpret Theorem \ref{thm:slocc_classification}: the fact that commutative Frobenius algebras on the qubit are either special or anti-special is, after all, a consequence of non-trivial operators on the qubit being either rank-1 or full rank. In general, structures on 2-dimensional vector spaces can be expected to come in two classes, which makes it more likely that they would be pairwise related, only for the correspondence to break down in higher dimensions.

More importantly, if the compositional approach led to a classification of entanglement for $n > 3$ qubits, what could it possibly look like? The one unambiguous discriminant that diagrammatic algebra has offered with regard to entanglement, so far, is the implication of separability from the disconnectedness of string diagrams, and it is a topological criterion. 

Suppose that we had a complete diagrammatic axiomatisation of $\cat{Qubit}$, and a normalisation strategy rewriting any diagram $x$ into a diagram $\hat{x}$, so that $\hat{x}$ is ``as disconnected as possible'': no connected component of $x'$ is interpreted as a product state. Then, we could associate to any $x$ the number $n(\hat{x})$ of connected components of its normal form. This number would be SLOCC-invariant: if $y$ is obtained by composing $x$ with an invertible single-qubit operator --- by universality, there is a corresponding binary string diagram --- then $n(\hat{y}) = n(\hat{x})$. 

This would be a reasonable, well-defined starting point, from which we could look for more refined SLOCC-invariants. At the same time, it is one which can obviously benefit from the ``string diagram rewriting as graph rewriting'' approach employed with ZX calculi; hence, from privileging symmetry in the choice of the generators of the theory, while keeping a lookout for algebraically motivated axioms. This was our perspective in developing the ZW calculus, the subject of the next and last chapter.

  \thispagestyle{empty} 
\chapter{A calculus of qubits} \label{chap:zwcalculus}
\thispagestyle{plain}

\noindent\emph{In this chapter:}
\begin{itemize}
	\item[$\triangleright$] The ZW calculus is a complete diagrammatic axiomatisation of the theory of qubits. We define it gradually, building up on its fragments: wire, even, pure, vanilla, each introducing a new generator, and making some of the earlier axioms redundant. --- \emph{Section \ref{sec:fragments}}
	\item[$\triangleright$] We derive some consequences and inductive generalisations of the axioms, leading to a more compact set of rules. Then, we prove the main theorem, completeness for $R$-bits, and discuss some minor variants. --- \emph{Section \ref{sec:completeness}}
	\item[$\triangleright$] Abstracting from an interpretation of the pure fragment as a theory of fermionic oscillators, through its anyonic generalisations, we introduce analogues of the ZW axioms and generators for higher-dimensional quantum systems. We prove their universality, and propose a route towards completeness. --- \emph{Section \ref{sec:qudits}}
\end{itemize}

\section{The ZW calculus and its fragments} \label{sec:fragments}

\begin{itemize}
\item[]	\textbf{Note.} This section and the following are partially based on material previously published in \cite{hadzihasanovic2015diagrammatic}.
\end{itemize}

Let us forget about qubits for a moment, and consider the $W$ Frobenius algebra abstractly. At this level, there is nothing special about complex vector spaces: for all commutative rings $R$, we can define its analogue on the free $R$-module on a set of two elements.

\begin{dfn} \index{R-bit@$R$-bit, $R\cat{bit}$}
Let $R$ be a commutative ring. The PROP of \emph{$R$-bits} is the full monoidal subcategory $R\cat{bit}$ of $R\cat{Mod}$ whose objects are tensor products of a finite number of copies of $R \oplus R$.
\end{dfn}

We will keep using bra-ket notation for the morphisms of $R\cat{bit}$, writing $\ket{0}, \ket{1}$ for the generators of $R \oplus R$, and $\ket{b_1\ldots b_n}$ for $\ket{b_1} \otimes \ldots \otimes \ket{b_n}$, $b_i \in \{0,1\}$.

The comonoid part of the $W$ Frobenius algebra is defined by
\begin{align*}
	c_\bbullet: \ket{0} \mapsto \ket{0} \otimes \ket{0}, \quad \quad \quad & \ket{1} \mapsto \ket{0} \otimes \ket{1} + \ket{1} \otimes \ket{0}, \\
	d_\bbullet: \ket{0} \mapsto 1, \quad \quad \quad & \ket{1} \mapsto 0.
\end{align*} 
This comonoid is also part of a well-known Hopf algebra, called the \emph{fermionic line} in the theory of quantum groups \cite[Example 14.6]{majid2002quantum}. The fermionic line has the remarkable property that its monoid part is the transpose of its comonoid part with respect to the canonical self-duality. For this to be a model of $\textit{Hopf}$, however, the braiding must be interpreted not as the usual swap map of $R$-modules, but as the map \index{fermionic line}
\begin{equation*}
	x: \ket{b_1} \otimes \ket{b_2} \mapsto (-1)^{b_1 b_2}\ket{b_2} \otimes \ket{b_1},
\end{equation*}
that we will depict as
\begin{equation} \label{eq:crossing}
\input{img/c5_1_crossing.tex}
\end{equation}

We call diagram (\ref{eq:crossing}) and its interpretation in $R\cat{bit}$ the \emph{crossing}. The crossing satisfies all the axioms of a braiding, including the naturality axioms $nat^f_{b,L}$, $nat^f_{b,R}$, for all the Hopf algebra operations, that is, $c$, $d$ and their transposes in $R\cat{bit}$, and even for the self-duality maps of $R \oplus R$. However, it does not satisfy naturality axioms for all maps of $R\cat{bit}$, so it cannot be chosen as a braiding for $R\cat{bit}$ as a PROB: for example, for the linear map $\ketbra{0}{1} + \ketbra{1}{0}$, as in diagram (\ref{eq:binary_w}), \index{crossing}
\begin{equation*}
\input{img/c5_1_nonnatural.tex}
\end{equation*} 

In order to ``make the crossing a braiding'', one can choose to work in the category $R\cat{Mod}[\mathbb{Z}_2]$ of $\mathbb{Z}_2$-graded $R$-modules, also known as supermodules in the theory of supersymmetry \cite{varadarajan2004supersymmetry}. The objects of $R\cat{Mod}[\mathbb{Z}_2]$ are $R$-modules $G$ together with a direct sum decomposition $G \simeq G_0 \oplus G_1$; we call $G_0$ the \emph{even} part, and $G_1$ the \emph{odd} part of $G$. Morphisms $f: G_0 \oplus G_1 \to H_0 \oplus H_1$ are required to preserve the decomposition, that is, map $G_i$ to $H_i$, for $i = 0,1$. \index{module!graded}

Then, $R\cat{Mod}[\mathbb{Z}_2]$ becomes a braided monoidal category with the tensor product of $R$-modules, decomposed as
\begin{equation*}
	(G \otimes H)_0 = (G_0 \otimes H_0) \oplus (G_1 \otimes H_1), \quad \quad 	(G \otimes H)_1 = (G_0 \otimes H_1) \oplus (G_1 \otimes H_0),
\end{equation*}
the tensor unit $R \oplus 1$, and the braiding $G \otimes H \to H \otimes G$ defined, on elements $g \in G_i, h \in H_j$, for $i, j \in \{0,1\}$, by
\begin{equation*}
	g \otimes h \mapsto (-1)^{ij} h \otimes g.
\end{equation*}
Restricted to the full monoidal subcategory of $R\cat{Mod}[\mathbb{Z}_2]$ whose objects are tensor products of a finite number of copies of $R \oplus R$, seen as a $\mathbb{Z}_2$-graded $R$-module, this is precisely our crossing. Notice that the map $\ketbra{0}{1} + \ketbra{1}{0}$ is not a $\mathbb{Z}_2$-graded operator on $R \oplus R$, since it maps the ``even'' generator $\ket{0}$ to the ``odd'' generator $\ket{1}$, and vice versa. \index{RModgraded@$R\cat{Mod}[\mathbb{Z}_2]$}

\begin{remark}
That the comonoid part of the $W$ Frobenius algebra is part of a bialgebra in the category of $\mathbb{Z}_2$-graded complex vector spaces had already been noticed by Vicary, as reported by Kissinger in \cite[Section 10.1.2]{kissinger2011pictures}.
\end{remark}

More simply, if we take the point of view that a braiding is a generator of an algebraic theory as much as any other, the fermionic line is a model of $\textit{Hopf}$ in $R\cat{bit}$ seen as a map $f: \textit{Hopf} \to R\cat{bit}$ of PROs, rather than a map of PROBs: we encountered a similar situation in Example \ref{exm:tensorpro}, where the definition of a distributive law of monads relied on the interpretation of a braiding-like operation in a 2-category that has no natural braidings.

We will adopt such a neutral approach, keeping the crossing as a generator of $R\cat{bit}$, alongside the swap operation that makes it a PROP.

\begin{dfn} \index{ZW calculus!wire fragment} \index{ZWwire@$ZW_\mathrm{wire}$}
The \emph{wire fragment} of the ZW calculus is the PROP $ZW_\mathrm{wire}$ containing $\textit{CSDual}$, together with the additional generator $x: [2] \to [2]$, depicted as in diagram (\ref{eq:crossing}), and satisfying the axioms (all labels left implicit)
\begin{equation*}
\input{img/c5_1_crossing_ax.tex}
\end{equation*} 
\begin{equation*}
\input{img/c5_1_crossing_nat.tex}
\end{equation*} 
\begin{equation*}
\input{img/c5_1_crossing_framed.tex}
\end{equation*} 
together with the usual axioms of PROPs.
\end{dfn}
Theories of two different braidings such as this are the subject of virtual knot theory \cite{kauffman2012virtual}. The axioms $frm$ and $frm'$ appear as a replacement of the first Reidemeister move in the theory of ``blackboard-framed'' knots; since, by the so-called Whitney trick, it is a consequence of $rei^x_2$ and $rei^x_3$ that
\begin{equation} \label{eq:whitney}
\input{img/c5_1_whitney.tex}
\end{equation} 
it follows from $frm$ and $frm'$ that
\begin{equation} \label{eq:crossing_ax2}
\input{img/c5_1_crossing_ax2.tex}
\end{equation} 
We call this operation $i: [1] \to [1]$; it is interpreted as the map $\ketbra{0}{0} - \ketbra{1}{1}$ in $R\cat{bit}$. Equation (\ref{eq:crossing_ax2}) can be used as a replacement of $frm$ and $frm'$ in the presentation of $ZW_\mathrm{wire}$. We refer, for example, to \cite{kauffman2001knots} for more details on these topics.

\begin{dfn} \index{ZW calculus!even fragment} \index{ZWeven@$ZW_\mathrm{even}$}
The \emph{even fragment} of the ZW calculus is the PROP $ZW_\mathrm{even}$ containing $ZW_\mathrm{wire}$, together with the additional generator $w: [1] \to [2]$, depicted as
\begin{equation} \label{diag:even}
\input{img/c5_1_even.tex}
\end{equation} 
subject to the following axioms:
\begin{enumerate}
	\item $w$ forms a cocommutative comonoid with the discard operation $\varepsilon w$, that is,
\begin{equation*}
\input{img/c5_1_even_ax1.tex}
\end{equation*} 
\begin{equation*}
\input{img/c5_1_even_ax1b.tex}
\end{equation*} 
	\item $w$ is natural and cocommutative with respect to the crossing, that is,
\begin{equation*}
\input{img/c5_1_even_ax2.tex}
\end{equation*} 
	\item defining $m: [2] \to [1]$ to be the transpose of $w$
\begin{equation*}
\input{img/c5_1_even_monoid.tex}
\end{equation*} 
	$(w, \varepsilon w, m, m \eta)$ form a Hopf algebra with $i$ as the antipode and $x$ as the braiding, that is,
\begin{equation*}
\input{img/c5_1_even_hopf1.tex}
\end{equation*} 
\begin{equation*}
\input{img/c5_1_even_hopf2.tex}
\end{equation*}  
\end{enumerate}
\end{dfn}
By self-duality, other equations can be derived as transposes of the ones we gave: for example, $m$ is also part of a monoid with unit $m \eta$, commutative both with respect to the swap and to the crossing. Moreover, by a standard diagrammatic proof, see for example \cite[Figure 4.6]{majid2013quantum}, the antipode of a Hopf algebra is provably a comonoid anti-homomorphism, hence, by cocommutativity, a homomorphism:
\begin{equation*}
\input{img/c5_1_antipode_homo.tex}
\end{equation*} 

\begin{cons}
Already in the even fragment, we can construct some interesting maps. The following ``triangle'' of $W$ operations
\begin{equation} \label{eq:wtriangle}
\input{img/c5_1_triangle.tex}
\end{equation} 
is interpreted in $R\cat{bit}$ as the $\mathbb{Z}_2$ group $R$-algebra
\begin{equation*}
	\ket{0} \otimes \ket{0}, \ket{1} \otimes \ket{1} \mapsto \ket{0}, \quad \quad \quad \ket{0} \otimes \ket{1}, \ket{1} \otimes \ket{0} \mapsto \ket{1}.
\end{equation*}
Precomposing it with a crossing, we obtain an algebra with unit $\ket{0}$, and such that $\ket{1} \otimes \ket{1} \mapsto -\ket{0}$. If we let $1 := \ket{0}$, $i := \ket{1}$, in $\mathbb{R}\cat{bit}$ this is the algebra of complex numbers, seen as a real algebra.

More in general, we can take a number $n$ of copies of the $\mathbb{Z}_2$ algebra, seen as an $R$-algebra with a single non-unit generator $\gamma_i$, with $i = 1,\ldots, n$. When we compose them together, we have a choice of using either the swap or the crossing, that is,
\begin{equation*}
\input{img/c5_1_composite.tex}
\end{equation*} 
The first choice corresponds to making $\gamma_1, \gamma_2$ commute in the product, that is, have $\gamma_1 \cdot \gamma_2 = \gamma_2 \cdot \gamma_1$, with the identification $\gamma_1 \sim \gamma_1 \otimes 1$, and $\gamma_2 \sim 1 \otimes \gamma_2$. The second choice corresponds to making them anti-commute: $\gamma_1 \cdot \gamma_2 = -\gamma_2 \cdot \gamma_1$. Then, we can see that
\begin{equation*}
\input{img/c5_1_clifford.tex}
\end{equation*} 
is an algebra on $n = p+q$ generators satisfying 
\begin{align*}
	\gamma_i \cdot \gamma_j & = - \gamma_j \cdot \gamma_i, \quad \quad i \neq j, \\
	\gamma_i \cdot \gamma_i & = \begin{cases}
		1, & i \leq p, \\
		-1, & i > p,
	\end{cases}
\end{align*}
for all $i, j = 0, \ldots, n$. In $\mathbb{R}\cat{bit}$, this is the real Clifford algebra with signature $(p,q)$; for example, the case $(0,2)$ corresponds to the algebra of quaternions, and the case $(1,3)$ to the algebra of Dirac matrices. \index{Clifford algebra}

Mixing swaps and braidings, different algebras based on commutation and anticommutation relations can be constructed. For example, 
\begin{equation*}
\input{img/c5_1_dualquaternion.tex}
\end{equation*} 
is the algebra of dual quaternions, used in mechanics to represent rigid motions in 3-dimensional space \cite{yang1964dual}.
\end{cons}

\begin{remark}
The ``triangle'' (\ref{eq:wtriangle}) is a scalar multiple of the $X$ spider with two inputs and one output in the ZX calculus. After transposition, as a state, this is SLOCC-equivalent to $\ghz$; it follows that already the even fragment of the ZW calculus contains representatives of both maximal SLOCC classes of three qubits. We will discuss the significance of this fact in Section \ref{sec:qudits}.
\end{remark} 

The even fragment of the ZW calculus only contains operations that preserve the $\mathbb{Z}_2$-grading of $R \oplus R$, so $ZW_\mathrm{even}$ can in fact be interpreted in $R\cat{Mod}[\mathbb{Z}_2]$. Equivalently, for all states in the image of $ZW_\mathrm{even}$, written as linear combinations
\begin{equation} \label{eq:linearcomb}
	\sum_{i=1}^q r_i \ket{b_{i1}\ldots b_{in}},
\end{equation}
such that $0 \neq r_i \in R$, and no pair of $n$-tuples $(b_{i1}\ldots b_{in})$ is equal, 
\begin{equation*}
	\sum_{j=1}^n b_{ij} \equiv 0\,\mathrm{mod}\,2,
\end{equation*}
for all $i = 1,\ldots,q$. The next step is to introduce an ``odd'' state, which can be either $\ket{1}$, or $\ket{01} + \ket{10}$; the latter choice seems to make more sense diagrammatically, since it allows us to decompose the generator (\ref{diag:even}) in the way that we described at the end of Chapter \ref{chap:quantum}.

\begin{dfn} \index{ZW calculus!pure fragment} \index{ZWpure@$ZW_\mathrm{pure}$}
The \emph{pure fragment} of the ZW calculus is the PROP $ZW_\mathrm{pure}$ containing $ZW_\mathrm{wire}$, together with the generators $w_3: [0] \to [3]$, $w_2: [0] \to [2]$,
\begin{equation*}
\input{img/c5_1_wgenerators.tex}
\end{equation*} 
satisfying the following axioms:
\begin{enumerate}
	\item both $w_3$ and $w_2$ are symmetric with respect to the swap:
\begin{equation*}
\input{img/c5_1_w3symmetric.tex}
\end{equation*} 
\begin{equation*}
\input{img/c5_1_w2symmetric.tex}
\end{equation*} 
	\item defining $n: [1] \to [1]$ and $w: [1] \to [2]$ by
	\begin{equation*}
\input{img/c5_1_w_definition.tex}
\end{equation*} 
	$w$ satisfies all the axioms of $ZW_\mathrm{even}$, and $n$ satisfies
\begin{equation*}
\input{img/c5_1_not_ax.tex}
\end{equation*} 
\end{enumerate}
\end{dfn}
The pure fragment is interpreted in $R\cat{bit}$ by the assignment
\begin{equation*}
	w_3 \mapsto \ket{001} + \ket{010} + \ket{100}, \quad \quad \quad
	w_2 \mapsto \ket{01} + \ket{10}.
\end{equation*}
By the symmetry axioms $sym_{3,L}$, $sym_{3,R}$, it suffices to impose only one out of $\coo{un}_{w,L}$ and $\coo{un}_{w,R}$, and replace $nat_w^w$ with a coassociativity axiom. It follows from the axioms that $w_2$ is also symmetric with respect to the crossing:
\begin{equation} \label{eq:w2crossing}
\input{img/c5_1_w2crossing.tex}
\end{equation} 
Moreover, the axiom $hopf$ is derivable from the others in $ZW_\mathrm{pure}$. The proof is essentially the same as the derivation of complementarity from strong complementarity in Section \ref{sec:frobenius}:
\begin{equation} \label{eq:purehopf}
\input{img/c5_1_pure_hopf.tex}
\end{equation} 
where the second equation comes from composing both sides of $\coo{un}_{w,R}$ with $n$ and using $inv$, and in the last diagram we used $inv$ again to simplify. By the transposes of $nat_w^{m\eta}$, the last diagram in (\ref{eq:purehopf}) is equal to
\begin{equation*}
\input{img/c5_1_pure_hopf2.tex}
\end{equation*} 
proving $hopf$.

It is apparent from our presentation that the pure fragment of the ZW calculus, much like the ZX calculus, can be seen as a calculus of undirected graphs, with the \emph{proviso} that, when there is a crossing of wires, we are not always allowed to slide nodes through it. We have used this liberty in the last couple of proofs, where we freely employed transposed versions of the axioms, without specifying what wires, in particular, were transposed.

In fact, it is possible to introduce a ``spider'' presentation for $ZW_\mathrm{pure}$. This has countable generators $w_n: [0] \to [n]$ for all $n \in \mathbb{N}$, depicted as \index{spider presentation}
\begin{equation*}
\input{img/c5_1_w_spider.tex}
\end{equation*} 
and symmetric both with respect to the crossing and to the swap. These are interpreted in $R\cat{bit}$ as the states
\begin{equation*}
	\sum_{k=1}^n \ket{\underbrace{0\ldots 0}_{k-1}1\underbrace{0\ldots 0}_{n-k}\,};
\end{equation*} 
in particular, $w_0 \mapsto 0$, $w_1 \mapsto \ket{1}$, and $w_2$ and $w_3$ correspond to the ones defined earlier.

Spiders $[k] \to [n-k]$ are obtained by transposing any $k$ wires of $w_n$ with the self-duality maps. They satisfy the axioms:
\begin{equation*}
\input{img/c5_1_spider_cut.tex}
\end{equation*} 
that is, two nodes can be merged as long as there is another binary node between them, and loops on a single node can be eliminated. Together with $inv$, $sym_{3,L}$, $sym_{3,R}$ and $sym_2$, they subsume the axioms that make (\ref{eq:w_algebra}) a commutative Frobenius algebra.

The pure fragment of the ZW calculus only contains operations that either preserve, or reverse the $\mathbb{Z}_2$-grading of $R \oplus R$; equivalently, all states in the image of $ZW_\mathrm{pure}$, written as in (\ref{eq:linearcomb}), are either ``purely even'' or ``purely odd'', that is,
\begin{equation*}
	\sum_{j=1}^n b_{ij} \equiv 0\,\mathrm{mod}\,2 \text{ for all $i$, or } 	\sum_{j=1}^n b_{ij} \equiv 1\,\mathrm{mod}\,2 \text{ for all $i$.}
\end{equation*}
It turns out that to obtain linear combinations of even and odd states, it suffices to add a single generator of mixed parity, namely the state $\ket{000} + \ket{111}$, proportional to $\ghz$.
\begin{dfn} \index{ZW calculus!vanilla} \index{ZWvanilla@$ZW$}
The \emph{vanilla ZW calculus} is the PROP $ZW$ containing $ZW_\mathrm{pure}$, together with the generator $z_3: [0] \to [3]$, depicted as
\begin{equation*}
\input{img/c5_1_zgenerator.tex}
\end{equation*} 
subject to the following axioms: 
\begin{enumerate}
	\item $z_3$ is symmetric with respect to the swap:
\begin{equation*}
\input{img/c5_1_z3symmetric.tex}
\end{equation*}
	\item defining $c: [1] \to [2]$ by
\begin{equation*}
\input{img/c5_1_z_definition.tex}
\end{equation*}
	$c$ forms a cocommutative comonoid with the discard operation $\varepsilon c$, that is,
\begin{equation*}
\input{img/c5_1_vanilla_ax1.tex}
\end{equation*}
	\item $i: [1] \to [1]$ behaves as a comodule homomorphism, and $n: [1] \to [1]$ as a comonoid homomorphism with respect to $c$, that is,
\begin{equation*}
\input{img/c5_1_vanilla_ax2.tex}
\end{equation*}
	\item $c$ is related to $w: [1] \to [2]$ and its transposes by
\begin{equation*}
\input{img/c5_1_vanilla_ax3.tex}
\end{equation*}
\begin{equation*}
\input{img/c5_1_vanilla_ax3b.tex}
\end{equation*}
\end{enumerate}
\end{dfn}
Of course, $z_3$ is the same as the white ternary $Z$ spider of the $ZX$ calculus, and forms a special commutative Frobenius algebra with its transpose almost by definition. Notice that specialness follows from the definition of the discard operation:
\begin{equation*}
\input{img/c5_1_coas_special.tex}
\end{equation*}
where the first equality uses coassociativity.

Again, we can define a spider presentation, this time with no surprises: it has countable generators $z_n: [0] \to [n]$ for all $n \in \mathbb{N}$, depicted as
\begin{equation*}
\input{img/c5_1_z_spider.tex}
\end{equation*}
and symmetric with respect to the swap, interpreted in $R\cat{bit}$ as the states
\begin{equation*}
	\ket{\underbrace{0\ldots 0}_n} + \ket{\underbrace{1\ldots 1}_n}.
\end{equation*}
In particular, $z_0 \mapsto 2$, $z_1 \mapsto \ket{0} + \ket{1}$, and $z_2$ has the same interpretation as the self-duality map $\eta$. Spiders $[k] \to [n-k]$ are obtained by transposition of $k$ wires, and satisfy
\begin{equation*}
\input{img/c5_1_z_spider_cut.tex}
\end{equation*}
\begin{equation*}
\input{img/c5_1_z_elim.tex}
\end{equation*}

So far, we have had no requirements on $R$ beside the structure of an abelian group, and we have got as far as we could: we will prove in Section \ref{sec:completeness} that $ZW$ is isomorphic to $\mathbb{Z}\cat{bit}$. To proceed further, we need to introduce generators labelled by elements of $R$.

\begin{dfn} \index{ZW calculus!$R$-labelled} \index{ZWlabel@$ZW_R$}
Let $R$ be a commutative ring. The \emph{$R$-labelled ZW calculus} is the PROP $ZW_R$ containing $ZW$, together with a family of generators $\{r: [1] \to [1]\}_{r \in R}$ labelled by elements of $R$, depicted as
\begin{equation*}
\begin{tikzpicture}[baseline={([yshift=-.5ex]current bounding box.center)}]
\begin{pgfonlayer}{bg}
	\path[fill, color=gray!10] (-1,-1) rectangle (1,1);
\end{pgfonlayer}
\begin{pgfonlayer}{mid}
	\draw[edge] (0,-1) to (0,1);
	\node[dotwhite] at (0,0) [label=right:$r$] {};
	\node[scale=1.25] at (1.5,-1) {,};
\end{pgfonlayer}
\end{tikzpicture}

\end{equation*}
and satisfying the following axioms:
\begin{enumerate}
	\item composition and convolution by the Hopf algebra $(w, \varepsilon w, m, m\eta)$ correspond to multiplication and addition in $R$, respectively, that is, \index{convolution!by a Hopf algebra}
\begin{equation*}
\input{img/c5_1_ring_ax1.tex}
\end{equation*}
\begin{equation*}
\input{img/c5_1_ring_ax2.tex}
\end{equation*}
	\item any $r: [1] \to [1]$ is a comonoid homomorphism for $w$, and a comodule homomorphism for $c$, that is,
\begin{equation*}
\input{img/c5_1_ring_ax3.tex}
\end{equation*}
\begin{equation*}
\input{img/c5_1_ring_ax3b.tex}
\end{equation*}
\end{enumerate}
\end{dfn}
The $R$-labelled ZW calculus is interpreted in $R\cat{bit}$ by the assignment
\begin{equation*}
	r \mapsto \ketbra{0}{0} + r \ketbra{1}{1},
\end{equation*}
supplementing those of the vanilla fragment. Thanks to the axioms $ph^r$, $ZW_R$ can be given an $R$-labelled spider presentation, similar to the theories of $G$-labelled Frobenius algebras of Section \ref{sec:frobenius}; that is, the $Z$ spiders carry a label $r \in R$, and the $cut_z$ axioms multiply labels in $R$, as in \index{spider presentation!labelled}
\begin{equation*}
\input{img/c5_1_z_labelled.tex}
\end{equation*}

To summarise, we have defined a sequence of maps of PROPs
\begin{equation*}
\begin{tikzpicture}[baseline={([yshift=-.5ex]current bounding box.center)}]
	\node[scale=1.25] (0) at (0,0) {$\textit{CSDual}$};
	\node[scale=1.25] (1) at (2.5,0) {$ZW_\mathrm{wire}$};
	\node[scale=1.25] (2) at (5,0) {$ZW_\mathrm{even}$};
	\node[scale=1.25] (3) at (7.5,0) {$ZW_\mathrm{pure}$};
	\node[scale=1.25] (4) at (9.8,0) {$ZW$};
	\node[scale=1.25] (5) at (11.8,0) {$ZW_R$,};
	\draw[1c] (0) to (1);
	\draw[1c] (1) to (2);
	\draw[1c] (2) to (3);
	\draw[1c] (3) to (4);
	\draw[1c] (4) to (5);
\end{tikzpicture}
\end{equation*}
each corresponding to the addition of a single generator, except the last one, which corresponds to the addition of a family of generators: in order,
\begin{equation*}
\input{img/c5_1_all_generators.tex}
\end{equation*}
the second of which can be further decomposed after adding the third, to obtain a calculus of undirected labelled graphs, similar to the ZX calculus. 

Up to the vanilla ZW calculus, the presentations we gave --- in the non-spider version --- are all finitary, having a small number of ``small'' axioms, usually only involving the interaction of two generators. In fact, for a finitely generated ring $R$ (which is, consequently, finitely presented), it is possible to add operations $r: [1] \to [1]$ only for a finite family of generators $r \in R$, and encode their relations directly in the diagrammatic language, to obtain a finitary presentation of $ZW_R$.  

While this is a conceptually interesting point, spider presentations are generally more convenient to work with. In the next section, we will introduce an equivalent presentation of $ZW_R$, consisting of ``spiderised'' axiom schemes, and use it to prove that $ZW_R$ is isomorphic to $R\cat{bit}$, for all commutative rings $R$; in particular, that $ZW_\mathbb{C}$ is complete for $\cat{Qubit}$.

\section{Derived rules and the completeness theorem} \label{sec:completeness}

Given a commutative ring $R$, we consider the PROP $ZW_R$ in the spider presentation, that is, with additional generators $w_n, z^r_n: [0] \to [n]$ for all $n \in \mathbb{N}$, $r \in R$, 
\begin{equation*}
\input{img/c5_2_generators.tex}
\end{equation*}
connected to the original ones by the axiom schemes
\begin{equation*}
\input{img/c5_1_spider_cut.tex}
\end{equation*} 
\begin{equation*}
\input{img/c5_2_ring_spider.tex}
\end{equation*} 
and the identification $z_3 \equiv z_3^1$; when $r = 1$, in general, we will avoid writing any labels. First of all, these can replace the axioms that make $w, c: [1] \to [2]$ comonoids with units $\varepsilon w, \varepsilon c: [1] \to [0]$, respectively, and also the $ph^r$ axioms.

We will prove that several axioms of $ZW_R$ admit inductive generalisations, that can be used to replace lengthy sequences of applications of $cut$ rules, followed by rules for ternary or binary spiders. The formal way of handling such axiom schemes in automated graph rewriting was studied in \cite{kissinger2014pattern}.

\begin{prop} The following are derived rules of $ZW_R$ for all $w_n$, $z_n^r$ with $n \geq 2$:
\begin{equation*}
\input{img/c5_2_w_symmetric.tex}
\end{equation*} 
\begin{equation*}
\input{img/c5_2_z_symmetric.tex}
\end{equation*} 
\end{prop}
\begin{proof}
For $sym_w$, $sym_z$, this is a simple induction from $sym_{3,L}$, $sym_{3,R}$, $sym_2$, $sym_{z,L}$, $sym_{z,R}$, using the $cut_w$ and $cut_z$ rules. 

For $sym_w^x$, start from $\coo{com}_w$ coupled with $inv$ for the ternary case, and from equation (\ref{eq:w2crossing}) for the binary case, then proceed by induction.
\end{proof}
This is what allows us to treat the diagrams of the ZW calculus as undirected graphs, a liberty that we will casually exploit.

\begin{prop} The following are derived rules of $ZW_R$ for all $n \geq 2$:
\begin{equation*}
\input{img/c5_2_crossnat.tex}
\end{equation*} 
\end{prop}
\begin{proof}
The case $n = 0$ is $nat_x^w$ combined with $nat_x^\varepsilon$, while $n = 1$ follows immediately from $inv$. The general case is proved by induction from $nat_x^w$, using the $cut_w$ rules.
\end{proof}

\begin{prop} \label{prop:wbialg}
The following are derived rules of $ZW_R$ for all $n, m \in \mathbb{N}$:
\begin{equation*}
\input{img/c5_2_w_bialgebra.tex}
\end{equation*} 
\end{prop}
\begin{proof}
The case $n = m = 0$ is $nat_{\varepsilon w}^{m\eta}$ combined with $tr_w$ and $inv$; the case $n = 0, m > 1$ is an inductive generalisation of $nat_w^{m\eta}$, and similarly $m = 0, n > 1$; the case $n = 1$ or $m = 1$ follows immediately from $inv$. Finally, $n = m = 2$ is $nat_w^m$, and from there we can proceed by double induction on $n, m$, using the $cut_w$ rules and the $xnat$ rules to slide black nodes through crossings.
\end{proof}

\begin{prop} The following are derived rules of $ZW_R$ for all $r \in R$, and all $n, m \in \mathbb{N}$ such that either $n = m = 0$, or $m > 0$:
\begin{equation*}
\input{img/c5_2_cw_bialgebra.tex}
\end{equation*} 
\end{prop}
\begin{proof}
For the case $m = 2, n = 0$,
\begin{equation*}
\input{img/c5_2_bialg_02.tex}
\end{equation*} 
merging the white spiders and using $sym_z$, we see that the latter is equal to
\begin{equation*}
\input{img/c5_2_bialg_02b.tex}
\end{equation*} 
The case $m = n = 0$ then follows:
\begin{equation*}
\input{img/c5_2_bialg_00.tex}
\end{equation*} 

The case $m = 1$, $n = 0$ is just $nat^r_{\varepsilon c}$, the cases $m = n = 1$ and $m = 2$, $n= 1$ follow from $inv$, the case $m = 1, n = 2$ is $nat^r_c$, and $m = n = 2$ is a combination of $nat^r_c$ and $nat^m_c$. From here, we proceed by double induction, as in Proposition \ref{prop:wbialg}.
\end{proof}

\begin{prop}
The following are derived rules of $ZW_R$ for all $n \in \mathbb{N}$:
\begin{equation*}
\input{img/c5_2_automorph.tex}
\end{equation*}
\end{prop}
\begin{proof}
The case $n = 1$ follows from $id$ (or $rng_1$), the case $n = 2$ is $nat^n_c$, and $n = 0$ is a combination of the latter with $inv$. The general case is a simple induction, starting from $nat^n_c$ and using the $cut_z$ rules.
\end{proof}

\begin{prop} \label{prop:zwloop}
The following are derived rules of $ZW_R$ for all $r \in R$, $n \geq 2$:
\begin{equation*}
\input{img/c5_2_zwloop.tex}
\end{equation*} 
\end{prop}
\begin{proof}
The case $n = 2$ is $loop$, together with $nat^r_{\varepsilon c}$ to get rid of $r$. The sequence
\begin{equation*}
\input{img/c5_2_zwloop_proof.tex}
\end{equation*} 
of equalities proves the general case.
\end{proof}

\begin{prop}
The following are derived rules of $ZW_R$ for all $n \in \mathbb{N}$, $r_i \in R$, $i = 1,\ldots, n$:
\begin{equation*}
\input{img/c5_2_sumrule.tex}
\end{equation*}
\end{prop}
\begin{proof}
The case $n = 0$ comes from
\begin{equation*}
\input{img/c5_2_sumrule_proof.tex}
\end{equation*}
The case $n = 1$ is just an application of $inv$, and $n = 2$ is $rng_+^{r,s}$. The general case is then proved by induction, similarly to Proposition \ref{prop:zwloop}.
\end{proof}

\begin{prop} \label{prop:crossminus}
The following are derived rules of $ZW_R$, for all $r \in R$:
\begin{equation*}
\input{img/c5_2_crossminus.tex} 
\end{equation*}
\end{prop} 
\begin{proof} Use $inv$ to introduce a pair of black nodes on the left hand side; then,
\begin{equation*}
\input{img/c5_2_crossminus_proof.tex} 
\end{equation*}
Transposing some wires and using $cut_w, cut_z$, and $id$, we can see that this is equal to
\begin{equation*}
\input{img/c5_2_crossminus_proof2.tex} 
\end{equation*}
The result follows by $inv$ and $sym_z$.
\end{proof}

Before proceeding, we recall how our interpretation map $zw_R: ZW_R \to R\cat{bit}$ is defined on the generators.
\begin{equation*}
\begin{tabular}{l l}
	$s: [2] \to [2]$ & $\mapsto \; \ketbra{00}{00} + \ketbra{01}{10} + \ketbra{10}{01} + \ketbra{11}{11}$ \\
	$\eta: [0] \to [2]$ & $\mapsto \; \ket{00} + \ket{11}$ \\
	$\varepsilon: [2] \to [0]$ & $\mapsto \; \bra{00} + \bra{11}$ \\
	$x: [2] \to [2]$ & $\mapsto \; \ketbra{00}{00} + \ketbra{01}{10} + \ketbra{10}{01} - \ketbra{11}{11}$ \\
	$w_n: [0] \to [n]$ & $\mapsto \; \sum\limits_{k=1}^n \ket{\underbrace{0\ldots 0}_{k-1}1\underbrace{0\ldots 0}_{n-k}\,}$ \\
	$z_n^r: [0] \to [n]$ & $\mapsto \; \ket{\underbrace{0\ldots 0}_n\,} + r\ket{\underbrace{1\ldots 1}_n\,}$.
\end{tabular}
\end{equation*}
It can be checked that all the axioms are sound for $R\cat{bit}$, that is, the map $zw_R$ is well-defined.

\begin{thm}[Universality of the ZW calculus] \label{thm:universality} \index{universality!of the ZW calculus}
The map $zw_R: ZW_R \to R\cat{bit}$ is full.
\end{thm}
\begin{proof}
Because of self-duality, it sufficies to check that every state in $R\cat{bit}$ is in the image of $zw_R$; any other morphism can be obtain by transposition. Write an arbitrary $n$-partite state $v$ in the form (\ref{eq:linearcomb}), that is, 
\begin{equation*} 
	\sum_{i=1}^m r_i \ket{b_{i1}\ldots b_{in}},
\end{equation*}
where $r_i \neq 0$, and no pair of $n$-tuples $(b_{i1},\ldots,b_{in})$ is equal. We claim that $v$ is the image through $zw_R$ of the diagram
\begin{equation} \label{eq:normalform}
\input{img/c5_2_normalform.tex}
\end{equation}
where, for $i = 1,\ldots, m$ and $j= 1,\ldots, n$, the dotted wire connecting the $i$-th white node to the $j$-th output is present if and only if $b_{ij} = 1$.

We can check this by a combination of diagrammatic and algebraic reasoning. By the interpretation of $w_m$, the diagram (\ref{eq:normalform}) is interpreted as the sum $\sum_{i=1}^m v_i$, where $v_i$ is the interpretation of
\begin{equation*}
\input{img/c5_2_normal_proof.tex}
\end{equation*}
with $\ket{1}$ plugged into the $i$-th white node, and $\ket{0}$ in all the others. By $ba_{zw}$, this is equal to
\begin{equation*}
\input{img/c5_2_normal_proof2.tex}
\end{equation*}
where the last equality uses either $cut_w$ or $ba_w$, depending on the wire being present or not, that is, on $b_{1j}$ being $1$ or $0$. We can do the same for all $i \neq j$, until we are left with the diagram
\begin{equation*}
\input{img/c5_2_normal_proof3.tex}
\end{equation*}
whose interpretation, by direct calculation, is $r_i \ket{b_{i1}\ldots b_{in}}$. 
\end{proof}

We say that a string diagram of $ZW_R$ is in \emph{normal form} if it is of the form (\ref{eq:normalform}) for some state $v$ of $R\cat{bit}$; this is unique up to a permutation of the white nodes, and can be made strictly unique by picking a specific ordering for the summands of (\ref{eq:linearcomb}), for example the one induced by the lexicographic ordering of the $n$-tuples $(b_{i1}, \ldots, b_{in})$. \index{normal form}

\begin{remark}
An embryo of this normal form appeared in \cite{bruni2006basic}, whose authors considered analogues of the $Z$ and $W$ monoids as the building blocks of a diagrammatic language for the category of finite sets and relations, seen as modules over the semiring of Booleans. 

While they achieved a certain completeness result, their axiomatisation included a large number of complicated axioms, including an axiom scheme not reducible to any finite set of equations. Nevertheless, it led to further work on algebras of connectors for the study of concurrent systems \cite{sobocinski2013connector}, which ended up crossing paths with ZX calculi \cite{bonchi2014interacting}.
\end{remark}

In order to prove that $ZW_R$ is complete for $R\cat{bit}$, it now suffices to show that any string diagram of $ZW_R$ can be rewritten in normal form using the axioms. In order to do that, we will prove in turn that
\begin{enumerate}
	\item any composite of two diagrams in normal form can be rewritten in normal form, and
	\item all generators can be rewritten in normal form.
\end{enumerate}
\begin{dfn} \index{pre-normal form}
A string diagram of $ZW_R$ is in \emph{pre-normal form} if it is of the form (\ref{eq:normalform}), where the following are also allowed:
\begin{itemize}
	\item two white nodes may be connected to the same outputs;
	\item $r_i$ may be $0$ for some $i$;
	\item a single white node may have more than one connection to the same black node.
\end{itemize}
\end{dfn}
\begin{exm}
The following string diagrams of $ZW_\mathbb{Z}$ are all in pre-normal form, but only the last one is in normal form:
\begin{equation*}
\input{img/c5_2_prenormal.tex}
\end{equation*}
In the first diagram, the second and third white node are both only connected to the second output; in the second diagram, the second node has two connections to the first output. Each of the diagrams, however, is interpreted in $\mathbb{Z}\cat{bit}$ as the state $\ket{00} + 2\ket{01} + \ket{11}$.
\end{exm}

In the next few proofs, we will often ``zoom in'' on a certain portion of a string diagram, which may require some reshuffling of nodes, using swapping or transposition of wires, with the implicit understanding that this can always be reversed later.

\begin{lem} \label{lem:prenormal}
Any string diagram in pre-normal form can be rewritten in normal form.
\end{lem}
\begin{proof}
Consider a diagram in pre-normal form, and suppose that two white nodes are connected to the same outputs. The relevant portion of the diagram looks like
\begin{equation*}
\input{img/c5_2_prenormal_proof.tex}
\end{equation*}
Now, the two input wires both lead to the bottom black node; zooming in on that, we find
\begin{equation} \label{eq:prenormal_proof}
\input{img/c5_2_prenormal_proof2.tex}
\end{equation}
which has a single white node, connected to the same outputs, replacing the two initial ones. 

Now, suppose that there is a white node labelled $0$. The relevant part of the diagram is
\begin{equation*}
\input{img/c5_2_prenormal_proof3.tex}
\end{equation*}
which after an application of $cut_w$ simply eliminates the white node.

Finally, suppose that there is a white node with two wires connecting it to the same black node, and let this be an output. The relevant part of the diagram is of the form
\begin{equation*}
\input{img/c5_2_prenormal_proof4.tex}
\end{equation*}
which, again, simply eliminates the white node. The case in which the black node is the bottom one is similar. 
\end{proof}
In the last proof, we have started using $cut_z$ and $cut_w$ rules tacitly, which we will do more and more often. 

\begin{lem}[Negation] \label{lem:negation}
The composition of one output of a string diagram in pre-normal form with $n: [1] \to [1]$ can be rewritten in normal form, and that has the effect of ``complementing'' the connections of the output to white nodes; that is, locally, 
\begin{equation*}
\input{img/c5_2_negation.tex}
\end{equation*}
\end{lem}
\begin{proof}
Using $cut_z$ and $cut_w$, we can rewrite the left hand side as
\begin{equation*}
\input{img/c5_2_negation_proof.tex}
\end{equation*}
which, using the $aut$ rule and transposing some wires, becomes
\begin{equation*}
\input{img/c5_2_negation_proof2.tex} 
\end{equation*}
From here, we can apply the same reasoning backwards, applying $ba_{zw}$ to the rightmost $m$ nodes, which leads us to the result. 
\end{proof}

In the following, and later statements, ``plugging one output of a string diagram into another'' means post-composition with the self-duality map $\varepsilon: [2] \to [0]$. 
\begin{lem}[Trace] \label{lem:trace}
The plugging of two outputs of a string diagram in pre-normal form into each other can be rewritten in normal form.
\end{lem}
\begin{proof}
While this can be proved directly using the $ba_w$ and $ba_{zw}$ rules, there is a simple proof making a repeated use of the negation lemma \ref{lem:negation}. Focus on the two relevant outputs, and subdivide the white nodes into four groups, based on their being connected to both outputs, to one of them, or neither of them:
\begin{equation*}
\input{img/c5_2_trace.tex} 
\end{equation*}
Using the negation lemma on the rightmost output, this becomes
\begin{equation*}
\input{img/c5_2_trace2.tex} 
\end{equation*}
after merging the black nodes with $cut_w$, the leftmost $n$ white nodes have two wires connecting them to the same black node. This means that $loop$ is applicable, leaving us with
\begin{equation*}
\input{img/c5_2_trace3.tex} 
\end{equation*}
The leftmost black nodes can be merged with any black node they are connected to, or eliminated with $ba_w$ if there is none, which rids us of the leftmost $n$ white nodes. Now, we can apply the negation lemma again, to find that the remaining portion of the diagram is equal to
\begin{equation*}
\input{img/c5_2_trace4.tex} 
\end{equation*}
which, focussing on the rightmost part, is equal, by $ba_w$, to
\begin{equation*}
\input{img/c5_2_trace5.tex} 
\end{equation*}
This rids us of the rightmost $q$ white nodes, and leaves us with a diagram in pre-normal form.
\end{proof}

A string diagram consisting of a single black node and no wires is interpreted as the scalar 0. The following lemma proves that it acts as an ``absorbing element'' for string diagrams in normal form.

\begin{lem}[Absorption] \label{lem:absorb}
For all diagrams in pre-normal form, the following is a derived rule of $ZW_R$:
\begin{equation} \label{eq:absorbed}
\input{img/c5_2_absorption.tex} 
\end{equation}
\end{lem}
\begin{proof}
Expanding the black node with $cut_w$, we can treat it as an output of a diagram in pre-normal form. Then, applying the negation lemma,
\begin{equation*}
\input{img/c5_2_absorption_proof.tex} 
\end{equation*}
where the new output is connected to all white nodes. From here, we can proceed as in the last part of the proof of Lemma \ref{lem:trace} in order to eliminate all the white nodes.
\end{proof}

We are now able to give the central proof of our completeness theorem.
\begin{thm} \label{thm:functoriality}
Any composition of two string diagrams in pre-normal form can be rewritten in normal form.
\end{thm}
\begin{proof}
We can factorise any composition of string diagrams in pre-normal form as a tensor product followed by a sequence of ``self-pluggings''; thus, by the trace lemma \ref{lem:trace}, it suffices to prove that a tensor product --- diagrammatically, the juxtaposition of two string diagrams in pre-normal form --- can be rewritten in normal form.

Given such a setup, we can create a pair of black nodes connected by a wire using $ba_w$, and apply the negation lemma on both sides, in order to obtain
\begin{equation*}
\input{img/c5_2_monoidality.tex} 
\end{equation*}
which is the plugging of two outputs connected to all the white nodes of their respective diagrams. The only case in which this still leaves the two diagrams disconnected is when one of the diagrams has no white nodes, that is, it looks like the right hand side of equation (\ref{eq:absorbed}). In this case, by Lemma \ref{lem:absorb}, we can use its isolated black node to ``absorb'' the other diagram, which produces a diagram in normal form.

So, suppose that $n, m > 0$. Focussing on the two outputs, 
\begin{equation} \label{eq:wbialgebra}
\input{img/c5_2_monoidality_2.tex} 
\end{equation}
then, we can use $ba_{zw}$ on each of the white nodes: for example, on the leftmost one, this leads to
\begin{equation*}
\input{img/c5_2_monoidality_3.tex} 
\end{equation*}
Each of the outgoing wires leads to a black node, so we can push the nodes indicated by arrows to the outside, and merge them using $cut_w$. Repeating this operation, we can push all black nodes to the outside, which leaves us with a tangle of $n\cdot m$ wires connecting white nodes, one for each pair $(r_i, r'_j)$, where $i = 1, \ldots, n$, and $j = 1, \ldots, m$.

This tangle is made of crossings, and not swaps, so we cannot push the white nodes through. However, each of the white nodes has at least one wire connecting it to a black node, which means that the axiom $unx$ is applicable: this allows us to turn all the crossings into swaps, and use $cut_z$ to merge each connected pair together.

After some rearranging, this leaves us with a string diagram of the form
\begin{equation*}
\input{img/c5_2_monoidality_4.tex} 
\end{equation*}
where we used the negation lemma again. A final application of $ba_w$ produces a string diagram in pre-normal form.
\end{proof} 

Theorem \ref{thm:functoriality} implies that the normal form is ``functorial'': that is, if we extend the definition of normal form from states to generic linear maps of $R\cat{bit}$, for example by choosing a specific way of transposing their inputs, we obtain a map of PROPs from $R\cat{bit}$ to $ZW_R$, that has $zw_R$ as a left inverse. It remains to prove that this is a two-sided inverse; for this, it is sufficient to prove that every generator of $ZW_R$ can be rewritten in normal form.

\begin{thm}[Completeness of the ZW calculus] \label{thm:completenesszw}
The map $zw_R: ZW_R \to R\cat{bit}$ is an isomorphism of PROPs. 
\end{thm} \index{completeness!of the ZW calculus}
\begin{proof}
Since the spider and non-spider presentations are equivalent, it suffices to show that the generators of the latter can be rewritten in normal form. For $w_3$ and $w_2$,
\begin{equation*}
\input{img/c5_2_w3_rewrite.tex} 
\end{equation*}
For $z_3^r$, with $r \in R$, 
\begin{equation*}
\input{img/c5_2_z2_rewrite.tex} 
\end{equation*}
and similarly for $z_2^r$. The self-duality maps can be seen as a special case of the latter, $r = 1$, using $id$ (or $rng_1$). 

For the crossing, we can start by rewriting the tensor product of two self-duality maps in normal form, which gives us
\begin{equation*}
\input{img/c5_2_dualities_rewrite2.tex} 
\end{equation*}
then,
\begin{equation*}
\input{img/c5_2_dualities_rewrite.tex} 
\end{equation*}
by Proposition \ref{prop:crossminus}, this is equal to
\begin{equation*}
\input{img/c5_2_dualities_rewrite3.tex} 
\end{equation*}
The case of the swap is similar, and easier. This concludes the proof.
\end{proof}

In particular, $ZW_\mathbb{C}$ is a complete axiomatisation of $\cat{Qubit}$. Observe that, in this case, we can also easily implement the dagger diagrammatically: the adjoint of a string diagram of $ZW_\mathbb{C}$ is the vertical reflection of the diagram, with labels $\lambda \in \mathbb{C}$ of white nodes turned into their complex conjugates $\overline{\lambda}$. For example,
\begin{equation*}
\input{img/c5_2_dagger_diagram.tex} 
\end{equation*}

We mentioned earlier that, instead of introducing generators $r: [1] \to [1]$ for all elements of a ring $R$, we can only introduce one for each of a family of generators of $R$, together with one axiom for each relation that they satisfy. Then, in the normal form, instead of having a white node labelled $r \in R$ at each end of the bottom black spider, we will need to have some canonical expression of $r$ by sums and products of generators. The completeness proof still goes through: essentially, we can work with diagrams in pre-normal form, where terms in a sum of products of generators are decomposed into different legs of the bottom spider, until the very end, then proceed as in Lemma \ref{lem:prenormal}, but stop before performing the steps in equation (\ref{eq:prenormal_proof}).

For example, the vanilla ZW calculus corresponds to having a single generator and no relations, so $ZW$ is isomorphic to $\mathbb{Z}\cat{bit}$, $\mathbb{Z}$ being the free commutative ring on one generator; this is what we originally proved in \cite{hadzihasanovic2015diagrammatic}. In this case, an integer $n$ is represented in the normal form as
\begin{equation*}
\input{img/c5_2_integers.tex} 
\end{equation*}
depending on $n$ being positive on negative. We can then obtain an axiomatisation of $\mathbb{Z}_n\cat{bit}$ for all $n \in \mathbb{N}$ simply by adding the axiom
\begin{equation*}
\input{img/c5_2_znloop.tex} 
\end{equation*}
The case $n = 2$ was considered in \cite{schumacher2012modal} under the name of modal quantum theory. The corresponding ZW calculus is particularly simple, for 
\begin{equation*}
\input{img/c5_2_mobits_swap.tex} 
\end{equation*}
becomes provable, making several axioms redundant.

In the case of $\cat{Qubit}$, a possibility is to only introduce generators $a: [1] \to [1]$ for positive real numbers, together with a single generator for the imaginary unit, satisfying
\begin{equation} \label{eq:complex_unit}
\input{img/c5_2_complex_unit.tex} 
\end{equation}
Then, the map $\lambda: [1] \to [1]$ for an arbitrary complex number $\lambda = a + bi$ can be represented as 
\begin{equation*}
\input{img/c5_2_complex_from_real.tex} 
\end{equation*}
or other variations, depending on $a$, $b$ being positive or negative. What is interesting, here, is that
\begin{equation*}
\input{img/c5_2_complex_loop.tex} 
\end{equation*}
a valid equation of ribbons, is strongly reminiscent of equation (\ref{eq:complex_unit}); moreover, the interpretation of a ribbon twist as the imaginary unit would be compatible with the interpretation of vertical reflection as the dagger of $\cat{Qubit}$.  

On the other hand, we could interpret the positive real generators $a: [1] \to [1]$ as ``wires with a length''; passing to logarithms in the labels, we can also make them compose additively, rather than multiplicatively. Overall, these are suggestions that there may be a further underlying geometry of the complex ZW calculus, which we are yet to uncover.

\subsection*{List of rules of the ZW calculus}
For ease of reference, we provide here the full list of axioms for $ZW_R$ in the spider presentation. Restricting to $r, s = 1$ and removing the last two axioms gives an axiomatisation of $ZW$. We provide inductive generalisations directly when they capture several axioms in one scheme. The axioms of PROPs are implied.

\begin{equation*}
\input{img/zw_dual_1.tex}
\end{equation*} 
\begin{equation*}
\input{img/zw_wire_1.tex}
\end{equation*} 
\begin{equation*}
\input{img/zw_wire_2.tex}
\end{equation*} 
\begin{equation*}
\input{img/zw_pure_1.tex}
\end{equation*} 
\begin{equation*}
\input{img/zw_pure_2.tex}
\end{equation*} 
\begin{equation*}
\input{img/zw_pure_3.tex}
\end{equation*} 
\begin{equation*}
\input{img/zw_full_1.tex}
\end{equation*} 
\begin{equation*}
\input{img/zw_full_2.tex}
\end{equation*} 
\begin{equation*}
\input{img/zw_full_3.tex}
\end{equation*} 
\begin{equation*}
\input{img/zw_full_4.tex}
\end{equation*} 
\begin{equation*}
\input{img/zw_full_5.tex}
\end{equation*}

\section{Anyonic oscillators. Towards qudits} \label{sec:qudits}

There are at least two directions in which we can push the $\cat{Qubit}$ completeness result:
\begin{itemize}
	\item up in dimension, trying to obtain analogous axiomatisations for $\cat{Qudit}$, the full monoidal subcategory of $\cat{FHilb}$ whose generating object is $\mathbb{C}^d \equiv \ell^2(\mathbb{Z}_d)$, for an arbitrary $d \in \mathbb{N}$;
	\item ``in depth'', trying to find a complete set of axioms for the subcategories spanned by fragments of the ZW calculus.
\end{itemize}
In fact, there are reasons to believe that a better understanding of the fragments may help solving the first problem as well. \index{qudit} \index{Quditcat@$\cat{Qudit}$}

Arguably, the crux of the completeness theorem lies in the possibility of merging together string diagrams in normal form, what we proved in Theorem \ref{thm:functoriality}. The ends of a diagram in normal form are all $n$-ary multiplications in the $W$ algebra, the generator of the even fragment; and the initial step of the proof, equation (\ref{eq:wbialgebra}), relies on our ability to slide them past their transposes, as granted by the $ba_w$ rules.

This indicates that the real potential of the $W$ algebra comes from its being part of a self-transpose Hopf algebra, rather than an anti-special Frobenius algebra, as suggested in \cite{coecke2010compositional}; and a theory of self-transpose Hopf algebras with an alternative symmetric braiding is what $ZW_\mathrm{even}$ describes. Hence, it is worth considering whether its model in $\cat{Qubit}$ generalises to $\cat{Qudit}$. 

Before answering this question, let us take a detour to consider the physical meaning of the even and pure fragments in isolation; for the basic notions involved, we refer to any textbook on many-particle quantum systems, such as \cite{boffi2004heisenberg}, and to \cite{blute1994fock, vicary2008oscillator} for a categorical approach. 

The maps of $\cat{Qubit}$ in the image of $ZW_\mathrm{pure}$ --- to which we can freely add the generators $\lambda: [1] \to [1]$ of $ZW_\mathbb{C}$ --- are valid transformations of a type of quantum system other than the qubit: the \emph{fermionic oscillator}. We translate between the two by viewing $\ket{0}, \ket{1}$ as the two states of the single-mode fermionic Fock space: $\ket{0}$ is the vacuum state, and $\ket{1}$ the state of one particle; no other states are possible, by the Pauli exclusion principle. \index{fermionic oscillator}

The $\mathbb{Z}_2$-grading of the maps in the even fragment can be interpreted as conservation of energy: if one particle comes in, one particle comes out. More concretely than in the qubit interpretation, we can see string diagrams as physical circuits transmitting fermions: for example, the generator of $ZW_\mathrm{even}$,
\begin{equation*}
\input{img/c5_1_even.tex}
\end{equation*}
can be seen as a beam splitter, sending the one-particle state $\ket{1}$ to the superposition $\ket{01} + \ket{10}$ of the particle being on the left branch, and the particle being on the right branch of a circuit. 

In the pure fragment, we are furthermore allowed to introduce new particles into the system, something that would violate conservation of energy in an isolated system --- whence the loss of $\mathbb{Z}_2$-grading, and naturality of the crossing --- yet makes perfect sense operationally. This enables us to define the operation
\begin{equation*}
\input{img/c5_3_creation.tex}
\end{equation*}
which introduces a new particle on a circuit; it is interpreted as the fermionic creation operator $\dagg{a}: \ket{0} \mapsto \ket{1}, \ket{1} \mapsto 0$. Dually, its adjoint \index{creation operator} \index{annihilation operator}
\begin{equation*}
\input{img/c5_3_annihilation.tex}
\end{equation*}
is interpreted as the annihilation operator $a: \ket{1} \mapsto \ket{0}, \ket{0} \mapsto 0$.

Using a mixture of string diagrams and linear algebra, we can even read the canonical anticommutation relations off the axioms: because, in the model, \index{commutation relation} \index{anticommutation relation}
\begin{equation*}
\input{img/c5_3_antipode_linear.tex}
\end{equation*}
we have
\begin{equation*}
\input{img/c5_3_anticommutation.tex}
\end{equation*}
which we can read as the equation $a\dagg{a} = 1 - \dagg{a}a$.

Fermionic oscillators support a model of computation that is computationally equivalent to the qubit circuit model \cite{bravyi2002fermionic}, yet operationally very different, with contrasting locality and entanglement properties \cite{dariano2014feynman}. If, as we strongly suspect, the pure fragment of the ZW calculus, together with the $\{\lambda: [1] \to [1]\}_{\lambda \in \mathbb{C}}$, is universal for fermionic quantum computation, it becomes interesting in its own right to find a complete set of axioms. 

Incidentally, one of the weaknesses of the ZW calculus with respect to ZX calculi is that it is less obvious whether a certain diagram represents a unitary transformation of qubits: in comparison, the generators of ZX calculi have a closer link to the standard building blocks of qubit circuits, namely the CNOT gate, phase gates and the Hadamard gate. Yet the weakness for one model of computation may turn out to be a strength for another, if diagrams in the pure fragment have a direct operational translation into fermionic circuits.

\begin{remark} 
To obtain the full ZW calculus from its pure fragment, all we need is to add
\begin{equation} \label{eq:copymap}
\input{img/c5_3_copymap.tex}
\end{equation}
In the fermionic interpretation, this can be seen as a non-physical operation that ``copies'' particles. There is an intriguing parallel with the way classical logic is obtained from linear logic, and cartesian categories from monoidal categories: is qubit computation, formally, the non-resource sensitive version of fermionic computation?
\end{remark} 

\begin{remark}
With a complete set of axioms at hand, we may want to solve the entanglement classification problem for fermionic quantum computation as a stepping stone towards the qubit case. Because nodes have a single colour in the pure fragment, the diagrams are more purely ``topological'': for example, rather than $\wstate$ and $\ghz$, we could take $\wstate$ and the ``triangle'',
\begin{equation*}
\input{img/c5_3_twoslocc.tex}
\end{equation*}
two topologically inequivalent networks, as representatives of the two 3-qubit SLOCC-maximal classes. In fact, 
\begin{equation*}
\input{img/c5_3_fourslocc.tex}
\end{equation*}
some of the non-obviously equivalent 4-qubit states that one can represent in the pure fragment, all belong to distinct SLOCC super-classes, as defined in \cite{lamata2007inductive}. However, we have not yet attempted to find representatives for each super-class, nor to establish a formal correspondence.
\end{remark} 

In generalising the pure fragment to an arbitrary dimension $d \in \mathbb{N}$, we want to keep these basic intuitions. The simplest $d$-level generalisation of a fermion is a special type of abelian anyon \cite{panangaden2010anyons}: a particle that, exchanged with another particle, receives a phase $q := e^\frac{2i\pi}{d}$, a primitive $d$-th root of the unit. This is modelled by the unitary operator on $\mathbb{C}^d \otimes \mathbb{C}^d$ defined by 
\begin{equation*}
	x: \ket{k} \otimes \ket{j} \mapsto q^{jk} \, \ket{j} \otimes \ket{k}, \quad \quad j, k = 0,\ldots,d-1,
\end{equation*}
with inverse
\begin{equation*}
	\dagg{x}: \ket{k} \otimes \ket{j} \mapsto q^{-jk} \, \ket{j} \otimes \ket{k}, \quad \quad j, k = 0,\ldots,d-1,
\end{equation*}
which we picture, respectively, as
\begin{equation*}
\input{img/c5_3_dcrossing.tex}
\end{equation*}
These satisfy the second and third Reidemeister moves, and $x^d = (\dagg{x})^d = 1$, or, equivalently, $\dagg{x} = x^{d-1}$. In the case $d = 2$, we recover the crossing of the ZW calculus, with $\dagg{x} = x$. 

Because we do not want to lose circuit-like undirectedness, we are inclined to keep the canonical self-duality maps,
\begin{equation*}
	\eta := \sum_{j=0}^{d-1} \ket{j} \otimes \ket{j}, \quad \quad \varepsilon := \sum_{j=0}^{d-1} \bra{j} \otimes \bra{j},
\end{equation*}
for all $d$. Unlike in the $d=2$ case, these do not respect, in general, the obvious $\mathbb{Z}_d$-grading, so they do not commute with the crossing; however, they do satisfy a kind of ``skew-naturality'',
\begin{equation*}
\input{img/c5_3_wire_skew.tex}
\end{equation*}
They also satisfy
\begin{equation*}
\input{img/c5_3_wire_framed.tex}
\end{equation*}
interpreted as the map $\sum_{j=0}^{d-1} q^{j^2} \ketbra{j}{j}$. We can see these as the generalisations of the axioms of the wire fragment.

In line with our discussion above, we expect the even fragment to be generalised by a cocommutative comonoid $(w,v)$ that
\begin{enumerate}
	\item is $\mathbb{Z}_d$-graded, or satisfies ``conservation of energy'', that is, maps $\ket{n}$ to a superposition of $\ket{k}\otimes \ket{n-k}$, $k = 0,\ldots, n$;
	\item forms a bialgebra with its transpose monoid, and $x$ as the braiding.
\end{enumerate}
We will depict this comonoid and its transpose as
\begin{equation*}
\input{img/c5_1_splitter.tex}
\end{equation*}
to avoid suggesting symmetries that are not present, and imposing in particular that 
\begin{equation} \label{eq:anyonic}
\input{img/c5_3_bialgebra.tex}
\end{equation}
The Hopf algebra that we will define is a ``symmetrisation'' of the anyonic line \cite[Chapter 16]{majid2002quantum}, a standard generalisation of the fermionic line, which satisfies the first requirement but is not self-transpose. To state the definition, we need some basic notions of $q$-arithmetic; see \cite[Chapter 7]{majid2002quantum} for more details. \index{anyonic line}

\begin{dfn} For all $q \in \mathbb{C}$, $n \in \mathbb{N}$, the \emph{$q$-integer} $\qint{n}$ is defined by \index{qinteger@$q$-integer}
\begin{equation*}
	\qint{n} := \sum_{k=0}^{n-1} q^k.
\end{equation*}
The \emph{$q$-factorial} of $n$ is then defined by \index{qfactorial@$q$-factorial}
\begin{equation*}
	\qint{n}! := \prod_{k=1}^n \qint{k},
\end{equation*}
and, for $k = 0,\ldots,n$, the \emph{$q$-binomial coefficients} are \index{qbinomial@$q$-binomial coefficient}
\begin{equation*}
	\binom{n}{k}_q := \frac{\qint{n}!}{\qint{k}!\qint{n-k}!}.
\end{equation*}
\end{dfn} 
We recover the standard integers, factorial, and binomial coefficients with the choice $q = 1$. When $q = e^\frac{2i\pi}{d}$, a primitive $d$-th root of the unit, we have in particular 
\begin{equation*}
	\qint{d} = 0,
\end{equation*}
hence $\qint{n}! = 0$ for all $n \geq d$. For our purposes, the most important thing about these numbers is that they satisfy the \emph{$q$-Vandermonde identities}
\begin{equation} \label{eq:vandermonde}
	\binom{n}{k}_q = \sum_{i=0}^k q^{(j-i)(k-i)}\,\binom{j}{i}_q \binom{n-j}{k-i}_q,
\end{equation} 
for all $n$, and $j, k = 0,\ldots,n$.

Having fixed $q = e^\frac{2i\pi}{d}$, we define a comonoid on $\mathbb{C}^d$ with copy map
\begin{equation*}
	w: \ket{n} \mapsto \sum_{k=0}^n \binom{n}{k}_q^\frac{1}{2} \ket{k} \otimes \ket{n-k}, 
\end{equation*}
for $n = 0,\ldots, d-1$, and discard map $v: \ket{0} \mapsto 1, \ket{n} \mapsto 0$ for $n=1,\ldots,d-1$; notice the exponent over the $q$-binomial coefficient. It is not hard to verify the coassociativity, cocommutativity and counitality of $(w,v)$. 
	
\begin{prop} \label{eq:selftranspose}
The comonoid $(w,v)$ forms a bialgebra with its transpose monoid, and $x$ as a braiding.
\end{prop}
\begin{proof}
We focus on equation (\ref{eq:anyonic}); the other three bialgebra axioms are easy to check. The left hand side can be calculated, in steps, to be
\begin{align*}
	\ket{j} & \otimes \ket{n-j} \; \mapsto \; \sum_{i=0}^j \sum_{l=0}^{n-j} \binom{j}{i}_q^\frac{1}{2} \binom{n-j}{l}_q^\frac{1}{2} \, \ket{i} \otimes \ket{j-i} \otimes \ket{l} \otimes \ket{n-j-l}  \\ 
	& \mapsto \; \sum_{i=0}^j \sum_{l=0}^{n-j} q^{l(j-i)} \binom{j}{i}_q^\frac{1}{2} \binom{n-j}{l}_q^\frac{1}{2} \, \ket{i} \otimes \ket{l} \otimes \ket{j-i} \otimes \ket{n-j-l}  \\
	& \mapsto \; \sum_{i=0}^j \sum_{l=0}^{n-j} q^{l(j-i)} \binom{j}{i}_q^\frac{1}{2} \binom{n-j}{l}_q^\frac{1}{2} \binom{i+l}{i}_q^\frac{1}{2} \binom{n-i-l}{j-i}_q^\frac{1}{2} \, \ket{i+l} \otimes \ket{n-i-l} = \\
	& = \; \sum_{k=0}^n \sum_{i=0}^{k} q^{(k-i)(j-i)} \binom{j}{i}_q^{\frac{1}{2}} \binom{n-j}{k-i}_q^{\frac{1}{2}} \binom{k}{i}_q^{\frac{1}{2}} \binom{n-k}{j-i}_q^{\frac{1}{2}} \, \ket{k} \otimes \ket{n-k},
\end{align*}
where in the last equality we substituted $k := i+l$ for $l$. The right hand side is
\begin{equation*}
	\ket{j} \otimes \ket{n-j} \; \mapsto \; \sum_{k=0}^n \binom{n}{j}_q^{\frac{1}{2}} \binom{n}{k}_q^{\frac{1}{2}} \, \ket{k} \otimes \ket{n-k};
\end{equation*}
comparing the coefficients, we see that it suffices to prove
\begin{equation} \label{eq:binomproof}
	\binom{n}{j}_q^{\frac{1}{2}} \binom{n}{k}_q^{\frac{1}{2}} = \sum_{i=0}^{k} q^{(k-i)(j-i)} \binom{j}{i}_q^{\frac{1}{2}} \binom{n-j}{k-i}_q^{\frac{1}{2}} \binom{k}{i}_q^{\frac{1}{2}} \binom{n-k}{j-i}_q^{\frac{1}{2}}.
\end{equation}
Suppose, without loss of generality, that $k = j + m$, with $m$ nonnegative. Then, expanding the $q$-binomial coefficients, we find that
\begin{equation*}
	\binom{n}{j}_q = \binom{n}{k}_q \, \prod_{l=1}^m\frac{\qint{j+l}}{\qint{n-(j+l-1)}},
\end{equation*}
so the left hand side of (\ref{eq:binomproof}) is equal to
\begin{equation*}
	\binom{n}{k}_q \, \left(\prod_{l=1}^m\frac{\qint{j+l}}{\qint{n-(j+l-1)}}\right)^\frac{1}{2}.
\end{equation*}
On the other hand,
\begin{equation*}
	\binom{k}{i}_q = \binom{j}{i}_q \, \prod_{l=1}^m\frac{\qint{j+l}}{\qint{j+l-i}}, \quad \quad \binom{n-k}{j-i}_q = \binom{n-j}{k-i} \, \prod_{l=1}^m \frac{\qint{j+l-i}}{\qint{n-(j+l-1)}},
\end{equation*}
so the right hand side of (\ref{eq:binomproof}) is equal to
\begin{equation*}
	\left( \sum_{i=0}^{k} q^{(k-i)(j-i)} \binom{j}{i}_q \binom{n-j}{k-i}_q \right) \left(\prod_{l=1}^m \frac{\qint{j+l}}{\qint{j+l-i}} \frac{\qint{j+l-i}}{\qint{n-(j+l-1)}} \right)^\frac{1}{2},
\end{equation*}
and the result follows by the $q$-Vandermonde identities (\ref{eq:vandermonde}).
\end{proof}

If $d = q = 1$, because the binomial coefficients are never 0, the comonoid $(w,v)$ can in fact be defined on the underlying infinite-dimensional vector space of $\ell^2(\mathbb{N})$, letting $n$ range over all natural numbers; by lack of self-duality maps, Proposition \ref{eq:selftranspose} still holds for the pair of $(w,v)$ and its transpose with respect to the fixed basis.

Except in the infinite-dimensional case, the bialgebra is in fact a Hopf algebra: depicting the state $\ket{d-1}$ as
\begin{equation*}
\begin{tikzpicture}[baseline={([yshift=-.5ex]current bounding box.center)}]
\begin{pgfonlayer}{bg}
	\path[fill, color=gray!10] (-1,-1) rectangle (1,1);
\end{pgfonlayer}
\begin{pgfonlayer}{mid}
	\draw[edge] (0,1) to (0,0);
	\node[dotdark] at (0,0) [label=below:$d-1$] {};
	\node[scale=1.25] at (1.5,-1) {,};
\end{pgfonlayer}
\end{tikzpicture}

\end{equation*}
we find that it satisfies
\begin{equation*}
\input{img/c5_3_disconnect.tex}
\end{equation*}
and the proof (\ref{eq:purehopf}) goes through, starting from equation \ref{eq:anyonic}, with the antipode
\begin{equation*}
\input{img/c5_3_antipode.tex}
\end{equation*}
This is the same as the antipode of the anyonic line, and can be explicitly calculated as in \cite[Example 14.6]{majid2002quantum} to be
\begin{equation*}
	\ket{n} \mapsto q^\frac{n(n-1)}{2} (-1)^n \ket{n}, \quad \quad n = 0,\ldots, d-1.
\end{equation*}

\begin{exm}
We give explicit expressions for the first few low-dimensional cases. With $d = 2$, we recover the pure fragment of the ZW calculus. For $d = 3$, the $w$ operation is defined by
\begin{equation*}
	\ket{0} \mapsto \ket{00}, \quad \quad  \ket{1} \mapsto \ket{01} + \ket{10}, \quad \quad
	\ket{2} \mapsto \ket{02} + e^\frac{i\pi}{6} \ket{11} + \ket{20},
\end{equation*}
and the antipode is
\begin{equation*}
	\ket{0} \mapsto \ket{0}, \quad \quad \ket{1} \mapsto -\ket{1}, \quad \quad \ket{2} \mapsto e^\frac{2i\pi}{3} \ket{2}.
\end{equation*}
For $d = 4$, the $w$ operation is
\begin{equation*}
	 \ket{0} \mapsto \ket{00}, \quad \quad \ket{1} \mapsto \ket{11}, \quad \quad \ket{2} \mapsto \ket{02} + \sqrt[4]{2} \, e^\frac{i\pi}{8} \ket{11} + \ket{20},
\end{equation*}
\begin{equation*}
	\ket{3} \mapsto \ket{03} + e^\frac{i\pi}{4} \ket{12} + e^\frac{i\pi}{4} \ket{21} + \ket{30},
\end{equation*}
and the antipode is
\begin{equation*}
	\ket{0} \mapsto \ket{0}, \quad \quad \ket{1} \mapsto -\ket{1}, \quad \quad \ket{2} \mapsto i\ket{2}, \quad \quad \ket{3} \mapsto -i \ket{3}.
\end{equation*}
The unit is always $\ket{0}$.
\end{exm}

In order to generalise the pure fragment, we want to enable particle creation and annihilation, which amounts to introducing the generator 
\begin{equation*}
\begin{tikzpicture}[baseline={([yshift=-.5ex]current bounding box.center)}]
\begin{pgfonlayer}{bg}
	\path[fill, color=gray!10] (-1,-1) rectangle (1,1);
\end{pgfonlayer}
\begin{pgfonlayer}{mid}
	\draw[edge] (0,1) to (0,0);
	\node[dotdark] at (0,0) [label=below:$1$] {};
	\node[scale=1.25] at (1.5,-1) {,};
\end{pgfonlayer}
\end{tikzpicture}

\end{equation*}
interpreted as the one-particle state $\ket{1}$. Then, we can construct creation and annihilation operators $\dagg{a}$, $a$,
\begin{equation*}
\input{img/c5_3_anyonic_creation.tex}
\end{equation*}
so that 
\begin{equation*}
\input{img/c5_3_almost_antipode.tex}
\end{equation*}
leading from the bialgebra axiom (\ref{eq:anyonic}) to the commutation relation
\begin{equation*}
\input{img/c5_3_anyon_commutation.tex}
\end{equation*}
that is, $a\dagg{a} = 1 + q\,\dagg{a}a$. In the infinite-dimensional case with $q = 1$, we recover the bosonic creation and annihilation operators
\begin{align*}
	\dagg{a}: \ket{n} \, \mapsto \; & \binom{n+1}{1}^\frac{1}{2} \ket{n+1} = \sqrt{n+1} \ket{n+1}, \\
	a: \ket{n} \, \mapsto \; & \binom{n}{1}^\frac{1}{2} \ket{n-1} = \sqrt{n} \ket{n-1},
\end{align*}
with their canonical commutation relation $a\dagg{a} - \dagg{a}a = 1$.

Next, we want to generalise the
\begin{equation*}
\input{img/c5_3_complex_number.tex}
\end{equation*}
For the completeness proof of the ZW calculus, the most important property that these satisfy is being comonoid homomorphisms for $w$, what allows us to have arbitrary labels in the $ba_{zw}$ axiom scheme. Imposing that $\lambda: \mathbb{C}^d \to \mathbb{C}^d$ be a $\mathbb{Z}_d$-graded map, which forces it to be of the form $\sum_{k=0}^{d-1} \lambda_k\, \ketbra{k}{k}$, and that 
\begin{equation*}
\input{img/c5_3_complex_commute.tex}
\end{equation*}
we find that necessarily $\lambda_n = \lambda_k \lambda_{n-k}$ for all $n$, $k = 0,\ldots,n$; by recursion, $\lambda_n = \lambda_1^n$, and we can identify $\lambda$ with the complex number $\lambda_1$. The convolution
\begin{equation} \label{eq:dconvolution}
\input{img/c5_3_convolution.tex}
\end{equation}
evaluates to 
\begin{equation*}
	\ket{n} \mapsto \left(\sum_{k=0}^n \binom{n}{k}_q \,\lambda^k \,\mu^{n-k}\right) \ket{n};
\end{equation*}
by the $q$-binomial theorem \cite[Lemma 7.1]{majid2002quantum}, the coefficient would be $(\lambda + \mu)^n$, if $\lambda, \mu$ were elements of a noncommutative algebra satisfying $\mu\lambda = q \, \lambda \mu$.

Now, we expect any generalisation of the copy generator (\ref{eq:copymap}) to be of the form
\begin{equation*}
	c: \ket{n} \mapsto c_n\,\ket{n} \otimes \ket{n}, \quad \quad n = 0,\ldots, d-1,
\end{equation*}
for some coefficients $c_n$. Imposing that $c$ form a bialgebra with the transpose of $w$, 
\begin{equation*}
\input{img/c5_3_bialgebra_copy.tex}
\end{equation*}
we obtain the recursive relation 
\begin{equation*}
	c_n = \binom{n}{k}_q^\frac{1}{2} c_k\,c_{n-k},
\end{equation*}
which, with $c_0 = c_1 = 1$, resolves as
\begin{equation*}
	c_n := \sqrt{\qint{n}!}
\end{equation*}
for all $n \in \mathbb{N}$. The following equation, generalising $loop$, also holds:
\begin{equation*}
\input{img/c5_3_dloop.tex}
\end{equation*}

We can verify that $c$ forms a cocommutative comonoid with the discard operation $\sum_{k=0}^{d-1} \frac{1}{c_k} \bra{k}$, depicted as
\begin{equation*}
\begin{tikzpicture}[baseline={([yshift=-.5ex]current bounding box.center)}]
\begin{pgfonlayer}{bg}
	\path[fill, color=gray!10] (-1,-1) rectangle (1,1);
\end{pgfonlayer}
\begin{pgfonlayer}{mid}
	\draw[edge] (0,-1) to (0,0.0);
	\node[dotwhite] at (0,0.0) {};
	\node[scale=1.25] at (1.5,-1) {,};
\end{pgfonlayer}
\end{tikzpicture}

\end{equation*}
and that the analogue of $ph^r$ holds; together, these allow us to introduce $\mathbb{C}$-labelled white spiders $z_n^\lambda: [0] \to [n]$, interpreted as
\begin{equation*}
	\sum_{k=0}^{d-1} c_k^{n-2} \lambda^k \ket{\underbrace{k\ldots k}_n}.
\end{equation*}

Then, we can prove that the generators \index{universality!of the ZW calculus}
\begin{equation*}
\input{img/c5_3_all_generators_d2.tex}
\end{equation*}
together with the swap and self-duality maps are universal for $\cat{Qudit}$. Write an arbitrary $n$-qudit state $v$ as
\begin{equation*}
	\sum_{i=1}^m \lambda_i \ket{k_{i1},\ldots, k_{in}},
\end{equation*}
for some $\lambda_i \in \mathbb{C}$, $k_{ij} = 0,\ldots,d-1$, and let
\begin{equation*}
	\tilde{\lambda}_i := \lambda_i \, \prod_{j=1}^n \sqrt{\frac{1}{\qint{k_{ij}}!}}.
\end{equation*}
We claim that $v$ is the image of 
\begin{equation*}
\input{img/c5_3_dnormalform.tex}
\end{equation*}
where we used spider notation for the $w$ comonoid and its transpose, and exactly $k_{ij}$ wires connect the $i$-th white node to the $j$-th output. The proof is the same as universality of the ZW calculus, Theorem \ref{thm:universality}, with the only catch that
\begin{equation*}
	\ket{\underbrace{1\ldots 1}_k} \mapsto \sqrt{\qint{k}!} \; \ket{k}, \quad \quad k < d,
\end{equation*}
with the $w$ multiplication, which is why we need to adjust the coefficients $\lambda_i$.  

This is what we know so far. The main building blocks of the qubit ZW calculus have qudit analogues, with similar relations, and so does its normal form. If we had an analogue of $unx$, allowing us to disentangle crossed wires --- the obvious generalisation does not hold --- we could directly transport the main part of the qubit completeness proof, Theorem \ref{thm:functoriality}, to the qudit case. Quite unexpectedly, with the choice $q = 1$, we even obtained a partial calculus for infinite-dimensional quantum systems.

We do not know, on the other hand, whether the general $w$ comonoid admits a decomposition into ``fully symmetric'' components, leading to a calculus of undirected graphs, in the style of ZX calculi and the qubit ZW calculus. We also do not know what algebraic structure plays, in the qudit case, the role that commutative rings played in the qubit ZW calculus, in the sense that a ``vanilla qudit ZW calculus'' would axiomatise ``the free such-and-such on $d^n$ generators''; the interpretation of the convolution (\ref{eq:dconvolution}) would suggest some kind of braided algebra, whose elements commute up to a factor $q$. 

All questions of completeness, for $\cat{Qudit}$ and its fragments, theories of anyonic oscillators, remain open.
  \thispagestyle{empty} 

\backmatter \pagestyle{plain}
\chapter{Conclusion}
\thispagestyle{plain}

I started the research that led to this thesis with many questions on the exact placement of string diagrams --- seen sometimes as a curiosity, at most a convenient calculational tool --- at the interface between algebra and geometry. It led me to recognise the role of polygraphs, with their dual nature between combinatorics of spaces and universal algebra, at the crux of the question, and to develop a blueprint for a compositional universal algebra, with topology-based interactions of polygraphs translating to interactions of algebraic theories, in Chapter \ref{chap:interacting}.

Although inspired by various sources, and based on intuitions shared by researchers in higher algebra, category theory, and rewriting theory, this is, to the best of my knowledge, an entirely novel framework, that did not come as an answer to any pre-existing technical questions. It is my opinion, moreover, that it has not reached a stable foundation yet: I did my best, in Chapter \ref{chap:polygraphs}, to connect various strands of research on strict $\omega$-categories and add some missing pieces, but the formalism itself seems to be inadequate for some of my purposes. In Chapter \ref{chap:directed}, I gave an indication of what an improved, usable foundation could look like.

As a consequence, at least by some of the usual standards, it is hard to assess my contribution at this time: we need to wait, and see whether it will produce new methods that will solve old problems. On the positive side, because it consists in a new connection between topology and higher-dimensional algebra, there is a wealth of ideas and techniques that only await translation from one context into the other. Moreover, there is a certain simplicity to it, a combinatorial nature that may be the mark of something fundamental that was overlooked; perhaps, if the question was not asked, it is because we did not expect that we could still learn something about, say, the notion of homomorphism.

By the same standards, it is easier to evaluate my other main contribution: the ZW calculus is the first diagrammatic calculus complete for the theory of qubits, solving a problem that had been open for almost a decade. My result has already been used to derive the strongest axiomatisations, to date, for ZX calculi, a major strand of categorical quantum mechanics. Furthermore, compared to previous partial completeness results, this relies on a normal form with a simple algebraic interpretation, rather than some \emph{ad hoc} factorisation, and lends itself to a variety of generalisations. While the ZW calculus still has to find serious applications in quantum information and computation, I am confident that they will come with time, as the focus in the community partially shifts from ZX to ZW.

On both sides, there is no shortage of directions for future research. The following is a brief survey, starting with a recap of the most immediate ones, which we have amply discussed in Subsection \ref{sec:special} and Section \ref{sec:qudits}, and progressing towards the long shots, in no particular thematic order.

\begin{enumerate}
\item \emph{Completeness of the qudit ZW calculi.} Of all developments, this seems most within reach: the qudit generators and equations I gave in Section \ref{sec:qudits} are tailored to support critical steps of the qubit completeness proof, so it may be only the ``minor'' axioms that need to be adapted. If this is brought to completion, a subsequent goal would be to connect all the qudit calculi, and obtain a full axiomatisation of $\cat{FHilb}$.

\item \emph{General theory of representable regular polygraphs.} Here, too, we have a fully developed theory in low dimension, and an intuition of how it should generalise to higher dimensions. Formulating the different divisibility conditions in arbitrary dimension is going to be the main challenge; an answer to the next item in this list may be helpful. There are various constructions that need to be subsequently implemented, in order to reach the same flexibility as $\omega$-categories: presumably, a ``representable completion'' functor, then truncations and skeleta, and so on.

\item \emph{Enumeration of globes.} If we hope to ever have a computational implementation of regular polygraphs beyond the lowest dimensions, we need a recursive procedure that generates all non-isomorphic globes. A simple characterisation of different ``divisibility configurations'' would be an added bonus. A related problem is connecting globular posets to the $\omega$-categorical combinatorics, and settling the conjecture on their relation to directed complexes.

\item \emph{Completeness of the pure fragments.} The pure fragment of the ZW calculus has the remarkable property, not shared by other quantum calculi, that any string diagram in its language represents a physically meaningful experimental setup, interpreted as a fermionic circuit diagram. However, we know little about it in isolation, not even that it is universal for fermionic quantum computation. The same holds for the anyonic fragments of the qudit ZW calculi.

\item \emph{Compositional algebra in Globular.} My approach shares with the \emph{Globular} proof assistant a focus on string diagrams and on freely generated higher categories; the main point of divergence is that \emph{Globular} has an algebraic approach to weakness, based on a ``quasistrict'' algebra of composition with weak interchange. As we improve the understanding of low-dimensional representable regular polygraphs, and especially of their strictifiability properties, it would be useful to make an explicit comparison, with the goal of implementing tensor products and other constructions of compositional universal algebra in the proof assistant.

\item \emph{Constructions on $\cat{Cosp}(\cat{Hilb}_1)$.} In Chapter \ref{chap:quantum}, we reconstructed the category $\cat{Hilb}_{\leq 1}$ of Hilbert spaces and contractions, together with its dagger, as the truncation of a bicategory of cospans. This led to the question: which dagger-preserving, functorial constructions on $\cat{Hilb}_{\leq 1}$ lift to this bicategory, and which ones are actually induced by functors on $\cat{Hilb}_1$, the category of Hilbert spaces and isometries? 

\item \emph{Compositional rewriting theory and directed topology.} Some of the fundamental methods of algebraic topology --- Seifert-van Kampen, Mayer-Vietoris, monoidality of homology functors --- are of a compositional nature; some of them have been given directed versions for other notions of directed space. It is plausible that, in conjunction with Squier's theorem and related methods, these could be used to obtain rewriting-theoretic results compositionally. Moreover, the particular notion of directed space embodied by (regular) polygraphs may admit more refined invariants, that still need to be defined.

\item \emph{Fermionic and qubit entanglement classification.} As we saw in Chapter \ref{chap:quantum}, the classification of entanglement is the purpose for which a calculus based on the GHZ and W states was initially conceived. In Chapter \ref{chap:zwcalculus}, I suggested that the case of fermionic entanglement, better reflected in the topology of diagrams, may be simpler to tackle first. Of course, this would depend on the completion of the fourth item in this list.

\item \emph{``Continuous'' directed spaces.} The definition of globes in an oriented thin poset demands, recursively, decompositions into pairs of mergeable globes, until an atomic globe is reached. It is conceivable that a coalgebraic version exists, where the bi-partitions continue indefinitely. In \cite{freyd2008algebraic}, Freyd reconstructed the real interval as a terminal coalgebra for a bi-partition operation in dimension 1; perhaps, coalgebraic $n$-globes may generalise this, leading to a continuous, ``point-set'' notion of directed space, in the spirit of coalgebraic real analysis, more refined than the $d$-spaces or partially ordered spaces used in directed algebraic topology.

\item \emph{Simpson's Conjecture.} In Chapter \ref{chap:directed}, I connected the notions of Saavedra and, indirectly, of Joyal-Kock weak unit to representability conditions on regular polygraphs. This opens the way to abstract generalisations of the concepts that were used to solve Simpson's conjecture, on the semistrictifiability of homotopy types, in the special case of dimension 3. If the theory of representable regular polygraphs reaches an advanced enough stage, we may hope to use it, eventually, to prove a version of the conjecture in arbitrary dimensions.
\end{enumerate}
These only concern the direct objects of this thesis: I could add the exploration of those connections that have only been suggested, or brushed by, from proof theory to topological quantum field theory. 

My expectation is that, in the future, ideas of diagrammatic algebra, distilled to their computational essence, will be an integral part of our understanding of higher-dimensional algebra, beyond the three dimensions where our ``geometric intuition'' supports us; string-diagrammatic equations will be seen as low-dimensional examples of the general principles at play. In my thesis, I focussed on individuating one component of this essence: the combinatorial-topological compositionality of higher algebraic theories. I hope it will be a part of the picture.
 \thispagestyle{empty} 
\chapter*{Acknowledgements}
\thispagestyle{plain}

Lists of names are not in a particular order.

I am grateful to my supervisor, Bob Coecke, for his feedback and encouragement throughout these years. I have immense respect for his anti-dogmatic attitude, and commitment to training and supporting independent researchers; what he has helped build in Oxford is one of a kind.

I am also grateful to the other architect of the Quantum Group, Samson Abramsky, for all the times that I have come to him with the vaguest idea, and left with the exact reference for which I should have been looking. Special thanks also to Jamie Vicary, my internal examiner, on whose initiative a substantial higher categories sub-group has grown in Oxford since my arrival, making it a great place for the more abstract-leaning of us; to Aleks Kissinger, for his help early in my doctorate, as he was passing the baton of research on the GHZ/W calculus; to Paul-Andr\'e Melli\`es, for conversations in Paris that turned out to be course-changing for my project; and to Joachim Kock, my external examiner, for precious feedback on an earlier version of this thesis. I have also benefitted from conversations with Chris Heunen, Kohei Kishida, Ross Duncan, Simon Perdrix, Dusko Pavlovic, Emmanuel Jeandel.

Among fellow doctoral students, past and present, for many discussions from which I have learnt something, or valuable comments on my work, special thanks to Alex Kavvos, David Reutter, Miriam Backens, Antonin Delpeuch, Stefano Gogioso, Brendan Fong, Linde Wester, Christoph Dorn, Matthijs V\'ak\'ar, Dominic Verdon, Dan Marsden, William Zeng.

On a more personal side, there are many people who, in these years, have made the good times great, and the rest bearable. In Oxford, these include, but are not limited to Samuel Meister, Marco Molteni, Shriya Misra, Paul Boes; outside of Oxford, Valentina Alfarano, Chiara Cotroneo, Francesco Fiero, Tommaso Sciotto, Martina Rosola, Lorenzo Mileti Nardo, Paola Codazzi, Daniele Malpetti, Davide Costanzo, Margherita Boselli, Giulio Bagattini. Thanks also to Marco Sgrignoli and Manfredi Pumo, for co-authoring \texttt{SIGARETTO}, which now retroactively encompasses this thesis.

A very special thanks to my partner Sophie Turner, for more than can be said, and to my family, always.

I am grateful to the EPSRC, the Department of Computer Science, and Wolfson College, Oxford, for their financial and practical support.
 \thispagestyle{empty} \cleardoublepage 
\cleardoublepage \addcontentsline{toc}{section}{Bibliography}
\bibliographystyle{alpha}
\small \bibliography{bibliography} \clearpage \thispagestyle{empty}
\cleardoublepage \addcontentsline{toc}{section}{Index}
\normalsize \printindex[main]

\end{document}